**МИНИСТЕРСТВО ОБРАЗОВАНИЯ И НАУКИ УКРАИНЫ**

**ЗАПОРОЖСКИЙ НАЦИОНАЛЬНЫЙ УНИВЕРСИТЕТ**

**С.В. Курапов**
**М.В. Давидовский**

# АЛГОРИТМИЧЕСКИЕ МЕТОДЫ КОНЕЧНЫХ ДИСКРЕТНЫХ СТРУКТУР

# ТОПОЛОГИЧЕСКИЙ РИСУНОК ГРАФА

## часть 2

**(монография)**

Визуализированный граф является мощным инструментом в задачах анализа и синтеза данных. При этом задача визуализации заключается не только в выводе вершин и ребер согласно представлению графа, но и в том, чтобы этот вывод был визуально прост для пользователя. Процесс визуализации включает решение нескольких задач, одной из которых является задача построения топологического рисунка плоской части непланарного графа с минимальным количеством удаленных ребер. В работе рассматривается математическая модель описания топологического рисунка графа, основанная на методах теории вращения вершин с индукцией простых циклов удовлетворяющих критерию планарности Маклейна. Показано, что топологический рисунок непланарного графа может быть построен на основе выделенной плоской части графа. Топологическая модель рисунка графа позволяет свести переборную задачу выделения плоского графа к дискретной задачи оптимизации – поиску подмножества изометрических циклов графа удовлетворяющих нулевому значению функционала Маклейна. Для выделения плоской части графа разработан вычислительный метод на основе линейной алгебры и алгебры структурных чисел, с полиномиальной вычислительной сложностью.

Для научных работников, студентов и аспирантов, специализирующихся на применение методов прикладной математики.

# Содержание

## Введение

Общая задача построения топологического рисунка несепарабельного графа состоит из четырех задач:

1. определение планарности графа и в случае положительного ответа, построение его топологического рисунка;

2. построение топологического рисунка максимально плоской части несепарабельного непланарного графа с минимальным числом удаленных ребер;

3. построениетопологического рисунка графа с минимальным числом пересечений ребер;

4. построение топологического рисунка с минимальной толщиной графа.

Решение второй задачи носит центральный характер, так как две последующие задачи осуществляют построения относительно плоской части с минимальным числом удаленных ребер (хотя это может быть и не определяющим фактом). Математической моделью для построения топологического рисунка плоского графа является теория вращения вершин графа, созданная Г.Рингелем и Дж.Янгсом [32]. Теория вращений вершин характеризует относительное расположение вершин графа на плоскости и одновременно индуцирует (порождает) простые циклы графа. Поэтому данная модель может быть положена в основу задачи построения топологического рисунка плоской части не планарного графа [14-22].

В работе [4] показано, что задача выделения максимально плоской части не планарного графа является NP-полной. То есть ее можно решить путем полного перебора вариантов и не существует полиномиального алгоритма для ее разрешимости.

В англоязычной литературе встречаются два определения плоских частей не планарного графа. Первое определение: *наибольший планарный подграф не планарного графа* (maximum planar subgraph) – это планарный граф с наибольшим количеством ребер среди всех подграфов графа G. Второе определение: *максимальный планарный подграф не планарного графа* G = (V,E) (maximal planar subgraph) – это планарный подграф P = (V,E\F) графа G, такой, что добавление любого ребра из F к P нарушает его планарность, то есть, $P \cup e$ не планарен для каждого $e \in F$.

Обычно рассматривается следующий способ построения максимально плоского подграфа [6-8,36]. Удаляют по отдельности каждое ребро из графа и проверяют после каждого удаления оставшуюся часть графа на планарность. Если планарная часть не найдена, то удаляют всевозможное количество пар ребер и оставшуюся часть графа снова проверяют на планарность. Если решение не найдено, то продолжают удалять по три ребра всевозможными способами. Если снова решение не найдено, то продолжают процесс исключения ребер до получения решения, каждый раз увеличивая количество удаляемых ребер. Такой метод решения будем называть ***методом последовательного исключения*** ***ребер***.

Данный подход предполагает применение на каждом шаге удаления ребер проверку графа на планарность. Алгоритм определения планарности давно известен - это полиномиальный алгоритм Хопкрофта-Тарьяна [38]. Алгоритм Хопкрофта-Тарьяна основан на выделение DFS-дерева графа, и по этому пути идут многие исследователи задачи выделения плоской части непланарного графа [6-8,36,38].

Но для того чтобы сказать, что путем удаления ребер выделена плоская часть графа, нужно получить и описать планарный рисунок. Естественно предположить, что такой рисунок можно построить на основании выделенного дерева в не планарном графе [29-31]. Хорошо, если дерево можно представить в линейном виде, но как правило, такие деревья в графе трудно найти. Процесс затрудняется и неоднозначным построением рисунка самого дерева.

Для представления рисунка плоского графа, возможно, применить критерий планарности Маклейна. Однако не ко всякому графу, можно применить критерий Маклейна [23] для проверки планарности или для выделения плоской части. Очевидно, что нельзя применить данный критерий к сепарабельным графам.

**Определение 1.** *Сепарабельным графом* будем называть связный граф, имеющий точки сочленения, или мосты, или петли, или кратные ребра, или вершины с валентность равной или меньше двум, или ориентированные ребра

Таким образом, требуется, как то обозначить вид графов, к которым можно применять понятие топологического рисунка графа и критерий планарности Маклейна. С этой целью расширим понятие несепарабельного графа [36].

**Определение 2.** *Несепарабельным графом* **G** будем называть связный неориентированный граф без петель и кратных ребер, без мостов и точек сочленения, без вершин с локальной степенью меньшей или равной двум.

Таким образом, возникает иной подход**,** *основанный на цикломатических свойствах* графа для решения задачи выделения максимально плоского подграфа.

В таком методе, основную роль играют уже не ребра, а простые циклы. И тогда, задача о построении максимально плоского суграфа, может быть сведена к задаче дискретной оптимизации [2,5]. При этом ищется приближенное решение, так как задача выделения максимально плоского суграфа остается NP-полной [4].

**Определение 3.** Максимально плоскую часть графа со всем множеством вершин, с удалением минимального количества ребер из несепарабельного графа, будем называть *максимально плоским суграфом*.

Будем рассматривать приближенное решение задачи в общем виде, и выделим только основные этапы. Для решения данной задачи применим методы диакоптики [2], то есть разобьем решение на части связанные между собой.

Рассмотрим граф G = (V,E) с пронумерованным множеством ребер E = $\{e_1, e_2, ..., e_m\}$ и пронумерованным множеством вершин V = $\{v_1, v_2, ..., v_n\}$, причем card(V) = $n$ и card(E) = $m$. Выделим несепарабельную часть графа.

Таким образом, диакоптика позволяет применить математическую модель, основанную на цикломатических свойствах графа, и тем самым подключить для решения задачи теорему Маклейна. Тогда процесс решения можно представить в виде, состоящим из двух последовательных этапов:

- построение максимально плоского суграфа для несепарабельного части графа;
- преобразовать выделенную плоскую часть, индуцировав до сепарабельного плоского рисунка (пусть даже с нулевым дополнением).

Здесь основную роль играют простые циклы графа. Следовательно, задача о построении максимально плоского суграфа может быть сведена к задаче комбинаторной оптимизации, то есть к поиску оптимального значения некоторого функционала.

**Задача комбинаторной оптимизации для построения топологического рисунка максимально плоского суграфа выглядит следующим образом:** *Найти на множестве простых циклов подмножество независимых циклов, описывающее плоский суграф и удовлетворяющее нулевому значению функционала Маклейна и уравнению Эйлера с максимальным числом ребер.*

Такая постановка позволяет свести перечислительную задачу к классу задач комбинаторной оптимизации. И это позволяет применить для её решения хорошо разработанный математический аппарат дискретной оптимизации. Параллельно с этим цикломатический подход позволяет строго и однозначно описывать топологический рисунок плоской части графа, так как полученная в результате решения независимая система циклов индуцирует (порождает) вращение вершин графа. Согласно теории вращений, вращение вершин создает топологический рисунок графа [32].

Множество суграфов в несепарабельном неориентированном графе можно рассматривать как линейное пространство $£_G$ над полем $GR^2$. Операция суммы $G_1 \oplus G_2$ суграфов $G_1$ и $G_2$ определена как суграф, множество рёбер которого является симметрической разностью множеств рёбер суграфов $G_1$ и $G_2$ [33].

В качестве базисной системы векторов пространства суграфов можно выбрать суграфы с единственным ребром. Размерность пространства суграфов графа G, состоящего из $n$ вершин равна $m$. В пространстве суграфов можно выделить два подпространства, которые называются подпространством разрезов и подпространством циклов графа. Размерность подпространства разрезов равна $n-1$, а элементы этого подпространства называются квалиразрезами. Размерность подпространства циклов равна $m-n+1$, а элементы этого подпространства называются квазициклами [7,8].

В любом несепарабельном неориентированном графе можно выделить подмножество простых циклов – множество $C_{\tau}$ изометрических циклов графа [15,16].

Будем искать приближённое решение задачи построения максимального планарного суграфа в виде подмножества изометрических циклов. Так как любой простой цикл есть объединение множества изометрических циклов. Очевидно, для построения плоской части графа, достаточно выделить максимальное множество изометрических циклов. И если необходимо, по законам диакоптики добавить недостающие простые циклы [13,16].

Применение методов дискретной оптимизации позволяют процесс выделения плоской части непланарного графа разбить на три последовательных этапа.

1 этап. Выделение базиса подпространства циклов состоящего из изометрических циклов графа с минимальным значением функционаола Маклейна [15].

2 этап. Выделение из базиса подмножества изометрических циклов с нулевым значением кубического функционала Маклейна путем удаления ребер с соблюдением условия Эйлера [18].

3 этап. Дополнение подмножества изометрических циклов с нулевым значением кубического функционала Маклейна простыми циклами с учетом ребер исключенных в процессе планаризации.

Следует заметить, сто построение базиса подпространства разрезов графа, можно просто осуществить с помощью карандаша и бумаги, взяв за основу инцидентность вершин и ребер графа. Осуществить вычисление базиса подпространства циклов графа уже не как просто. Вначале требуется выделить все множество изометрических циклов несепарабельного неориентированного графа, проведя необходимые построения для каждого ребра. Затем выбрать независимую систему циклов из всего множества изометрических циклов. С этой целью нужно создать и развить алгоритмические методы для обработки информации на ЭВМ.

Перейдем к последовательному описанию этапов выделения максимально плоского суграфа из не планарного неориентированного несепарабельного графа.

## Глава 6. МЕТОД ГАУССА И ПЛОСКАЯ ЧАСТЬ ГРАФА

### 6.1. Линейная независимость множества изометрических циклов

В процессе решения задачи выделения плоской части для непланарного графа, могут появляться зависимые системы изометрических циклов. Для проверки условия независимости можно применить методы вычисления ранга матрицы циклов. Известно, что ранг матрицы будет равен числу ненулевых строк в матрице после приведения её к ступенчатой форме при помощи элементарных преобразований над строками матрицы. Обычно, для этих преобразований применяют метод Гаусса [12].

Применение метода Гаусса для определения ранга матрицы изометрических циклов имеет свои особенности.

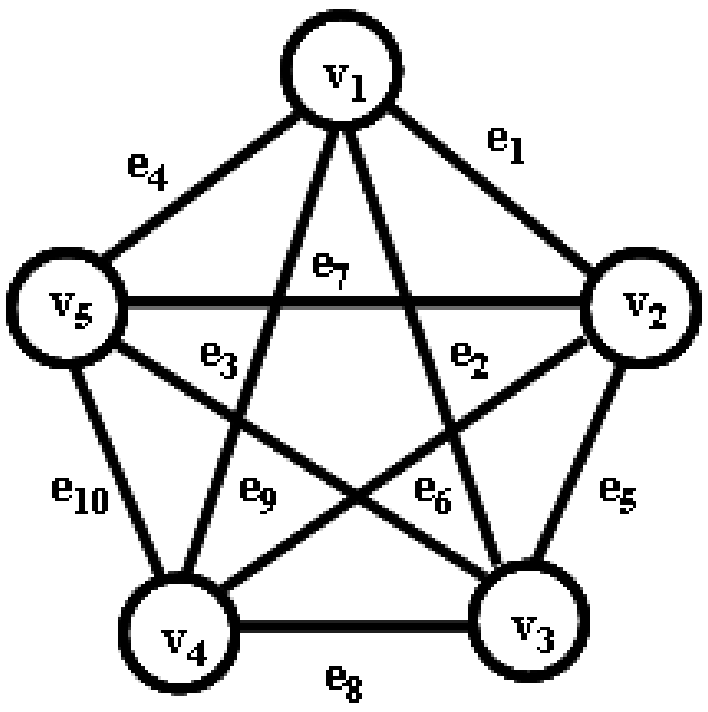

Рис. 6.1. Граф $K_5$.

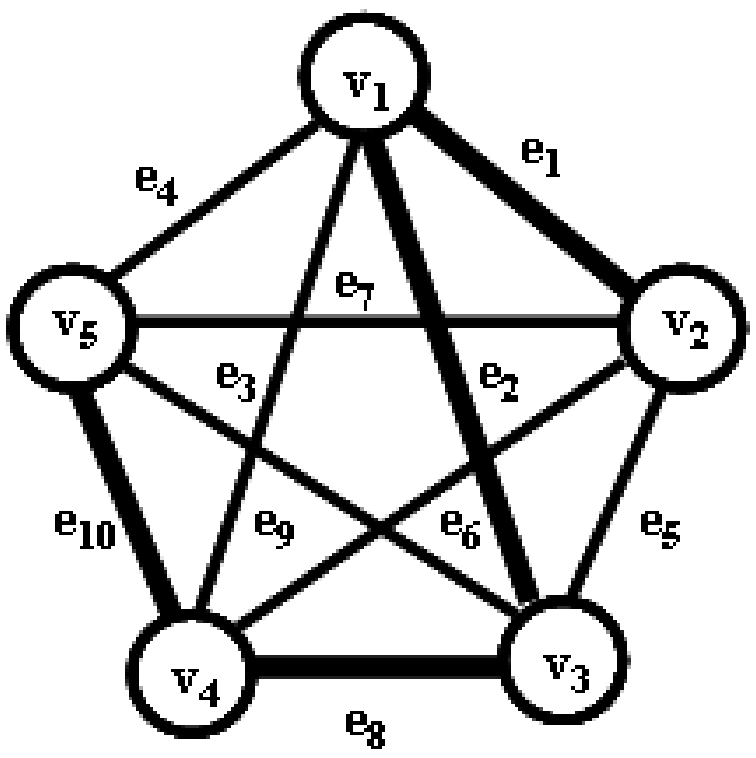

Рис. 6.2. Дерево графа $K_5$.

В результате выделения дерева в графе **G**, множество ребер разбивается на подмножество ребер дерева и подмножество хорд [7,8]. Тогда матрицу независимых циклов графа с количеством строк равным цикломатическому числу графа $\nu(G)$ = m - n + 1, можно разбить на две подматрицы - подматрицу хорд и подматрицу ветвей дерева. Если произвести линейное преобразование над строками, то возможно получить единичную подматрицу хорд. Основная идея метода Гаусса для вычисления ранга матрицы, состоит в построении единичной подматрицы хорд путем линейной комбинации строк, и получения пустых строк для зависимых циклов в случае, когда количество циклов больше цикломатического числа. Если имеются нулевые строки, то система циклов зависима. В процессе проведения линейного преобразования на месте пустых строк образуются простые циклы как результат кольцевого сложения циклов. Тогда процесс формирования нулевых строк индуцирует определение циклов необходимых для исключения из независимой системы циклов. Сказанное рассмотрим на примере полного графа $K_5$ (см. рис. 6.1).

***Пример 6.1.*** В качестве примера рассмотрим полный граф $K_5$.

Выберем дерево в графе **T** = $\{e_1,e_2,e_8,e_{10}\}$. Тогда хорды образуют множество **H** = $\{e_3,e_4,e_5,e_6,e_7,e_9\}$ (см. рис. 6.2). Множество изометрических циклов имеет вид:

$c_1 = \{e_1,e_2,e_5\}$; $c_2 = \{e_1,e_3,e_6\}$; $c_3 = \{e_1,e_4,e_7\}$; $c_4 = \{e_2,e_3,e_8\}$;

$c_5 = \{e_2,e_4,e_9\}$; $c_6 = \{e_3,e_4,e_{10}\}$; $c_7 = \{e_5,e_6,e_8\}$; $c_8 = \{e_5,e_7,e_9\}$;
$c_9 = \{e_6,e_7,e_{10}\}$; $c_{10} = \{e_8,e_9,e_{10}\}$.

Выберем следующее подмножество изометрических циклов $\{c_1,c_2,c_3,c_4,c_7,c_{10}\}$ с мощностью равной цикломатическому числу графа $\nu(\mathrm{G}) = m - n + 1$. Проверим, является ли данное подмножество циклов базисом подпространства циклов. Составим матрицу циклов и разобьем ее на две части, подматрицу ребер графа принадлежащие хордам и подматрицу ребер принадлежащих к ветвям дерева:

| | $e_3$ | $e_4$ | $e_5$ | $e_6$ | $e_7$ | $e_9$ | $e_1$ | $e_2$ | $e_8$ | $e_{10}$ |
|---|---|---|---|---|---|---|---|---|---|---|
| $c_1$ | | | 1 | | | | 1 | 1 | | |
| $c_2$ | 1 | | | 1 | | | 1 | | | |
| $c_3$ | | 1 | | | 1 | | 1 | | | |
| $c_4$ | 1 | | | | | | | 1 | 1 | |
| $c_7$ | | | 1 | 1 | | | | | 1 | |
| $c_{10}$ | | | | | | 1 | | | 1 | 1 |

Последовательно преобразуем строки матрицы, приводя матрицу к диагональному виду, используя элементарные преобразования строк. Выбираем в качестве главного элемента строки только ребра принадлежащие хордам графа:

| | $e_3$ | $e_4$ | $e_5$ | $e_6$ | $e_7$ | $e_9$ | $e_1$ | $e_2$ | $e_8$ | $e_{10}$ |
|---|---|---|---|---|---|---|---|---|---|---|
| $c_1$ | | | 1 | | | | 1 | 1 | | |
| $c_2$ | 1 | | | | | | | 1 | 1 | |
| $c_3$ | | 1 | | | 1 | | 1 | | | |
| $c_4$ | | | | 1 | | | 1 | 1 | 1 | |
| $c_7$ | | | | | | | | | | |
| $c_{10}$ | | | | | | 1 | | | 1 | 1 |

Результат окончательного преобразования для подматрицы хорд представим в ступенчатом виде:

| | $e_3$ | $e_4$ | $e_5$ | $e_6$ | $e_7$ | $e_9$ | $e_1$ | $e_2$ | $e_8$ | $e_{10}$ |
|---|---|---|---|---|---|---|---|---|---|---|
| $c_2$ | 1 | | | | | | | 1 | 1 | |
| $c_3$ | | 1 | | | 1 | | 1 | 1 | | |
| $c_1$ | | | 1 | | | | 1 | 1 | | |
| $c_4$ | | | | 1 | | | 1 | 1 | 1 | |
| $c_7$ | | | | | | | | | | |
| $c_{10}$ | | | | | | 1 | | | 1 | 1 |

Таким образом, выбранное множество изометрических циклов не является базисом подпространства циклов, так как имеются строки с нулевым количеством элементов.

В данном случае, базис может быть образован, если цикл $c_7$ заменить циклом $c_8$ или циклом $c_9$ с выбором хорды $e_7$ в качестве главного элемента строки.

Базис, состоящий из изометрических циклов $\mathrm{b}_1 = \{c_1,c_2,c_3,c_4,c_9,c_{10}\}$:

| | $e_3$ | $e_4$ | $e_5$ | $e_6$ | $e_7$ | $e_9$ | $e_1$ | $e_2$ | $e_8$ | $e_{10}$ |
|---|---|---|---|---|---|---|---|---|---|---|
| $c_2$ | 1 | | | | | | | 1 | 1 | |
| $c_3$ | | 1 | | | 1 | | 1 | 1 | | |
| $c_1$ | | | 1 | | | | 1 | 1 | | |
| $c_4$ | | | | 1 | | | 1 | 1 | 1 | |
| $c_9$ | | | | 1 | 1 | | | | | 1 |
| $c_{10}$ | | | | | | 1 | | | 1 | 1 |

Или после окончательного приведения подматрицы хорд к ступенчатому виду:

| | $e_3$ | $e_4$ | $e_5$ | $e_6$ | $e_7$ | $e_9$ | $e_1$ | $e_2$ | $e_8$ | $e_{10}$ |
|---|---|---|---|---|---|---|---|---|---|---|
| $c_2$ | 1 | | | | | | | 1 | 1 | |
| $c_3$ | | 1 | | | | | | | 1 | 1 |
| $c_1$ | | | 1 | | | | 1 | 1 | | |
| $c_4$ | | | | 1 | | | 1 | 1 | 1 | |
| $c_9$ | | | | | 1 | | 1 | 1 | 1 | 1 |
| $c_{10}$ | | | | | | 1 | | | 1 | 1 |

После преобразования, строк с нулевым содержанием элементов в матрице циклов не имеется. Следовательно, система циклов независима.

## 6.2. Изометрические циклы и плоские конфигурации

Введем следующее определение.

**Определение 6.1.** Подграф, состоящий из подмножества простых циклов, совместно с ободом, для которых значение функционала Маклейна равно нулю и кольцевая сумма элементов есть пустое множество, будем называть *плоской конфигурацией*.

Следует заметить, что плоская конфигурация может состоять полностью из изометрических циклов, а может состоять совместно как из простых циклов, так и из изометрических циклов.

***Пример 6.2.*** Рассмотрим граф, представленный на рис. 6.3.

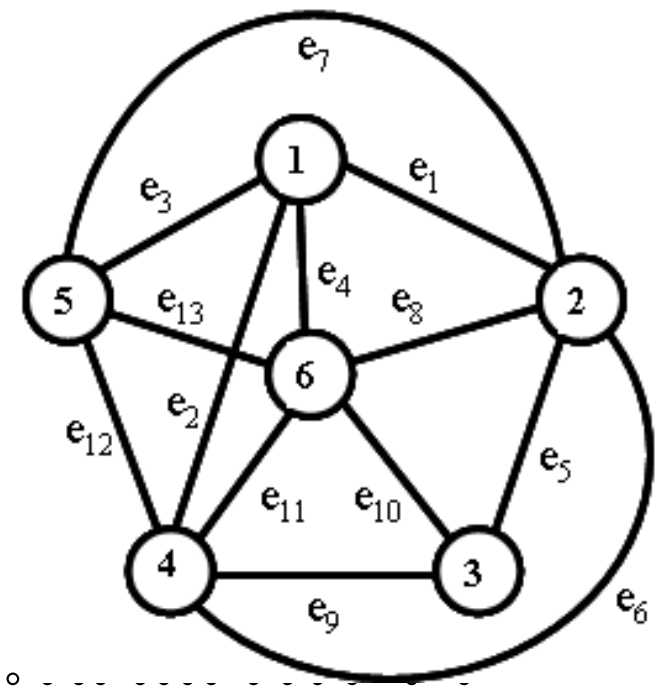

Рис. 6.3. Граф $G_1$.

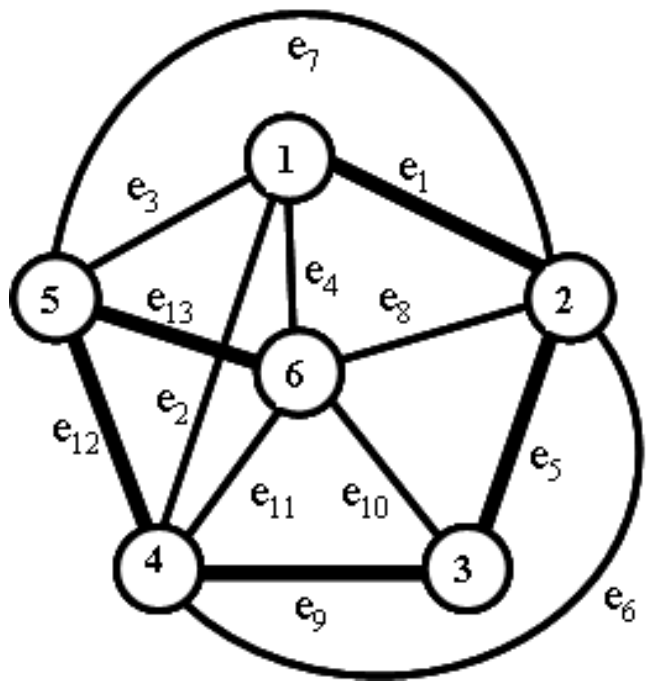

Рис. 6.4. Дерево графа $G_1$.

Множество изометрических циклов графа $G_1$ [18,40]:

| | Циклы в ребрах | Циклы в вершинах |
|---|---|---|
| цикл 1: | $\{e_1,e_2,e_6\} \rightarrow$ | $\{v_1,v_2,v_4\}$ |
| цикл 2: | $\{e_1,e_3,e_7\} \to$ | $\{v_1,v_2,v_5\}$ |
| цикл 3: | $\{e_1,e_4,e_8\} \to$ | $\{v_1,v_2,v_6\}$ |
| цикл 4: | $\{e_2,e_3,e_{11}\} \to$ | $\{v_1,v_4,v_5\}$ |
| цикл 5: | $\{e_2,e_4,e_{12}\} \to$ | $\{v_1,v_4,v_6\}$ |
| цикл 6: | $\{e_3,e_4,e_{13}\} \to$ | $\{v_1,v_5,v_6\}$ |
| цикл 7: | $\{e_5,e_6,e_9\} \to$ | $\{v_2,v_3,v_4\}$ |
| цикл 8: | $\{e_5,e_8,e_{10}\} \to$ | $\{v_2,v_3,v_6\}$ |
| цикл 9: | $\{e_6,e_7,e_{11}\} \to$ | $\{v_2,v_4,v_5\}$ |
| цикл 10: | $\{e_6,e_8,e_{12}\} \to$ | $\{v_2,v_4,v_6\}$ |
| цикл 11: | $\{e_7,e_8,e_{13}\} \to$ | $\{v_2,v_5,v_6\}$ |
| цикл 12: | $\{e_9,e_{10},e_{12}\} \to$ | $\{v_3,v_4,v_6\}$ |
| цикл 13: | $\{e_{11},e_{12},e_{13}\} \to$ | $\{v_4,v_5,v_6\}$ |

Плоские конфигурации, состоящие из изометрических циклов:

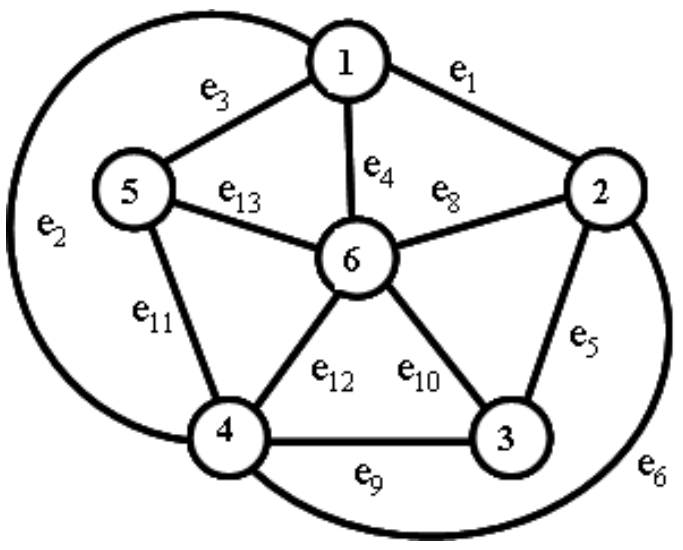

Рис. 6.5. Конфигурация $k_1$.

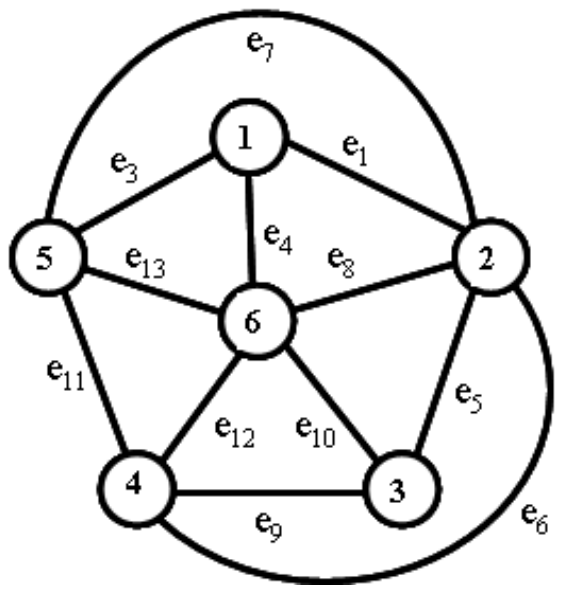

Рис. 6.5. Конфигурация $k_2$.

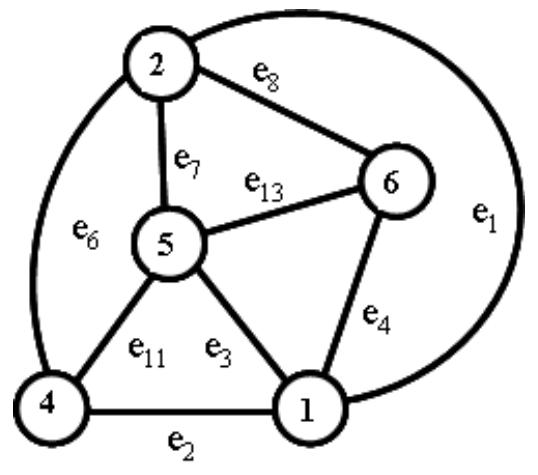

Рис. 6.7. Конфигурация $k_4$.

$k_1 = \{c_{13} \oplus c_{12} \oplus c_8 \oplus c_7 \oplus c_6 \oplus c_4 \oplus c_3 \oplus c_1\} = \varnothing;$
$k_2 = \{c_{13} \oplus c_{12} \oplus c_9 \oplus c_8 \oplus c_7 \oplus c_6 \oplus c_3 \oplus c_2\} = \varnothing;$
$k_3 = \{c_{13} \oplus c_{12} \oplus c_{11} \oplus c_8 \oplus c_7 \oplus c_4 \oplus c_2 \oplus c_1\} = \varnothing;$
$k_4 = \{c_{11} \oplus c_9 \oplus c_6 \oplus c_4 \oplus c_3 \oplus c_1\} = \varnothing;$
$k_5 = \{c_{12} \oplus c_8 \oplus c_7 \oplus c_5 \oplus c_3 \oplus c_1\} = \varnothing;$
$k_6 = \{c_{13} \oplus c_{10} \oplus c_6 \oplus c_4 \oplus c_3 \oplus c_1\} = \varnothing;$
$k_7 = \{c_{10} \oplus c_5 \oplus c_3 \oplus c_1\} = \varnothing;$
$k_8 = \{c_{13} \oplus c_6 \oplus c_5 \oplus c_4\} = \varnothing;$
$k_9 = \{c_{12} \oplus c_{10} \oplus c_8 \oplus c_7\} = \varnothing;$
$k_{10} = \{c_{11} \oplus c_6 \oplus c_3 \oplus c_2\} = \varnothing;$
$k_{11} = \{c_9 \oplus c_4 \oplus c_2 \oplus c_1\} = \varnothing.$

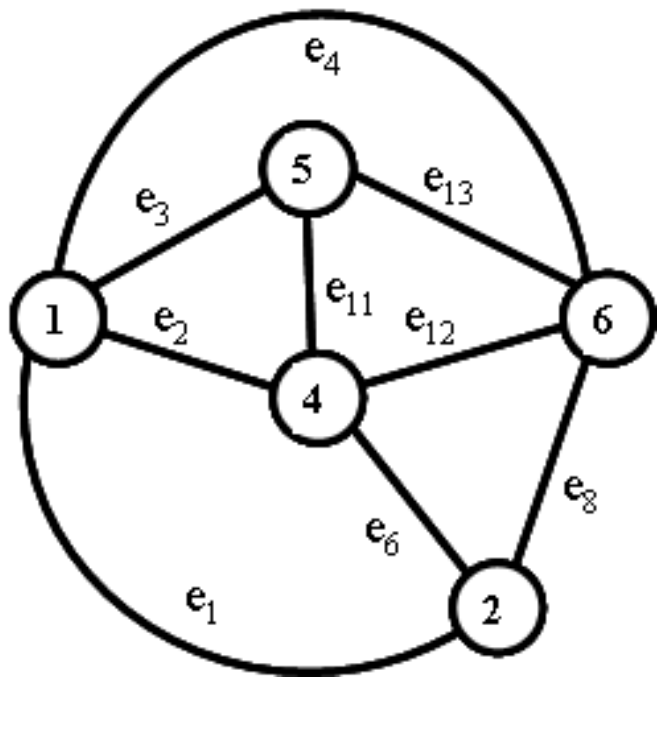

Рис. 6.8. Конфигурация $k_5$.

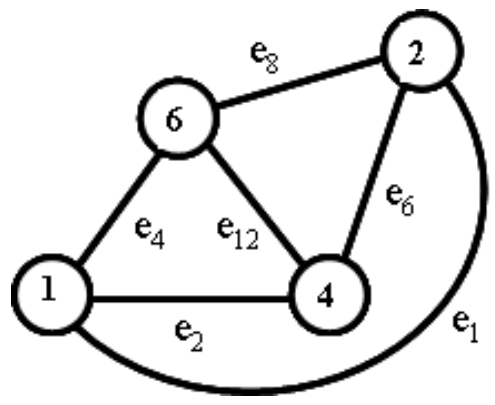

Рис. 6.9. Конфигурация $k_7$.

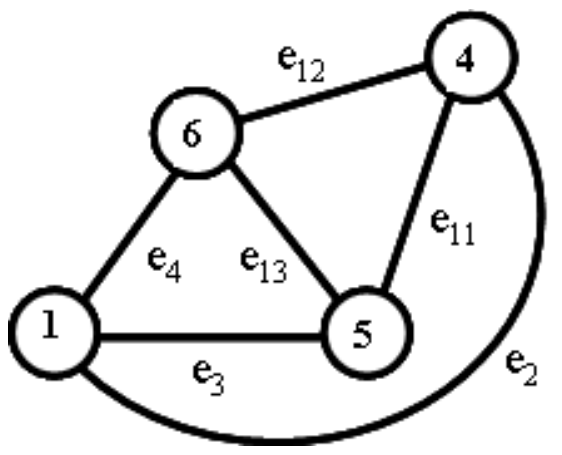

Рис. 6.10. Конфигурация $k_8$.

Для выделения плоских конфигураций модифицируем метод Гаусса для определения ранга матрицы. С этой целью выделим дерево графа, и определим хорды (см. рис. 6.4). В нашем примере, хорды - это следующие ребра графа $H = \{e_2,e_3,e_4,e_6,e_7,e_8,e_{10},e_{11}\}$. Рассмотрим модификацию алгоритма Гаусса для построения плоских конфигураций.

### 6.3. Модифицированный алгоритм Гаусса

Матрица изометрических циклов графа $G_1$ (см. рис. 6.4) имеет вид:

| | $e_1$ | **$e_2$** | **$e_3$** | **$e_4$** | $e_5$ | **$e_6$** | **$e_7$** | **$e_8$** | $e_9$ | **$e_{10}$** | **$e_{11}$** | $e_{12}$ | $e_{13}$ | номер цикла |
|---|---|---|---|---|---|---|---|---|---|---|---|---|---|---|
| $c_1$ | 1 | 1 | | | | 1 | | | | | | | | |
| $c_2$ | 1 | | 1 | | | | 1 | | | | | | | |
| $c_3$ | 1 | | | 1 | | | | 1 | | | | | | |
| $c_4$ | | 1 | 1 | | | | | | | | 1 | | | |
| $c_5$ | | 1 | | 1 | | | | | | | | 1 | | |
| $c_6$ | | | 1 | 1 | | | | | | | | | 1 | |
| $c_7$ | | | | | 1 | 1 | | | 1 | | | | | |
| $c_8$ | | | | | 1 | | | 1 | | 1 | | | | |
| $c_9$ | | | | | | 1 | 1 | | | | 1 | | | |
| $c_{10}$ | | | | | | 1 | | 1 | | | | 1 | | |
| $c_{11}$ | | | | | | | 1 | 1 | | | | | 1 | |
| $c_{12}$ | | | | | | | | | 1 | 1 | | 1 | | |
| $c_{13}$ | | | | | | | | | | | 1 | 1 | 1 | |

Выбираем в первой строке главный элемент, в данном случае это ребро $e_2$ (это обязательно должна быть хорда). Строку, в которой выбирается главный элемент, будем называть основной строкой. Находим в матрице циклы содержащие ребро $e_2$ и переставляет эти циклы в конец матрицы.

| | $e_1$ | **$e_2$** | **$e_3$** | **$e_4$** | $e_5$ | **$e_6$** | **$e_7$** | **$e_8$** | $e_9$ | **$e_{10}$** | **$e_{11}$** | $e_{12}$ | $e_{13}$ | номер цикла |
|---|---|---|---|---|---|---|---|---|---|---|---|---|---|---|
| $c_1$ | 1 | **1** | | | | 1 | | | | | | | | |
| $c_2$ | 1 | | 1 | | | | 1 | | | | | | | |
| $c_3$ | 1 | | | 1 | | | | 1 | | | | | | |
| $c_4$ | | 1 | 1 | | | | | | | | 1 | | | |
| $c_5$ | | 1 | | 1 | | | | | | | | 1 | | |
| $c_6$ | | | 1 | 1 | | | | | | | | | 1 | |
| $c_7$ | | | | | 1 | 1 | | | 1 | | | | | |
| $c_8$ | | | | | 1 | | | 1 | | 1 | | | | |
| $c_9$ | | | | | | 1 | 1 | | | | 1 | | | |
| $c_{10}$ | | | | | | 1 | | 1 | | | | 1 | | |
| $c_{11}$ | | | | | | | 1 | 1 | | | | | 1 | |
| $c_{12}$ | | | | | | | | | 1 | 1 | | 1 | | |
| $c_{13}$ | | | | | | | | | | | 1 | 1 | 1 | |

Будем называть *строки пересекающимися,* если они расположены, ниже основной строки и в них присутствует главный элемент. Объединяем пересекающиеся строки содержащие элемент $e_2$ с первой строкой и помещаем в конец матрицы циклов. Отмечаем пересекающиеся строки. В результате строим новую матрицу циклов с помеченными номером выбранного цикла для пересекающихся строк .

| | $e_1$ | **$e_2$** | **$e_3$** | **$e_4$** | $e_5$ | **$e_6$** | **$e_7$** | **$e_8$** | $e_9$ | **$e_{10}$** | **$e_{11}$** | $e_{12}$ | $e_{13}$ | номер цикла |
|---|---|---|---|---|---|---|---|---|---|---|---|---|---|---|
| $c_1$ | 1 | **1** | | | | 1 | | | | | | | | |
| $c_2$ | 1 | | 1 | | | | 1 | | | | | | | |
| $c_3$ | 1 | | | 1 | | | | 1 | | | | | | |
| $c_6$ | | | 1 | 1 | | | | | | | | | 1 | |
| $c_7$ | | | | | 1 | 1 | | | 1 | | | | | |
| $c_8$ | | | | | 1 | | | 1 | | 1 | | | | |
| $c_9$ | | | | | | 1 | 1 | | | | 1 | | | |
| $c_{10}$ | | | | | | 1 | | 1 | | | | 1 | | |
| $c_{11}$ | | | | | | | 1 | 1 | | | | | 1 | |
| $c_{12}$ | | | | | | | | | 1 | 1 | | 1 | | |
| $c_{13}$ | | | | | | | | | | | 1 | 1 | 1 | |
| $c_4$ | 1 | | 1 | | | 1 | | | | | 1 | | 1 | 1 |
| $c_5$ | 1 | | | 1 | | 1 | | | | | | 1 | | 1 |

Выбираем главный элемент во второй строке, это хорда $e_3$. Находим в нижних строках матрицы циклов, циклы содержащие ребро $e_3$ и помещаем их в конец матрицы. Объединяем пересекающие строки с выбранной строкой с учетом кольцевой операции суммирования. Помечаем выделенные строки номерами циклов.

| | $e_1$ | **$e_2$** | **$e_3$** | **$e_4$** | $e_5$ | **$e_6$** | **$e_7$** | **$e_8$** | $e_9$ | **$e_{10}$** | **$e_{11}$** | $e_{12}$ | $e_{13}$ | номер цикла |
|---|---|---|---|---|---|---|---|---|---|---|---|---|---|---|
| $c_1$ | 1 | **1** | | | | 1 | | | | | | | | |
| $c_2$ | 1 | | **1** | | | | 1 | | | | | | | |
| $c_3$ | 1 | | | 1 | | | | 1 | | | | | | |
| $c_7$ | | | | | 1 | 1 | | | 1 | | | | | |
| $c_8$ | | | | | 1 | | | 1 | | 1 | | | | |
| $c_9$ | | | | | | 1 | 1 | | | | 1 | | | |
| $c_{10}$ | | | | | | 1 | | 1 | | | | 1 | | |
| $c_{11}$ | | | | | | | 1 | 1 | | | | | 1 | |
| $c_{12}$ | | | | | | | | | 1 | 1 | | 1 | | |
| $c_{13}$ | | | | | | | | | | | 1 | 1 | 1 | |
| $c_4$ | | | | | | 1 | 1 | | | | 1 | | 1 | 1,2 |
| $c_5$ | 1 | | | 1 | | 1 | | | | | | 1 | | 1 |
| $c_6$ | 1 | | | 1 | | | 1 | | | | | | 1 | 2 |

Выбираем главный элемент в третьей строке, это хорда $e_4$. Объединяем пересекающие строки и помечаем их номерами циклов.

| | $e_1$ | **$e_2$** | **$e_3$** | **$e_4$** | $e_5$ | **$e_6$** | **$e_7$** | **$e_8$** | $e_9$ | **$e_{10}$** | **$e_{11}$** | $e_{12}$ | $e_{13}$ | номер цикла |
|---|---|---|---|---|---|---|---|---|---|---|---|---|---|---|
| $c_1$ | 1 | **1** | | | | 1 | | | | | | | | |
| $c_2$ | 1 | | **1** | | | | 1 | | | | | | | |
| $c_3$ | 1 | | | **1** | | | | 1 | | | | | | |
| $c_7$ | | | | | 1 | 1 | | | 1 | | | | | |
| $c_8$ | | | | | 1 | | | 1 | | 1 | | | | |
| $c_9$ | | | | | | 1 | 1 | | | | 1 | | | |
| $c_{10}$ | | | | | | 1 | | 1 | | | | 1 | | |
| $c_{11}$ | | | | | | | 1 | 1 | | | | | 1 | |
| $c_{12}$ | | | | | | | | | 1 | 1 | | 1 | | |
| $c_{13}$ | | | | | | | | | | | 1 | 1 | 1 | |
| $c_4$ | | | | | | 1 | 1 | | | | 1 | | 1 | 1,2 |
| $c_5$ | | | | | | 1 | | 1 | | | | 1 | | 1,3 |
| $c_6$ | | | | | | | 1 | 1 | | | | | 1 | 2,3 |

Выбираем главный элемент в четвертой строке, это хорда $c_5$. Находим пересекающисся строки, помещаем их в конец матрицы. Объединяем пересекающиеся строки.

| | $e_1$ | **$e_2$** | **$e_3$** | **$e_4$** | $e_5$ | **$e_6$** | **$e_7$** | **$e_8$** | $e_9$ | **$e_{10}$** | **$e_{11}$** | $e_{12}$ | $e_{13}$ | номер |
|---|---|---|---|---|---|---|---|---|---|---|---|---|---|---|
| $c_1$ | 1 | **1** | | | | 1 | | | | | | | | |
| $c_2$ | 1 | | **1** | | | | 1 | | | | | | | |
| $c_3$ | 1 | | | **1** | | | | 1 | | | | | | |
| $c_7$ | | | | | 1 | **1** | | | 1 | | | | | |
| $c_8$ | | | | | 1 | | | 1 | | 1 | | | | |
| $c_9$ | | | | | | 1 | 1 | | | | 1 | | | |
| $c_{10}$ | | | | | | 1 | | 1 | | | | 1 | | |
| $c_{11}$ | | | | | | | 1 | 1 | | | | | 1 | |
| $c_{12}$ | | | | | | | | | 1 | 1 | | 1 | | |
| $c_{13}$ | | | | | | | | | | | 1 | 1 | 1 | |
| $c_4$ | | | | | 1 | | 1 | | 1 | | 1 | | | 1,2,7 |
| $c_5$ | | | | | 1 | | | 1 | 1 | | | 1 | | 1,3,7 |
| $c_6$ | | | | | | | 1 | 1 | | | | | 1 | 2,3 |
| $c_9$ | | | | | 1 | | 1 | | 1 | | | 1 | | 7 |
| $c_{10}$ | | | | | 1 | | | 1 | 1 | | | 1 | | 7 |

Выбираем главный элемент в пятой строке, это хорда $e_8$. Находим пересекающиеся строки и помещаем их в конец матрицы. Объединяем с основной строкой.

| | $e_1$ | **$e_2$** | **$e_3$** | **$e_4$** | $e_5$ | **$e_6$** | **$e_7$** | **$e_8$** | $e_9$ | **$e_{10}$** | **$e_{11}$** | $e_{12}$ | $e_{13}$ | номер цикла |
|---|---|---|---|---|---|---|---|---|---|---|---|---|---|---|
| $c_1$ | 1 | **1** | | | | 1 | | | | | | | | |
| $c_2$ | 1 | | **1** | | | | 1 | | | | | | | |
| $c_3$ | 1 | | | **1** | | | | 1 | | | | | | |
| $c_7$ | | | | | 1 | **1** | | | 1 | | | | | |
| $c_8$ | | | | | 1 | | | **1** | | 1 | | | | |
| $c_{12}$ | | | | | | | | | 1 | 1 | | 1 | | |
| $c_{13}$ | | | | | | | | | | | 1 | 1 | 1 | |
| $c_4$ | | | | | 1 | | 1 | | 1 | | 1 | | | 1,2,7 |
| $c_5$ | | | | | | | | | 1 | 1 | | 1 | | 1,3,7,8 |
| $c_6$ | | | | | 1 | | 1 | | | 1 | | | 1 | 2,3,8 |
| $c_9$ | | | | | 1 | | 1 | | 1 | | 1 | | | 7 |
| $c_{10}$ | | | | | | | | | 1 | 1 | | 1 | | 7,8 |
| $c_{11}$ | | | | | 1 | | 1 | | | 1 | | | 1 | 8 |

Выбираем главный элемент в шестой строке, это хорда $e_{10}$. Находим пересекающиеся строки и помещаем их в конец матрицы. Объединяем с основной строкой.

| | $e_1$ | **$e_2$** | **$e_3$** | **$e_4$** | $e_5$ | **$e_6$** | **$e_7$** | **$e_8$** | $e_9$ | **$e_{10}$** | **$e_{11}$** | $e_{12}$ | $e_{13}$ | номер цикла |
|---|---|---|---|---|---|---|---|---|---|---|---|---|---|---|
| $c_1$ | 1 | **1** | | | | 1 | | | | | | | | |
| $c_2$ | 1 | | **1** | | | | 1 | | | | | | | |
| $c_3$ | 1 | | | **1** | | | | 1 | | | | | | |
| $c_7$ | | | | | 1 | **1** | | | 1 | | | | | |
| $c_8$ | | | | | 1 | | | **1** | | 1 | | | | |
| $c_{12}$ | | | | | | | | | 1 | **1** | | 1 | | |
| $c_{13}$ | | | | | | | | | | | 1 | 1 | 1 | |
| $c_4$ | | | | | 1 | | 1 | | 1 | | 1 | | | 1,2,7 |
| $c_5$ | | | | | | | | | | | | | | 1,3,7,8,12 |

| | | | | | | | | | | | | | |
|---|---|---|---|---|---|---|---|---|---|---|---|---|---|
| $c_6$ | | | | | 1 | | 1 | | 1 | | | 1 | 1 | 2,3,8,12 |
| $c_9$ | | | | | 1 | | 1 | | 1 | | 1 | | | 7 |
| $c_{10}$ | | | | | | | | | | | | | | 7,8,12 |
| $c_{11}$ | | | | | 1 | | 1 | | 1 | | | 1 | 1 | 8,12 |

Выбираем главный элемент в седьмой строке, это хорда $e_{11}$. Находим пересекающиеся строки и помещаем их в конец матрицы. Объединяем с основной строкой.

| | $e_1$ | **$e_2$** | **$e_3$** | **$e_4$** | $e_5$ | **$e_6$** | **$e_7$** | **$e_8$** | $e_9$ | **$e_{10}$** | **$e_{11}$** | $e_{12}$ | $e_{13}$ | номер цикла |
|---|---|---|---|---|---|---|---|---|---|---|---|---|---|---|
| $c_1$ | 1 | **1** | | | | 1 | | | | | | | | |
| $c_2$ | 1 | | **1** | | | | 1 | | | | | | | |
| $c_3$ | 1 | | | **1** | | | | 1 | | | | | | |
| $c_7$ | | | | | 1 | **1** | | | 1 | | | | | |
| $c_8$ | | | | | 1 | | | **1** | | 1 | | | | |
| $c_{12}$ | | | | | | | | | 1 | **1** | | 1 | | |
| $c_{13}$ | | | | | | | | | | | **1** | 1 | 1 | |
| $c_5$ | | | | | | | | | | | | | | 1,3,7,8,12 |
| $c_6$ | | | | | 1 | | 1 | | 1 | | | 1 | 1 | 2,3,8,12 |
| $c_9$ | | | | | 1 | | 1 | | 1 | | | 1 | 1 | 7,13 |
| $c_{10}$ | | | | | | | | | | | | | | 7,8,12 |
| $c_{11}$ | | | | | 1 | | 1 | | 1 | | | 1 | 1 | 8,12 |
| $c_4$ | | | | | 1 | | 1 | | 1 | | | 1 | 1 | 1,2,7,13 |

Выбираем главный элемент в следующей непустой строке, это хорда $e_7$. Находим пересекающиеся строки и помещаем их в конец матрицы. Объединяем с основной строкой.

| | $e_1$ | **$e_2$** | **$e_3$** | **$e_4$** | $e_5$ | **$e_6$** | **$e_7$** | **$e_8$** | $e_9$ | **$e_{10}$** | **$e_{11}$** | $e_{12}$ | $e_{13}$ | номер цикла |
|---|---|---|---|---|---|---|---|---|---|---|---|---|---|---|
| $c_1$ | 1 | **1** | | | | 1 | | | | | | | | |
| $c_2$ | 1 | | **1** | | | | 1 | | | | | | | |
| $c_3$ | 1 | | | **1** | | | | 1 | | | | | | |
| $c_7$ | | | | | 1 | **1** | | | 1 | | | | | |
| $c_8$ | | | | | 1 | | | **1** | | 1 | | | | |
| $c_{12}$ | | | | | | | | | 1 | **1** | | 1 | | |
| $c_{13}$ | | | | | | | | | | | **1** | 1 | 1 | |
| $c_6$ | | | | | 1 | | **1** | | 1 | | | 1 | 1 | 2,3,8,12 |
| $c_9$ | | | | | 1 | | 1 | | 1 | | | 1 | 1 | 7,13 |
| $c_{11}$ | | | | | 1 | | 1 | | 1 | | | 1 | 1 | 8,12 |
| $c_4$ | | | | | 1 | | 1 | | 1 | | | 1 | 1 | 1,2,7,13 |
| $c_5$ | | | | | | | | | | | | | | 1,3,7,8,12 |
| $c_{10}$ | | | | | | | | | | | | | | 7,8,12 |

Окончательно, имеем:

| | $e_1$ | **$e_2$** | **$e_3$** | **$e_4$** | $e_5$ | **$e_6$** | **$e_7$** | **$e_8$** | $e_9$ | **$e_{10}$** | **$e_{11}$** | $e_{12}$ | $e_{13}$ | номер цикла |
|---|---|---|---|---|---|---|---|---|---|---|---|---|---|---|
| $c_1$ | 1 | **1** | | | | 1 | | | | | | | | |
| $c_2$ | 1 | | **1** | | | | 1 | | | | | | | |
| $c_3$ | 1 | | | **1** | | | | 1 | | | | | | |
| $c_7$ | | | | | 1 | **1** | | | 1 | | | | | |
| $c_8$ | | | | | 1 | | | **1** | | 1 | | | | |
| $c_{12}$ | | | | | | | | | 1 | **1** | | 1 | | |
| $c_{13}$ | | | | | | | | | | | **1** | 1 | 1 | |
| $c_6$ | | | | | 1 | | **1** | | 1 | | | 1 | 1 | 2,3,8,12 |
| $c_9$ | | | | | | | | | | | | | | 7,13,6 |

| | | | | | | | | | | | | | |
|---|---|---|---|---|---|---|---|---|---|---|---|---|---|
| $c_{11}$ | | | | | | | | | | | | | 8,12,6 |
| $c_4$ | | | | | | | | | | | | | 1,2,7,13,6 |
| $c_5$ | | | | | | | | | | | | | 1,3,7,8,12 |
| $c_{10}$ | | | | | | | | | | | | | 7,8,12 |

Строим столбец взаимодействия циклов, заполняя ячейки соответствующими ранее определенными номерами циклов. Например, строка цикла $c_5$ заполняется номерами циклов 1,3,7,8,12.

| | $c_1$ | $c_2$ | $c_3$ | $c_7$ | $c_8$ | $c_{12}$ | $c_{13}$ | $c_6$ | $c_9$ | $c_{10}$ | $c_{11}$ | $c_4$ | $c_5$ |
|---|---|---|---|---|---|---|---|---|---|---|---|---|---|
| $c_1$ | | | | | | | | | | | | | |
| $c_2$ | | | | | | | | | | | | | |
| $c_3$ | | | | | | | | | | | | | |
| $c_7$ | | | | | | | | | | | | | |
| $c_8$ | | | | | | | | | | | | | |
| $c_{12}$ | | | | | | | | | | | | | |
| $c_{13}$ | | | | | | | | | | | | | |
| $c_6$ | | 2 | 3 | | 8 | 12 | | | | | | | |
| $c_9$ | | | | 7 | | | 12 | 6 | | | | | |
| $c_{10}$ | | | | 7 | 8 | 12 | | | | | | | |
| $c_{11}$ | | | | | 8 | 12 | | 6 | | | | | |
| $c_4$ | 1 | 2 | | 7 | | | 13 | 6 | | | | | |
| $c_5$ | 1 | | 3 | 7 | 8 | 12 | | | | | | | |

Просматриваем строки до диагонального элемента. Если ячейки до диагонального элемента пусты, то в диагональ записываем номер цикла. Если элементы строк до диагонального элемента не пусты, то записываем объединение номеров циклов с учетом кольцевого суммирования и номера самого цикла. Полученная запись распространяется на нижние строки, где записан данный цикл. Например, для строки с циклом $c_6$ до диагонали записаны номера циклов $c_2$,$c_3$,$c_8$,$c_{12}$ и добавляем номер цикла $c_5$. Диагональная запись имеет вид: 2,3,6,8,12. Помещаем данную запись в нижние строки матрицы вместо номера цикла $c_5$.

| | $c_1$ | $c_2$ | $c_3$ | $c_7$ | $c_8$ | $c_{12}$ | $c_{13}$ | $c_6$ | $c_9$ | $c_{10}$ | $c_{11}$ | $c_4$ | $c_5$ |
|---|---|---|---|---|---|---|---|---|---|---|---|---|---|
| $c_1$ | 1 | | | | | | | | | | | | |
| $c_2$ | | 2 | | | | | | | | | | | |
| $c_3$ | | | 3 | | | | | | | | | | |
| $c_7$ | | | | 7 | | | | | | | | | |
| $c_8$ | | | | | 8 | | | | | | | | |
| $c_{12}$ | | | | | | 12 | | | | | | | |
| $c_{13}$ | | | | | | | 13 | | | | | | |
| $c_6$ | | 2 | 3 | | 8 | 12 | | 2,3,6, 8,12 | | | | | |
| $c_9$ | | | | 7 | | | 13 | 2,3,6, 8,12 | 2,3,6,7, 8,9,12,13 | | | | |
| $c_{10}$ | | | | 7 | 8 | 12 | | | | 7,8,10, 12 | | | |
| $c_{11}$ | | | | | 8 | 12 | | 2,3,6, 8,12 | | | 2,3, 6,11 | | |
| $c_4$ | 1 | 2 | | 7 | | | 13 | 2,3,6, | | | | 1,3,4,6, | |
| | | | | | | | | 8,12 | | | | 7,8,12,13 | |
| $c_5$ | 1 | | 3 | 7 | 8 | 12 | | | | | | | 1,3,5,<br>7,8,12 |

В записи строки цикла $c_{11}$ номера циклов $c_8$ и $c_{12}$ повторяются дважды, поэтому они не включаются в диагональную запись.

Просматриваем диагонали последних строк (выделены зеленным цветом) для которых порядок больше чем цикломатическое число m-n+1 и выписываем подмножества циклов определяющих плоские конфигурации.

$k_1 = c_{13} \oplus c_{12} \oplus c_9 \oplus c_8 \oplus c_7 \oplus c_6 \oplus c_3 \oplus c_2 = \varnothing$;
$k_2 = c_{12} \oplus c_{10} \oplus c_8 \oplus c_7 = \varnothing$;
$k_3 = c_{11} \oplus c_6 \oplus c_3 \oplus c_2 = \varnothing$;
$k_4 = c_{13} \oplus c_{12} \oplus c_8 \oplus c_7 \oplus c_6 \oplus c_4 \oplus c_3 \oplus c_1 = \varnothing$;
$k_5 = c_{12} \oplus c_8 \oplus c_7 \oplus c_5 \oplus c_3 \oplus c_1 = \varnothing$.

Мы специально отправляли пересекающиеся строки в конец матрицы циклов с целью получения конфигураций с большей мощностью. Так как в общем случае, подмножество плоских конфигураций зависит от порядка расположения изометрических циклов. Например, располагая изометрические циклы в следующем порядке и в процессе расчета не переставляя строки, получим:

| | Циклы в ребрах | Циклы в вершинах |
|---|---|---|
| цикл 1: | $\{e_1,e_2,e_6\}$ | $\{v_1,v_2,v_4\}$ |
| цикл 2: | $\{e_1,e_3,e_7\}$ | $\{v_1,v_2,v_5\}$ |
| цикл 3: | $\{e_1,e_4,e_8\}$ | $\{v_1,v_2,v_6\}$ |
| цикл 4: | $\{e_2,e_3,e_{11}\}$ | $\{v_1,v_4,v_5\}$ |
| цикл 5: | $\{e_2,e_4,e_{12}\}$ | $\{v_1,v_4,v_6\}$ |
| цикл 6: | $\{e_3,e_4,e_{13}\}$ | $\{v_1,v_5,v_6\}$ |
| цикл 7: | $\{e_5,e_6,e_9\}$ | $\{v_2,v_3,v_4\}$ |
| цикл 8: | $\{e_5,e_8,e_{10}\}$ | $\{v_2,v_3,v_6\}$ |
| цикл 9: | $\{e_6,e_7,e_{11}\}$ | $\{v_2,v_4,v_5\}$ |
| цикл 10: | $\{e_6,e_8,e_{12}\}$ | $\{v_2,v_4,v_6\}$ |
| цикл 11: | $\{e_7,e_8,e_{13}\}$ | $\{v_2,v_5,v_6\}$ |
| цикл 12: | $\{e_9,e_{10},e_{12}\}$ | $\{v_3,v_4,v_6\}$ |
| цикл 13: | $\{e_{11},e_{12},e_{13}\}$ | $\{v_4,v_5,v_6\}$ |

Строим столбец взаимодействия циклов, заполняя ячейки соответствующими ранее определенными номерами циклов.

| | $c_1$ | $c_2$ | $c_3$ | $c_4$ | $c_5$ | $c_6$ | $c_7$ | $c_8$ | $c_9$ | $c_{10}$ | $c_{11}$ | $c_{12}$ | $c_{13}$ |
|---|---|---|---|---|---|---|---|---|---|---|---|---|---|
| $c_1$ | | | | | | | | | | | | | |
| $c_2$ | | | | | | | | | | | | | |
| $c_3$ | | | | | | | | | | | | | |
| $c_4$ | 1 | 2 | | | | | | | | | | | |
| $c_5$ | 1 | | 3 | | | | | | | | | | |
| $c_6$ | | 2 | 3 | 4 | 5 | | | | | | | | |
| $c_7$ | | | | | | | | | | | | | |
| $c_8$ | | | | 4 | 5 | 6 | 7 | | | | | | |
| $c_9$ | | | | 4 | | | | | | | | | |
| $c_{10}$ | | | | | 5 | | | | | | | | |
| $c_{11}$ | | | | 4 | 5 | 6 | | | | | | | |
| $c_{12}$ | | | | | | 6 | | 8 | | | | | |
| $c_{13}$ | | | | | | 6 | | | | | | | |

Вычисляем диагональные элементы.

| | $c_1$ | $c_2$ | $c_3$ | $c_4$ | $c_5$ | $c_6$ | $c_7$ | $c_8$ | $c_9$ | $c_{10}$ | $c_{11}$ | $c_{12}$ | $c_{13}$ |
|---|---|---|---|---|---|---|---|---|---|---|---|---|---|
| $c_1$ | 1 | | | | | | | | | | | | |
| $c_2$ | | 2 | | | | | | | | | | | |
| $c_3$ | | | 3 | | | | | | | | | | |
| $c_4$ | 1 | 2 | | 1,2,4 | | | | | | | | | |
| $c_5$ | 1 | | 3 | | 1,3,5 | | | | | | | | |
| $c_6$ | | 2 | 3 | 1,2,4 | 1,3,5 | 4,5,6 | | | | | | | |
| $c_7$ | | | | | | | 7 | | | | | | |
| $c_8$ | | | | 1,2,4 | 2,3,4,5 | 4,5,6 | 7 | 1,3,4, 6,7,8 | | | | | |
| $c_9$ | | | | 1,2,4 | | | | | 1,2,4,9 | | | | |
| $c_{10}$ | | | | | 1,3,5 | | | | | 1,3,5, 10 | | | |
| $c_{11}$ | | | | 1,2,4 | 1,3,5 | 4,5,6 | | | | | 2,3,6, 11 | | |
| $c_{12}$ | | | | | | 4,5,6 | | 1,3,4, 6,7,8 | | | | 1,3,5,7, 8,12 | |
| $c_{13}$ | | | | | | 4,5,6 | | | | | | | 4,5,6,13 |

Выделена следующая система зависимых изометрических циклов:

$$c_9 \oplus c_4 \oplus c_2 \oplus c_1 = \varnothing;$$
$$c_{10} \oplus c_5 \oplus c_3 \oplus c_1 = \varnothing;$$
$$c_{11} \oplus c_6 \oplus c_3 \oplus c_2 = \varnothing;$$
$$c_{13} \oplus c_6 \oplus c_5 \oplus c_4 = \varnothing;$$
$$c_{12} \oplus c_8 \oplus c_7 \oplus c_5 \oplus c_3 \oplus c_1 = \varnothing.$$

Попутно образуются простые циклы:

$$c_1 \oplus c_2 \oplus c_4 = \{e_6,e_7,e_{11}\};$$
$$c_1 \oplus c_3 \oplus c_5 = \{e_6,e_8,e_{12}\};$$
$$c_4 \oplus c_5 \oplus c_6 = \{e_{11},e_{12},e_{13}\};$$
$$c_1 \oplus c_3 \oplus c_4 \oplus c_6 \oplus c_7 \oplus c_8 = \{e_9,e_{10},e_{11},e_{12}\}.$$

Расставляем изометрические циклы в другой последовательности:

Множество хорд: $\{e_2,e_3,e_4,e_6,e_7,e_8,e_{10},e_{13}\}$.

| Новый номер | Старый номер | | |
|---|---|---|---|
| цикл 1: | цикл 2: | $\{e_1,e_3,e_7\}$ | $\{v_1,v_2,v_5\}$ |
| цикл 2: | цикл 3: | $\{e_1,e_4,e_8\}$ | $\{v_1,v_2,v_6\}$ |
| цикл 3: | цикл 6: | $\{e_3,e_4,e_{13}\}$ | $\{v_1,v_5,v_6\}$ |
| цикл 4: | цикл 7: | $\{e_5,e_6,e_9\}$ | $\{v_2,v_3,v_4\}$ |
| цикл 5: | цикл 8: | $\{e_5,e_8,e_{10}\}$ | $\{v_2,v_3,v_6\}$ |
| цикл 6: | цикл 12: | $\{e_9,e_{10},e_{12}\}$ | $\{v_3,v_4,v_6\}$ |
| цикл 7: | цикл 13: | $\{e_{11},e_{12},e_{13}\}$ | $\{v_4,v_5,v_6\}$ |
| цикл 8: | цикл 1: | $\{e_1,e_2,e_6\}$ | $\{v_1,v_2,v_4\}$ |
| цикл 9: | цикл 4: | $\{e_2,e_3,e_{11}\}$ | $\{v_1,v_4,v_5\}$ |
| цикл 10: | цикл 5: | $\{e_2,e_4,e_{12}\}$ | $\{v_1,v_4,v_6\}$ |
| цикл 11: | цикл 9: | $\{e_6,e_7,e_{11}\}$ | $\{v_2,v_4,v_5\}$ |
| цикл 12: | цикл 10: | $\{e_6,e_8,e_{12}\}$ | $\{v_2,v_4,v_6\}$ |
| цикл 13: | цикл 11: | $\{e_7,e_8,e_{13}\}$ | $\{v_2,v_5,v_6\}$ |

Строим столбец взаимодействия циклов, заполняя ячейки соответствующими ранее определенными номерами циклов.

| | $c_1$ | $c_2$ | $c_3$ | $c_4$ | $c_5$ | $c_6$ | $c_7$ | $c_8$ | $c_9$ | $c_{10}$ | $c_{11}$ | $c_{12}$ | $c_{13}$ |
|---|---|---|---|---|---|---|---|---|---|---|---|---|---|
| $c_1$ | | | | | | | | | | | | | |
| $c_2$ | | | | | | | | | | | | | |
| $c_3$ | 1 | 2 | | | | | | | | | | | |
| $c_4$ | | | | | | | | | | | | | |
| $c_5$ | | | | | | | | | | | | | |
| $c_6$ | | | | | | | | | | | | | |
| $c_7$ | | | | | | | | | | | | | |
| $c_8$ | | | | 4 | | | | | | | | | |
| $c_9$ | 1 | | 3 | | 5 | 6 | 7 | 8 | | | | | |
| $c_{10}$ | | 2 | | | 5 | 6 | | 8 | | | | | |
| $c_{11}$ | | | 3 | 4 | 5 | 6 | 7 | | | | | | |
| $c_{12}$ | | | | 4 | 5 | 6 | | | | | | | |
| $c_{13}$ | | | 3 | | | | | | | | | | |

Вычисляем диагональные элементы.

| | $c_1$ | $c_2$ | $c_3$ | $c_4$ | $c_5$ | $c_6$ | $c_7$ | $c_8$ | $c_9$ | $c_{10}$ | $c_{11}$ | $c_{12}$ | $c_{13}$ |
|---|---|---|---|---|---|---|---|---|---|---|---|---|---|
| $c_1$ | 1 | | | | | | | | | | | | |
| $c_2$ | | 2 | | | | | | | | | | | |
| $c_3$ | 1 | 2 | 1,2,3 | | | | | | | | | | |
| $c_4$ | | | | 4 | | | | | | | | | |
| $c_5$ | | | | | 5 | | | | | | | | |
| $c_6$ | | | | | | 6 | | | | | | | |
| $c_7$ | | | | | | | 7 | | | | | | |
| $c_8$ | | | | 4 | | | | 4,8 | | | | | |
| $c_9$ | 1 | | 1,2,3 | | 5 | 6 | 7 | 4,8 | 2,3,4,5,6,7,8,9 | | | | |
| $c_{10}$ | | 2 | | | 5 | 6 | | 4,8 | | 2,4,5,6,8,10 | | | |
| $c_{11}$ | | | 1,2,3 | 4 | 5 | 6 | 7 | | | | 1,2,3,4,5,6,7,11 | | |
| $c_{12}$ | | | | 4 | 5 | 6 | | | | | | 4,5,6,12 | |
| $c_{13}$ | | | 1,2,3 | | | | | | | | | | 1,2,3,13 |

Множество плоских конфигураций относительно первоначального расположения циклов имеет вид:

$$c_{13} \oplus c_{12} \oplus c_8 \oplus c_7 \oplus c_6 \oplus c_4 \oplus c_3 \oplus c_1 = \varnothing;$$
$$c_{12} \oplus c_8 \oplus c_7 \oplus c_5 \oplus c_3 \oplus c_1 = \varnothing;$$
$$c_{13} \oplus c_{12} \oplus c_9 \oplus c_8 \oplus c_7 \oplus c_6 \oplus c_3 \oplus c_2 = \varnothing;$$
$$c_{12} \oplus c_{10} \oplus c_8 \oplus c_7 = \varnothing;$$
$$c_{11} \oplus c_6 \oplus c_3 \oplus c_2 = \varnothing.$$

Попутно образуются простые циклы:

$c_2 \oplus c_3 \oplus c_6 = \{e_7, e_8, e_{13}\}$;

$c_1 \oplus c_7 = \{e_1, e_2, e_5, e_9\}$.

## 6.4. Построение плоских конфигураций с заданным ободом

Выделим в графа простой цикл и будем считать его ободом для совокупности изометрических циклов. Поместим его в качестве последнего элемента в матрице изометрических циклов. Например, для нашего графа G выбираем цикл, состоящий из следующих ребер $\{e_1,e_2,e_5,c_{10},c_{11},c_{13}\}$. Применим метод Гаусса и определим состав изометрических циклов для заданного обода. Матрица циклов имеет вид:

цикл 1: $\{e_1,e_2,e_6\}$
цикл 2: $\{e_1,e_3,e_7\}$
цикл 3: $\{e_1,e_4,e_8\}$
цикл 4: $\{e_2,e_3,e_{11}\}$
цикл 5: $\{e_2,e_4,e_{12}\}$
цикл 6: $\{e_3,e_4,e_{13}\}$
цикл 7: $\{e_5,e_6,e_9\}$
цикл 8: $\{e_5,e_8,e_{10}\}$
цикл 9: $\{e_6,e_7,e_{11}\}$
цикл 10: $\{e_6,e_8,e_{12}\}$
цикл 11: $\{e_7,e_8,e_{13}\}$
цикл 12: $\{e_9,e_{10},e_{12}\}$
цикл 13: $\{e_{11},e_{12},e_{13}\}$
цикл 14: $\{e_1,e_2,e_5,c_{10},c_{11},c_{13}\}$

Строим столбец взаимодействия циклов, заполняя ячейки соответствующими ранее определенными номерами циклов.

| | $c_1$ | $c_2$ | $c_3$ | $c_4$ | $c_5$ | $c_6$ | $c_7$ | $c_8$ | $c_9$ | $c_{10}$ | $c_{11}$ | $c_{12}$ | $c_{13}$ | $c_{14}$ |
|---|---|---|---|---|---|---|---|---|---|---|---|---|---|---|
| $c_1$ | | | | | | | | | | | | | | |
| $c_2$ | | | | | | | | | | | | | | |
| $c_3$ | | | | | | | | | | | | | | |
| $c_4$ | 1 | | | | | | | | | | | | | |
| $c_5$ | 1 | | 3 | 4 | | | | | | | | | | |
| $c_6$ | | | 3 | | 5 | | | | | | | | | |
| $c_7$ | | | | 4 | | | | | | | | | | |
| $c_8$ | | | | | 5 | | 7 | | | | | | | |
| $c_9$ | | 2 | | 4 | | | | | | | | | | |
| $c_{10}$ | | | | 4 | 5 | | | | | | | | | |
| $c_{11}$ | | 2 | | | 5 | 6 | | | | | | | | |
| $c_{12}$ | | | | | | | | 8 | | | | | | |
| $c_{13}$ | | | | | | 6 | 7 | | | | | | | |
| $c_{14}$ | 1 | | | 4 | | 6 | 7 | 8 | | | | | | |

Вычисляем диагональные элементы.

| | $c_1$ | $c_2$ | $c_3$ | $c_4$ | $c_5$ | $c_6$ | $c_7$ | $c_8$ | $c_9$ | $c_{10}$ | $c_{11}$ | $c_{12}$ | $c_{13}$ | $c_{14}$ |
|---|---|---|---|---|---|---|---|---|---|---|---|---|---|---|
| $c_1$ | 1 | | | | | | | | | | | | | |
| $c_2$ | | 2 | | | | | | | | | | | | |
| $c_3$ | | | 3 | | | | | | | | | | | |
| $c_4$ | 1 | | | 1,4 | | | | | | | | | | |
| $c_5$ | 1 | | 3 | 1,4 | 3,4,5 | | | | | | | | | |
| $c_6$ | | | 3 | | 3,4,5 | 4,5,6 | | | | | | | | |
| $c_7$ | | | | 1,4 | | | 1,4,7 | | | | | | | |
| $c_8$ | | | | | 3,4,5 | | 1,4,7 | 1,3,5, 7,8 | | | | | | |
| $c_9$ | | 2 | | 1,4 | | | | | 1,2, 4,9 | | | | | |
| $c_{10}$ | | | | 1,4 | 3,4,5 | | | | | 1,3, 5,10 | | | | |
| $c_{11}$ | | 2 | | | 3,4,5 | 4,5,6 | | | | | 2,3, 6,11 | | | |
| $c_{12}$ | | | | | | | | 1,3,5, 7,8 | | | | 1,3,5, 7,8,12 | | |
| $c_{13}$ | | | | | | 4,5,6 | 1,4,7 | | | | | | 5,6, 7,13 | |
| $c_{14}$ | 1 | | | 1,4 | | 4,5,6 | 1,4,7 | 1,3,5, 7,8 | | | | | | 3,4,6, 8,14 |

Откуда получаем (см. рис. 6.11):

$c_{14} \oplus c_8 \oplus c_6 \oplus c_4 \oplus c_3 = \varnothing$.

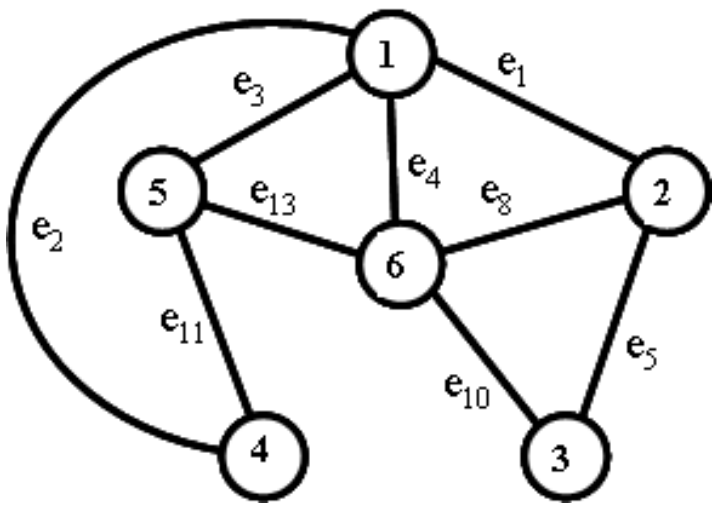

Рис. 6.11. Плоская конфигурация с заданным ободом.

Естественно, результат вычислений зависит от порядка расположения элементов в матрице циклов. Следует заметить, что одновременно можно располагать не один цикл, а несколько простых циклов, считая каждый из них ободом графа, располагая их в конце матрицы циклов.

### 6.5. Дополнение плоских конфигураций

Плоские конфигурации можно дополнять ребрами и вершинами.

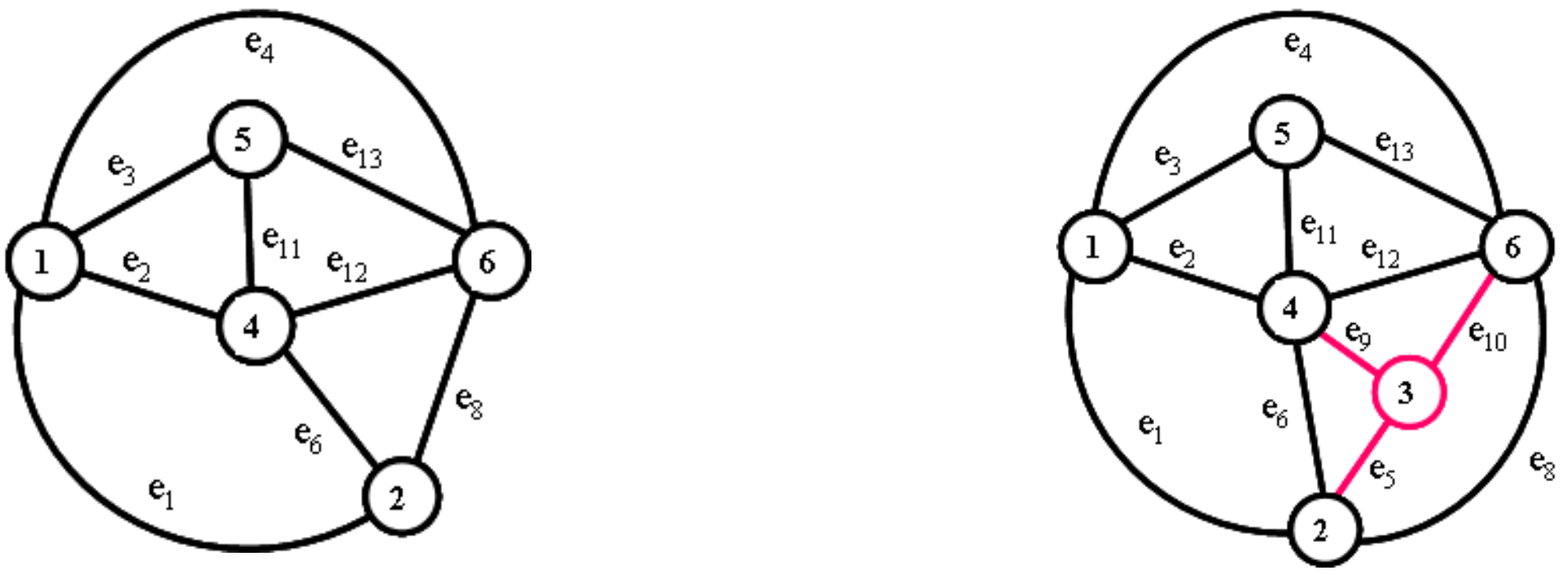

Рис. 6.12. Добавление вершины $v_3$ с ребрами.

Например, плоскую конфигурацию $\{c_{13},c_{10},c_6,c_4,c_3,c_1\}$ можно дополнить включив вершину $v_3$ вместе с ребрами $e_5,e_9,e_{10}$ и построить рисунок для максимально плоского суграфа.

Рисунок максимально плоского суграфа можно достроить, подключив к ободу ребра.

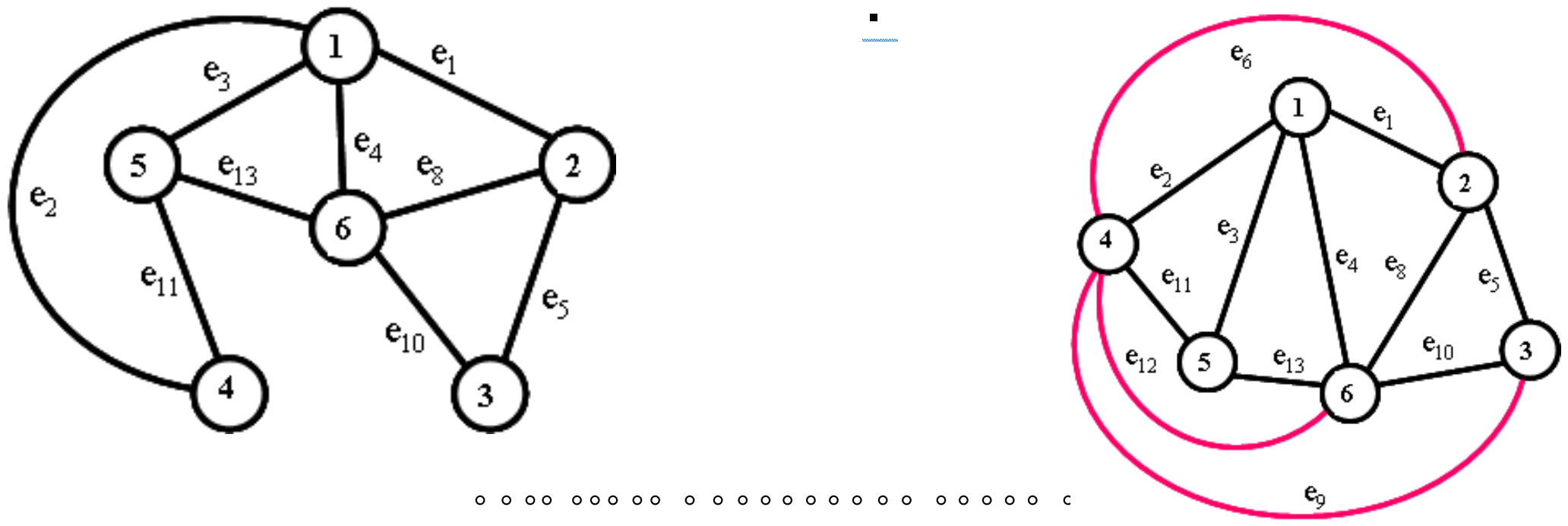

Рис. 6.13. Добавление ребер к ободу.

Для некоторых графов, построения максимально плоского суграфа может быть определено нахождением максимальной плоской конфигурации.

Например, рассмотрим граф представленного на рис. 6.14. Множество изометрических циклов графа имеет вид:

| | Множество изометрических циклов графа в виде ребер: | Множество изометрических циклов графа в виде вершин: |
|---|---|---|
| множество 1 | $\{e_1,e_2,e_{10},e_{12}\}$; | $\{v_1,v_3,v_7,v_9\}$; |
| множество 2 | $\{e_1,e_3,e_{10},e_{13}\}$; | $\{v_1,v_3,v_7,v_{10}\}$; |
| множество 3 | $\{e_1,e_3,e_{15},e_{16}\}$; | $\{v_1,v_4,v_7,v_{10}\}$; |
| множество 4 | $\{e_1,e_4,e_{15},e_{17}\}$; | $\{v_1,v_4,v_7,v_{11}\}$; |
| множество 5 | $\{e_2,e_3,e_6,e_7\}$; | $\{v_1,v_2,v_9,v_{10}\}$; |
| множество 6 | $\{e_2,e_3,e_{12},e_{13}\}$; | $\{v_1,v_3,v_9,v_{10}\}$; |
| множество 7 | $\{e_2,e_4,e_{21},e_{22}\}$; | $\{v_1,v_5,v_9,v_{11}\}$; |
| множество 8 | $\{e_3,e_4,e_{16},e_{17}\}$; | $\{v_1,v_4,v_{10},v_{11}\}$; |
| множество 9 | $\{e_5,e_6,e_{11},e_{12}\}$; | $\{v_2,v_3,v_8,v_9\}$; |
| множество 10 | $\{e_5,e_6,e_{20},e_{21}\}$; | $\{v_2,v_5,v_8,v_9\}$; |
| множество 11 | $\{e_5,e_7,e_{11},e_{13}\}$; | $\{v_2,v_3,v_8,v_{10}\}$; |
| множество 12 | $\{e_5,e_8,e_{20},e_{23}\}$; | $\{v_2,v_5,v_8,v_{12}\}$; |
| множество 13 | $\{e_6,e_7,e_{12},e_{13}\}$; | $\{v_2,v_3,v_9,v_{10}\}$; |
| множество 14 | $\{e_6,e_8,e_{21},e_{23}\}$; | $\{v_2,v_5,v_9,v_{12}\}$; |
| множество 15 | $\{e_7,e_8,e_{16},e_{18}\}$; | $\{v_2,v_4,v_{10},v_{12}\}$; |
| множество 16 | $\{e_9,e_{10},e_{14},e_{15}\}$; | $\{v_3,v_4,v_6,v_7\}$; |
| множество 17 | $\{e_9,e_{11},e_{19},e_{20}\}$; | $\{v_3,v_5,v_6,v_8\}$; |
| множество 18 | $\{e_9,e_{12},e_{19},e_{21}\}$; | $\{v_3,v_5,v_6,v_9\}$; |
| множество 19 | $\{e_9,e_{13},e_{14},e_{16}\}$; | $\{v_3,v_4,v_6,v_{10}\}$; |
| множество 20 | $\{e_{10},e_{13},e_{15},e_{16}\}$; | $\{v_3,v_4,v_7,v_{10}\}$; |
| множество 21 | $\{e_{11},e_{12},e_{20},e_{21}\}$; | $\{v_3,v_5,v_8,v_9\}$; |
| множество 22 | $\{e_{14},e_{17},e_{19},e_{22}\}$; | $\{v_4,v_5,v_6,v_{11}\}$; |
| множество 23 | $\{e_{14},e_{18},e_{19},e_{23}\}$; | $\{v_4,v_5,v_6,v_{12}\}$; |
| множество 24 | $\{e_{17},e_{18},e_{22},e_{23}\}$; | $\{v_4,v_5,v_{11},v_{12}\}$; |

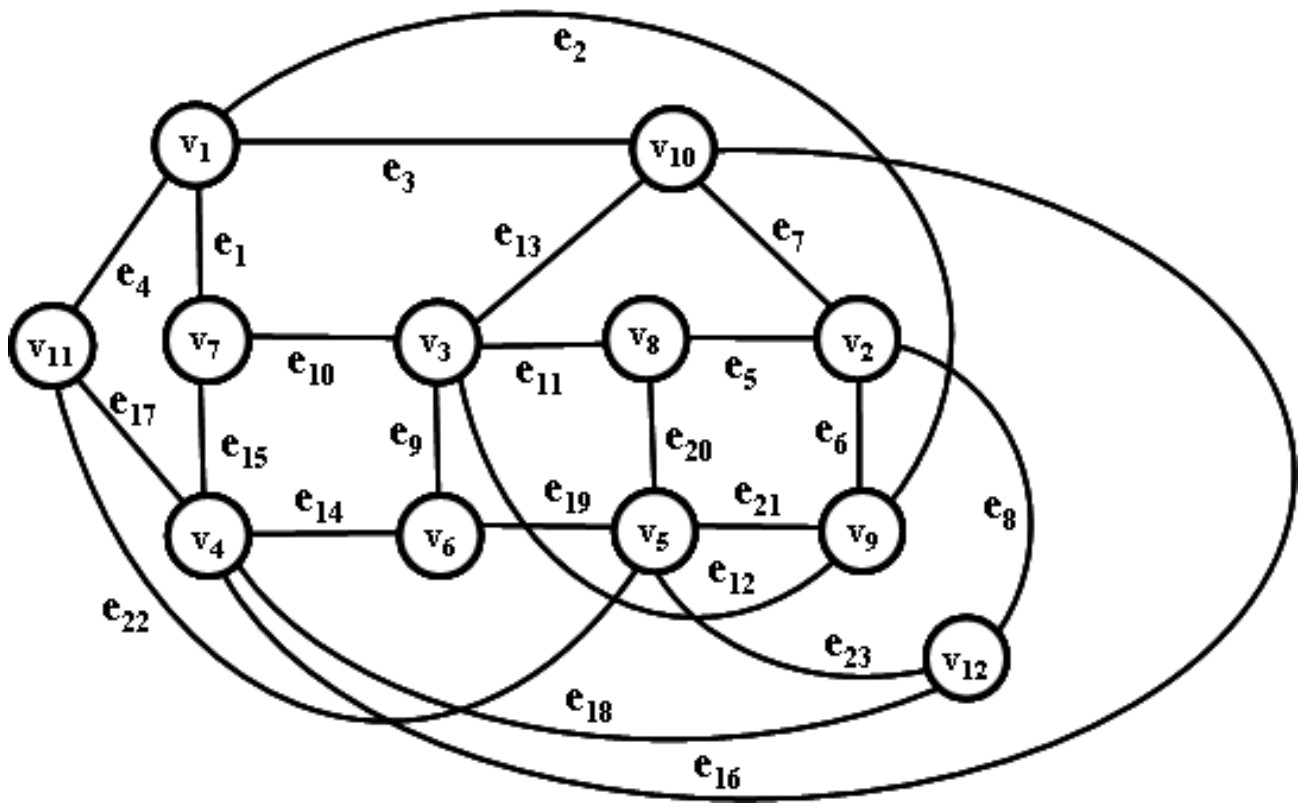

Рис. 6.14. . Граф G.

Здесь плоская конфигурация сразу порождает рисунок максимально плоского суграфа с удаленными ребрами $\{e_2, e_{12}, e_{22}\}$.

$$\{c_2 \oplus c_4 \oplus c_8 \oplus c_{10} \oplus c_{11} \oplus c_{14} \oplus c_{15} \oplus c_{16} \oplus c_{17} \oplus c_{23}\} = \varnothing.$$

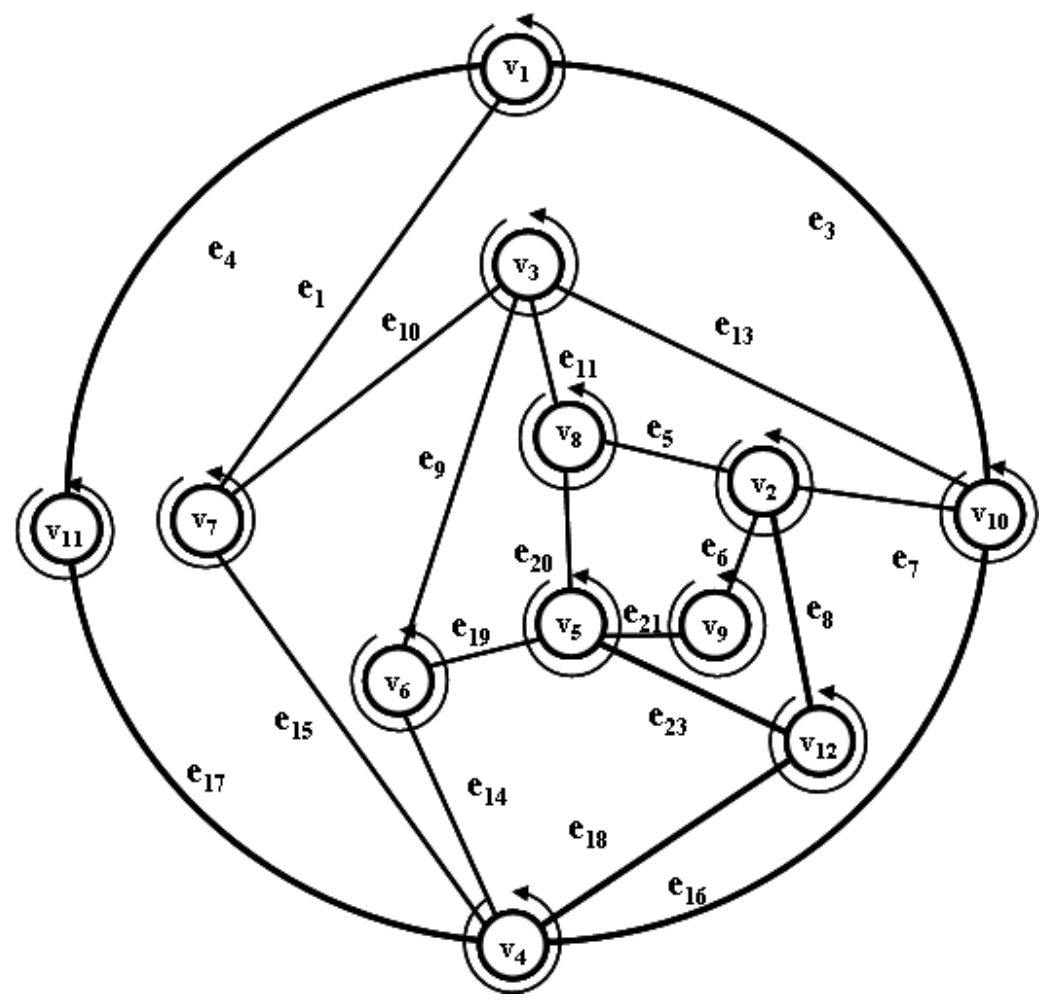

Рис. 6.15. Максимально плоский суграф не планарного графа.

### 6.5. Фрагментарный алгоритм (подход снизу)

Введем понятие фрагментарной структуры и фрагментарного алгоритма [10].

**Определение 6.2**. Фрагментарной структурой $(Y, E)$ на конечном множестве Y называется семейство его подмножеств $E = \{E_1, E_2, \ldots, E_n\}$, такое, что $\forall E_i \in E, E \neq \varnothing$ $\exists e \in E_i : E_i \setminus \{e\} \in E$.

Элементы из множества E будем называть допустимыми фрагментами.

Таким образом, для любого допустимого фрагмента $E_i$ существует нумерация его элементов $E_i = \{e_{i1}, e_{i2}, \ldots, e_{is}\}$, такая, что $\{e_{i1}, e_{i2}, \ldots, e_{is}\} \in E$ для всех $k=1,2,\ldots,s$.

**Определение 6.3.** Одноэлементные множества, которые являются допустимыми фрагментами, будем называть элементарными фрагментами.

**Определение 6.4.** Фрагмент называется максимальным, если он не является подмножеством никакого другого фрагмента.

Фрагменты обладают следующими свойствами:

**Свойство 6.1.** Пустое множество является фрагментом; $\varnothing \in \mathrm{E}$.

**Свойство 6.2.** Пусть $\max\limits_{i=1,n} |\mathrm{E}_i| = \mathrm{M}$. Тогда для любого целого числа $m$ в интервале $0 \le m \le \mathrm{M}$ найдется элемент в множестве E, мощность которого равна $m$.

**Теорема 6.1.** Если $(\mathrm{Y}, \mathrm{E})$ – фрагментарная структура на множестве Y, то для любого непустого множества $\mathrm{A} \in \mathrm{E}$ существует нумерация его элементов $A=\{x_1,x_2,\ldots,x_m\}$ такая, что для всех $k = \overline{1,m}$ множество $\{x_1,x_2,\ldots,x_k\} \in \mathrm{E}$.

Из теоремы вытекает, что всякий допустимый фрагмент можно построить из пустого множества, последовательно добавляя к нему элементы так, чтобы на каждом шаге такой процедуры полученное подмножество было допустимым фрагментом.

Максимальный фрагмент может быть построен с помощью следующего «жадного» алгоритма [10]:

- элементы множества Y линейно упорядочиваются;

- на начальном шаге выбирается пустое множество $\mathrm{Y}_0 = \varnothing$;

- на шаге с номером $k$+1 выбирается первый по порядку элемент $y \in \mathrm{Y} \setminus \mathrm{Y}_k$, такой, что $\mathrm{Y}_k \cup \{y\} \in \mathrm{E}$. Строится множество $\mathrm{Y}_{k+1} = \mathrm{Y}_k \cup \{y\}$;

алгоритм заканчивает работу, если на очередном шаге не удалось найти элемент $y \in \mathrm{Y} \setminus \mathrm{Y}_k$ с требуемым свойством.

Приведённый алгоритм будем называть фрагментарным алгоритмом. Результат применения фрагментарного алгоритма определяется заданным линейным порядком на множестве Y. Таким образом, любой максимальный фрагмент может быть задан некоторой перестановкой элементов множества Y. Пусть $\mathrm{A} \in \mathrm{E}$. Условие для элемента $y \in \mathrm{Y}$, при котором $\mathrm{A} \cup \{y\} \in \mathrm{E}$, будем называть *условием присоединения элемента y.*

С множеством изометрических циклов в графе G(X,U) связан ряд инвариантов графа [10]. Одним из таких инвариантов является мощность множества изометрических циклов card $\mathrm{C}_\tau$. Другим инвариантом может служить вектор количества изометрических циклов, упорядоченный по возрастанию их длин. Следующим инвариантом, по аналогии с вектором локальных степеней, является вектор количества изометрических циклов, проходящих по рёбрам графа, – будем называть его вектором циклов по рёбрам:

$$\mathrm{P}_u = < a_{u_1}, a_{u_2}, a_{u_3}, \ldots, a_{u_m} >,$$

где $a_{u_i}$ – количество изометрических циклов суграфа, проходящих по ребру $u_i$.

Инвариантом является также вектор количества изометрических циклов, проходящих по вершинам графа. Будем называть его вектором циклов по вершинам:

$$\mathrm{P}_x =< a_{x_1}, a_{x_2}, a_{x_3}, ..., a_{u_n} >,$$

где $a_{x_j}$ – количество циклов, проходящих через вершину $x_j$.

Для любого подмножества изометрических циклов $\mathrm{C}_\tau$ мощностью $k$ можно определить значение функционала Маклейна:

$$\mathrm{F}(\mathrm{C}_\tau)_k = \sum_{i=1}^{m} a_{u_i}^2 - 3\sum_{i=1}^{m} a_{u_i} + 2m \qquad (6.1)$$

где $a_{u_i}$ – количество циклов, проходящих по ребру с номером $i$, а $m$ –количество рёбер графа (суграфа).

Следует заметить, что функционал Маклейна принимает нулевое значение тогда, когда выбранное подмножество изометрических циклов определяет планарный суграф в графе G:

$$\mathrm{F}(\mathrm{C}_\tau)_k = 0 \qquad (6.2)$$

Между количеством выбранных изометрических циклов $k$, количеством рёбер $m$ и количеством вершин $n$ для данного подмножества изометрических циклов существует связь (формула Эйлера):

$$k\text{-}m\text{+}n\text{-}1\text{=}0 \qquad (6.3)$$

**Определение 6.5.** Сумму элементов подмножества изометрических циклов мощностью $k$ будем называть ободом:

$$\mathrm{C}_{0k} = c_1 \oplus c_2 \oplus ... \oplus c_k$$

Будем рассматривать множество циклов графа как фрагментарную структуру, причём изометрические циклы представляют элементарные фрагменты этой структуры. Любой планарный суграф графа G можно представить в виде суммы изометрических циклов (элементарных фрагментов).

Построение максимального фрагмента для задачи построения максимальной плоской части графа опишем с помощью приведённого выше фрагментарного алгоритма с соблюдением следующих свойств:

а) элементы множества изометрических циклов $\mathrm{C}_\tau$ представляем в виде кортежа (упорядочиваем и нумеруем);

б) на начальном шаге выбираем пустое множество изометрических циклов $\mathrm{C}_p = \varnothing$;

в) на шаге с номером $k$ выбираем следующий по порядку элемент кортежа, имеющий общее ребро с внешним простым циклом (ободом), построенным как сумма циклов предыдущих шагов $c_{0k}=c_{0(k-1)} \oplus c_{0k}$, при условии, что обод не является пустым множеством;

г) если выполнены условия (6.2) и (6.3), то выбранный изометрический цикл помещается в множество $C_p$;

д) алгоритм заканчивает работу, если на очередном шаге не удалось найти элемент $c_i \in C_\tau$ с требуемым свойством.

Рассмотрим работу данного алгоритма на примере.

***Пример 6.3.*** Выделим плоскую часть графа K5, представленного на рис. 6.15.

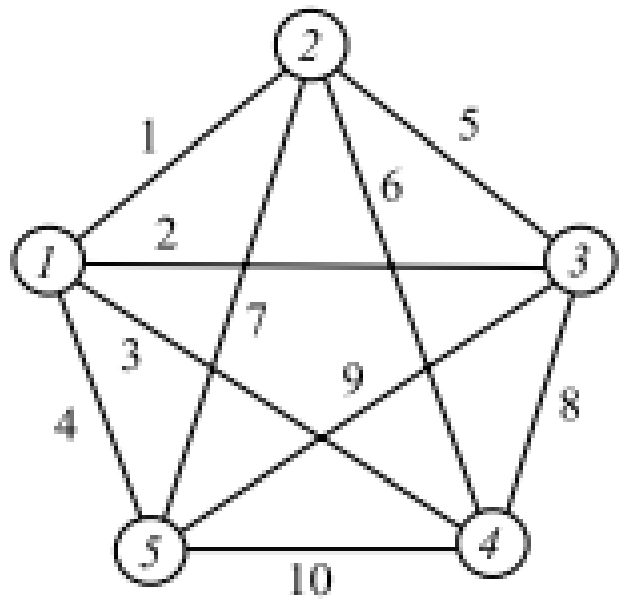

Рисунок 6.15. Граф K5.

Множество изометрических циклов графа G представляем в виде кортежа:

$$C_\tau = < c_1, c_2, c_3, c_4, c_5, c_6, c_7, c_8, c_9, c_{10} >:$$

$c_1 = \{e_1, e_2, e_5\} \to \{v_1, v_2, v_3\}$; $c_2 = \{e_1, e_3, e_6\} \to \{v_1, v_2, v_4\}$;

$c_3 = \{e_1, e_4, e_7\} \to \{v_1, v_2, v_5\}$; $c_4 = \{e_2, e_3, e_8\} \to \{v_1, v_3, v_4\}$;

$c_5 = \{e_2, e_4, e_9\} \to \{v_1, v_3, v_5\}$; $c_6 = \{e_3, e_{10}, e_7\} \to \{v_1, v_4, v_5\}$;

$c_7 = \{e_3, e_6, e_8\} \to \{v_2, v_3, v_4\}$; $c_8 = \{e_5, e_7, e_9\} \to \{v_2, v_3, v_5\}$;

$c_9 = \{e_6, e_7, e_{10}\} \to \{v_2, v_4, v_5\}$; $c_{10} = \{e_8, e_9, e_{10}\} \to \{v_3, v_4, v_5\}$.

Для упрощения записи иногда рёбра будем обозначать числами.

Выбираем в множестве изометрических циклов первый по порядку цикл и помещаем его в список циклов, характеризующий плоскую часть графа:

$$c_1 = \{e_1, e_2, e_5\} \to \{v_1, v_2, v_3\},\ C_p = \{c_1\}.$$

Строим векторы циклов по рёбрам $P_e = < a_{e_1}, a_{e_2}, a_{e_3}, ..., a_{e_m} >$ и по вершинам $P_v = < a_{v_1}, a_{v_2}, a_{v_3}, ..., a_{v_n} >$. Вычисляем значение функционала Маклейна (6.1) для выбранной системы циклов: $F(C_p) = 0$; ограничение (6.3) также выполнено: $1 - 3 + 3 - 1 = 0$.

Так как цикл единственный, его внешний простой цикл (обод) определяется самим циклом:

$$c_0 = \varnothing \oplus c_1 = \{e_1, e_2, e_5\}.$$

Выбираем во множестве изометрических циклов следующий по порядку пересекающийся (в множественном смысле) цикл $c_2 = \{e_1, e_3, e_6\} \to \{c_1, c_2, c_4\}$, $|c_0 \cap c_2| = 1$. Строим векторы $P_e$ = <2,1,1,0,1,1,0,0,0,0> и $P_v$ = <2,2,1,1,0>. Отсюда определяем значение функционала Маклейна (6.1) для выбранной системы циклов: $F(C_p)=0$; ограничение (6.3) также выполнено: $2-5+4-1=0$. Обод определяется суммированием $c_1$ и $c_2$:

$$c_0 = c_1 \oplus c_2 = \{e_1, e_2, e_5\} \oplus \{e_1, e_3, e_6\} = \{e_2, e_3, e_5, e_6\}.$$

Обод $c_0 \neq 0$. Все условия выполнены. Помещаем цикл $c_2$ во множество, характеризующее плоскую часть графа: $C_p = \{c_1, c_2\}$.

Выбираем во множестве изометрических циклов следующий пересекающийся цикл:

$$c_3 = \{e_1, e_4, e_7\} \to \{v_1, v_2, v_5\},\ |c_0 \cap c_3| = 1.$$

Но мы не можем поместить его в список циклов, характеризующий плоскую часть графа, так как значение функционала Маклейна для выбранной совокупности циклов больше нуля, согласно $P_e$ = <3,1,1,1,1,1,1,0,0,0>.

Выбираем во множестве изометрических циклов следующий пересекающийся цикл:

$$c_4 = \{e_2, e_3, e_8\} \to \{v_1, v_3, v_4\},\ |c_0 \cap c_4| = 1.$$

Строим векторы $P_e$=<2,2,2,0,1,1,0,1,0,0> и $P_v$=<3,2,2,2,0>. Отсюда определяем значение функционала Маклейна (7.1) для выбранной системы циклов: $F(C_p)=0$; ограничение (6.3) также выполняется: $3-6+3-1=0$. Обод определяется как сумма предыдущего обода $c_0$ и $c_4$:

$$c_0 = \{e_2, e_3, e_5, e_6\} \oplus \{e_2, e_3, e_8\} = \{e_5, e_6, e_8\}.$$

Все условия выполнены, помещаем цикл $c_4$ во множество, характеризующее плоскую часть графа: $C_p=\{c_1,c_2,c_4\}$.

Циклы $c_5$ и $c_6$ не имеют общих рёбер с ободом и исключаются из рассмотрения. Следующим циклом, имеющим с предыдущим ободом общее ребро, является цикл $c_7 = \{e_3, e_6, e_8\} \to \{v_2, v_3, v_4\}$. Но его присоединение нарушает ограничение (6.3): $4-6+4-1 \neq 0$. Кроме того, внешний простой цикл (обод) в этом случае есть пустое множество: $c_0 = \{e_5, e_6, e_8\} \oplus \{e_5, e_6, e_8\} = \varnothing$. Следовательно, цикл не может быть добавлен в множество $C_p$.

Просматривая далее список, последовательно присоединим цикл $c_8$, а затем $c_9$:

$$c_8 = \{e_5, e_7, e_9\} \to \{v_2, v_3, v_5\}; \quad c_9 = \{e_6, e_7, e_{10}\} \to \{v_2, v_4, v_5\}.$$

Ограничения (6.2) и (6.3) выполняются и циклы $c_8$ и $c_9$ могут быть помещены в множество $C_p$.

В конечном итоге для множества $C_p=\{c_1,c_2,c_4,c_8,c_9\}$ получаем $P_e$=<2,2,2,0,2,2,2,1,2,2> и $P_v$=<3,4,3,3,2>. Значение функционала Маклейна (6.1) для выбранной системы циклов равно нулю: $F(C_p)=0$; ограничение (6.3) также выполнено: $5-9+5-1=0$. Обод определяется как сумма предыдущего обода и $c_8$ и $c_9$:

$$c_0 = \{e_5, e_6, e_8\} \oplus \{e_5, e_7, e_9\} \oplus \{e_6, e_7, e_{10}\} = \{e_8, e_9, e_{10}\}.$$

На последующих шагах не удается найти цикл, удовлетворяющий условиям (6.2) и (6.3); конец работы алгоритма выделения фрагментов (см. рис. 6.17).

Таким образом, мы выделили множество изометрических циклов, характеризующее плоскую часть графа, описываемую перестановкой (1,2,3,4,5,6,7,8,9,10). Количество рёбер в данной плоской части $k_{\text{u}} = 9$, а количество вершин $k_{\text{x}} = 5$.

Отметим в заключение следующее свойство построенной фрагментарной модели задачи: любой максимальный планарный суграф в несепарабельном неориентированном графе G может быть построен предлагаемым методом при надлежащем выборе перестановки элементарных фрагментов. Это свойство фрагментарной модели называется свойством достижимости.

$$(\mathbf{1},\mathbf{2},3,\mathbf{4},5,6,7,\mathbf{8},\mathbf{9},10) \rightarrow \mathrm{C}_p = \{c_1, c_2, c_4, c_8, c_9\}$$

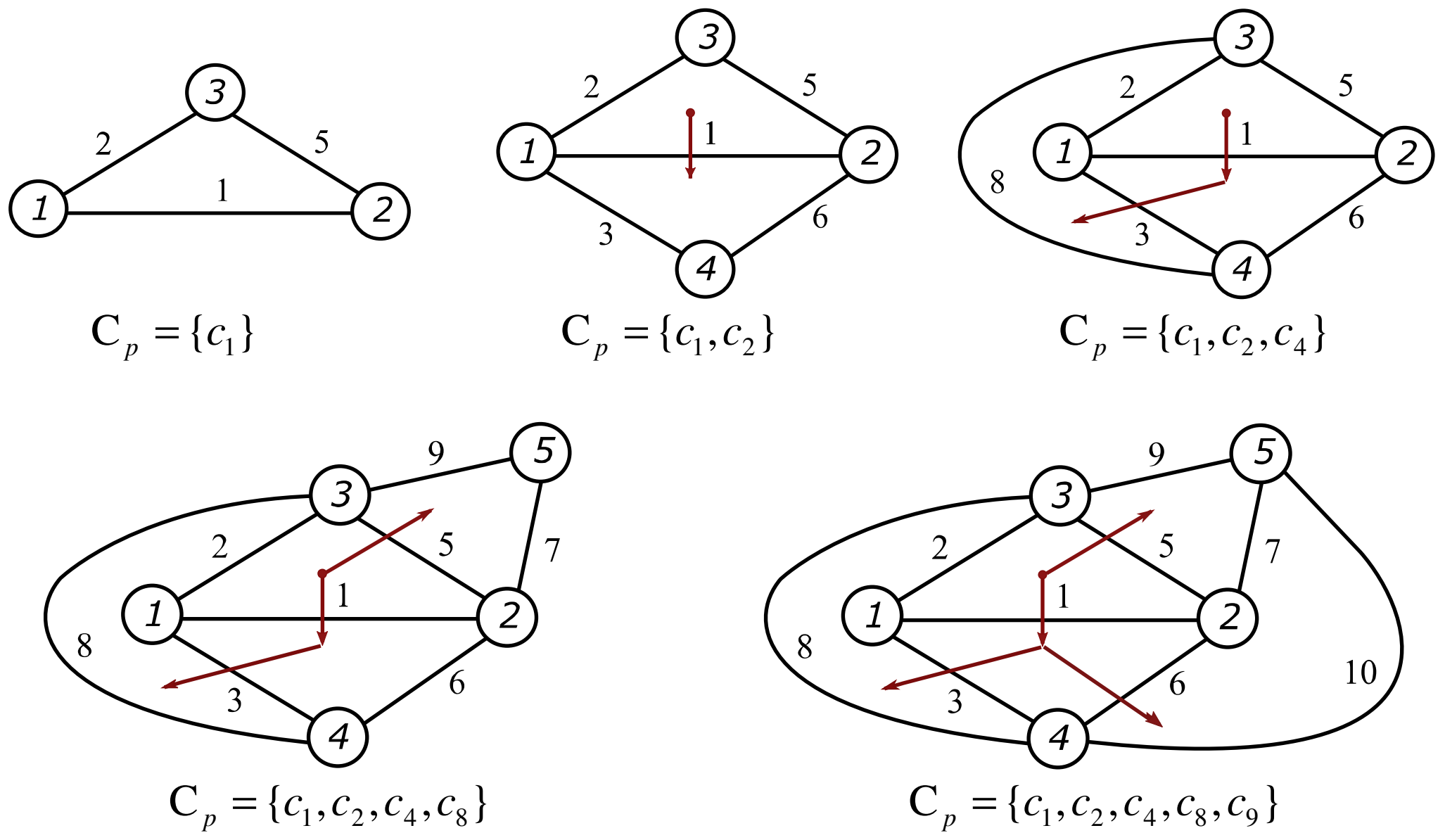

Рисунок 6.17. Плоская часть графа К$_5$.

Из сказанного вытекает, что задача поиска планарного суграфа с максимальным числом рёбер сводится к поиску максимального фрагмента во фрагментарной структуре, т.е. к комбинаторной задаче $\mathrm{F}(s) \underset{s \in S_n}{\rightarrow} \max$ $F(s) \underset{s \in S_n}{\rightarrow} \max$ отыскания некоторой экстремальной перестановки на множестве перестановок изометрических циклов ($N$ – мощность множества изометрических циклов). При этом любая перестановка является допустимой [11].

Поставим в соответствие каждой перестановке значение критерия на соответствующем ей максимальном фрагменте, а именно число рёбер планарного суграфа. На множестве перестановок можно предложить несколько простых алгоритмов приближённого поиска оптимальной перестановки. Опишем некоторые из них.

а) Метод локального поиска. Множество перестановок может рассматриваться как метрическое пространство с некоторой метрикой $-\rho : S_n \times S_n \rightarrow R_+^1$ [10]. Выбирается начальная перестановка и ищется оптимальная перестановка в её $\varepsilon$-окрестности. Эта найденная перестановка является базовой для следующего шага. Условие остановки алгоритма: в $\varepsilon$-окрестности не существует перестановки с лучшим критерием.

В качестве $\varepsilon$ можно брать любое натуральное число 1,2,...

б) Хорошие результаты даёт эволюционная модель на фрагментарной структуре [15]. В этой модели в качестве пространства поиска выбирается множество $S_n=\{i_1,i_2,\dots,i_N\}$ всех перестановок чисел 1,2,…,$N$.

Оператор построения начальной популяции выделяет произвольное подмножество заданной мощности $Q$ из множества Y. Правило вычисления критерия селекции устроено следующим образом: по заданной перестановке фрагментов с помощью фрагментарного алгоритма строится максимальный допустимый фрагмент. Вычисляется значение целевой функции для этого фрагмента.

**Комментарии**

Одной из промежуточных задач используемой для построения математической модели выделения плоской части графа, является задача определения независимой системы циклов. Как правило, для проверки независимости системы применяется алгоритм Гаусса для проверки ранга системы. Учитывая специфику выделения независимой системы простых циклов, требуется произвести модификацию алгоритма с целью построения плоских конфигураций. Подробно разбирается процесс построения независимой системы простых циклов. Показано, что процесс построения независимой системы циклов зависит от последовательности расположения элементов в системе циклов. Приведены методы построения плоских конфигураций в задаче выделения плоских частей графа.

Развитие методов построения плоских конфигураций приводит к методу решения задачи выделения плоской части называемым *подходом снизу*. Данный подход основан на последовательном подключении простых циклов с соблюдением условия Эйлера. Представлен фрагментарный алгоритм для проведения последовательности действий и приведен пример построения плоской части графа.

# Глава 7. ВЫДЕЛЕНИЕ ПЛОСКОЙ ЧАСТИ ГРАФА

## 7.1. Выделение множества циклов непланарного графа

В главе 6 мы рассмотрели построение плоской конфигурации применяя подход снизу, то есть когда исходное множество циклов для описания рисунка графа пусто. Рассмотрим подход построения плоской конфигурации сверху, то есть когда начальное множество не пусто и состоит из множества изометрических циклов графа.

Такой способ построения алгоритма решения задачи выделения плоской части графа с минимальным количеством удаленных ребер и построение его топологического рисунка, требует предварительного рассмотрения некоторых способов решения и их описание. Очевидно, что метод должен предпологать последовательное удаление циклов из заданного множества циклов до получения подмножества циклов с заданным значением функционала Маклейна. При этом предпологается удаление минимального количества ребер. Есстественно, что множество исходных циклов должно состоять из изометрических циклов, так как любой простой цикл в графе есть линейная комбинация независимых изометрических циклов графа [33].

Будем рассматривать алгоритмы выделения множества изометрических циклов с минимальным значением функционала Маклейна в подпространстве циклов C.

В работах [3,35] показано, что структурное число **W**, полученное как произведение однострочных структурных чисел, характеризующих простые циклы, проходящие по хордам графа на выделенном дереве, определяет все множество базисов, состоящих из простых циклов, и его построение не зависит от выбора дерева графа.

Так как любой столбец (не повторяющийся четное число раз) структурного числа C характеризует систему независимых циклов, то соответствующему столбцу можно поставить в соответствие значение функционала Маклейна. Однако, ввиду того, что величина функционала Маклейна является интегральной характеристикой множества (подмножества) циклов, не удается для каждого цикла графа поставить в соответствие число $\lambda_i$ так, чтобы сумма всех $\lambda_i$ ($i = 1,2,...,\nu(G)$) определяла какое-то числовое значение функционала Маклейна. Однако можно поставить подмножеству изометрических циклов в соответствие значение функционала Маклейна [24]. При этом следует заметить, что подмножество выделяемых циклов может быть как линейно независимо, так и линейно зависимо.

Напомним, что алгебраической обратной производной структурного числа называется структурное число $\delta W/\delta\alpha$, равное:

$$\frac{\delta W}{\delta \alpha} = W \,\Big|\, \text{ столбцы, содержащие элемент } \alpha \text{, опущены.} \qquad (7.1)$$

Воспользовавшись способом записи структурного числа в виде семейства множеств, можно записать обратную производную как:

$$\frac{\delta W}{\delta \alpha} = \{c_i \mid \alpha \notin c_i, c_i \in C\} \tag{7.2}$$

Алгебраическая обратная производная структурного числа W эквивалентна рассмотрению множеств базисов, не содержащих базисы с номером цикла $\alpha$.

Таким образом, множество базисов, состоящих из простых циклов (или суграфов) можно записывать и хранить в компактном виде как результат произведения однострочных структурных чисел, состоящих из номеров суграфов (циклов), включающих хорды графа (т.е. в виде структурного числа C).

И так, из множества изометрических циклов можно удалять циклы с использованием операции взятия алгебраической обратной производной структурного числа графа.

Рассмотрим задачу построения топологического рисунка планарного графа как задачу дискретной оптимизации [23,26,28]. Будем рассматривать один из методов дискретной оптимизации - метод наискорейшего спуска.

### 7.2. Метод наискорейшего спуска

В нашем случае метод наискорейшего спуска выглядит следующим образом. Задано дискретное пространство $C_\tau$ и целевая функция в виде функционала Маклейна F(c), $c \in C_\tau$. Требуется найти $c'$: $F(c') = \min F(c)$, $c \in C_\tau$, card p = m-n+1, в задаче глобальной оптимизации или $F(c') = \min F(c)$, card p = m-n+1,$\in Q \subset C_\tau$ , в задаче локальной оптимизации. Где p – количество циклов в выделенном подмножестве.

Метод наискорейшего спуска является локально оптимальным и реализует некоторую стратегию частичного перебора. Он заключается в следующем.

Задается некоторое число равное значению фукционала Маклейна, определенное для всего множества изометрических циклов $C_\tau$. Затем путем последовательного исключения каждого цикла, определяем значение функционала Маклейна для каждого случая. Определяем цикл $c_1$ максимально изменяющий значение функционала Маклейна. Для подмножества изометрических циклов без выбранного цикла (обозначим его $Q_1 \in C_\tau \setminus \{c_1\}$), будем иметь значение функционала $F(Q_1)$ Вновь путем последовательного исключения каждого цикла, определяем значение функционала Маклейна для каждого случая. Определяем цикл $c_2$ максимально изменяющий значение функционала Маклейна. Для подмножества изометрических циклов без выбранных двух циклов (обозначим его $Q_2 \in C_\tau \setminus \{c_1, c_2\}$), будем иметь значение функционала $F(Q_2)$ и т.д. Поиск решения заканчивается, когда количество оставшихся изометрических циклов в подмножестве равно цикломатическому числу графа.

В результате формируется последовательность локальных минимумов $F(C_\tau) > F(Q_1) > F(Q_2) > \dots F(Q_{p-m+n-1})$. Выделенное подмножество изометрических циклов $Q_{p-m+n-1}$ является решением задачи в случае выделения подмножества изометрических циклов с минимальным значением функционала Маклейна. Если значение функционала Маклейна для выделенного подмножества изометрических циклов равно нулю, то граф планарен.

Рассмотрим вопрос подробнее. Пусть имеется некоторое множество изометрических циклов $C_\tau$ мощностью card $(C_\tau) = p$ и, определенным образом, связанное с ним семейство его подмножеств $Q(C_\tau)$. Будем считать, что множество $C_\tau$ упорядоченно.

Составим характеристический двоичный вектор $\mathbf{Z}$ состоящий из нулей и единиц. $z_i = 1$ означает, что элемент (цикл) с номером *i* выбирается в качестве элемента подмножества, а $z_i = 0$ означает, что элемент с номером *i* не входит в подмножество. Любое подмножество $\mathbf{Z} \in C_\tau$ однозначно задается своим характеристическим вектором.

Каждому характеристическому вектору имеющего не менее $\nu(G) = m-n+1$ единиц можно поставить в соответствии значение функционала Маклейна. И тогда по условию задачи, необходимо найти характеристический вектор размерностью $\nu(G) = m-n+1$ с минимальным значением функционала Маклейна. Если граф планарен, то такой двоичный характеристический вектор существует, и значение функцирнала Маклейна равно нулю.

***Пример 7..1.*** Для графа, представленного на рис. 7..1, множество $C_\tau$ изометрических циклов имеет вид:

$c_1 = \{e_1,e_3,e_5\}$;
$c_2 = \{e_2,e_3,e_4\}$;
$c_3 = \{e_4,e_8,e_9\}$;
$c_4 = \{e_5,e_6,e_{12}\}$;
$c_5 = \{e_5,e_8,e_{13}\}$;
$c_6 = \{e_6,e_7,e_{11}\}$;
$c_7 = \{e_7,e_8,e_{10}\}$;
$c_8 = \{e_1,e_2,e_9,e_{13}\}$;
$c_9 = \{e_{10},e_{11},e_{12},e_{13}\}$.

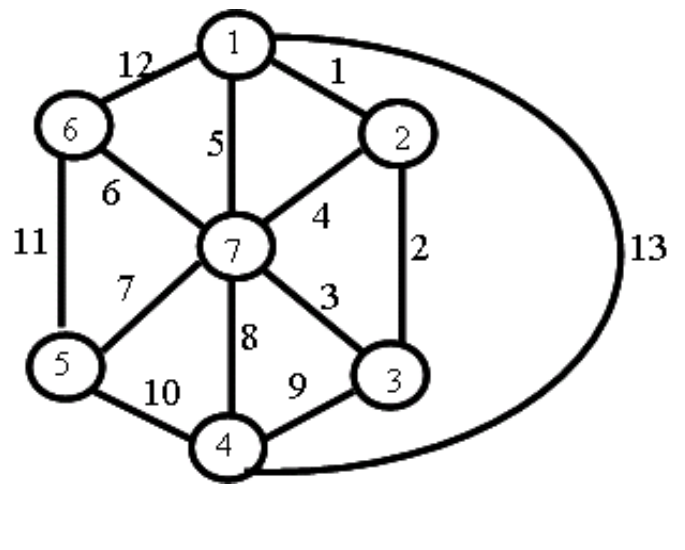

Рис. 7.1. Граф $G_2$.

И тогда можно записать точки векторного подпространства циклов в виде двоичных

векторов, где:

$q_0 = \{c_1,c_2,c_3,c_4,c_5,c_6,c_7,c_8,c_9\} \rightarrow$ кортеж $z_0 = <1,1,1,1,1,1,1,1,1>$ – описывает подмножество состоящее из всех изометрических циклов графа;

$q_1 = \{c_2,c_3,c_4,c_5,c_6,c_7,c_8,c_9\} \rightarrow$ кортеж $z_1 = <0,1,1,1,1,1,1,1,1>$ – описывает подмножество изометрических циклов графа за исключением цикла $c_1$;

………………………………………………………………..…………………

$q_9 = \{c_1,c_2,c_3,c_4,c_5,c_6,c_7,c_8,c_9\} \rightarrow$ кортеж $z_9 = <1,1,1,1,1,1,1,1,0>$ – описывает подмножество изометрических циклов графа за исключением цикла $c_9$;

$q_{10} = \{c_3,c_4,c_5,c_6,c_7,c_8,c_9\} \rightarrow$ кортеж $z_{10} = <0,0,1,1,1,1,1,1,1>$ – описывает подмножество изометрических циклов графа за исключением циклов $c_1$ и $c_2$;

…………………………………………………………………..………………

$q_{45} = \{c_1,c_2,c_3,c_4,c_5,c_6,c_7\} \rightarrow$ кортеж $z_{45} = <1,1,1,1,1,1,1,0,0>$ – описывает подмножество изометрических циклов графа за исключением циклов $c_8$ и $c_9$. И так далее.

Множество, состоящее из изометрических циклов, можно рассматривать как однострочное структурное число Q (кортеж изометрических циклов). Наличие операции взятия обратной производной структурного числа, позволяет пошагово записать процесс перехода от подмножеств, имеющих характеристические вектора, состоящие из g единиц (то есть мощностью g), к подмножествам, имеющим характеристические вектора, состоящие из g-1 единиц (то есть мощностью g-1). Таким образом, имеется возможность построить алгоритм последовательного исключения циклов из множества изометрических циклов до получения подмножества с минимальным значением функционала Маклейна. Причем направление изменения значения функционала Маклейна, по аналогии с методами оптимизации непрерывных функций, задает операция взятия обратной производной структурного числа Q.

Свяжем метод наискорейшего спуска с методами алгебры структурных чисел для определения поставленной цели [32].

Требуется из множества изометрических циклов графа $Q = \{c_1,c_2,\dots c_p\}$, $Q \in C_\tau$ найти элемент $q_k$ с определенным значением функционала Маклейна.

$$F(Q_p) = F\frac{\delta Q^h}{\delta c_1 \delta c_2 \dots \delta c_h} \rightarrow const, \qquad (7.3)$$

здесь $h = \mathbf{card}(C_\tau) - (m-n+1)$ – количество удаленных циклов, $c_i$ - удаляемый цикл, $F(Q_i)$ – значение функционала Маклейна для подмножества циклов $Q_i \subset C_\tau$.

И тогда подмножество, состоящее из определенного набора двоичных векторов, можно описывать структурным числом. Например:

Q – однострочное структурное число, описывающее кортеж состоящий из всего множества изометрических циклов графа;

$\frac{\delta Q}{\delta c_1}$ - однострочное структурное число описывающее кортеж изометрических циклов графа, с исключенным циклом $c_1$. Данный процесс соответствует характеристическому вектору $z_1$ =<0,1,1,1,1,1,1,1,1>.

$\frac{\delta^2 C_b}{\delta c_5 \delta c_9}$ - структурное число описывающее подмножество циклов, за исключением циклов $c_5$ и $c_9$. Данный процесс соответствует характеристическому вектору z =<1,1,1,1,0,1,1,1,0>.

При применении градиентных методов, в случае решения оптимизационных задач для непрерывных функций, замечено, что они сходятся к минимуму с высокой скоростью (со скоростью геометрической прогрессии) для гладких выпуклых функций. У таких функций наибольшее M и наименьшее $m$ (собственные значения матрицы вторых производных - матрицы Гессе H(z)) мало отличаются друг от друга, т. е. матрица H(z) хорошо обусловлена.

В нашем случае, целевая функция и множество циклов выпуклы, и для каждого подмножества Z, состоящего из изометрических циклов с мощностью большей цикломатического числа, можно составить матрицу со свойствами адекватными матрице Гесса:

$$H(Z) = \begin{Vmatrix} \varnothing & \frac{\delta^2 Z}{\delta\alpha_1 \delta\alpha_2} & \dots & \frac{\delta^2 Z}{\delta\alpha_1 \delta\alpha_t} \\ \frac{\delta^2 Z}{\delta\alpha_2 \delta\alpha_1} & \varnothing & \dots & \frac{\delta^2 Z}{\delta\alpha_2 \delta\alpha_t} \\ \dots & \dots & \dots & \dots \\ \frac{\delta^2 Z}{\delta\alpha_t \delta\alpha_1} & \frac{\delta^2 Z}{\delta\alpha_t \delta\alpha_2} & \dots & \varnothing \end{Vmatrix}. \qquad (7.4)$$

где $t$ - количество изометрических циклов в данном подмножестве.

Множество изометрических циклов графа будем обозначать как $C_\tau$. Элементы линейного подпространства циклов состоящего из изометрических циклов мощностью равной цикломатическому числу графа $\nu(G) = m - n + 1$ будем обозначать, как $c_b$. В свою очередь, $c_b \in C_\tau$.

С множеством изометрических циклов в графе связан ряд инвариантов графа []. Одним из таких инвариантов является мощность множества изометрических циклов $C_\tau$. Другим инвариантом может служить вектор изометрических циклов, упорядоченный по возрастанию их длин. Следующим инвариантом, по аналогии с вектором локальных степеней, является вектор количества изометрических циклов проходящих по ребрам графа. Будем называть его *вектором циклов по ребрам*:

$$P_e = < a_1, a_2, ..., a_m >, \quad (7.5)$$

где $a_i$- количество изометрических циклов графа, проходящих по ребру $e_i$,$i$=(1,2,…,$m$).

Инвариантом является также вектор количества изометрических циклов проходящих по вершинам графа. Будем называть его вектором циклов по вершинам:

$$P_v = < b_1, b_2, ..., b_n >, \quad (7.6)$$

где $b_j$ - количество изометрических циклов графа, проходящих по вершине $v_j$,$j$=(1,2,…,$n$).

Теперь можно описать алгоритм выделения множества, состоящего из изометрических циклов графа методом наискорейшего спуска.

**Шаг 1**. [**Взятие обратной производной**].

К каждому элементу подмножества изометрических циклов применяем операцию взятия обратной производной. В результате получим подмножество из оставшихся элементов. Идем на шаг 2.

**Шаг 2**. [**Выделение подмножества с минимальным значением функционала Маклейна**].

Выделяем подмножество изометрических циклов, имеющее минимальное значение функционала Маклейна. Если таких подмножеств несколько, то среди них выбираем то, у которого в строке имеется минимальный элемент матрицы H(Z). Идем на шаг 3.

**Шаг 3**. [**Определение количества элементов**].

Если количество элементов в подмножестве больше цикломатического числа, то идем на шаг 1. Если количество элементов равно цикломатическому числу графа, то конец работы алгоритма. При значении функционала Маклейна равного нулю - граф планарен. В противном случае - непланарен.

Сказанное поясним на следующем примере.

***Пример 7..2.*** Определить является ли граф $G_3$ представленный на рис. 7..2 планарным, и если он планарен построить плоский рисунок для графа $G_3$.

$|V|=11$; $|E|=20$; $\nu(G)=10$.

Множество изометрических циклов графа $C_\tau$:

$c_1 = \{e_1,e_3,e_5,e_{11}\}$; $c_2 = \{e_1,e_4,e_6,e_{16}\}$; $c_3 = \{e_2,e_3,e_7,e_9\}$; $c_4 = \{e_2,e_3,e_8,e_{13}\}$;
$c_5 = \{e_2,e_4,e_8,e_{14}\}$; $c_6 = \{e_3,e_4,e_{13},e_{14}\}$; $c_7 = \{e_5,e_6,e_{12},e_{14},e_{16}\}$;
$c_8 = \{e_5,e_6,e_{12},e_{15},e_{19}\}$; $c_9 = \{e_7,e_8,e_{10}\}$; $c_{10} = \{e_9,e_{10},e_{13}\}$; $c_{11} = \{e_{11},e_{12},e_{13}\}$;
$c_{12} = \{e_1,e_2,e_5,e_8,e_{12}\}$; $c_{13} = \{e_1,e_4,e_5,e_{12},e_{14}\}$; $c_{14} = \{e_{14},e_{15},e_{16},e_{19}\}$;
$c_{15} = \{e_{14},e_{15},e_{17},e_{20}\}$; $c_{16} = \{e_{16},e_{17},e_{18}\}$; $c_{17} = \{e_{18},e_{19},e_{20}\}$.

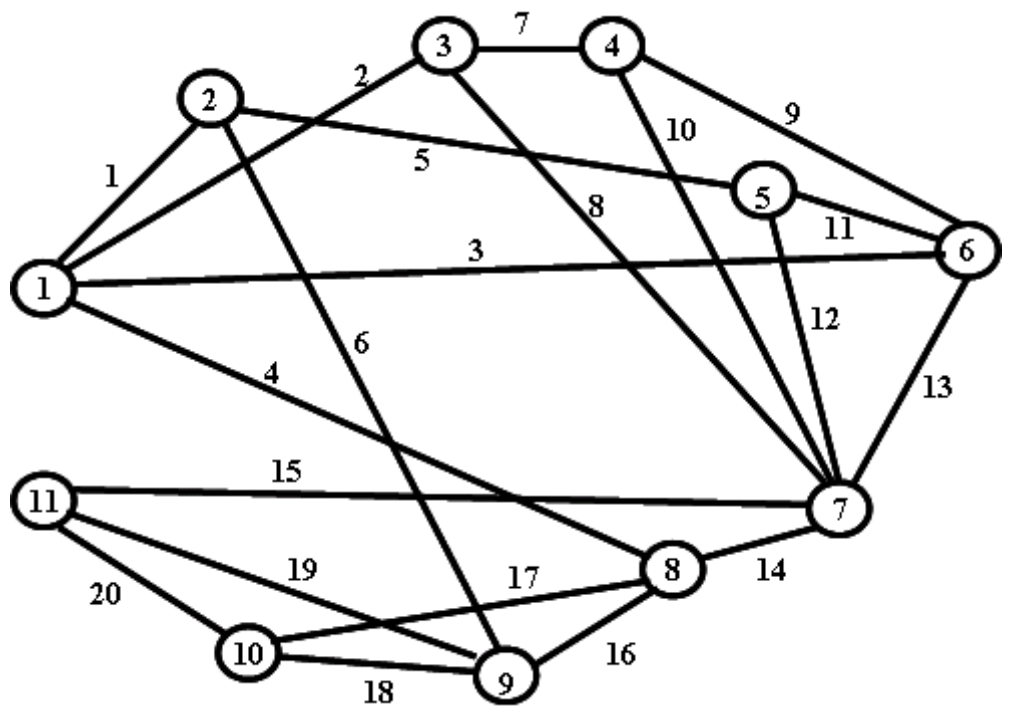

Рис. 7.2. Граф $G_3$.

Функционал Маклейна для множества изометрических циклов равен: $F(G_\tau) = 92$;

$$F\left(\frac{\delta C_b}{\delta c_1}\right) = 78;\ F\left(\frac{\delta C_b}{\delta c_2}\right) = 78;\ F\left(\frac{\delta C_b}{\delta c_3}\right) = 84;\ F\left(\frac{\delta C_b}{\delta c_4}\right) = 76;\ F\left(\frac{\delta C_b}{\delta c_5}\right) = 72;$$

$$F\left(\frac{\delta C_b}{\delta c_6}\right) = 72;\ F\left(\frac{\delta C_b}{\delta c_7}\right) = 66;\ F\left(\frac{\delta C_b}{\delta c_8}\right) = 74;\ F\left(\frac{\delta C_b}{\delta c_9}\right) = 88;\ F\left(\frac{\delta C_b}{\delta c_{10}}\right) = 88;$$

$$F\left(\frac{\delta C_b}{\delta c_{11}}\right) = 82;\ F\left(\frac{\delta C_b}{\delta c_{12}}\right) = 68;\ F\left(\frac{\delta C_b}{\delta c_{13}}\right) = 64;\ F\left(\frac{\delta C_b}{\delta c_{14}}\right) = 76;\ F\left(\frac{\delta C_b}{\delta c_{15}}\right) = 82;$$

$$F\left(\frac{\delta C_b}{\delta c_{16}}\right) = 88;\ F\left(\frac{\delta C_b}{\delta c_{17}}\right) = 90.$$

Максимальное изменение значения функционала Маклейна получается после удаления цикла $c_{13}$. Удаляем цикл $c_{13}$. Продолжая процесс, получим:

$$F\left(\frac{\delta^2 C_b}{\delta c_{13}\delta c_1}\right) = 54;\ F\left(\frac{\delta^2 C_b}{\delta c_{13}\delta c_2}\right) = 54;\ F\left(\frac{\delta^2 C_b}{\delta c_{13}\delta c_3}\right) = 56;\ F\left(\frac{\delta^2 C_b}{\delta c_{13}\delta c_4}\right) = 48;$$

$$F\left(\frac{\delta^2 C_b}{\delta c_{13}\delta c_5}\right) = 48;\ F\left(\frac{\delta^2 C_b}{\delta c_{13}\delta c_6}\right) = 48;\ F\left(\frac{\delta^2 C_b}{\delta c_{13}\delta c_7}\right) = 44;\ F\left(\frac{\delta^2 C_b}{\delta c_{13}\delta c_8}\right) = 50;$$

$$F\left(\frac{\delta^2 C_b}{\delta c_{13}\delta c_9}\right) = 60;\ F\left(\frac{\delta^2 C_b}{\delta c_{13}\delta c_{10}}\right) = 60;\ F\left(\frac{\delta^2 C_b}{\delta c_{13}\delta c_{11}}\right) = 56;\ F\left(\frac{\delta^2 C_b}{\delta c_{13}\delta c_{12}}\right) = 46;$$

$$F\left(\frac{\delta^2 C_b}{\delta c_{13}\delta c_{14}}\right) = 50;\ F\left(\frac{\delta^2 C_b}{\delta c_{13}\delta c_{15}}\right) = 56;\ F\left(\frac{\delta^2 C_b}{\delta c_{13}\delta c_{16}}\right) = 60;\ F\left(\frac{\delta^2 C_b}{\delta c_{13}\delta c_{17}}\right) = 62.$$

Максимальное изменение значения функционала Маклейна получается после удаления цикла $c_7$. Удаляем $c_7$. Продолжая процесс, получим:

$$F\left(\frac{\delta^3 C_b}{\delta c_{13}\delta c_7\delta c_1}\right) = 36;\ F\left(\frac{\delta^3 C_b}{\delta c_{13}\delta c_7\delta c_2}\right) = 38;\ F\left(\frac{\delta^3 C_b}{\delta c_{13}\delta c_7\delta c_3}\right) = 36;$$

$$F\left(\frac{\delta^3 C_b}{\delta c_{13}\delta c_7\delta c_4}\right) = 28;\ F\left(\frac{\delta^3 C_b}{\delta c_{13}\delta c_7\delta c_5}\right) = 30;\ F\left(\frac{\delta^3 C_b}{\delta c_{13}\delta c_7\delta c_6}\right) = 30;$$

$$F\left(\frac{\delta^3 C_b}{\delta c_{13}\delta c_7\delta c_8}\right) = 36;\quad F\left(\frac{\delta^3 C_b}{\delta c_{13}\delta c_7\delta c_9}\right) = 40;\quad F\left(\frac{\delta^3 C_b}{\delta c_{13}\delta c_7\delta c_{10}}\right) = 40;$$

$$F\left(\frac{\delta^3 C_b}{\delta c_{13}\delta c_7\delta c_{11}}\right) = 38;\quad F\left(\frac{\delta^3 C_b}{\delta c_{13}\delta c_7\delta c_{12}}\right) = 30;\quad F\left(\frac{\delta^3 C_b}{\delta c_{13}\delta c_7\delta c_{14}}\right) = 34;$$

$$F\left(\frac{\delta^3 C_b}{\delta c_{13}\delta c_7\delta c_{15}}\right) = 38;\quad F\left(\frac{\delta^3 C_b}{\delta c_{13}\delta c_7\delta c_{16}}\right) = 42;\quad F\left(\frac{\delta^3 C_b}{\delta c_{13}\delta c_7\delta c_{17}}\right) = 42.$$

Максимальное изменение значения функционала Маклейна получается после удаления цикла $c_4$. Удаляем $c_4$. Продолжая процесс, получим:

$$F\left(\frac{\delta^4 C_b}{\delta c_{13}\delta c_7\delta c_4\delta c_1}\right) = 22;\quad F\left(\frac{\delta^4 C_b}{\delta c_{13}\delta c_7\delta c_4\delta c_2}\right) = 22;\quad F\left(\frac{\delta^4 C_b}{\delta c_{13}\delta c_7\delta c_4\delta c_3}\right) = 24;$$

$$F\left(\frac{\delta^4 C_b}{\delta c_{13}\delta c_7\delta c_4\delta c_5}\right) = 18;\quad F\left(\frac{\delta^4 C_b}{\delta c_{13}\delta c_7\delta c_4\delta c_6}\right) = 18;\quad F\left(\frac{\delta^4 C_b}{\delta c_{13}\delta c_7\delta c_4\delta c_8}\right) = 20;$$

$$F\left(\frac{\delta^4 C_b}{\delta c_{13}\delta c_7\delta c_4\delta c_9}\right) = 26;\quad F\left(\frac{\delta^4 C_b}{\delta c_{13}\delta c_7\delta c_4\delta c_{10}}\right) = 26;\quad F\left(\frac{\delta^4 C_b}{\delta c_{13}\delta c_7\delta c_4\delta c_{11}}\right) = 24;$$

$$F\left(\frac{\delta^4 C_b}{\delta c_{13}\delta c_7\delta c_4\delta c_{12}}\right) = 18;\quad F\left(\frac{\delta^4 C_b}{\delta c_{13}\delta c_7\delta c_4\delta c_{14}}\right) = 18;\quad F\left(\frac{\delta^4 C_b}{\delta c_{13}\delta c_7\delta c_4\delta c_{15}}\right) = 22;$$

$$F\left(\frac{\delta^4 C_b}{\delta c_{13}\delta c_7\delta c_4\delta c_{16}}\right) = 26;\quad F\left(\frac{\delta^4 C_b}{\delta c_{13}\delta c_7\delta c_4\delta c_{17}}\right) = 25.$$

Максимальное изменение значения функционала Маклейна получается после удаления циклов $c_5$,$c_6$,$c_{12}$ и $c_{14}$.

Рассмотрим вариант их попарного удаления:

$$F\left(\frac{\delta^5 C_b}{\delta c_{13}\delta c_7\delta c_4\delta c_5\delta c_6}\right) = 12;\quad F\left(\frac{\delta^5 C_b}{\delta c_{13}\delta c_7\delta c_4\delta c_5\delta c_{12}}\right) = 12;$$

$$F\left(\frac{\delta^5 C_b}{\delta c_{13}\delta c_7\delta c_4\delta c_5\delta c_{14}}\right) = 10;\quad F\left(\frac{\delta^5 C_b}{\delta c_{13}\delta c_7\delta c_4\delta c_6\delta c_{12}}\right) = 8;$$

$$F\left(\frac{\delta^5 C_b}{\delta c_{13}\delta c_7\delta c_4\delta c_6\delta c_{14}}\right) = 10;\quad F\left(\frac{\delta^5 C_b}{\delta c_{13}\delta c_7\delta c_4\delta c_{12}\delta c_{14}}\right) = 8.$$

Рассмотрим пересечение циклов:

$c_5 \cap c_6 = \{e_2,e_4,e_8,e_{14}\} \cap \{e_3,e_4,e_{13},e_{14}\}=\{e_4,e_{14}\}$;

$c_5 \cap c_{12} = \{e_2,e_4,e_8,e_{14}\} \cap \{e_1,e_2,e_5,e_8,e_{12}\}=\{e_2,e_8\}$;

$c_5 \cap c_{14} = \{e_2,e_4,e_8,e_{14}\} \cap \{e_{14},e_{15},e_{16},e_{19}\}=\{e_{14}\}$;

$c_6 \cap c_{12} = \{e_3,e_4,e_{13},e_{14}\} \cap \{e_1,e_2,e_5,e_8,e_{12}\}=\varnothing$;

$c_6 \cap c_{14} = \{e_2,e_4,e_8,e_{14}\} \cap \{e_{14},e_{15},e_{16},e_{19}\}=\{e_{14}\}$;

$c_{12} \cap c_{14} = \{e_1,e_2,e_5,e_8,e_{12}\} \cap \{e_{14},e_{15},e_{16},e_{19}\}=\varnothing$.

Здесь перспективно выбрать удаление одного цикла $c_{12}$ (как участвующего и в пересечении $c_6 \cap c_{12}$ и в пересечении $c_{12} \cap c_{14}$).

Удаляем $c_{12}$. Продолжаем процесс:

$$F(\frac{\delta^5 C_b}{\delta c_{13}\delta c_7\delta c_4\delta c_{12}\delta c_1}) = 12;\ F(\frac{\delta^5 C_b}{\delta c_{13}\delta c_7\delta c_4\delta c_{12}\delta c_2}) = 14;$$

$$F(\frac{\delta^5 C_b}{\delta c_{13}\delta c_7\delta c_4\delta c_{12}\delta c_3}) = 14;\ F(\frac{\delta^5 C_b}{\delta c_{13}\delta c_7\delta c_4\delta c_{12}\delta c_5}) = 10;$$

$$F(\frac{\delta^5 C_b}{\delta c_{13}\delta c_7\delta c_4\delta c_{12}\delta c_6}) = 10;\ F(\frac{\delta^5 C_b}{\delta c_{13}\delta c_7\delta c_4\delta c_{12}\delta c_8}) = 14;\ F(\frac{\delta^5 C_b}{\delta c_{13}\delta c_7\delta c_4\delta c_{12}\delta c_9}) = 16;$$

$$F(\frac{\delta^5 C_b}{\delta c_{13}\delta c_7\delta c_4\delta c_{12}\delta c_{10}}) = 16;\ F(\frac{\delta^5 C_b}{\delta c_{13}\delta c_7\delta c_4\delta c_{12}\delta c_{11}}) = 14;\ F(\frac{\delta^5 C_b}{\delta c_{13}\delta c_7\delta c_4\delta c_{14}2\delta c_{14}}) = 8;$$

$$F(\frac{\delta^5 C_b}{\delta c_{13}\delta c_7\delta c_4\delta c_{12}\delta c_{15}}) = 16;\ F(\frac{\delta^5 C_b}{\delta c_{13}\delta c_7\delta c_4\delta c_{12}\delta c_{16}}) = 18;\ F(\frac{\delta^5 C_b}{\delta c_{13}\delta c_7\delta c_4\delta c_{12}\delta c_{17}}) = 18.$$

Максимальное изменение значения функционала Маклейна получается после удаления цикла $c_{14}$. Удаляем $c_{14}$. Продолжаем процесс:

$$F(\frac{\delta^6 C_b}{\delta c_{13}\delta c_7\delta c_4\delta c_{12}\delta c_{14}\delta c_1}) = 6;\ F(\frac{\delta^6 C_b}{\delta c_{13}\delta c_7\delta c_4\delta c_{12}\delta c_{14}\delta c_2}) = 6;$$

$$F(\frac{\delta^6 C_b}{\delta c_{13}\delta c_7\delta c_4\delta c_{12}\delta c_{14}\delta c_3}) = 6;\ F(\frac{\delta^6 C_b}{\delta c_{13}\delta c_7\delta c_4\delta c_{12}\delta c_{14}\delta c_5}) = 4;$$

$$F(\frac{\delta^6 C_b}{\delta c_{13}\delta c_7\delta c_4\delta c_{12}\delta c_{14}\delta c_6}) = 0;\ F(\frac{\delta^6 C_b}{\delta c_{13}\delta c_7\delta c_4\delta c_{12}\delta c_{14}\delta c_8}) = 8;$$

$$F(\frac{\delta^6 C_b}{\delta c_{13}\delta c_{17}\delta c_4\delta c_{12}\delta c_{14}\delta c_9}) = 8;\ F(\frac{\delta^6 C_b}{\delta c_{13}\delta c_7\delta c_4\delta c_{12}\delta c_{14}\delta c_{10}}) = 6;$$

$$F(\frac{\delta^6 C_b}{\delta c_{13}\delta c_7\delta c_4\delta c_{12}\delta c_{14}\delta c_{11}}) = 6;\ F(\frac{\delta^6 C_b}{\delta c_{13}\delta c_7\delta c_4\delta c_{12}\delta c_{14}\delta c_{15}}) = 6;$$

$$F(\frac{\delta^6 C_b}{\delta c_{13}\delta c_7\delta c_4\delta c_{12}\delta c_{14}\delta c_{16}}) = 8;\ F(\frac{\delta^6 C_b}{\delta c_{13}\delta c_7\delta c_4\delta c_{12}\delta c_{14}\delta c_{17}}) = 8.$$

Максимальное изменение значения функционала Маклейна получается после удаления цикла $c_5$. Удаляем $c_5$. Продолжаем процесс:

$$F(\frac{\delta^7 C_b}{\delta c_{13}\delta c_7\delta c_4\delta c_{12}\delta c_{14}\delta c_6\delta c_1}) = 0;\quad F(\frac{\delta^7 C_b}{\delta c_{13}\delta c_7\delta c_4\delta c_{12}\delta c_{14}\delta c_6\delta c_2}) = 0;$$

$$F(\frac{\delta^7 C_b}{\delta c_{13}\delta c_7\delta c_4\delta c_{12}\delta c_{14}\delta c_6\delta c_3}) = 0;\quad F(\frac{\delta^7 C_b}{\delta c_{13}\delta c_7\delta c_4\delta c_{12}\delta c_{14}\delta c_6\delta c_5}) = 0;$$

$$F(\frac{\delta^7 C_b}{\delta c_{13}\delta c_7\delta c_4\delta c_{12}\delta c_{14}\delta c_6\delta c_8}) = 0;\quad F(\frac{\delta^7 C_b}{\delta c_{13}\delta c_7\delta c_4\delta c_{12}\delta c_{14}\delta c_6\delta c_9}) = 0;$$

$$F\left(\frac{\delta^7 C_b}{\delta c_{13}\delta c_7\delta c_4\delta c_{12}\delta c_{14}\delta c_6\delta c_{10}}\right)=0;\quad F\left(\frac{\delta^7 C_b}{\delta c_{13}\delta c_7\delta c_4\delta c_{12}\delta c_{14}\delta c_6\delta c_{11}}\right)=0;$$

$$F\left(\frac{\delta^7 C_b}{\delta c_{13}\delta c_7\delta c_4\delta c_{12}\delta c_{14}\delta c_6\delta c_{15}}\right)=0;\quad F\left(\frac{\delta^7 C_b}{\delta c_{13}\delta c_7\delta c_4\delta c_{12}\delta c_{14}\delta c_6\delta c_{16}}\right)=0;$$

$$F\left(\frac{\delta^7 C_b}{\delta c_{13}\delta c_7\delta c_4\delta c_{12}\delta c_{14}\delta c_6\delta c_{17}}\right)=0.$$

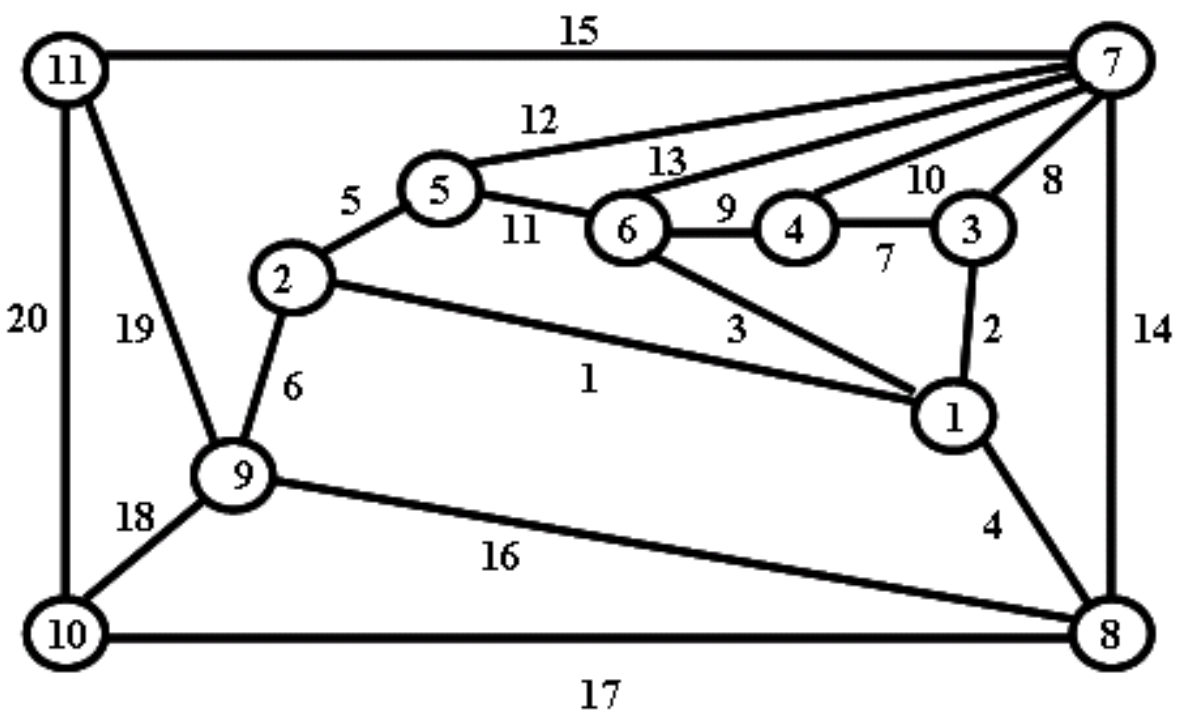

Рис. 7.3. Рисунок планарного графа $G_3$.

Удаляем $c_{15}$. После удаления изометрических циклов $c_{13},c_7,c_4,c_{12},c_{14},c_6,c_{15}$ цикломатическое число графа $\nu(G_3)$ =17-7=10 и функционал Маклейна имеет нулевое значение, следовательно граф $G_3$ планарный (см. рис. 7..3).

На рис. 7..4 представлен график изменения значения целевой функции Маклейна при пошаговом удалении циклов. Показано, что в случае удаления цикла $c_5$, вместо цикла $c_{12}$ (на четвертом шаге), можно и не выделить базис с нулевым значением функционала Маклейна.

Для выделенной системы изометрических циклов алгоритмом (описанным в данной главе) строится вращение вершин графа, что и описывает топологический рисунок графа:

вращение вершины $\sigma_1$: 3 8 2 6
вращение вершины $\sigma_2$: 5 1 9
вращение вершины $\sigma_3$: 4 7 1
вращение вершины $\sigma_4$: 6 7 3
вращение вершины $\sigma_5$: 2 7 6
вращение вершины $\sigma_6$: 5 7 4 1
вращение вершины $\sigma_7$: 8 3 4 6 5 11
вращение вершины $\sigma_8$: 9 1 7 10
вращение вершины $\sigma_9$: 2 8 10 11
вращение вершины $\sigma_{10}$: 9 8 11
вращение вершины $\sigma_{11}$: 7 9 10

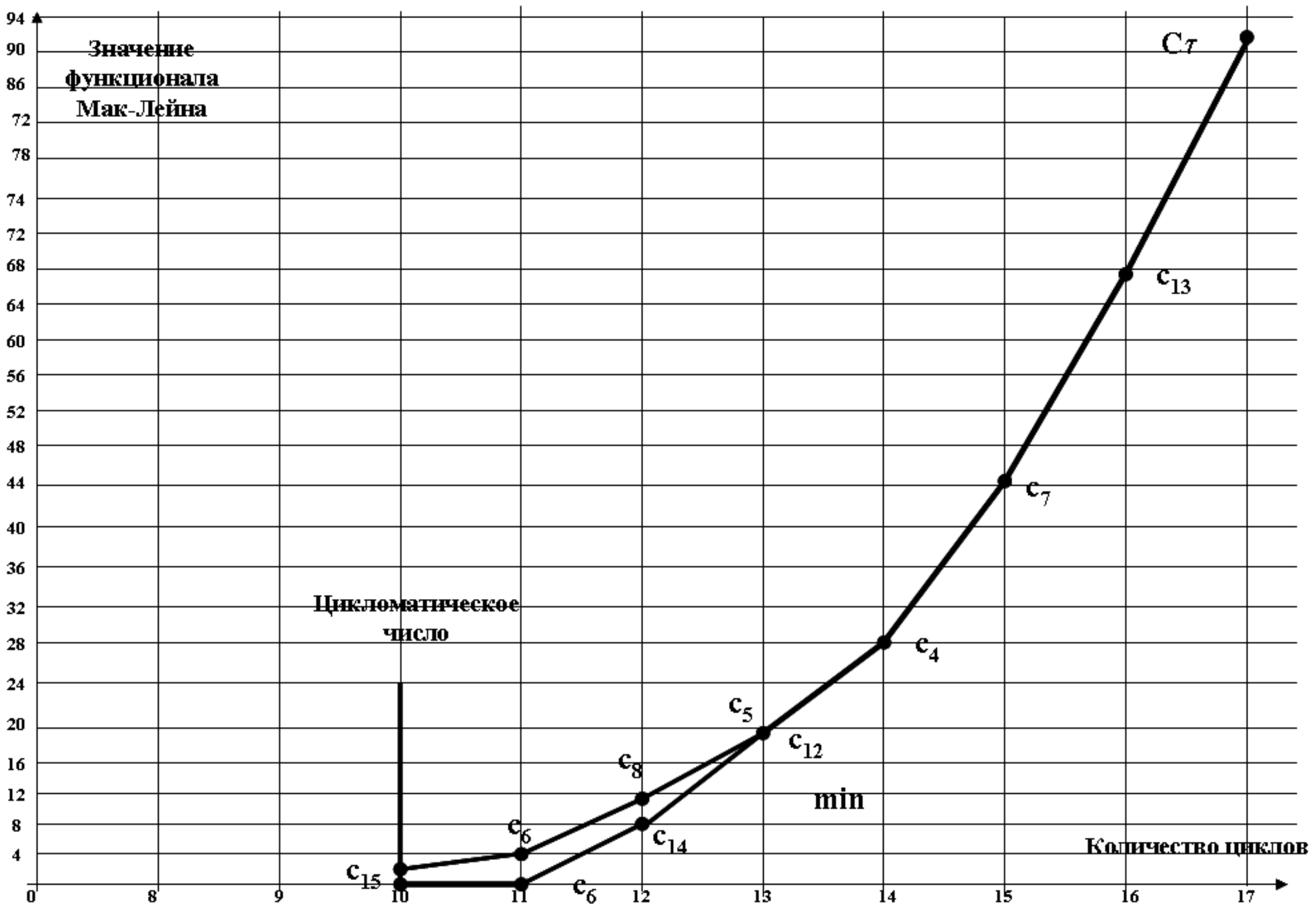

Рис. 7.4. График изменения значений целевой функции.

Удаленные циклы с базисными циклами составляют 0-подмножества:

$c_{13} = \{e_1,e_4,e_5,e_{12},e_{14}\} = c_1 \oplus c_3 \oplus c_5 \oplus c_9 \oplus c_{10} \oplus c_{11} = \{e_1,e_3,e_5,e_{11}\} \oplus$
$\oplus \{e_2,e_3,e_7,e_9\} \oplus \{e_2,e_4,e_8,e_{14}\} \oplus \{e_7,e_8,e_{10}\} \oplus \{e_9,e_{10},e_{13}\} \oplus \{e_{11},e_{12},e_{13}\}$;

$c_7 = \{e_5,e_6,e_{12},e_{14},e_{16}\} = c_1 \oplus c_2 \oplus c_3 \oplus c_5 \oplus c_9 \oplus c_{10} \oplus c_{11} = \{e_1,e_3,e_5,e_{11}\} \oplus$
$\{e_1,e_4,e_6,e_{16}\} \oplus \{e_2,e_3,e_7,e_9\} \oplus \{e_2,e_4,e_8,e_{14}\} \oplus \{e_7,e_8,e_{10}\} \oplus$
$\oplus \{e_9,e_{10},e_{\ 3}\} \oplus \{e_{11},e_{12},e_{13}\}$;

$c_4 = \{e_2,e_3,e_8,e_{13}\} = c_3 \oplus c_9 \oplus c_{10} = \{e_2,e_3,e_7,e_9\} \oplus \{e_7,e_8,e_{10}\} \oplus \{e_9,e_{10},e_{13}\}$;

$c_{12}= \{e_1,e_2,e_5,e_8,e_{12}\} = c_1 \oplus c_3 \oplus c_9 \oplus c_{10} \oplus c_{11} = \{e_1,e_3,e_5,e_{11}\} \oplus \{e_2,e_3,e_7,e_9\} \oplus$
$\oplus \{e_7,e_8,e_{10}\} \oplus \{e_9,e_{10},e_{13}\} \oplus \{e_{11},e_{12},e_{13}\}$;

$c_{14} = \{e_{14},e_{15},e_{16},e_{19}\} = c_1 \oplus c_2 \oplus c_3 \oplus c_5 \oplus c_8 \oplus c_9 \oplus c_{10} \oplus c_{11} = \{e_1,e_3,e_5,e_{11}\} \oplus$
$\oplus \{e_1,e_4,e_6,e_{16}\} \oplus \{e_2,e_3,e_7,e_9\} \oplus \{e_2,e_4,e_8,e_{14}\} \oplus \{e_5,e_6,e_{12},e_{15},e_{19}\}\oplus$
$\oplus \{e_7,e_8,e_{10}\} \oplus \{e_9,e_{10},e_{13}\} \oplus \{e_{11},e_{12},e_{13}\}$;

$c_6 = \{e_3,e_4,e_{13},e_{14}\} = c_3 \oplus c_5 \oplus c_9 \oplus c_{10} = \{e_2,e_3,e_7,e_9\} \oplus \{e_2,e_4,e_8,e_{14}\} \oplus$
$\oplus \{e_7,e_8,e_{10}\} \oplus \{e_9,e_{10},e_{13}\}$;

$c_{15} = \{e_{14},e_{15},e_{17},e_{20}\} = c_1 \oplus c_2 \oplus c_3 \oplus c_5 \oplus c_8 \oplus c_9 \oplus c_{10} \oplus c_{11} \oplus c_{16} \oplus c_{17} =$
$= \{e_1,e_3,e_5,e_{11}\} \oplus \{e_1,e_4,e_6,e_{16}\} \oplus \{e_2,e_3,e_7,e_9\} \oplus \{e_2,e_4,e_8,e_{14}\} \oplus$
$\oplus \{e_5,e_6,e_{12},e_{15},e_{19}\} \oplus \{e_7,e_8,e_{10}\} \oplus \{e_9,e_{10},e_{13}\} \oplus \{e_{11},e_{12},e_{13}\} \oplus$
$\oplus \{e_{16},e_{17},e_{18}\} \oplus \{e_{18},e_{19},e_{20}\}$.

Таким образом, применение операции взятия обратной производной структурного числа W, относительно отдельно взятого единичного цикла, позволяет поставить в соответствие оставшемуся подмножеству циклов определенное целое положительное число. Данное число определяется значением функционала Маклейна, полученное для подмножества изометрических циклов без учета удаленных.

Естественно, что в приведенном алгоритме, после каждого удаления цикла проверяется вектор количества циклов проходящих по вершинам $V_x$ на отсутствия нулевых координат. В случае появления нулевых значений, выбранный цикл возвращается в систему и выбирается другой цикл, удаление которого не приводит к появлению нулевых значений.

Следует заметить, что мы рассмотрели алгоритм наискорейшего спуска для установления планарности графа, который раскрывает механизм построения плоского топологического рисунка графа G. Однако данный процесс неоднозначен (см. рис.7..4) и имеет большую вычислительную сложность, чем алгоритм Хопкрофта-Тарьяна [28] позволяющий распознать планарность графа за линейное время.

### 7.3. Выделение базиса подпространства циклов методом структурных чисел

Все множество базисов состоящих из изометрических циклов **W** можно получить как произведение однострочных структурных чисел (исключая столбцы, повторяющиеся четное число раз), характеризующих подмножество изометрических циклов с ребром, принадлежащим хорде графа для выбранного дерева. Таким образом:

$$W = q_1 \times q_2 \times \ldots \times q_{n-m+1}. \qquad (7..7)$$

Множество баз матроида изометрических циклов также можно представлять как структурное число W, где каждый столбец (элемент) характеризует базис подпространства циклов. С другой стороны, структурное число можно представлять как произведение однострочных структурных чисел характеризующих изометрические циклы, проходящие по выбранной хорде. Причем согласно правилам алгебры структурных чисел [3] выбор дерева графа не влияят на конечный результат.

В свою очередь, подмножеству циклов можно поставить в соответствие два вектора. Вектор циклов по ребрам $P_e$, определяющий количество изометрических циклов, проходящих по ребрам данного подмножества, и вектор циклов по вершинам $P_v$, определяющий количество изометрических циклов, проходящих по вершинам данного подмножества. Базис линейного подпространства циклов не планарного графа должен удовлетворять следующим требованиям:

1. Значение функционала Маклейна должно стремиться к минимальному значению.
2. В векторе циклов по ребрам не должно быть нулевых элементов.
3. В векторе циклов по ребрам обязательно должны содержаться единичные элементы.
4. В векторе циклов по вершинам не должно быть нулевых элементов.

Естественно, что значение функционала Маклейна для конфигурации непланарного графа больше нуля. Следовательно, для достижения нулевого значения кубического функционала Маклейна нужно удалить часть циклов.

***Пример 7..3.*** Рассмотрим процесс выделения плоской части графа, используя в качестве базиса элемент базы матроида изометрических циклов [23] на примере следующего неориентированного графа.

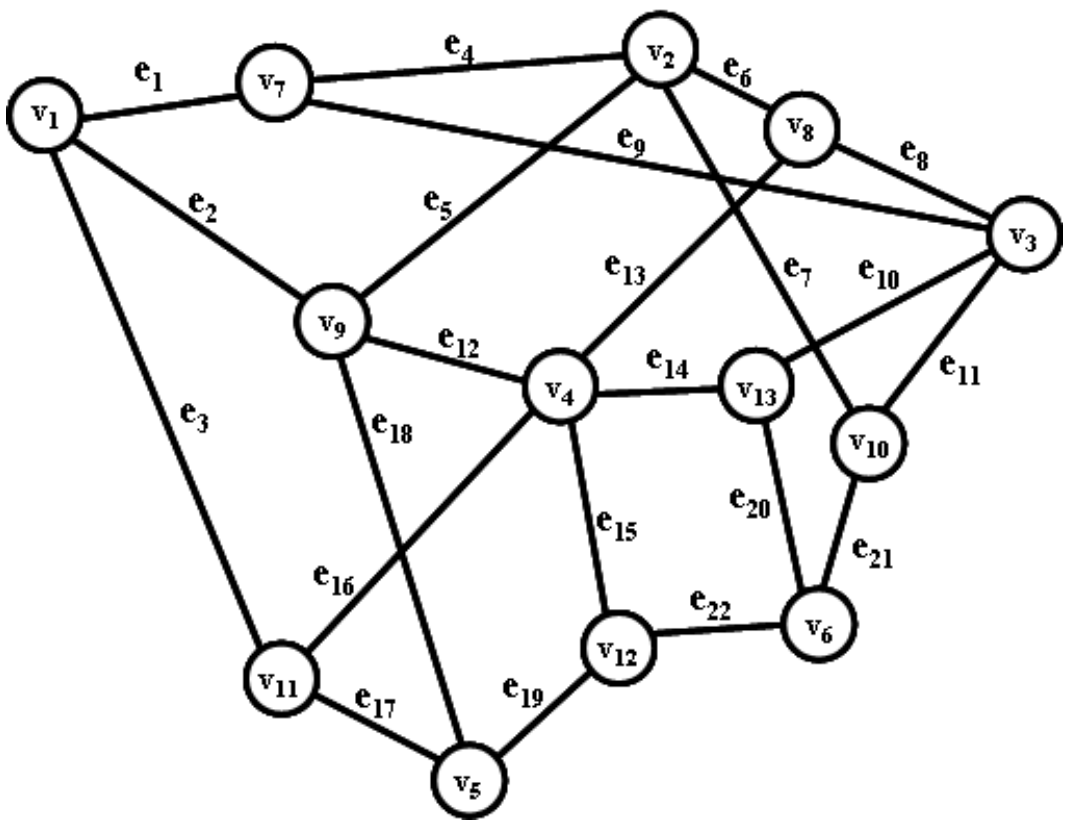

Рис. 7.5. Граф $G_4$.

Множество изометрических циклов графа:

| | Множество изометрических циклов графа в виде ребер: | Множество изометрических циклов графа в виде вершин: |
|---|---|---|
| $c_1$ | $\{e_1,e_2,e_4,e_5\}$; | $\{v_1,v_2,v_{7.},v_9\}$; |
| $c_2$ | $\{e_1,e_3,e_8,e_9,e_{13},e_{16}\}$; | $\{v_1,v_3,v_4,v_7,v_8,v_{11}\}$; |
| $c_3$ | $\{e_1,e_3,e_9,e_{10},e_{14},e_{16}\}$; | $\{v_1,v_3,v_4,v_7,v_{11},v_{13}\}$; |
| $c_4$ | $\{e_2,e_3,e_{12},e_{16}\}$; | $\{v_1,v_4,v_9,v_{11}\}$; |
| $c_5$ | $\{e_2,e_3,e_{17},e_{18}\}$; | $\{v_1,v_5,v_9,v_{11}\}$; |
| $c_6$ | $\{e_4,e_6,e_8,e_9\}$; | $\{v_2,v_3,v_7,v_8\}$; |
| $c_7$ | $\{e_4,e_7,e_9,e_{11}\}$; | $\{v_2,v_3,v_7,v_{10}\}$; |
| $c_8$ | $\{e_5,e_6,e_{12},e_{13}\}$; | $\{v_2,v_4,v_8,v_9\}$; |
| $c_9$ | $\{e_5,e_7,e_{18},e_{19},e_{21},e_{22}\}$; | $\{v_2,v_5,v_6,v_9,v_{10},v_{12}\}$; |
| $c_{10}$ | $\{e_6,e_7,e_8,e_{11}\}$; | $\{v_2,v_3,v_8,v_{10}\}$; |
| $c_{11}$ | $\{e_5,e_7,e_{12},e_{15},e_{21},e_{22}\}$; | $\{v_2,v_4,v_6,v_9,v_{10},v_{12}\}$; |
| $c_{12}$ | $\{e_5,e_7,e_{12},e_{14},e_{20},e_{21}\}$; | $\{v_2,v_4,v_6,v_9,v_{10},v_{13}\}$; |
| $c_{13}$ | $\{e_8,e_{10},e_{13},e_{14}\}$; | $\{v_3,v_4,v_8,v_{13}\}$; |
| $c_{14}$ | $\{e_1,e_2,e_9,e_{10},e_{12},e_{14}\}$; | $\{v_1,v_3,v_4,v_7,v_9,v_{13}\}$; |
| $c_{15}$ | $\{e_{10},e_{11},e_{20},e_{21}\}$; | $\{v_3,v_6,v_{10},v_{13}\}$; |
| $c_{16}$ | $\{e_{12},e_{16},e_{17},e_{18}\}$; | $\{v_4,v_5,v_9,v_{11}\}$; |
| $c_{17}$ | $\{e_{12},e_{15},e_{18},e_{19}\}$; | $\{v_4,v_5,v_9,v_{12}\}$; |
| $c_{18}$ | $\{e_{14},e_{15},e_{20},e_{22}\}$; | $\{v_4,v_6,v_{12},v_{13}\}$; |
| $c_{19}$ | $\{e_{15},e_{16},e_{17},e_{19}\}$; | $\{v_4,v_5,v_{11},v_{12}\}$; |
| $c_{20}$ | $\{e_6,e_7,e_{13},e_{15},e_{21},e_{22}\}$. | $\{v_2,v_4,v_6,v_8,v_{10},v_{12}\}$. |

Циклы матроида $\zeta$ ( $C_\tau$ ) :

1-й цикл матроида $\{c_8,c_6,c_4,c_1,c_2\}$;
2-й цикл матроида $\{c_{10},c_6,c_7\}$;
3-й цикл матроида $\{c_{13},c_3,c_2\}$;
4-й цикл матроида $\{c_{14},c_3,c_4\}$;
5-й цикл матроида $\{c_{15},c_{12},c_3,c_7,c_4,c_1\}$;
6-й цикл матроида $\{c_{16},c_5,c_4\}$;
7-й цикл матроида $\{c_{17},c_{11},c_9\}$;
8-й цикл матроида $\{c_{18},c_{12},c_{11}\}$;

9-й цикл матроида $\{c_{19},c_5,c_{11},c_9,c_4\}$;
10-й цикл матроида $\{c_{20},c_6,c_{11},c_4,c_1,c_2\}$;
11-й цикл матроида $\{c_{14},c_1,c_8,c_6,c_{13}\}$;
12-й цикл матроида $\{c_{15},c_8,c_7,c_6,c_{12},c_{13}\}$;
13-й цикл матроида $\{c_{16},c_5,c_1,c_2,c_8,c_6\}$;
14-й цикл матроида $\{c_{18},c_{17},c_{12},c_9\}$;
15-й цикл матроида $\{c_{19},c_5,c_1,c_2,c_8,c_6,c_{17}\}$;
16-й цикл матроида $\{c_{20},c_8,c_{17},c_9\}$;
17-й цикл матроида $\{c_{16},c_5,c_2,c_{14},c_{13}\}$;
18-й цикл матроида $\{c_{18},c_{17},c_{14},c_1,c_7,c_9,c_{15}\}$;
19-й цикл матроида $\{c_{19},c_5,c_2,c_{17},c_{14},c_{13}\}$;
20-й цикл матроида $\{c_{20},c_6,c_{17},c_{14},c_1,c_9,c_{13}\}$.

Рассмотрим 18-ый цикл матроида изометрических циклов:

$$c_{18} \oplus c_{17} \oplus c_{14} \oplus c_1 \oplus c_7 \oplus c_9 \oplus c_{15} = \{e_{14},e_{15},e_{20},e_{22}\} \oplus \{e_{12},e_{15},e_{18},e_{19}\} \oplus$$
$$\oplus \{e_1,e_2,e_9,e_{10},e_{12},e_{14}\} \oplus \{e_1,e_2,e_4,e_5\} \oplus \{e_4,e_7,e_9,e_{11}\} \oplus \{e_5,e_7,e_{18},e_{19},e_{21},e_{22}\} \oplus$$
$$\oplus \{e_1,e_{11},e_{20},e_{21}\} = \varnothing$$

Если удалить из множества изометрических циклов следующие циклы: $\{c_6,c_9,c_{14},c_{12},c_{16},c_{13}\}$, то получим элемент базы матроида изометрических циклов состоящий из 14 циклов. Следует заметить, что в данном примере, циклы матроида – это плоские конфигурации.

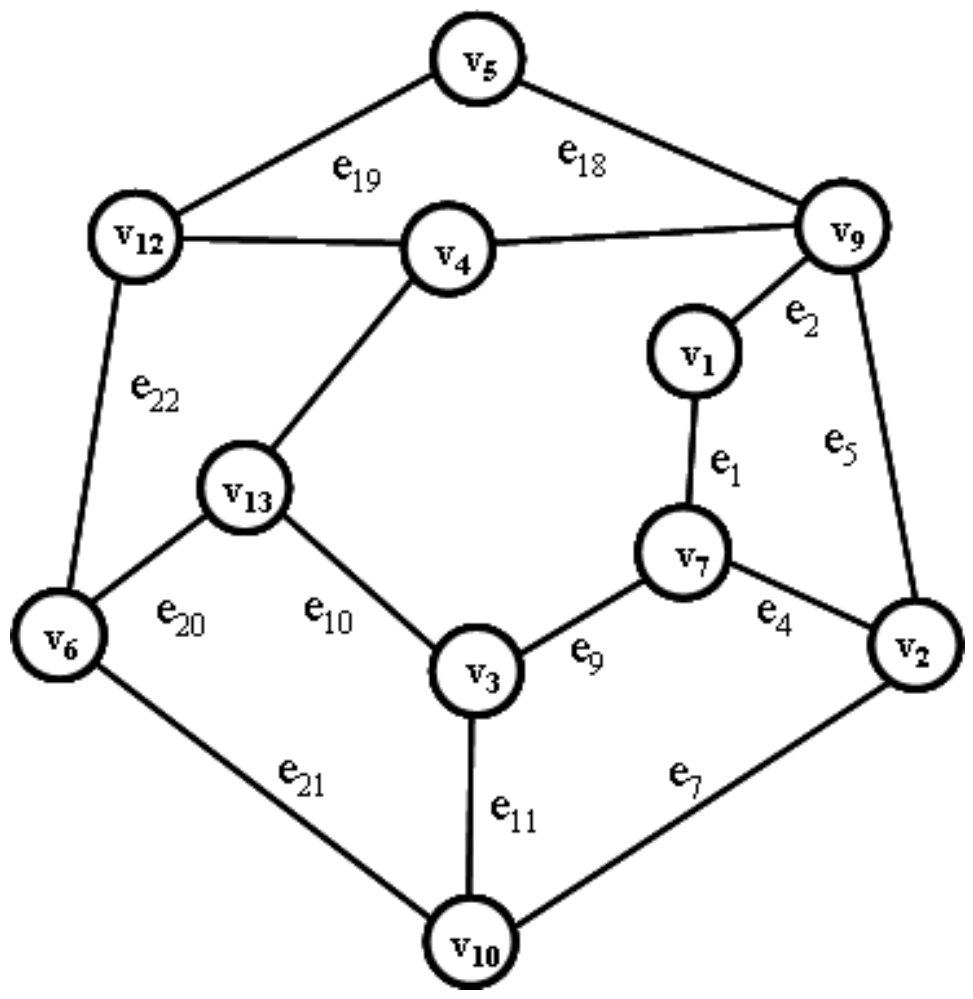

Рис. 7.6. Топологический рисунок 18-го цикла матроида.

Сформируем однострочные структурные числа для множества изометрических циклов. Для выделенного дерева графа $T = \{e_1,e_4,e_5,e_8,e_{10},e_{12},e_{13},e_{17},e_{19},e_{20},e_{21},e_{22}\}$ множество хорд $H = \{e_2,e_3,e_6,e_7,e_9,e_{11},e_{14},e_{15},e_{16},e_{18}\}$. И тогда однострочные структурные числа имеют вид:

по хорде $e_{14}$ проходят циклы: $[c_3,c_{12},c_{13},c_{14},c_{18}]$;
по хорде $e_{18}$ проходят циклы: $[c_5,c_9,c_{16},c_{17}]$;
по хорде $e_2$ проходят циклы: $[c_1,c_4,c_5,c_{14}]$;
по хорде $e_6$ проходят циклы: $[c_6,c_8,c_{10},c_{20}]$;
по хорде $e_9$ проходят циклы: $[c_2,c_3,c_6,c_7,c_{14}]$;
по хорде $e_{11}$ проходят циклы: $[c_7,c_{10},c_{15}]$;
по хорде $e_3$ проходят циклы: $[c_2,c_3,c_4,c_5]$;
по хорде $e_7$ проходят циклы: $[c_7,c_9,c_{10},c_{11},c_{12},c_{20}]$;

по хорде $e_{16}$ проходят циклы: [$c_2$,$c_3$,$c_4$,$c_{16}$,$c_{19}$];
по хорде $e_{15}$ проходят циклы: [$c_{11}$,$c_{17}$,$c_{18}$,$c_{19}$,$c_{20}$].

Построим произведение однострочных структурных чисел для множества изометрических циклов. Длина элемента структурного числа всегда равна количеству хорд графа, в данном случае – десяти.

Алгоритмом «бегущая строка» выделим все множество элементов произведения однострочных структурных чисел и определим их количество.

Элементы произведения однострочных структурных чисел имеют вид:

1-й - элемент структурного числа = {$c_3$,$c_5$,$c_1$,$c_8$,$c_2$,$c_7$,$c_4$,$c_{10}$,$c_{19}$,$c_{11}$};
2-й - элемент структурного числа = {$c_3$,$c_5$,$c_1$,$c_8$,$c_2$,$c_7$,$c_4$,$c_{10}$,$c_{19}$,$c_{17}$};
3-й - элемент структурного числа = {$c_3$,$c_5$,$c_1$,$c_8$,$c_2$,$c_7$,$c_4$,$c_{10}$,$c_{19}$,$c_{18}$};
4-й - элемент структурного числа = {$c_3$,$c_5$,$c_1$,$c_8$,$c_2$,$c_7$,$c_4$,$c_{10}$,$c_{19}$,$c_{20}$};
5-й - элемент структурного числа = {$c_3$,$c_5$,$c_1$,$c_8$,$c_2$,$c_7$,$c_4$,$c_{11}$,$c_{19}$,$c_{17}$};
6-й - элемент структурного числа = {$c_3$,$c_5$,$c_1$,$c_8$,$c_2$,$c_7$,$c_4$,$c_{11}$,$c_{19}$,$c_{18}$};
……………………………………………………………………
3122-й - элемент структурного числа = {$c_{18}$,$c_{17}$,$c_5$,$c_{20}$,$c_7$,$c_{15}$,$c_4$,$c_{10}$,$c_2$,$c_{19}$};
3123-й - элемент структурного числа = {$c_{18}$,$c_{17}$,$c_5$,$c_{20}$,$c_7$,$c_{15}$,$c_4$,$c_{10}$,$c_3$,$c_{11}$};
3124-й - элемент структурного числа = {$c_{18}$,$c_{17}$,$c_5$,$c_{20}$,$c_7$,$c_{15}$,$c_4$,$c_{10}$,$c_3$,$c_{19}$};
3125-й - элемент структурного числа = {$c_{18}$,$c_{17}$,$c_5$,$c_{20}$,$c_7$,$c_{15}$,$c_4$,$c_{10}$,$c_{19}$,$c_{11}$};
3126-й - элемент структурного числа = {$c_{18}$,$c_{17}$,$c_5$,$c_{20}$,$c_7$,$c_{15}$,$c_4$,$c_{11}$,$c_2$,$c_{19}$};
3127-й - элемент структурного числа = {$c_{18}$,$c_{17}$,$c_5$,$c_{20}$,$c_7$,$c_{15}$,$c_4$,$c_{11}$,$c_3$,$c_{19}$}.

Применение метода «бегущая строка» требует проверки на независимость, так как возможно четное число раз повторений данной конфигурации циклов.

***Пример 7.4.*** В качестве примера, выделим следующий элемент произведения однострочных структурных чисел {$c_3$,$c_5$,$c_1$,$c_8$,$c_2$,$c_7$,$c_4$,$c_{10}$,$c_{19}$,$c_{11}$}. Построим множество однострочных структурных чисел, исключив множество циклов, не вошедших в рассматриваемый элемент:

по хорде $e_{14}$ проходят циклы: [$c_3$];
по хорде $e_{18}$ проходят циклы: [$c_5$];
по хорде $e_2$ проходят циклы: [$c_1$,$c_4$,$c_5$];
по хорде $e_6$ проходят циклы: [$c_8$,$c_{10}$];
по хорде $e_9$ проходят циклы: [$c_2$,$c_3$,$c_7$];
по хорде $e_{11}$ проходят циклы: [$c_7$,$c_{10}$];
по хорде $e_3$ проходят циклы: [$c_2$,$c_3$,$c_4$,$c_5$];
по хорде $e_7$ проходят циклы: [$c_7$,$c_{10}$,$c_{11}$];
по хорде $e_{16}$ проходят циклы: [$c_2$,$c_3$,$c_4$,$c_{19}$];
по хорде $e_{15}$ проходят циклы: [$c_{11}$,$c_{19}$].

В произведении однострочных структурных чисел данная комбинация встречается четыре раза. Следовательно, данная конфигурация циклов линейно зависима:

1 - элемент структурного числа = {$c_3$,$c_5$,$c_8$,$c_7$,$c_{11}$,$c_1$,$c_{10}$,$c_2$,$c_4$,$c_{19}$};
2 - элемент структурного числа = {$c_3$,$c_5$,$c_8$,$c_{10}$,$c_{11}$,$c_1$,$c_7$,$c_2$,$c_4$,$c_{19}$};
3 - элемент структурного числа = {$c_3$,$c_5$,$c_8$,$c_{10}$,$c_{19}$,$c_1$,$c_{11}$,$c_7$,$c_2$,$c_4$};
4 - элемент структурного числа = {$c_3$,$c_5$,$c_8$,$c_{10}$,$c_{19}$,$c_1$,$c_{11}$,$c_7$,$c_4$,$c_2$}.

Мы рассмотрели три метода выделения базисов подпространства циклов, состоящих из изометрических циклов при подходе сверху.

1. Базис подпространства циклов $C_b$, состоящий из изометрических циклов может быть выделен из случайной последовательности расположения элементов в кортеже изометрических циклов. Недостаток метода состоит в возможности получения не оптимального значения функционала Маклейна для выделенного базиса.

2. Недостаток метода наискорейшего спуска заключается в возможности получения зависимой системы циклов, состоящей из множества циклов с мощностью равной цикломатическому числу графа и с минимальным значением функционала Маклейна, при случайном расположении элементов в кортеже изометрических циклов.

3. Базис линейного подпространства циклом может быть выбран алгоритмом «бегущая строка» из структурного числа изометрических циклов. Однако выделенный элемент структурного числа может состоять из зависимой системы циклов, так как может повторяться в структурном числе четное число раз. Кроме того, значение функционала Маклейна может быть не оптимальным.

Имея несколько методов выделения базиса подпространства циклов, возможно создание комплексных подходов к решению поставленной задачи [25,27,28]. Очевидно, что при создании практических систем приближенного решения для выделения базисов подпространства циклов можно воспользоваться методом Монте-Карло. То есть произвести случайным образом расположение элементов в кортежах или в однострочных структурных числах.

В свою очередь, подмножеству циклов можно поставить в соответствие два вектора. Вектор циклов по ребрам $P_e$, определяющий количество изометрических циклов, проходящих по ребрам данного подмножества, и вектор циклов по вершинам $P_v$, определяющий количество изометрических циклов, проходящих по вершинам данного подмножества. Базис линейного подпространства циклов не планарного графа должен удовлетворять следующим требованиям:

- значение функционала Маклейна должно стремиться к минимальному значению.
- в векторе циклов по ребрам не должно быть нулевых элементов.
- в векторе циклов по ребрам обязательно должны содержаться единичные элементы.
- в векторе циклов по вершинам не должно быть нулевых элементов.

## 7.4. Кубический функционал Маклейна

Для решения задачи построения топологического рисунка максимально плоского суграфа непланарного графа, необходимо описывать процесс удаления цикла с одновременным удалением ребра из выделенного базиса циклов с минимальным значением функционала Маклейна для оставшихся циклов. Применение правила удаления цикла исходит из удовлетворения уравнения Эйлера. Обозначим цикломатическое число графа как $\nu(\mathbf{G})$. Тогда формула Эйлера для подмножества независимых циклов графа имеет вид:

$$m - n - \nu(\mathbf{G}) + \mathbf{1} = \mathbf{0}. \tag{7.8}$$

Здесь $m$ - количество ребер графа, а $n$ - количество вершин графа.

Если мы удаляем из графа ребро, то для сохранения выполнения уравнения Эйлера должны изменить и значение цикломатического числа нового суграфа на единицу:

$$(m-1) - n - (\nu(\mathbf{G}) - \mathbf{1}) + \mathbf{1} = m - n - \nu(\mathbf{G}) + 1 = \mathbf{0}. \tag{7.9}$$

Если мы удаляем из графа два ребра, то для сохранения выполнения уравнения Эйлера должны изменить значение цикломатического числа нового суграфа на единицу и уменьшить количество компонент связности графа тоже на единицу.

$$(m-2) - (n-1) - (\nu(\mathbf{G}) - \mathbf{1}) + \mathbf{1} = m - n - \nu(\mathbf{G}) + \mathbf{1} = \mathbf{0} \tag{7.10}$$

Для выполнения этих целей, воспользуемся понятием дифференцирования структурного числа. Напомним, что *дифференцирование структурного числа* называется структурное число $\partial C/\partial \alpha$, равное:

$$\frac{\partial C}{\partial \alpha} = C \,\Big|\, \text{столбцы, не содержащие элемент } \alpha\text{, исключены.} \tag{7.11}$$

Очевидно также, что при удалении цикла и одновременного удаления ребра должно уменьшиться и цикломатическое число графа на единицу, при условии, что оставшаяся система циклов покрывает все множество вершин графа.

В этом случае функционал Маклейна уже не будет отвечать процессу построения топологического рисунка графа, так как он строился в предположении отсутствия нулевых значений для вектора $P_e$. Для более точного описания процесса удаления ребер из базиса, введем понятие кубического функционал Маклейна:

$$\mathrm{FP(C)} = \sum_{\mathrm{i}=1}^{\mathrm{m}} a_{\mathrm{i}}(a_{\mathrm{i}}-1)(a_{\mathrm{i}}-2) = \sum_{i=1}^{m} a_{\mathrm{i}}^3 - 3\sum_{i=1}^{m} a_{\mathrm{i}}^2 + 2\sum_{i=1}^{m} a_{\mathrm{i}}, \tag{7.12}$$

где $a_{\mathrm{i}}$ – количество изометрических циклов проходящих по ребру $i$ для заданного подмножества циклов.

Данный функционал подобен функционалу Маклейна, но позволяет учитывать нулевые значения в векторе $P_e$ в процессе удаления ребер из базиса.

Процесс выделения плоской части графа из базиса может быть описан как процедура дифференцирования элемента структурного числа:

$$C_b' = \frac{\partial^q C_b}{\partial c_1 \partial c_2 \partial c_3 ... \partial c_{q-1} \partial c_q}, \quad (7.13)$$

где $q$ – количество удаленных ребер. Структурное число $\mathbf{C_b}$ описывает базис изометрических циклов.

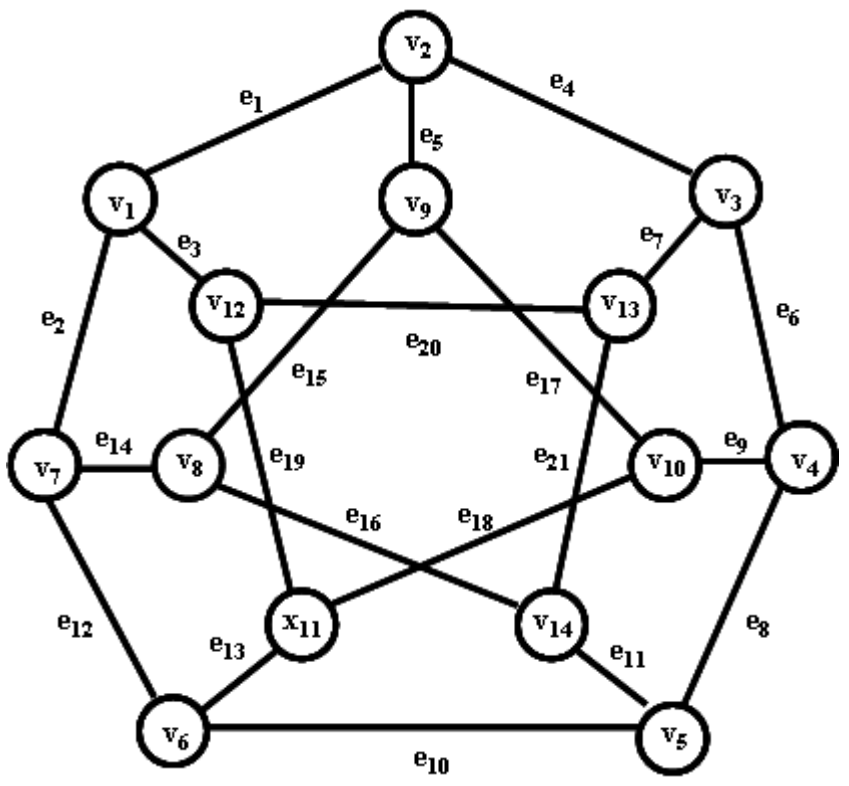

Рис. 7.7. Граф $G_5$.

Рассмотрим пример решения задачи выделение плоской части не планарного графа.

***Пример 7..5.*** Выделить максимально плоский суграф для графа, представленного на рис. 7..7 и продемонстрировать применение кубического функционала Маклейна.

Множество $C_\tau$ изометрических циклов графа $G_5$ представленного на рис. 7..7 записывается в виде:

$c_1 = \{e_1,e_2,e_5,e_{14},e_{15}\} \to \{v_1,v_2,v_7,v_8,v_9\}$;
$c_2 = \{e_1,e_3,e_4,e_7,e_{20}\} \to \{v_1,v_2,v_3,v_{12},v_{13}\}$;
$c_3 = \{e_2,e_3,e_{12},e_{13},e_{19}\} \to \{v_1,v_6,v_7,v_{11},v_{12}\}$;
$c_4 = \{e_4,e_5,e_6,e_9,e_{17}\} \to \{v_2,v_3,v_4,v_9,v_{10}\}$;
$c_5 = \{e_6,e_7,e_8,e_{11},e_{21}\} \to \{v_3,v_4,v_5,v_{13},v_{14}\}$;
$c_6 = \{e_8,e_9,e_{10},e_{13},e_{18}\} \to \{v_4,v_5,v_6,v_{10},v_{11}\}$;
$c_7 = \{e_{10},e_{11},e_{12},e_{14},e_{16}\} \to \{v_5,v_6,v_7,v_8,v_{14}\}$;
$c_8 = \{e_8,e_9,e_{11},e_{15},e_{16},e_{17}\} \to \{v_4,v_5,v_8,v_9,v_{10},v_{14}\}$;
$c_9 = \{e_4,e_5,e_7,e_{15},e_{16},e_{21}\} \to \{v_2,v_3,v_8,v_9,v_{13},v_{14}\}$;
$c_{10} = \{e_{12},e_{13},e_{14},e_{15},e_{17},e_{18}\} \to \{v_6,v_7,v_8,v_9,v_{10},v_{11}\}$;
$c_{11} = \{e_1,e_3,e_5,e_{17},e_{18},e_{19}\} \to \{v_1,v_2,v_9,v_{10},v_{11},v_{12}\}$;
$c_{12} = \{e_6,e_7,e_9,e_{18},e_{19},e_{20}\} \to \{v_3,v_4,v_{10},v_{11},v_{12},v_{13}\}$;
$c_{13} = \{e_{10},e_{11},e_{13},e_{19},e_{20},e_{21}\} \to \{v_5,v_6,v_{11},v_{12},v_{13},v_{14}\}$;
$c_{14} = \{e_2,e_3,e_{14},e_{16},e_{20},e_{21}\} \to \{v_1,v_7,v_8,v_{12},v_{13},v_{14}\}$.

Значение квадратичного функционала Маклейна для множества изометрических циклов равно: $F(G_\tau) = 100$.

Множество изометрических циклов графа $G_5$ позволяет построить вектор циклов $P_e$ - определяющий количество изометрических циклов, обозначаемых латинской буквой s, проходящих по ребру, где порядок записи в кортеже соответствует номеру ребра:

$P_e = \langle a_1,a_2,a_3,a_4,a_5,a_6,a_7,a_8,a_9,a_{10},a_{11},a_{12},a_{13},a_{14},a_{15},a_{16},a_{17},a_{18},a_{19},a_{20},a_{21}\rangle =$

= < 3,3,4,3,4,3,5,3,4,3,4,3,4,4,4,4,3,4,4,4,4 >.

Выделим базис подпространства циклов, состоящий из подмножества циклов $C_b = \{c_1,c_2,c_3,c_4,c_5,c_6,c_7,c_8\}$. Для выделенного базиса определим значение кубического функционала Маклейна: $FP(C_b) = 18$. Для удаления цикла воспользуемся понятием дифференцирования структурного числа:

$$FP(\frac{\partial C_b}{\partial c_1}) = 18;\quad FP(\frac{\partial C_b}{\partial c_2}) = 18;\quad FP(\frac{\partial C_b}{\partial c_3}) = 18;\quad FP(\frac{\partial C_b}{\partial c_4}) = 12;$$
$$FP(\frac{\partial C_b}{\partial c_5}) = 6;\quad FP(\frac{\partial C_b}{\partial c_6}) = 6;\quad FP(\frac{\partial C_b}{\partial c_7}) = 12;\quad FP(\frac{\partial C_b}{\partial c_8}) = 0.$$

И хотя при удалении $c_8$ достигается минимальное нулевое значение кубического функционала Маклейна, удаление ребра не производится, так как не удовлетворяется условие: при удалении цикла – удаляется одно и только одно ребро. При удалении $c_6$ – удаляется ребро $e_{18}$, а при удалении $c_5$ - удаляется ребро $e_{21}$. В этих случаях, значение кубического функционала Маклейна меньше из всех кандидатов на удаление и равно 5. Выбираем для удаления – $c_5$. Продолжаем процесс удаления циклов с одновременным удалением одного и только одного ребра:

$$FP(\frac{\partial^2 C_b}{\partial c_6 \partial c_1}) = 6;\quad FP(\frac{\partial^2 C_b}{\partial c_6 \partial c_2}) = 6;\quad FP(\frac{\partial^2 C_b}{\partial c_6 \partial c_3}) = 6;\quad FP(\frac{\partial^2 C_b}{\partial c_6 \partial c_4}) = 6;$$
$$FP(\frac{\partial^2 C_b}{\partial c_6 \partial c_5}) = 0;\quad FP(\frac{\partial^2 C_b}{\partial c_6 \partial c_7}) = 0;\quad FP(\frac{\partial^2 C_b}{\partial c_6 \partial c_8}) = 0.$$

И снова, при удалении цикла $c_8$ получаем нулевое значение кубического функционала Маклейна, но удаление ребра не может быть произведено, так как не выполняется условие совместного удаления ребра и цикла. При удалении $c_7$ – удаляется ребро $e_{10}$, а при удалении $c_5$ – удаляется ребро $e_{21}$. В последних случаях, значение кубического функционала Маклейна равно 0 и, при удалении его из базиса, удаляется одно и только одно ребро. Выбираем для удаления – $c_7$.

В результате мы получили плоский суграф (см. рис. 7..8). Шесть циклов $c_1,c_2,c_3,c_4,c_5,c_8$ и обод образуют плоский суграф из непланарного графа при исключении ребер $e_{18}$ и $e_{10}$.

Полученное множество циклов индуцирует вращение вершин плоского суграфа, тем самым, описывая следующей диаграммой его топологический рисунок [32]:

вращение вершины $\sigma_1$: $v_2$ $v_{12}$ $v_7$
вращение вершины $\sigma_2$: $v_9$ $v_3$ $v_1$
вращение вершины $\sigma_3$: $v_4$ $v_{13}$ $v_2$
вращение вершины $\sigma_4$: $v_5$ $v_3$ $v_{10}$
вращение вершины $\sigma_5$: $v_4$ $v_{14}$
вращение вершины $\sigma_6$: $v_7$ $v_{11}$
вращение вершины $\sigma_7$: $v_8$ $v_1$ $v_6$
вращение вершины $\sigma_8$: $v_9$ $v_7$ $v_{14}$

вращение вершины $\sigma_9$: $v_8$ $v_{10}$ $v_2$
вращение вершины $\sigma_{10}$: $v_9$ $v_4$
вращение вершины $\sigma_{11}$: $v_6$ $v_{12}$
вращение вершины $\sigma_{12}$: $v_1$ $v_{13}$ $v_{11}$
вращение вершины $\sigma_{13}$: $v_3$ $v_{14}$ $v_{12}$
вращение вершины $\sigma_{14}$: $v_8$ $v_{13}$ $v_5$.

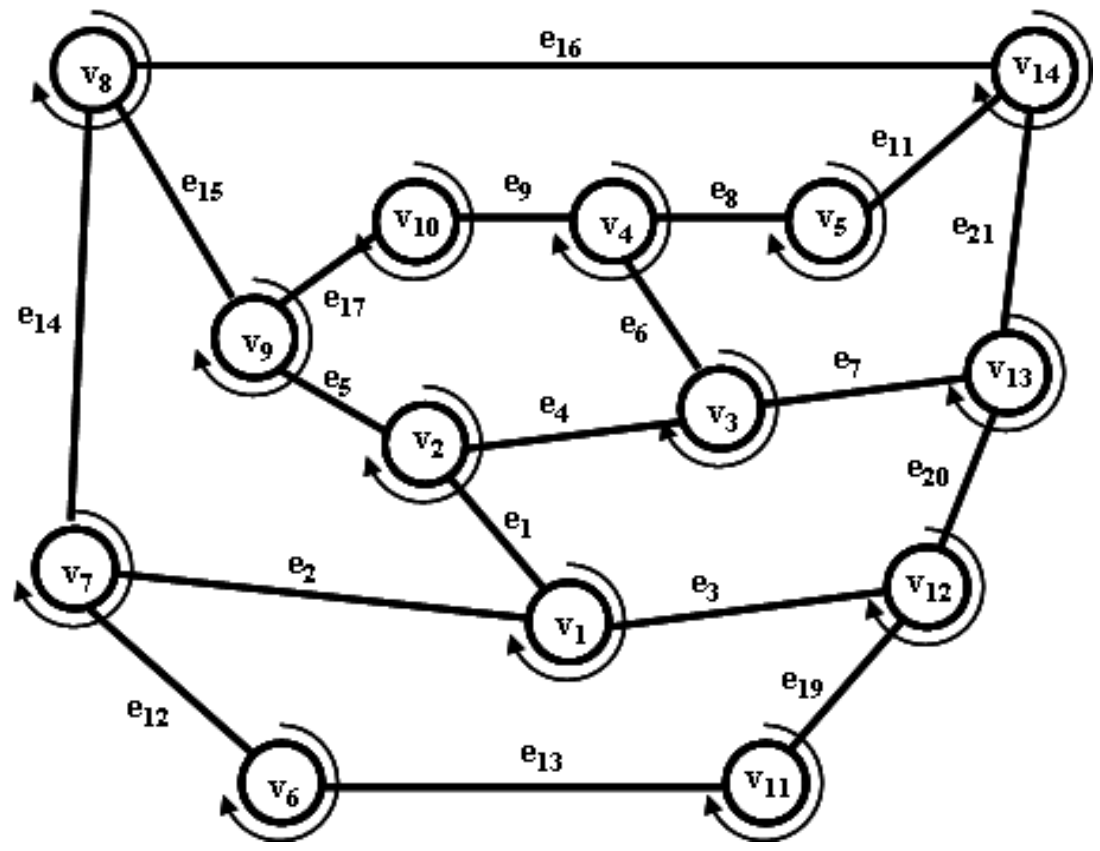

Рис. 7.8. Плоский суграф графа $G_5$.

Рассмотрим следующий пример.

***Пример 7..5.*** Выделить максимально плоский суграф графа представленного на рис. 7..9 и сформировать его топологический рисунок.

Множество изометрических циклов графа $G_6$ имеет вид:

$c_1 = \{e_1,e_2,e_6\}$; $c_2 = \{e_1,e_3,e_7\}$; $c_3 = \{e_1,e_4,e_8,e_{19}\}$;
$c_4 = \{e_2,e_3,e_{10}\}$; $c_5 = \{e_2,e_5,e_{11},e_{20}\}$; $c_6 = \{e_3,e_4,e_{15}\}$;
$c_7 = \{e_3,e_5,e_{17}\}$; $c_8 = \{e_4,e_5,e_{18}\}$; $c_9 = \{e_6,e_7,e_{10}\}$;
$c_{10} = \{e_6,e_8,e_9,e_{14}\}$; $c_{11} = \{e_7,e_8,e_{12},e_{14}\}$; $c_{12} = \{e_7,e_8,e_{15},e_{19}\}$;
$c_{13} = \{e_9,e_{10},e_{12}\}$; $c_{14} = \{e_9,e_{11},e_{13}\}$; $c_{15} = \{e_{10},e_{11},e_{16}\}$;
$c_{16} = \{e_{12},e_{13},e_{16}\}$; $c_{17} = \{e_{12},e_{14},e_{15},e_{19}\}$; $c_{18} = \{e_{15},e_{17},e_{18}\}$;
$c_{19} = \{e_{16},e_{17},e_{20}\}$.

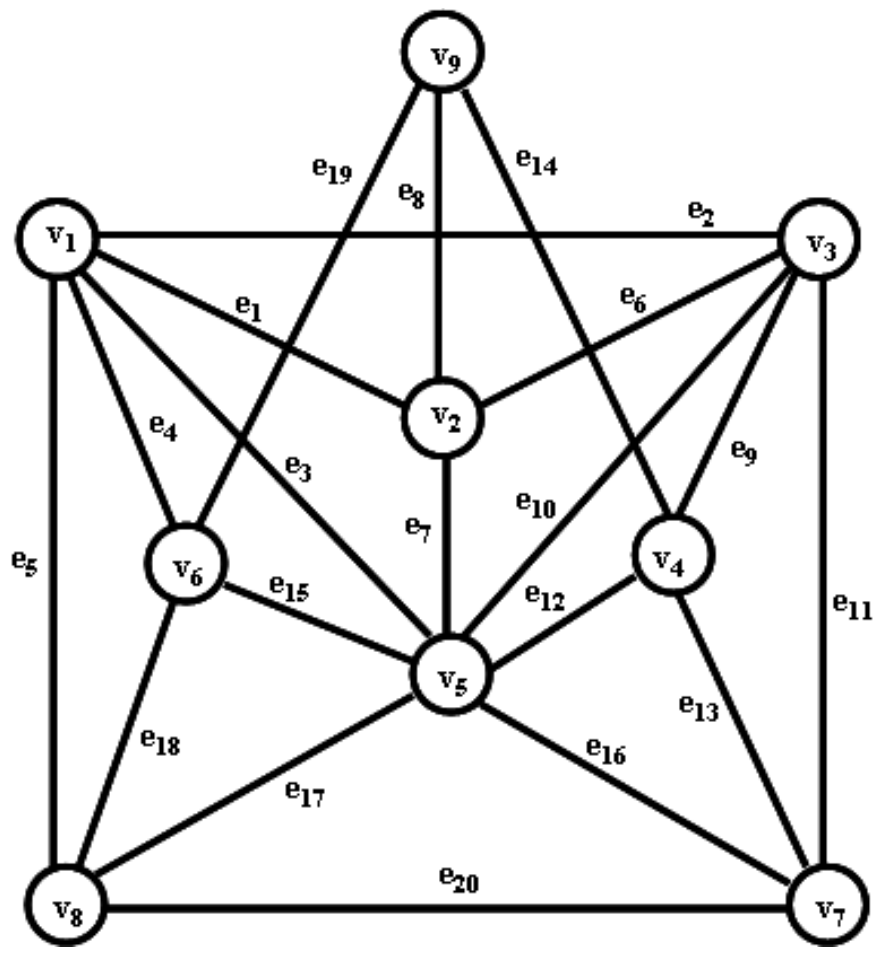

Рис. 7.9. Граф $G_5$.

Методом наискорейшего спуска выделим базис, состоящий из следующих циклов $\mathbf{C}_b$ = $\{c_2,c_3,c_5,c_8,c_9,c_{10},c_{14},c_{15},c_{16},c_{17},c_{18},c_{19}\}$. Значение кубического функционала Маклейна для данного независимого базиса равно 12.

Для выделения максимально плоского суграфа, воспользуемся понятием дифференцирования структурного числа:

$$FP(\frac{\partial C_b}{\partial c_2}) = 12; \quad FP(\frac{\partial C_b}{\partial c_3}) = 12; \quad FP(\frac{\partial C_b}{\partial c_5}) = 6; \quad FP(\frac{\partial C_b}{\partial c_8}) = 12;$$

$$FP(\frac{\partial C_b}{\partial c_9}) = 12; \quad FP(\frac{\partial C_b}{\partial c_{10}}) = 12; \quad FP(\frac{\partial C_b}{\partial c_{14}}) = 6; \quad FP(\frac{\partial C_b}{\partial c_{15}}) = 0;$$

$$FP(\frac{\partial C_b}{\partial c_{16}}) = 6; \quad FP(\frac{\partial C_b}{\partial c_{17}}) = 12; \quad FP(\frac{\partial C_b}{\partial c_{18}}) = 12; \quad FP(\frac{\partial C_b}{\partial c_{19}}) = 5.$$

Удалить цикл $c_{15}$ нельзя, так как в векторе $P_e$ для его ребер не содержится единичных значений. А вот в цикле $c_5$ содержится ребро $e_2$ со значением равным единице в векторе $P_e$. Удаляем цикл $c_5$ и ребро $e_2$, продолжаем процесс дальше:

$$FP(\frac{\partial^2 C_b}{\partial c_5 \partial c_2}) = 6; \quad FP(\frac{\partial^2 C_b}{\partial c_5 \partial c_3}) = 6; \quad FP(\frac{\partial^2 C_b}{\partial c_5 \partial c_8}) = 6; \quad FP(\frac{\partial^2 C_b}{\partial c_5 \partial c_9}) = 6;$$

$$FP(\frac{\partial^2 C_b}{\partial c_5 \partial c_{10}}) = 6; \quad FP(\frac{\partial^2 C_b}{\partial c_5 \partial c_{14}}) = 6; \quad FP(\frac{\partial^2 C_b}{\partial c_5 \partial c_{15}}) = 0; \quad FP(\frac{\partial^2 C_b}{\partial c_5 \partial c_{16}}) = 0;$$

$$FP(\frac{\partial^2 C_b}{\partial c_5 \partial c_{17}}) = 6; \quad FP(\frac{\partial^2 C_b}{\partial c_5 \partial c_{18}}) = 6; \quad FP(\frac{\partial^2 C_b}{\partial c_5 \partial c_{19}}) = 0.$$

Среди циклов, удаление которых приводит к минимальному значению кубического функционала Маклейна, только цикл $c_{19}$ содержит единицу для ребра $e_{20}$ в векторе $P_e$. Удаляем цикл $c_{19}$ и ребро $e_{20}$, получаем максимально плоский суграф (см. рис. 7.10).

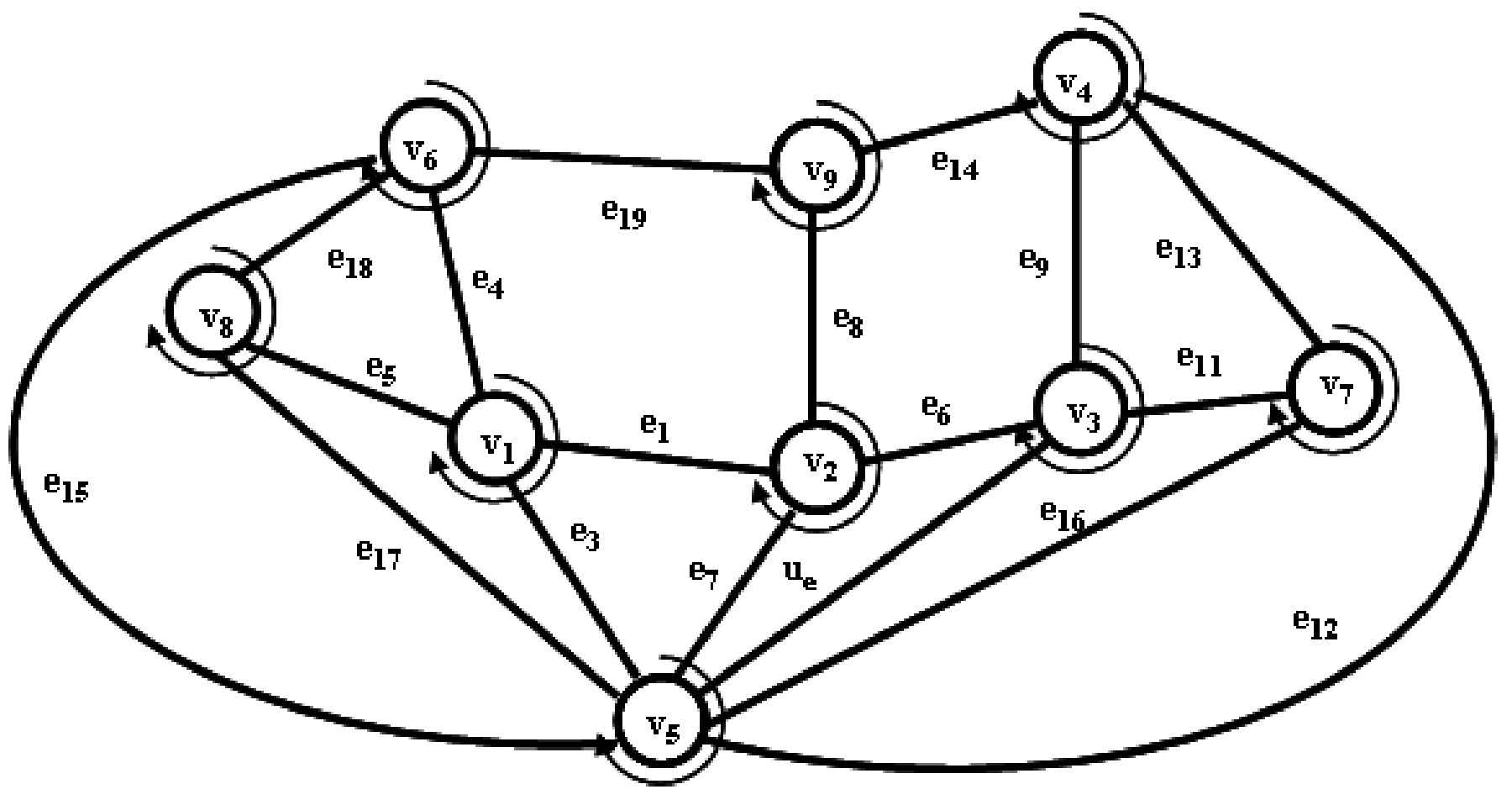

Рис. 7.10. Максимально плоский суграф графа $G_5$.

Произвольным образом выделим базис, состоящий из следующих циклов $\mathbf{C_b} = \{c_1,c_2,c_3,c_6,c_8,c_9,c_{11},c_{13},c_{14},c_{16},c_{18},c_{19}\}$. Значение кубического функционала Маклейна для данного базиса равно 24.

Для выделения максимально плоского суграфа, воспользуемся понятием дифференцирования структурного числа:

$$\mathrm{FP}(\frac{\partial C_b}{\partial c_1}) = 18;\quad \mathrm{FP}(\frac{\partial C_b}{\partial c_2}) = 12;\quad \mathrm{FP}(\frac{\partial C_b}{\partial c_3}) = 12;\quad \mathrm{FP}(\frac{\partial C_b}{\partial c_6}) = 18;$$

$$\mathrm{FP}(\frac{\partial C_b}{\partial c_8}) = 18;\quad \mathrm{FP}(\frac{\partial C_b}{\partial c_9}) = 18;\quad \mathrm{FP}(\frac{\partial C_b}{\partial c_{11}}) = 12;\quad \mathrm{FP}(\frac{\partial C_b}{\partial c_{13}}) = 18;$$

$$\mathrm{FP}(\frac{\partial C_b}{\partial c_{14}}) = 24;\quad \mathrm{FP}(\frac{\partial C_b}{\partial c_{16}}) = 18;\quad \mathrm{FP}(\frac{\partial C_b}{\partial c_{18}}) = 24;\quad \mathrm{FP}(\frac{\partial C_b}{\partial c_{19}}) = 24.$$

Удаляем цикл $c_3$, содержащий ребро $e_{19}$ со значением равным единице в векторе $\mathbf{P}_u$. Продолжаем процесс дальше:

$$\mathrm{FP}(\frac{\partial^2 C_b}{\partial c_3 \partial c_1}) = 12;\quad \mathrm{FP}(\frac{\partial^2 C_b}{\partial c_3 \partial c_2}) = 6;\quad \mathrm{FP}(\frac{\partial^2 C_b}{\partial c_3 \partial c_6}) = 12;\quad \mathrm{FP}(\frac{\partial^2 C_b}{\partial c_3 \partial c_8}) = 12;$$

$$\mathrm{FP}(\frac{\partial^2 C_b}{\partial c_3 \partial c_9}) = 6;\quad \mathrm{FP}(\frac{\partial^2 C_b}{\partial c_3 \partial c_{11}}) = 0;\quad \mathrm{FP}(\frac{\partial^2 C_b}{\partial c_3 \partial c_{13}}) = 6;\quad \mathrm{FP}(\frac{\partial^2 C_b}{\partial c_3 \partial c_{14}}) = 12;$$

$$\mathrm{FP}(\frac{\partial^2 C_b}{\partial c_3 \partial c_{16}}) = 6;\quad \mathrm{FP}(\frac{\partial^2 C_b}{\partial c_3 \partial c_{18}}) = 12;\quad \mathrm{FP}(\frac{\partial^2 C_b}{\partial c_3 \partial c_{19}}) = 12.$$

Среди циклов, содержится цикл $c_{11}$ удаление, которого приводит к нулевому значению кубического функционала Маклейна, но содержащий несколько единиц для ребер $e_8$ и $e_{14}$ в векторе $P_e$. Удаление этого цикла приводит к одновременному удалению двух ребер и вершины $v_9$. В результате получаем максимально плоский суграф (см. рис. 7.11).

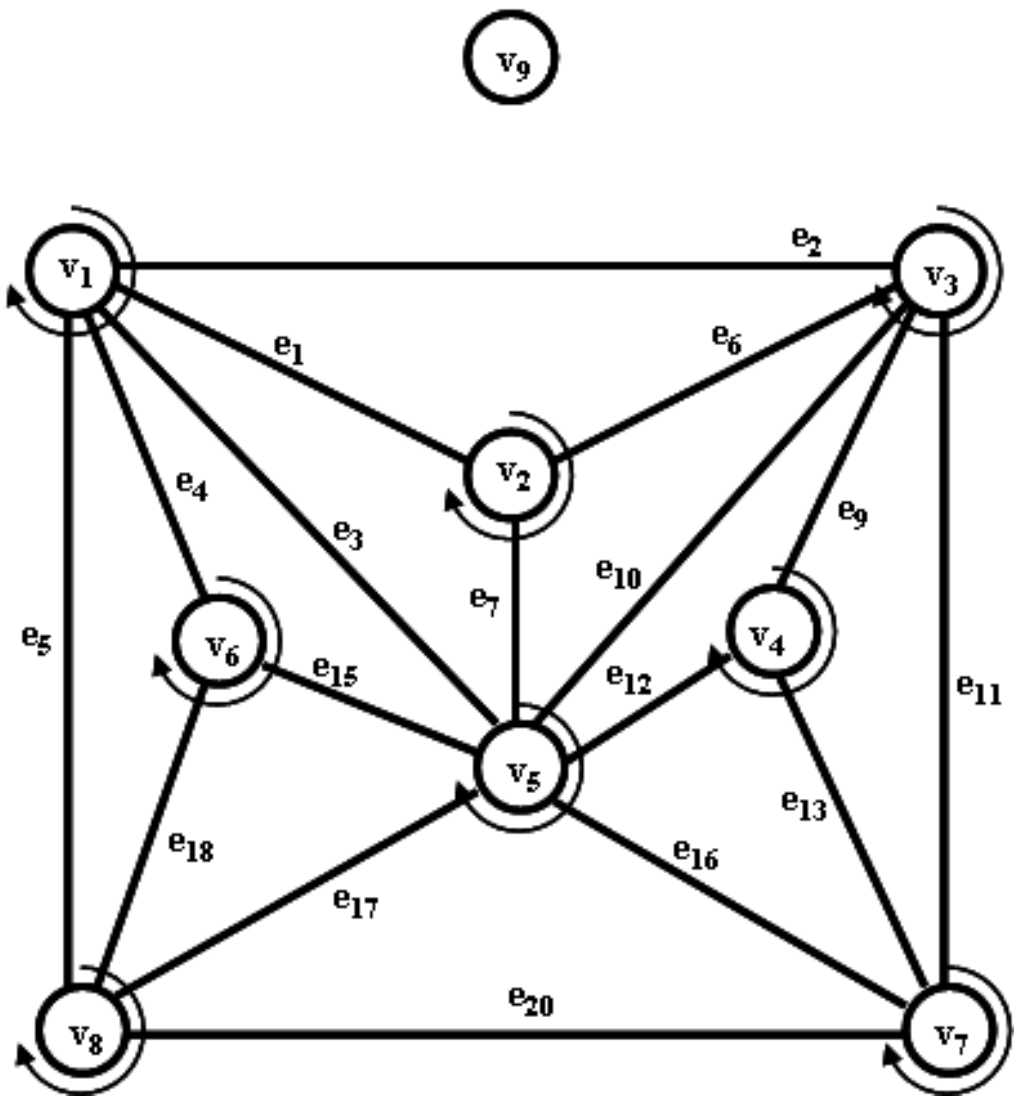

Рис. 7.11. Плоский суграф с тремя удаленными ребрами.

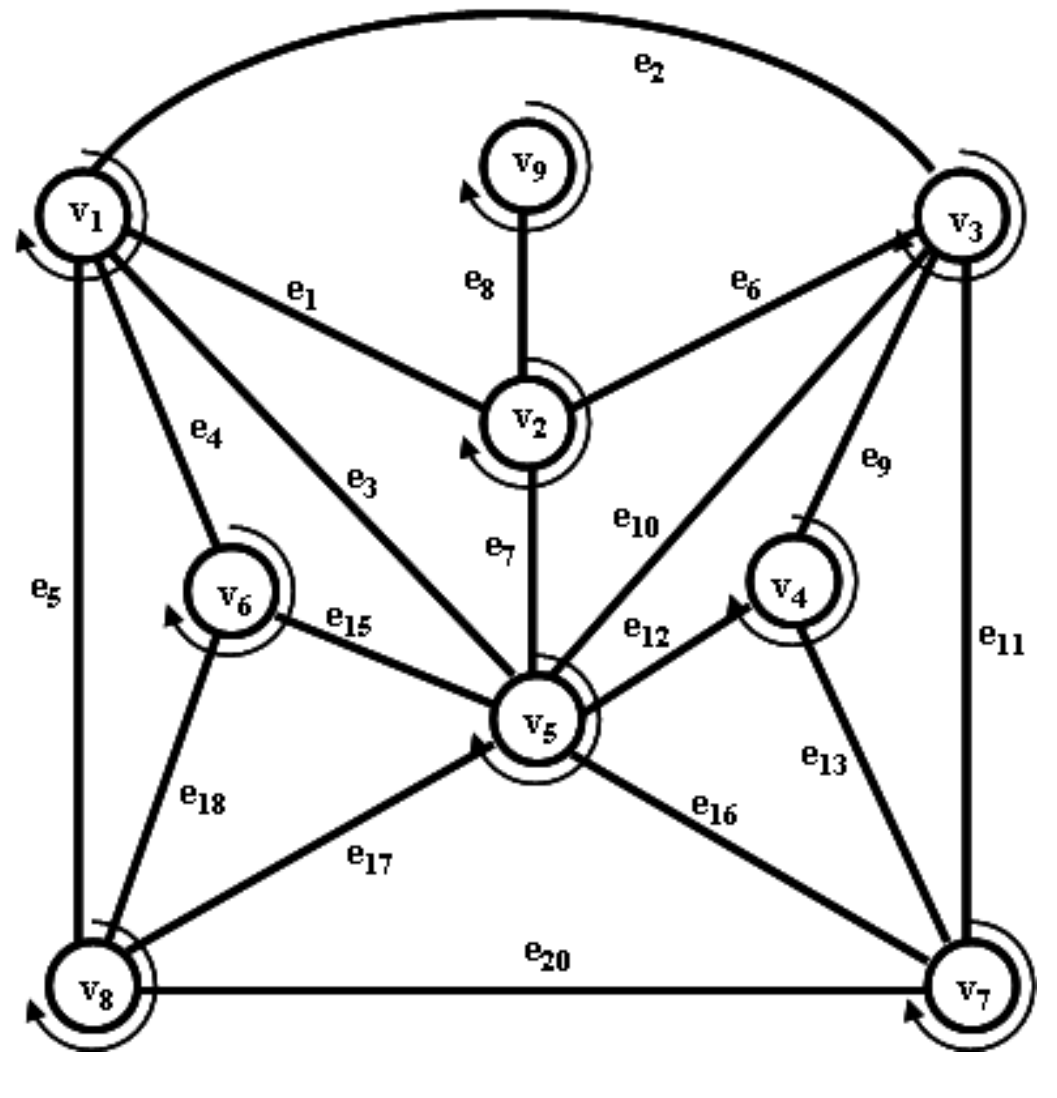

Рис. 7.12. Плоский подграф с двумя удаленными ребрами.

В случае нашего графа, получается, что у нас удалено три ребра и одна вершина. Однако всегда можно найти цикл и расположить в нем удаленную вершину соединив непересекающимся ребром смежную к ней вершину (см. рис. 7..12). Тогда количество удаленных ребер в подграфе будет равно двум для всех накрывающих вершин.

Рассматривая пример 7..6, можно заметить, что при удалении цикла возможен случай удаление двух ребер и вершины. Таким образом, возможно удаление нескольких ребер при удалении одного цикла, не забывая при этом, что одновременно происходит удаление вершин согласно выражению Эйлера.

### 7.5. Добавление простых циклов к топологическому рисунку

Третьей этап осуществляет подключение дополнительных ребер к топологическому рисунку плоской части непланарного графа, созданного двумя первыми этапами. Данный этап характерен применением методов векторной алгебры пересечений для построения топологического рисунка плоской части непланарного графа [29-31].

***Пример 7..7.*** Третьй этап рассмотрим на примере построения следующего графа $G_7$, заданного инцидентором P.

Элементы матрицы инциденций:

ребро 1: $(v_1,v_{19})$ или $(v_{19},v_1)$; ребро 2: $(v_1,v_{28})$ или $(v_{28},v_1)$;
ребро 3: $(v_1,v_{30})$ или $(v_{30},v_1)$; ребро 4: $(v_1,v_{35})$ или $(v_{35},v_1)$;
ребро 5: $(v_1,v_{44})$ или $(v_{44},v_1)$; ребро 6: $(v_2,v_{23})$ или $(v_{23},v_2)$;
ребро 7: $(v_2,v_{24})$ или $(v_{24},v_2)$; ребро 8: $(v_2,v_{26})$ или $(v_{26},v_2)$;
ребро 9: $(v_2,v_{33})$ или $(v_{33},v_2)$; ребро 10: $(v_2,v_{41})$ или $(v_{41},v_2)$;
ребро 11: $(v_3,v_{18})$ или $(v_{18},v_3)$; ребро 12: $(v_3,v_{19})$ или $(v_{19},v_3)$;
ребро 13: $(v_3,v_{23})$ или $(v_{23},v_3)$; ребро 14: $(v_3,v_{25})$ или $(v_{25},v_3)$;
ребро 15: $(v_3,v_{42})$ или $(v_{42},v_3)$; ребро 16: $(v_3,v_{45})$ или $(v_{45},v_3)$;
ребро 17: $(v_4,v_{18})$ или $(v_{18},v_4)$; ребро 18: $(v_4,v_{24})$ или $(v_{24},v_4)$;
ребро 19: $(v_4,v_{30})$ или $(v_{30},v_4)$; ребро 20: $(v_5,v_{20})$ или $(v_{20},v_5)$;
ребро 21: $(v_5,v_{21})$ или $(v_{21},v_5)$; ребро 22: $(v_5,v_{25})$ или $(v_{25},v_5)$;
ребро 23: $(v_5,v_{26})$ или $(v_{26},v_5)$; ребро 24: $(v_6,v_{13})$ или $(v_{13},v_6)$;
ребро 25: $(v_6,v_{18})$ или $(v_{18},v_6)$; ребро 26: $(v_6,v_{21})$ или $(v_{24},v_6)$;
ребро 27: $(v_6,v_{24})$ или $(v_{24},v_6)$; ребро 28: $(v_7,v_{27})$ или $(v_{27},v_7)$;
ребро 29: $(v_7,v_{30})$ или $(v_{30},v_7)$; ребро 30: $(v_7,v_{31})$ или $(v_{31},v_7)$;
ребро 31: $(v_8,v_{29})$ или $(v_{29},v_8)$; ребро 32: $(v_8,v_{32})$ или $(v_{32},v_8)$;
ребро 33: $(v_8,v_{35})$ или $(v_{35},v_8)$; ребро 34: $(v_9,v_{11})$ или $(v_{11},v_9)$;
ребро 35: $(v_9,v_{26})$ или $(v_{26},v_9)$; ребро 36: $(v_9,v_{27})$ или $(v_{27},v_9)$;
ребро 37: $(v_{10},v_{31})$ или $(v_{31},v_{10})$; ребро 38: $(v_{10},v_{32})$ или $(v_{32},v_{10})$;
ребро 39: $(v_{10},v_{34})$ или $(v_{34},v_{10})$; ребро 40: $(v_{11},v_{35})$ или $(v_{35},v_{11})$;
ребро 41: $(v_{11},v_{43})$ или $(v_{43},v_{11})$; ребро 42: $(v_{12},v_{29})$ или $(v_{29},v_{12})$;
ребро 43: $(v_{12},v_{31})$ или $(v_{31},v_{12})$; ребро 44: $(v_{12},v_{33})$ или $(v_{33},v_{12})$;
ребро 45: $(v_{13},v_{32})$ или $(v_{32},v_{13})$; ребро 46: $(v_{13},v_{39})$ или $(v_{39},v_{13})$;
ребро 47: $(v_{14},v_{28})$ или $(v_{28},v_{14})$; ребро 48: $(v_{14},v_{36})$ или $(v_{36},v_{14})$;
ребро 49: $(v_{14},v_{39})$ или $(v_{39},v_{14})$; ребро 50: $(v_{15},v_{38})$ или $(v_{38},v_{15})$;
ребро 51: $(v_{15},v_{42})$ или $(v_{42},v_{15})$; ребро 52: $(v_{15},v_{44})$ или $(v_{44},v_{15})$;
ребро 53: $(v_{16},v_{37})$ или $(v_{37},v_{16})$; ребро 54: $(v_{16},v_{39})$ или $(v_{39},v_{16})$;
ребро 55: $(v_{16},v_{45})$ или $(v_{45},v_{16})$; ребро 56: $(v_{17},v_{37})$ или $(v_{37},v_{17})$;
ребро 57: $(v_{17},v_{41})$ или $(v_{41},v_{17})$; ребро 58: $(v_{17},v_{43})$ или $(v_{43},v_{17})$;
ребро 59: $(v_{19},v_{41})$ или $(v_{41},v_{19})$; ребро 60: $(v_{20},v_{40})$ или $(v_{40},v_{20})$;
ребро 61: $(v_{20},v_{43})$ или $(v_{43},v_{20})$; ребро 62: $(v_{21},v_{36})$ или $(v_{36},v_{21})$;
ребро 63: $(v_{22},v_{36})$ или $(v_{36},v_{22})$; ребро 64: $(v_{22},v_{37})$ или $(v_{37},v_{22})$;
ребро 65: $(v_{22},v_{40})$ или $(v_{40},v_{22})$; ребро 66: $(v_{23},v_{27})$ или $(v_{27},v_{23})$;
ребро 67: $(v_{25},v_{44})$ или $(v_{44},v_{25})$; ребро 68: $(v_{28},v_{40})$ или $(v_{40},v_{28})$;
ребро 69: $(v_{29},v_{34})$ или $(v_{34},v_{29})$; ребро 70: $(v_{33},v_{38})$ или $(v_{38},v_{33})$;
ребро 71: $(v_{34},v_{38})$ или $(v_{38},v_{34})$; ребро 72: $(v_{42},v_{45})$ или $(v_{45},v_{42})$.

Выделим в данном графе $G_7$ плоскую часть графа, например следующую:

| | | | |
|---|---|---|---|
| $c_1$ | $\{e_{42},e_{44},e_{69},e_{70},e_{71}\}$ | $\{v_{29},v_{34},v_{38},v_{33},v_{12}\}$ | $<v_{29},v_{34},v_{38},v_{33},v_{12}>$ |
| $c_2$ | $\{e_{17},e_{18},e_{25},e_{27}\}$ | $\{v_{18},v_6,v_{24},v_4\}$ | $<\{v_4,v_{24},v_6,v_{18}>$ |
| $c_3$ | $\{e_3,e_4,e_{29},e_{30},e_{31},e_{33},e_{42},e_{43}\}$ | $\{v_{30},v_7,v_{31},v_{12},v_{29},v_8,v_{35},v_1\}$ | $<\{v_1,v_{35},v_8,v_{29},v_{12},v_{31},v_7,v_{30}>$ |
| $c_4$ | $\{e_{15},e_{16},e_{72}\}$ | $\{v_{42},v_{45},v_3\}$ | $<v_{42},v_{45},v_3>$ |
| $c_5$ | $\{e_{31},e_{32},e_{38},e_{39},e_{69}\}$ | $\{v_{29},v_{34},v_{10},v_{32},v_8\}$ | $<v_8,v_{32},v_{10},v_{34},v_{29}>$ |
| $c_6$ | $\{e_1,e_3,e_{12},e_{13},e_{28},e_{29},e_{66}\}$ | $\{v_{19},v_3,v_{23},v_{27},v_7,v_{30},v_1\}$ | $<v_1,v_{30},v_7,v_{27},v_{23},v_3,v_{19}>$ |
| $c_7$ | $\{e_{24},e_{26},e_{46},e_{48},e_{49},e_{62}\}$ | $\{v_{13},v_{39},v_{14},v_{36},v_{21},v_6\}$ | $<v_6,v_{21},v_{36},v_{14},v_{39},v_{13}>$ |
| $c_8$ | $\{e_6,e_8,e_{35},e_{36},e_{66}\}$ | $\{v_{23},v_{27},v_9,v_{26},v_2\}$ | $<v_{23},v_{27},v_9,v_{26},v_2>$ |
| $c_9$ | $\{e_{12},e_{16},e_{53},e_{55},e_{56},e_{57},e_{59}\}$ | $\{v_{19},v_{41},v_{17},v_{37},v_{16},v_{45},v_3\}$ | $<v_3,v_{45},v_{16},v_{37},v_{17},v_{41},v_{19}>$ |
| $c_{10}$ | $\{e_{14},e_{15},e_{51},e_{52},e_{67}\}$ | $\{v_{25},v_{44},v_{15},v_{42},v_3\}$ | $<v_{25},v_{44},v_{15},v_{42},v_3>$ |
| $c_{11}$ | $\{e_1,e_4,e_{40},e_{41},e_{57},e_{58},e_{59}\}$ | $\{v_{19},v_{41},v_{17},v_{43},v_{11},v_{35},v_1\}$ | $<v_{19},v_{41},v_{17},v_{43},v_{11},v_{35},v_1>$ |
| $c_{12}$ | $\{e_{20},e_{21},e_{60},e_{62},e_{63},e_{65}\}$ | $\{v_{20},v_{40},v_{22},v_{36},v_{21},v_5\}$ | $<v_{20},v_{40},v_{22},v_{36},v_{21},v_5>$ |
| $c_{13}$ | $\{e_7,e_8,e_{21},e_{23},e_{26},e_{27}\}$ | $\{v_{24},v_6,v_{21},v_5,v_{26},v_2\}$ | $<v_2,v_{26},v_5,v_{21},v_6,v_{24}>$ |
| $c_{14}$ | $\{e_{47},e_{48},e_{63},e_{65},e_{68}\}$ | $\{v_{28},v_{40},v_{22},v_{36},v_{14}\}$ | $<v_{14},v_{36},v_{22},v_{40},v_{28}>$ |
| обод | $\{e_{44},e_{70},e_{71},e_{17},e_{18},e_{25},e_{30},$ $e_{33},e_{43},e_{72},e_{32},e_{38},e_{39},e_{13},$ $e_{28},e_{24},e_{46},e_{49},e_6,e_{35},e_{36},$ $e_{53},e_{55},e_{56},e_{14},e_{51},e_{52},e_{67},$ $e_{40},e_{41},e_{58},e_{20},e_{60},e_7,e_{23},$ $e_{47},e_{68}\}$ | $\{v_{33},v_{38},v_{34},v_{10},v_{32},v_8,$ $v_{35},v_{11},v_{43},v_{17},v_{37},v_{16},v_{45},$ $v_{42},v_{15},v_{44},v_{25},v_3,v_{23},v_2,v_{24},$ $v_4,v_{18},v_6,v_{13},v_{39},v_{14},v_{28},v_{40},$ $v_{20},v_5,v_{26},v_9,v_{27},v_7,v_{31},v_{12}\}$ | $<v_{33},v_{38},v_{34},v_{10},v_{32},v_8,$ $v_{35},v_{11},v_{43},v_{17},v_{37},v_{16},v_{45},v_{42},$ $v_{15},v_{44},v_{25},v_3,v_{23},v_2,v_{24},v_4,v_{18},$ $v_6,v_{13},v_{39},v_{14},v_{28},v_{40},v_{20},v_5,$ $v_{26},v_9,v_{27},v_7,v_{31},v_{12}>$ |

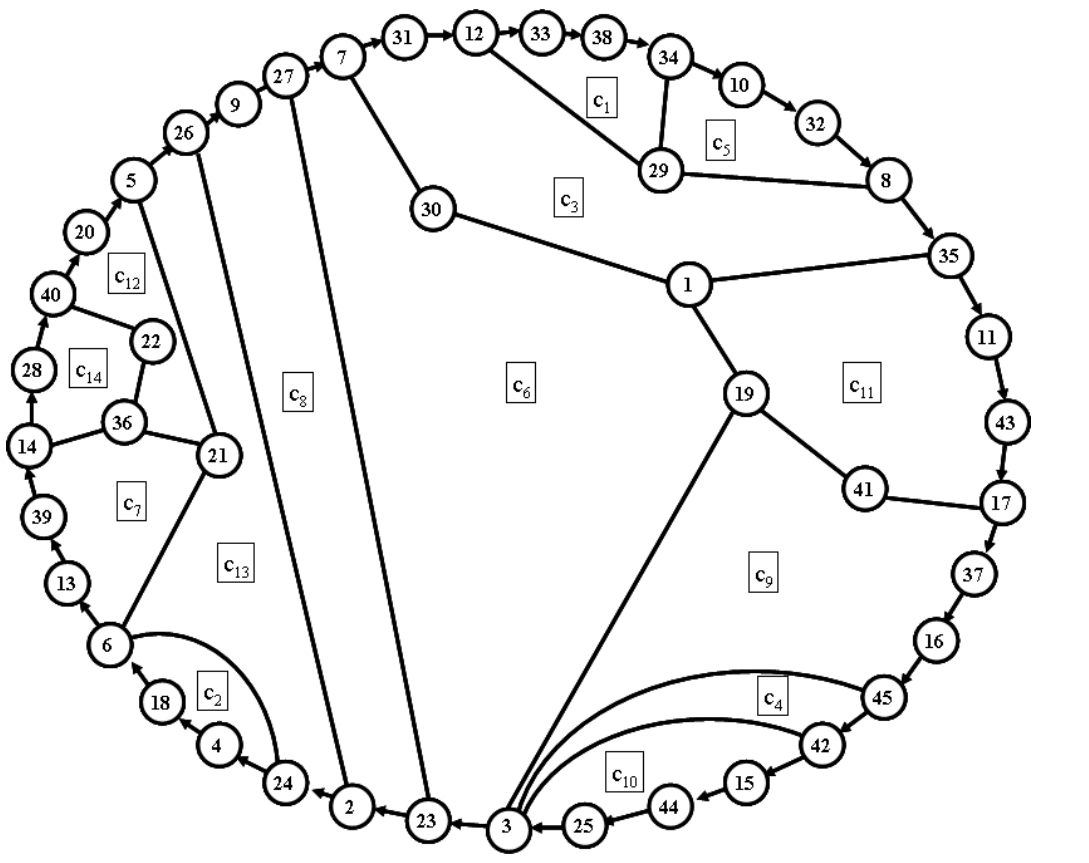

Рис. 7..13. Топологический рисунок плоской части графа $G_7$.

Рассмотрим обод данной плоской конфигурации как замкнутую последовательность ориентированных ребер. Назовём такое построение координатно-базисной системой КБС [10,11]. Проведем ребра, удаленные в процессе планаризации, которые имеют две концевые вершины совместимые с вершинами обода (см. рис. 7..14).

Запишем КБС в виде замкнутого кортежа, состоящего из ориентированных ребер $<e_{44},e_{70},e_{71},e_{39},e_{38},e_{32},e_{33},e_{40},e_{41},e_{58},e_{56},e_{53},e_{55},e_{72},e_{51},e_{52},e_{67},e_{14},e_{13},e_6,e_7,e_{18},e_{17},e_{25},e_{24},e_{46},e_{49},e_{47},e_{68},$ $e_{60},e_{20},e_{23},e_{35},e_{36},e_{28},e_{30},e_{43}>$. Каждое неориентированное ребро можно представить в виде двух разнонаправленных ориентированных дуг. Тогда для каждого ориентированного ребра существует проекция (в теоретико-множественном смысле) на координатно-базисную систему.

Например:

проекция $e_9(v_2,v_{33}) = <e_7,e_{18},e_{17},e_{25},e_{24},e_{46},e_{49},e_{47},e_{68},e_{60},e_{20},e_{23},e_{35},e_{36},e_{28},e_{30},e_{43},e_{44}>$;
проекция $e_9(v_{33},v_2) = <e_{70},e_{71},e_{39},e_{38},e_{32},e_{33},e_{40},e_{41},e_{58},e_{56},e_{53},e_{55},e_{72},e_{51},e_{52},e_{67},e_{14},e_{13},e_6>$;
проекция $e_{50}(v_{38},v_{15}) = <e_{71},e_{39},e_{38},e_{32},e_{33},e_{40},e_{41},e_{58},e_{56},e_{53},e_{55},e_{72},e_{51}>$;
проекция $e_{50}(v_{15},v_{38}) = <e_{52},e_{67},e_{14},e_{13},e_6,e_7,e_{18},e_{17},e_{25},e_{24},e_{46},e_{49},e_{47},e_{68},e_{60},e_{20},e_{23},e_{35},e_{36},e_{28},$ $e_{30},e_{43},e_{44},e_{70}>$;
проекция $e_{61}(v_{20},v_{43}) = <e_{20},e_{23},e_{35},e_{36},e_{28},e_{30},e_{43},e_{44},e_{70},e_{71},e_{39},e_{38},e_{32},e_{33},e_{40},e_{41}>$;
проекция $e_{61}(v_{43},v_{20}) = <e_{58},e_{56},e_{53},e_{55},e_{72},e_{51},e_{52},e_{67},e_{14},e_{13},e_6,e_7,e_{18},e_{17},e_{25},e_{24},e_{46},e_{49},$ $e_{47},e_{68},e_{60}>$;
проекция $e_{45}(v_{13},v_{32}) = <e_{46},e_{49},e_{47},e_{68},e_{60},e_{20},e_{23},e_{35},e_{36},e_{28},e_{30},e_{43},e_{44},e_{70},e_{71},e_{39},e_{38}>$;
проекция $e_{45}(v_{32},v_{13}) = <e_{32},e_{33},e_{40},e_{41},e_{58},e_{56},e_{53},e_{55},e_{72},e_{51},e_{52},e_{67},e_{14},e_{13},e_6,e_7,e_{18},e_{17},e_{25},e_{24}>$;
проекция $e_{54}(v_{16},v_{39}) = <e_{55},e_{72},e_{51},e_{52},e_{67},e_{14},e_{13},e_6,e_7,e_{18},e_{17},e_{25},e_{24},e_{46}>$;
проекция $e_{54}(v_{39}.v_{16}) = <e_{49},e_{47},e_{68},e_{60},e_{20},e_{23},e_{35},e_{36},e_{28},e_{30},e_{43},e_{44},e_{70},e_{71},e_{39},e_{38},e_{32},e_{33},e_{40},$ $e_{41},e_{58},e_{56},e_{53}>$;

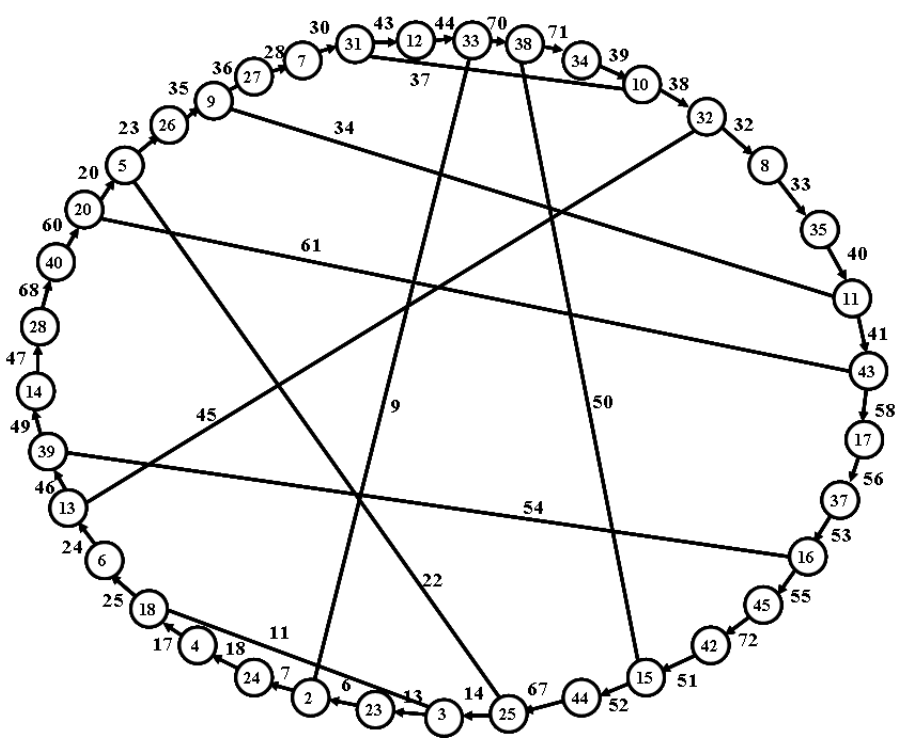

Рис. 7.14. Координатно-базисная система и проведенные ребра.

проекция $e_{22}(v_{25},v_5) = <e_{14},e_{13},e_6,e_7,e_{18},e_{17},e_{25},e_{24},e_{46},e_{49},e_{47},e_{68},e_{60},e_{20}>$;
проекция $e_{22}(v_5,v_{25}) = <e_{23},e_{35},e_{36},e_{28},e_{30},e_{43},e_{44},e_{70},e_{71},e_{39},e_{38},e_{32},e_{33},e_{40},e_{41},e_{58},e_{56},e_{53},e_{55},$ $e_{72},e_{51},e_{52},e_{67}>$;
проекция $e_{11}(v_3,v_{18}) = <e_{13},e_6,e_7,e_{18},e_{17}>$;
проекция $e_{11}(v_{18},v_3) = <e_{25},e_{24},e_{46},e_{49},e_{47},e_{68},e_{60},e_{20},e_{23},e_{35},e_{36},e_{28},e_{30},e_{43},e_{44},e_{70},e_{71},e_{39},e_{38},$ $e_{32},e_{33},e_{40},e_{41},e_{58},e_{56},e_{53},e_{55},e_{72},e_{51},e_{52},e_{67},e_{14}>$;
проекция $e_{37}(v_{31},v_{10}) = <e_{43},e_{44},e_{70},e_{71},e_{39}>$;
проекция $e_{37}(v_{10},v_{31}) = <e_{38},e_{32},e_{33},e_{40},e_{41},e_{58},e_{56},e_{53},e_{55},e_{72},e_{51},e_{52},e_{67},e_{14},e_{13},e_6,e_7,e_{18},e_{17},$ $e_{25},e_{24},e_{46},e_{49},e_{47},e_{68},e_{60},e_{20},e_{23},e_{35},e_{36},e_{28},e_{30}>$;
проекция $e_{34}(v_9,v_{11}) = <e_{36},e_{28},e_{30},e_{43},e_{44},e_{70},e_{71},e_{39},e_{38},e_{32},e_{33},e_{40}>$;
проекция $e_{34}(v_{11},v_9) = <e_{41},e_{58},e_{56},e_{53},e_{55},e_{72},e_{51},e_{52},e_{67},e_{14},e_{13},e_6,e_7,e_{18},e_{17},e_{25},e_{24},e_{46},e_{49},$ $e_{47},e_{68},e_{60},e_{20},e_{23},e_{35}>$.

Запишем КБС в виде замкнутого кортежа, состоящего из ориентированных ребер $<e_{44},e_{70},e_{71},e_{39},e_{38},e_{32},e_{33},e_{40},e_{41},e_{58},e_{56},e_{53},e_{55},e_{72},e_{51},e_{52},e_{67},e_{14},e_{13},e_6,e_7,e_{18},e_{17},e_{25},e_{24},e_{46},e_{49},e_{47},e_{68},$ $e_{60},e_{20},e_{23},e_{35},e_{36},e_{28},e_{30},e_{43}>$. Тогда неориентированное ребро можно представить в виде двух направленных в разные стороны ориентированных дуг. Для каждого ориентированного ребра существует проекция (в множественном смысле) на координатно-базисную систему. Например:

Из двух проекций ребра, выделим одну проекцию, но минимальную по длине.

проекция $e_9(v_2,v_{33}) = <e_7,e_{18},e_{17},e_{25},e_{24},e_{46},e_{49},e_{47},e_{68},e_{60},e_{20},e_{23},e_{35},e_{36},e_{28},e_{30},e_{43},e_{44}>$;
проекция $e_{50}(v_{38},v_{15}) = <e_{71},e_{39},e_{38},e_{32},e_{33},e_{40},e_{41},e_{58},e_{56},e_{53},e_{55},e_{72},e_{51}>$;
проекция $e_{61}(v_{20},v_{43}) = <e_{20},e_{23},e_{35},e_{36},e_{28},e_{30},e_{43},e_{44},e_{70},e_{71},e_{39},e_{38},e_{32},e_{33},e_{40},e_{41}>$;
проекция $e_{45}(v_{13},v_{32}) = <e_{46},e_{49},e_{47},e_{68},e_{60},e_{20},e_{23},e_{35},e_{36},e_{28},e_{30},e_{43},e_{44},e_{70},e_{71},e_{39},e_{38}>$;

проекция $e_{54}(v_{16},v_{39}) = \langle e_{55},e_{72},e_{51},e_{52},e_{67},e_{14},e_{13},e_{6},e_{7},e_{18},e_{17},e_{25},e_{24},e_{46} \rangle$;
проекция $e_{22}(v_{25},v_{5}) = \langle e_{14},e_{13},e_{6},e_{7},e_{18},e_{17},e_{25},e_{24},e_{46},e_{49},e_{47},e_{68},e_{60},e_{20} \rangle$;
проекция $e_{11}(v_{3},v_{18}) = \langle e_{13},e_{6},e_{7},e_{18},e_{17} \rangle$;
проекция $e_{37}(v_{31},v_{10}) = \langle e_{43},e_{44},e_{70},e_{71},e_{39} \rangle$;
проекция $e_{34}(v_{9},v_{11}) = \langle e_{36},e_{28},e_{30},e_{43},e_{44},e_{70},e_{71},e_{39},e_{38},e_{32},e_{33},e_{40} \rangle$;

В соответствии с законами векторной алгебры пересечений будем считать, что ребра пересекаются, если пересекаются их проекции (в терминах пересечения теории множеств). Ребра не пересекаются, если результат пересечения проекций есть пустое множество или одна проекция полностью включается в другую.

Рассмотрим пересечения ребра $e_9$. Пересечение ребер будем обозначать символом $\perp$. Проекцию ребра на координатно-базисную систему будем обозначать двумя символами **пр**.

$e_9 \perp e_{50} = \textbf{пр}(e_9) \cap \textbf{пр}(e_{50}) = \{e_7,e_{18},e_{17},e_{25},e_{24},e_{46},e_{49},e_{47},e_{68},e_{60},e_{20},e_{23},e_{35},e_{36},e_{28},e_{30},e_{43},e_{44}\} \cap$
$\cap \{e_{38},e_{32},e_{33},e_{40},e_{41},e_{58},e_{56},e_{53},e_{55},e_{72},e_{51}\} = \varnothing$;
$e_9 \perp e_{61} = \textbf{пр}(e_9) \cap \textbf{пр}(e_{61}) = \{e_7,e_{18},e_{17},e_{25},e_{24},e_{46},e_{49},e_{47},e_{68},e_{60},e_{20},e_{23},e_{35},e_{36},e_{28},e_{30},e_{43},e_{44}\} \cap$
$\cap \{e_{20},e_{23},e_{35},e_{36},e_{28},e_{30},e_{43},e_{44},e_{70},e_{71},e_{39},e_{38},e_{32},e_{33},e_{40},e_{41}\} = \{e_{20},e_{23},e_{35},e_{36},e_{28},e_{30},e_{43},e_{44}\}$;
$e_9 \perp e_{45} = \textbf{пр}(e_9) \cap \textbf{пр}(e_{45}) = \{e_7,e_{18},e_{17},e_{25},e_{24},e_{46},e_{49},e_{47},e_{68},e_{60},e_{20},e_{23},e_{35},e_{36},e_{28},e_{30},e_{43},e_{44}\} \cap$
$\cap \{e_{46},e_{49},e_{47},e_{68},e_{60},e_{20},e_{23},e_{35},e_{36},e_{28},e_{30},e_{43},e_{44},e_{70},e_{71},e_{39},e_{38}\} =$
$= \{e_{49},e_{47},e_{68},e_{60},e_{20},e_{23},e_{35},e_{36},e_{28},e_{30},e_{43},e_{44}\}$;
$e_9 \perp e_{54} = \textbf{пр}(e_9) \cap \textbf{пр}(e_{54}) = \{e_7,e_{18},e_{17},e_{25},e_{24},e_{46},e_{49},e_{47},e_{68},e_{60},e_{20},e_{23},e_{35},e_{36},e_{28},e_{30},e_{43},e_{44}\} \cap$
$\cap \{e_{55},e_{72},e_{51},e_{52},e_{67},e_{14},e_{13},e_{6},e_{7},e_{18},e_{17},e_{25},e_{24},e_{46}\} = \{e_7,e_{18},e_{17},e_{25},e_{24},e_{46}\}$;
$e_9 \perp e_{22} = \textbf{пр}(e_9) \cap \textbf{пр}(e_{22}) = \{e_7,e_{18},e_{17},e_{25},e_{24},e_{46},e_{49},e_{47},e_{68},e_{60},e_{20},e_{23},e_{35},e_{36},e_{28},e_{30},e_{43},e_{44}\} \cap$
$\cap \{e_{14},e_{13},e_{6},e_{7},e_{18},e_{17},e_{25},e_{24},e_{46},e_{49},e_{47},e_{68},e_{60},e_{20}\} = \{e_7,e_{18},e_{17},e_{25},e_{24},e_{46},e_{49},e_{47},e_{68},e_{60},e_{20}\}$;
$e_9 \perp e_{11} = \textbf{пр}(e_9) \cap \textbf{пр}(e_{11}) = \{e_7,e_{18},e_{17},e_{25},e_{24},e_{46},e_{49},e_{47},e_{68},e_{60},e_{20},e_{23},e_{35},e_{36},e_{28},e_{30},e_{43},e_{44}\} \cap$
$\cap \{e_{13},e_{6},e_{7},e_{18},e_{17}\} = \{e_7,e_{18},e_{17}\}$;
$e_9 \perp e_{37} = \textbf{пр}(e_9) \cap \textbf{пр}(e_{37}) = \{e_7,e_{18},e_{17},e_{25},e_{24},e_{46},e_{49},e_{47},e_{68},e_{60},e_{20},e_{23},e_{35},e_{36},e_{28},e_{30},e_{43},e_{44}\} \cap$
$\cap \{e_{43},e_{44},e_{70},e_{71},e_{39}\} = \{e_{43},e_{44}\}$;
$e_9 \perp e_{34} = \textbf{пр}(e_9) \cap \textbf{пр}(e_{34}) = \{e_7,e_{18},e_{17},e_{25},e_{24},e_{46},e_{49},e_{47},e_{68},e_{60},e_{20},e_{23},e_{35},e_{36},e_{28},e_{30},e_{43},e_{44}\} \cap$
$\cap \{e_{36},e_{28},e_{30},e_{43},e_{44},e_{70},e_{71},e_{39},e_{38},e_{32},e_{33},e_{40}\} = \{e_{36},e_{28},e_{30},e_{43},e_{44}\}$.

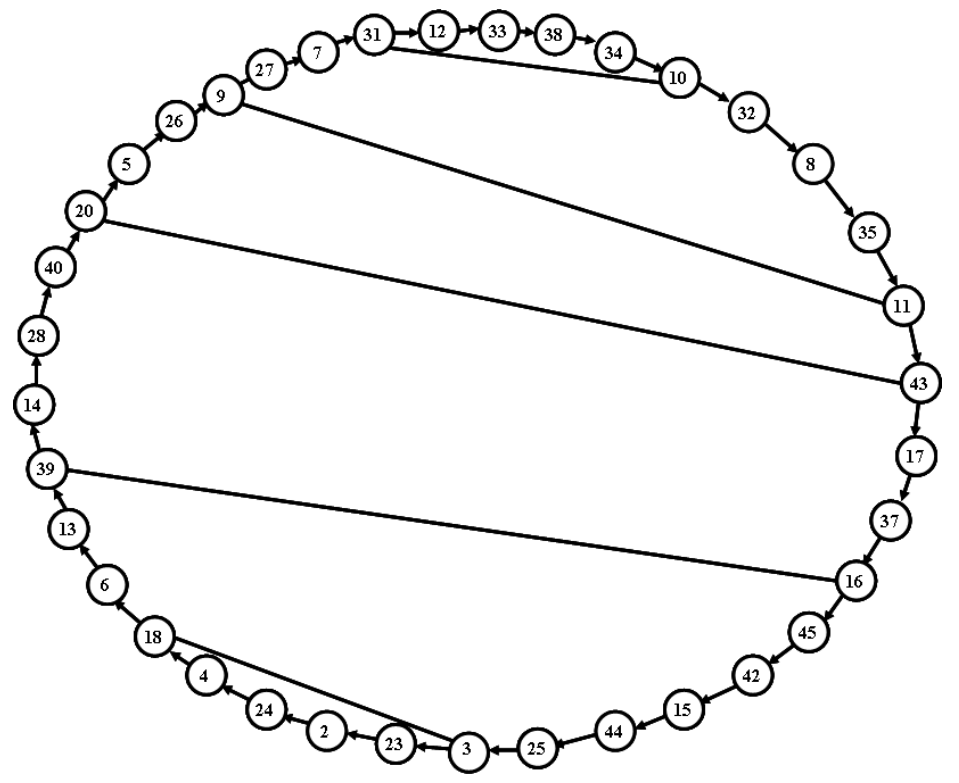

Рис. 7.15. Непересекающиеся соединения графа $G_7$.

Таким образом, ребро $e_9$ пересекается со всеми ребрами кроме ребра $e_{50}$.

В результате выделяются новые простые циклы, которые можно добавить к существующим циклам.

$c_{15} = \{e_{37},e_{43},e_{41},e_{70},e_{71},e_{39}\}$;

$c_{16} = \{e_{37},e_{38},e_{32},e_{33},e_{40},e_{34},e_{36},e_{28},e_{30}\}$;
$c_{17} = \{e_{34},e_{41},e_{61},e_{20},e_{23},e_{35}\}$;
$c_{18} = \{e_{54},e_{55},e_{72},e_{51},e_{52},e_{67},e_{11},e_{25},e_{24},e_{46}\}$;
$c_{19} = \{e_{11},e_{17},e_{18},e_{6},e_{7},e_{15}\}$;
$c_{0} = \{e_{61},e_{58},e_{56},e_{53},e_{54},e_{49},e_{47},e_{68},e_{60}\}$.

Для определения удаления ребер, необходимо рассмотреть все случаи парного пересечения ребер. После этого будем последовательно удалять ребра, максимально пересекающиеся с другими. В конце процесса выделится множество непересекающихся ребер. В нашем случае непересекающиеся соединения представлены на рис. 7.15.

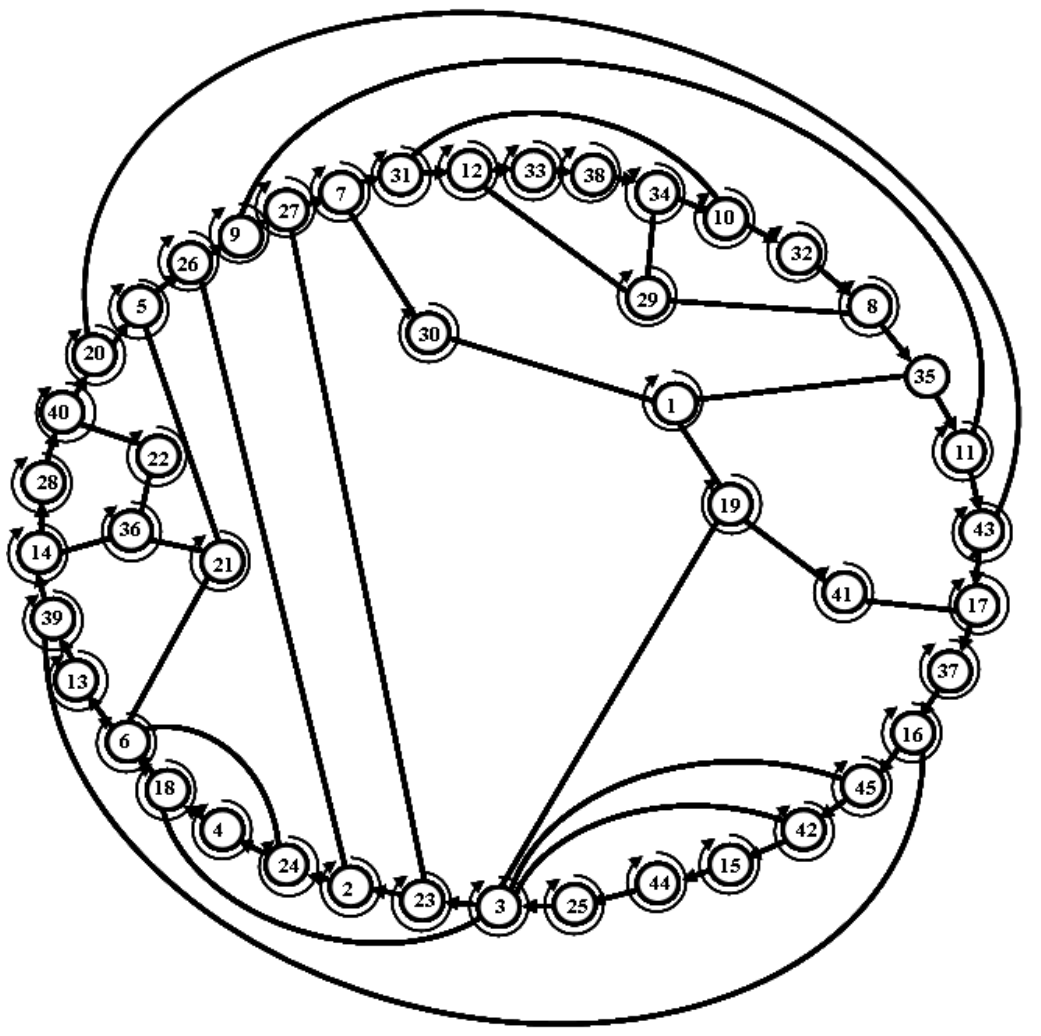

Рис. 7.16. Топологический рисунок плоского суграфа графа $G_7$.

Множество простых и изометрических циклов образует топологический рисунок плоского суграфа (см. рис. 7.16).

Возможен путь решения задачи выделения плоской части начиная с этапа выделения базиса подпространства циклов. В этом случае выделенный базис не обязательно обладает минимальным значением функционала Маклейна. Рассмотрим данный процесс на следующем примере.

***Пример 7.8.*** Пусть задан граф $G_8$ с количеством вершин = 21 и количеством ребер = 33. Граф задан инцидентором P:

ребро 1: $(v_1,v_5)$ или $(v_5,v_1)$;
ребро 2: $(v_1,v_8)$ или $(v_8,v_1)$;
ребро 3: $(v_1,v_{11})$ или $(v_{11},v_1)$;
ребро 4: $(v_1,v_{21})$ или $(v_{21},v_1)$;
ребро 5: $(v_2,v_{12})$ или $(v_{12},v_2)$;
ребро 6: $(v_2,v_{16})$ или $(v_{16},v_2)$;
ребро 7: $(v_2,v_{19})$ или $(v_{19},v_2)$;
ребро 8: $(v_2,v_{21})$ или $(v_{21},v_2)$;
ребро 9: $(v_3,v_7)$ или $(v_7,v_3)$;
ребро 10: $(v_3,v_{12})$ или $(v_{12},v_3)$;
ребро 11: $(v_3,v_{16})$ или $(v_{16},v_3)$;
ребро 12: $(v_4,v_{12})$ или $(v_{12},v_4)$;
ребро 13: $(v_4,v_{13})$ или $(v_{13},v_4)$;
ребро 14: $(v_4,v_{14})$ или $(v_{14},v_4)$;
ребро 15: $(v_4,v_{20})$ или $(v_{20},v_4)$;
ребро 16: $(v_5,v_{17})$ или $(v_{17},v_5)$;
ребро 17: $(v_5,v_{20})$ или $(v_{20},v_5)$;
ребро 18: $(v_6,v_{10})$ или $(v_{10},v_6)$;
ребро 19: $(v_6,v_{14})$ или $(v_{14},v_6)$;
ребро 20: $(v_6,v_{15})$ или $(v_{15},v_6)$;
ребро 21: $(v_7,v_{10})$ или $(v_{10},v_7)$;
ребро 22: $(v_7,v_{18})$ или $(v_{18},v_7)$;
ребро 23: $(v_8,v_{15})$ или $(v_{15},v_8)$;
ребро 24: $(v_8,v_{18})$ или $(v_{18},v_8)$;
ребро 25: $(v_9,v_{13})$ или $(v_{13},v_9)$;
ребро 26: $(v_9,v_{17})$ или $(v_{17},v_9)$;
ребро 27: $(v_9,v_{18})$ или $(v_{18},v_9)$;
ребро 28: $(v_{10},v_{20})$ или $(v_{20},v_{10})$;

ребро 29: $(v_{11},v_{15})$ или $(v_{15},v_{11})$; ребро 30: $(v_{11},v_{17})$ или $(v_{17},v_{11})$;
ребро 31: $(v_{13},v_{21})$ или $(v_{21},v_{13})$; ребро 32: $(v_{14},v_{19})$ или $(v_{19},v_{14})$;
ребро 33: $(v_{16},v_{19})$ или $(v_{19},v_{16})$.

Множество изометрических циклов:

| | ребра | вершины |
|---|---|---|
| цикл 1: | $\{e_1,e_3,e_{16},e_{30}\}$ | $\{v_1,v_5,v_{11},v_{17}\}$ |
| цикл 2: | $\{e_1,e_4,e_{13},e_{15},e_{17},e_{31}\}$ | $\{v_1,v_4,v_5,v_{13},v_{20},v_{21}\}\{$ |
| цикл 3: | $\{e_1,e_2,e_{16},e_{24},e_{26},e_{27}\}$ | $\{v_1,v_5,v_8,v_9,v_{17},v_{18}\}\{$ |
| цикл 4: | $\{e_2,e_3,e_{23},e_{29}\}$ | $\{v_1,v_8,v_{11},v_{15}\}$ |
| цикл 5: | $\{e_2,e_4,e_{24},e_{25},e_{27},e_{31}\}$ | $\{v_1,v_8,v_9,v_{13},v_{18},v_{21}\}$ |
| цикл 6: | $\{e_1,e_4,e_{16},e_{25},e_{26},e_{31}\}$ | $\{v_1,v_5,v_9,v_{13},v_{17},v_{21}\}$ |
| цикл 7: | $\{e_3,e_4,e_{25},e_{26},e_{30},e_{31}\}$ | $\{v_1,v_9,v_{11},v_{13},v_{17},v_{21}\}$ |
| цикл 8: | $\{e_5,e_6,e_{10},e_{11}\}$ | $\{v_2,v_3,v_{12},v_{16}\}$ |
| цикл 9: | $\{e_5,e_7,e_{12},e_{14},e_{32}\}$ | $\{v_2,v_4,v_{12},v_{14},v_{19}\}$ |
| цикл 10: | $\{e_5,e_8,e_{12},e_{13},e_{31}\}$ | $\{v_2,v_4,v_{12},v_{13},v_{21}\}$ |
| цикл 11: | $\{e_6,e_7,e_{33}\}$ | $\{v_2,v_{16},v_{19}\}$ |
| цикл 12: | $\{e_9,e_{10},e_{12},e_{15},e_{21},e_{28}\}$ | $\{v_3,v_4,v_7,v_{10},v_{12},v_{20}\}$ |
| цикл 13: | $\{e_{13},e_{15},e_{16},e_{17},e_{25},e_{26}\}$ | $\{v_4,v_5,v_9,v_{13},v_{17},v_{20}\}$ |
| цикл 14: | $\{e_{14},e_{15},e_{18},e_{19},e_{28}\}$ | $\{v_4,v_6,v_{10},v_{14},v_{20}\}$ |
| цикл 15: | $\{e_{18},e_{20},e_{21},e_{22},e_{23},e_{24}\}$ | $\{v_6,v_7,v_8,v_{10},v_{15},v_{18}\}$ |
| цикл 16: | $\{e_9,e_{10},e_{12},e_{13},e_{22},e_{25},e_{27}\}$ | $\{v_3,v_4,v_7,v_9,v_{12},v_{13},v_{18}\}$ |
| цикл 17: | $\{e_2,e_3,e_{24},e_{26},e_{27},e_{30}\}$ | $\{v_1,v_8,v_9,v_{11},v_{17},v_{18}\}$ |
| цикл 18: | $\{e_{23},e_{24},e_{26},e_{27},e_{29},e_{30}\}$ | $\{v_8,v_9,v_{11},v_{15},v_{17},v_{18}\}$ |

Случайным образом выделим множество изометрических циклов с мощностью равным цикломатическому числу графа.

| | ребра | вершины |
|---|---|---|
| цикл 1: | $\{e_1,e_3,e_{16},e_{30}\}$ | $\{v_1,v_5,v_{11},v_{17}\}$ |
| цикл 2: | $\{e_1,e_4,e_{13},e_{15},e_{17},e_{31}\}$ | $\{v_1,v_4,v_5,v_{13},v_{20},v_{21}\}\{$ |
| цикл 3: | $\{e_1,e_2,e_{16},e_{24},e_{26},e_{27}\}$ | $\{v_1,v_5,v_8,v_9,v_{17},v_{18}\}\{$ |
| цикл 4: | $\{e_2,e_3,e_{23},e_{29}\}$ | $\{v_1,v_8,v_{11},v_{15}\}$ |
| цикл 5: | $\{e_2,e_4,e_{24},e_{25},e_{27},e_{31}\}$ | $\{v_1,v_8,v_9,v_{13},v_{18},v_{21}\}$ |
| цикл 8: | $\{e_5,e_6,e_{10},e_{11}\}$ | $\{v_2,v_3,v_{12},v_{16}\}$ |
| цикл 9: | $\{e_5,e_7,e_{12},e_{14},e_{32}\}$ | $\{v_2,v_4,v_{12},v_{14},v_{19}\}$ |
| цикл 10: | $\{e_5,e_8,e_{12},e_{13},e_{31}\}$ | $\{v_2,v_4,v_{12},v_{13},v_{21}\}$ |
| цикл 11: | $\{e_6,e_7,e_{33}\}$ | $\{v_2,v_{16},v_{19}\}$ |
| цикл 12: | $\{e_9,e_{10},e_{12},e_{15},e_21,e_{28}\}$ | $\{v_3,v_4,v_7,v_{10},v_{12},v_{20}\}$ |
| цикл 14: | $\{e_{14},e_{15},e_{18},e_{19},e_{28}\}$ | $\{v_4,v_6,v_{10},v_{14},v_{20}\}$ |
| цикл 15: | $\{e_{18},e_{20},e_{21},e_{22},e_{23},e_{24}\}$ | $\{v_6,v_7,v_8,v_{10},v_{15},v_{18}\}$ |
| цикл 16: | $\{e_9,e_{10},e_{12},e_{13},e_{22},e_{25},e_{27}\}$ | $\{v_3,v_4,v_7,v_9,v_{12},v_{13},v_{18}\}$ |

Длина элемента структурного числа = 13.
Количество элементов структурного числа = 1.
Усеченные однострочные структурные числа:
циклы, проходящие по хорде $e_{17}$: $[c_2]$;
циклы, проходящие по хорде $e_{30}$: $[c_1]$;
циклы, проходящие по хорде $e_{33}$: $[c_{11}]$;
циклы, проходящие по хорде $e_3$: $[c_1,c_4]$;
циклы, проходящие по хорде $e_4$: $[c_2,c_5]$;
циклы, проходящие по хорде $e_7$: $[c_9,c_{11}]$;
циклы, проходящие по хорде $e_{14}$: $[c_9,c_{14}]$;

циклы, проходящие по хорде $e_{22}$: $[c_{15},c_{16}]$;
циклы, проходящие по хорде $e_2$: $[c_3,c_4,c_5]$;
циклы, проходящие по хорде $e_{10}$: $[c_8,c_{12},c_{16}]$;
циклы, проходящие по хорде $e_{15}$: $[c_2,c_{12},c_{14}]$;
циклы, проходящие по хорде $e_{27}$: $[c_3,c_5,c_{16}]$;
циклы, проходящие по хорде $e_{31}$: $[c_2,c_5,c_{10}]$;

1 - элемент структурного числа = $\{c_2,c_1,c_{11},c_4,c_5,c_9,c_{14},c_{15},c_3,c_8,c_{12},c_{16},c_{10}\}$.

Количество элементов в произведении одностросных структурных чисел нечетно, следовательно, выбранная система циклов не зависима.

Алгоритмом наискорейшего спуска выделим плоскую часть графа

| | ребра | вершины |
|---|---|---|
| цикл 1: | $\{e_1,e_3,e_{16},e_{30}\}$ | $\{v_1,v_5,v_{11},v_{17}\}$ |
| цикл 4: | $\{e_2,e_3,e_{23},e_{29}\}$ | $\{v_1,v_8,v_{11},v_{15}\}$ |
| цикл 5: | $\{e_2,e_4,e_{24},e_{25},e_{27},e_{31}\}$ | $\{v_1,v_8,v_9,v_{13},v_{18},v_{21}\}$ |
| цикл 9: | $\{e_5,e_7,e_{12},e_{14},e_{32}\}$ | $\{v_2,v_4,v_{12},v_{14},v_{19}\}$ |
| цикл 11: | $\{e_6,e_7,e_{33}\}$ | $\{v_2,v_{16},v_{19}\}$ |
| цикл 12: | $\{e_9,e_{10},e_{12},e_{15},e_21,e_{28}\}$ | $\{v_3,v_4,v_7,v_{10},v_{12},v_{20}\}$ |
| цикл 14: | $\{e_{14},e_{15},e_{18},e_{19},e_{28}\}$ | $\{v_4,v_6,v_{10},v_{14},v_{20}\}$ |
| цикл 15: | $\{e_{18},e_{20},e_{21},e_{22},e_{23},e_{24}\}$ | $\{v_6,v_7,v_8,v_{10},v_{15},v_{18}\}$ |

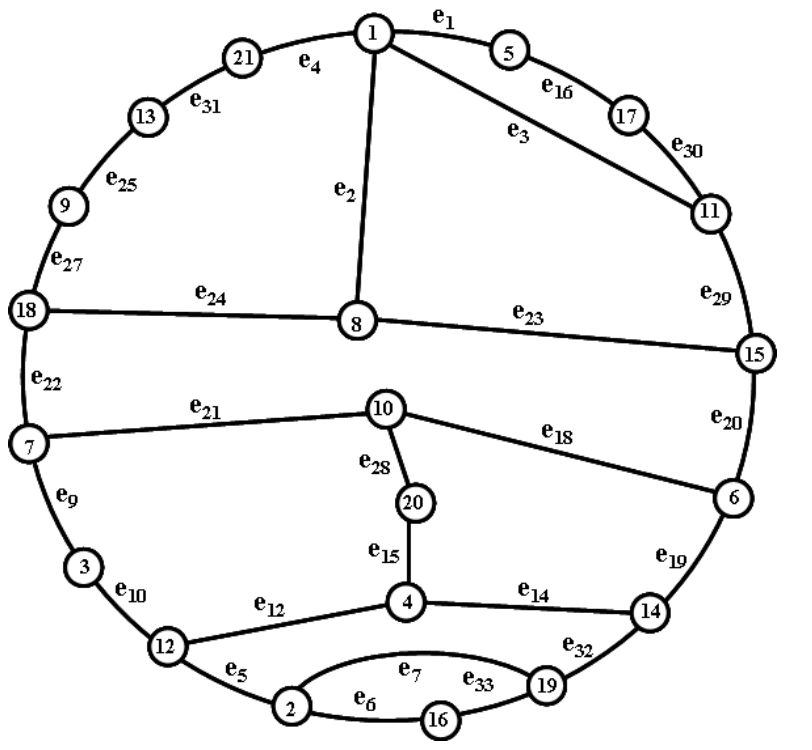

Рис. 7.17. Рисунок плоской части графа $G_8$.

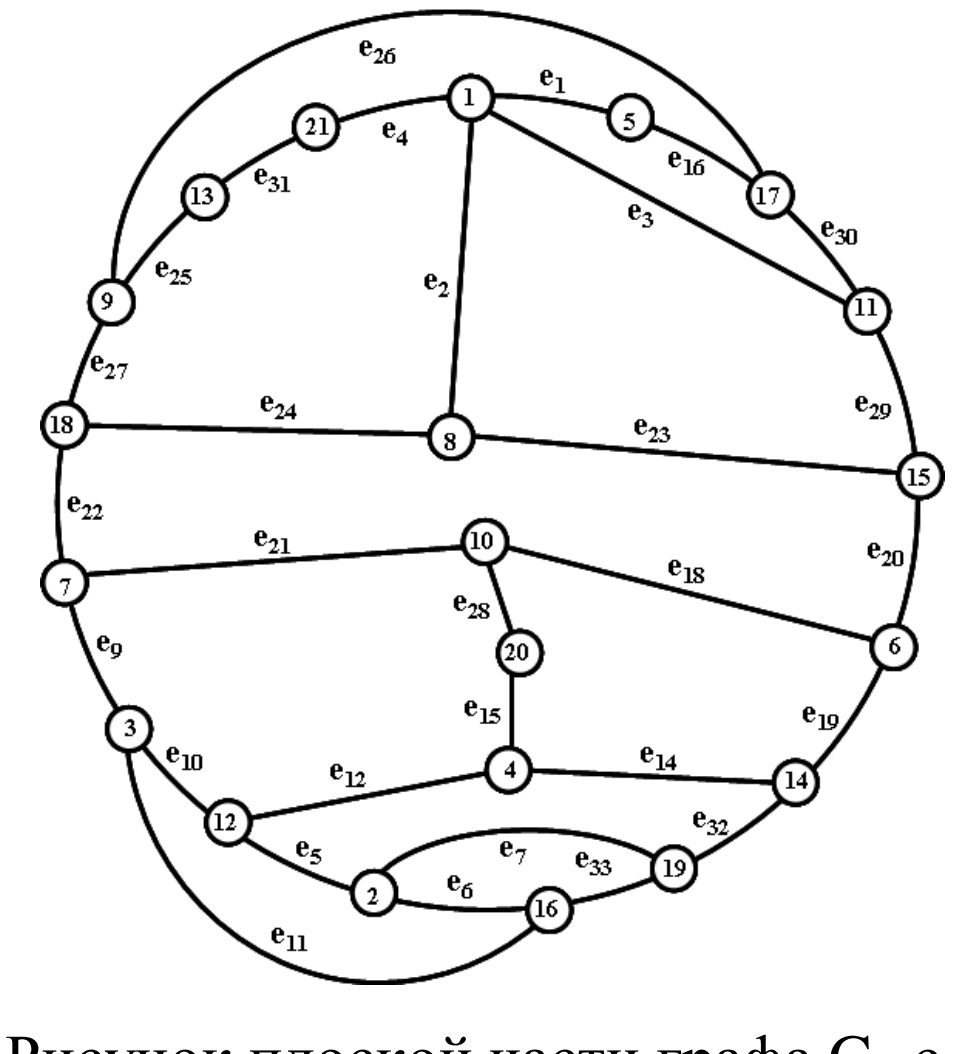

Рис. 7.18. Рисунок плоской части графа $G_8$ с добавлением ребер.

Окончанием процесса служит определение подмножества простых циклов для вновь введенных ребер. Далее производится включение данного подмножества циклов во множество циклов описывающих плоский сугpаф и определяется обод сугpафа.

Последним этапом является построение вращения вершин графа для выбранного множества циклов.

**Комментарии**

В данной главе рассматриваются методы построения топологического рисунка плоской части непланарного графа *подходом сверху*. Процесс решения состоит из четырех этапов.

На первом этапе из множества изометрических циклов графа выделяется подмножество циклов с мощностью равному цикломатическому числу графа. Данное подмножество должно быть базисом подпространства циклов и удовлетворять некоторым критериям. Если для графа G, существует базис подпространства циклов со значением функционала Маклейна равному нулю, то несепарабельный граф планарен. Показано, что для описания процесса выделения базиса изометрических циклов можно применять методы алгебры структурных чисел [3].

Выделение подмножества циклов с мощностью равному цикломатическому числу графа методом наискорейшего спуска с большой долей вероятности приводит к появления зависимого подмножества циклов

Выделение подмножества циклов с мощностью равному цикломатическому числу графа алгоритмом «прыгающая строка» требует проверки на линейную независимость. С этой целью строится усеченное число однострочных структурных чисел для определения четности повторения элементов, или применяется модифицированный алгоритм Гаусса

Второй этап удаления минимального количества циклов из базиса, для непланарных графов, предназначен для построения плоской части графа. Этот этап осуществляется путем выделения из базиса циклов подмножества, имеющего нулевое значение кубического функционала Маклейна. Причем удаление циклов из базиса циклов (элементов), должно производиться с выполнением условия Эйлера. С целью выделения такого подмножества, обосновывается введение кубического функционала Маклейна. Для описания процесса удаления циклов применяется метод наискорейшего спуска с вычислением кубического функционала Маклейна каждом шаге.

Третий этап предназначен для введения в топологический рисунок ребер ранее удаленных в процессе планаризации. Для определения пересечения ребер применяются методы алгебры пересечений соединений [29-31]. Показаны способы построения координатной-базисной системы векторов (ориентированных ребер) графа. Методами векторной алгебры пересечений строится дополнительная система простых циклов графа относительно обода.

Четвертый этап предназначен для построения вращения вершин графа выбранной системы циклов, тем самым определяя топологический рисунок плоской части.

Очевидно, что существует множество способов выделения плоской части графа, с минимальным количеством удаленных ребер. В зависимости от комбинаций способов, возможно построение различных практических систем программного обеспечения.

# Глава 8. ПОСТРОЕНИЕ ТОПОЛОГИЧЕСКОГО РИСУНКА ПЛОСКОЙ ЧАСТИ НЕПЛАНАРНОГО ГРАФА

## 8.1. Топологический рисунок плоской части непланарного графа

Будем рассматривать несепарабельный неориентированный непланарный граф G(V,E) с множеством вершин card V = *n* и множеством ребер card E=*m*. Выделим множество изометрических циклов графа $C_\tau$ [15]. Из множества изометрических циклов $C_\tau$ выделим множество базисов принадлежащих пространству суграфов £(G) $B_\tau = \{b_{\tau 1}, b_{\tau 2}, ..., b_{\tau p}\}$, $\operatorname{card} B_\tau = p$, $b_{\tau i} \subset C_\tau$, $card\,(b_{\tau i}) = m - n + 1$ для всех $i$ = 1,2,…,p. Построение базиса будем осуществлять модифицированным алгоритмом Гаусса из множества изометрических циклов графа со случайным порядком расположения элементов в кортеже циклов. Порядок расположения элементов будем описывать с помощью фрагментарного представления перестановок циклов [10,11].

Для выделения базиса пространства суграфов £(G), состоящего из множества изометрических циклов, воспользуемся модифицированным алгоритмом Гаусса для определения ранга системы [12].

Так как процесс выделения циклов модифицированным методом Гаусса очень чувствителен к порядку расположению элементов во множестве $C_\tau$, то базисы формируются в результате случайного порядка расположения изометрических циклов в $C_\tau$.

Следующим шагом является выделение циклов из базиса методом наискорейшего спуска до достижения нулевого значения кубического функционала Маклейна [18]

$$\text{FP}(b_\tau) = \sum_{\text{i}=1}^{\text{m}} a_\text{i}(a_\text{i} - 1)(a_\text{i} - 2) = \sum_{i=1}^{m} a_\text{i}^3 - 3\sum_{i=1}^{m} a_\text{i}^2 + 2\sum_{i=1}^{m} a_\text{i} \to 0 \cdot \qquad (8.1)$$

Здесь коэффициент $a_\text{i}$ характеризует количество изометрических циклов проходящих по ребру $e_\text{i}$ в подмножестве циклов. Каждый раз, исключая цикл из базиса необходимо соблюдать выполнение условия Эйлера для количества оставшихся циклов $k_c$ в системе [16]

$$k_c - m + n - 1 = 0\,. \qquad (8.2)$$

Среди множества плоских конфигураций, с заданной мощностью *p*, выбираются те, которые описывают рисунок с минимальным количеством удаленных ребер.

Применим метод локальной оптимизации для формирования множества топологического рисунка плоской части непланарного графа с минимальным числом удаленных ребер.

В качестве примера рассмотрим граф $G_1$ представленный на рис. 8.1.

Рассмотрим несколько вариантов выделения плоской части графа $G_1$. Для графа $G_1$ цикломатическое число $\nu(G_7) = m - n + 1 = 10$.

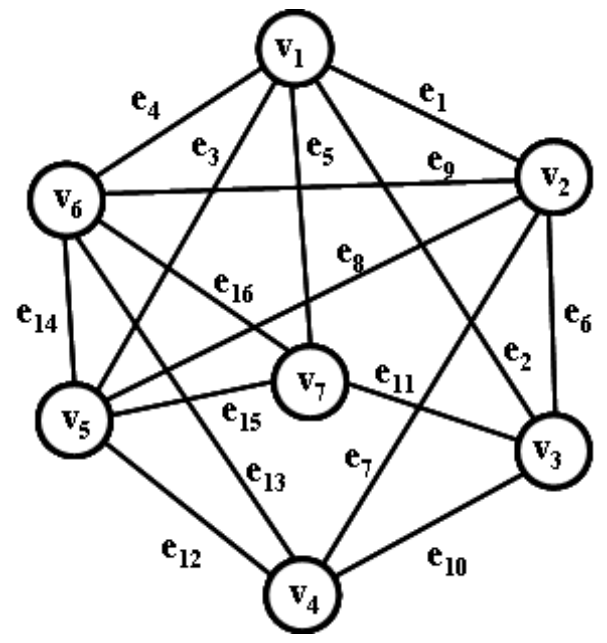

Рис. 8.1. Граф $G_9$.

Множество изометрических циклов графа $G_1$ имеет вид:

$c_1 = \{e_1,e_2,e_6\} \rightarrow \{v_1,v_2,v_3\}$;
$c_2 = \{e_1,e_3,e_8\} \rightarrow \{v_1,v_2,v_5\}$;
$c_3 = \{e_1,e_4,e_9\} \rightarrow \{v_1,v_2,v_6\}$;
$c_4 = \{e_2,e_3,e_{10},e_{12}\} \rightarrow \{v_1,v_3,v_4,v_5\}$;
$c_5 = \{e_2,e_4,e_{10},e_{13}\} \rightarrow \{v_1,v_3,v_4,v_6\}$;
$c_6 = \{e_2,e_5,e_{11}\} \rightarrow \{v_1,v_3,v_7\}$;
$c_7 = \{e_3,e_4,e_{14}\} \rightarrow \{v_1,v_5,v_6\}$;
$c_8 = \{e_3,e_5,e_{15}\} \rightarrow \{v_1,v_5,v_7\}$;
$c_9 = \{e_4,e_5,e_{16}\} \rightarrow \{v_1,v_6,v_7\}$;
$c_{10} = \{e_6,e_7,e_{10}\} \rightarrow \{v_2,v_3,v_4\}$;
$c_{11} = \{e_6,e_8,e_{11},e_{15}\} \rightarrow \{v_2,v_3,v_5,v_7\}$;
$c_{12} = \{e_6,e_9,e_{11},e_{16}\} \rightarrow \{v_2,v_3,v_6,v_7\}$;
$c_{13} = \{e_7,e_8,e_{12}\} \rightarrow \{v_2,v_4,v_5\}$;
$c_{14} = \{e_7,e_9,e_{13}\} \rightarrow \{v_2,v_4,v_6\}$;
$c_{15} = \{e_8,e_9,e_{14}\} \rightarrow \{v_2,v_5,v_6\}$;
$c_{16} = \{e_{10},e_{11},e_{12},e_{15}\} \rightarrow \{v_3,v_4,v_5,v_7\}$;
$c_{17} = \{e_{10},e_{11},e_{13},e_{16}\} \rightarrow \{v_3,v_4,v_6,v_7\}$;
$c_{18} = \{e_{12},e_{13},e_{14}\} \rightarrow \{v_4,v_5,v_6\}$;
$c_{19} = \{e_{14},e_{15},e_{16}\} \rightarrow \{v_5,v_6,v_7\}$.

**Вариант 1**. Модифицированным алгоритмом Гаусса выделим базис изометрических циклов из следующей фрагментарной перестановки изометрических циклов (циклы базиса подпространства циклов выделены красным цветом).

$f_1 = <$**$c_3$,$c_{17}$,$c_{18}$,$c_{15}$,$c_{16}$,$c_{11}$**,$c_{19}$,**$c_5$,$c_2$,$c_6$**,$c_8$,$c_{12}$,$c_9$,$c_4$,**$c_{13}$**,$c_1$,$c_{10}$,$c_{14}$,$c_7>$

$b_{\tau 1} = \{c_3,c_{17},c_{18},c_{15},c_{16},c_{11},c_5,c_2,c_6,c_{13}\}$.

$c_3 = \{e_1,e_4,e_9\}$;
$c_{17} = \{e_{10},e_{11},e_{13},e_{16}\}$;
$c_{18} = \{e_{12},e_{13},e_{14}\}$;
$c_{15} = \{e_8,e_9,e_{14}\}$;
$c_{16} = \{e_{10},e_{11},e_{12},e_{15}\}$;
$c_{11} = \{e_6,e_8,e_{11},e_{15}\}$;
$c_5 = \{e_2,e_4,e_{10},e_{13}\}$;
$c_2 = \{e_1,e_3,e_8\}$;
$c_6 = \{e_2,e_5,e_{11}\}$;
$c_{13} = \{e_7,e_8,e_{12}\}$.

Методом градиентного спуска выделяем систему циклов с нулевым значением кубического функционала Маклейна

$$FP(\frac{\partial b_{\tau 1}^3}{\partial c_{17} \partial c_{11} \partial c_{13}} = 0 .$$

В результате выделена следующая плоская конфигурация:

$c_3 = \{e_1,e_4,e_9\} \rightarrow \{v_2,v_6,v_1\}$;
$c_{18} = \{e_{12},e_{13},e_{14}\} \rightarrow \{v_5,v_6,v_4\}$;
$c_{15} = \{e_8,e_9,e_{14}\} \rightarrow \{v_5,v_6,v_2\}$;

$c_{16} = \{e_{10},e_{11},e_{12},e_{15}\} \rightarrow \{v_4,v_5,v_7,v_3\}$;
$c_5 = \{e_2,e_4,e_{10},e_{13}\} \rightarrow \{v_3,v_4,v_6,v_1\}$;
$c_2 = \{e_1,e_3,e_8\} \rightarrow \{v_2,v_5,v_1\}$;
$c_6 = \{e_2,e_5,e_{11}\} \rightarrow \{v_3,v_7,v_1\}$.

Топологический рисунок представлен на рис. 8.2. Удалены ребра $e_6,e_7,e_{15}$.

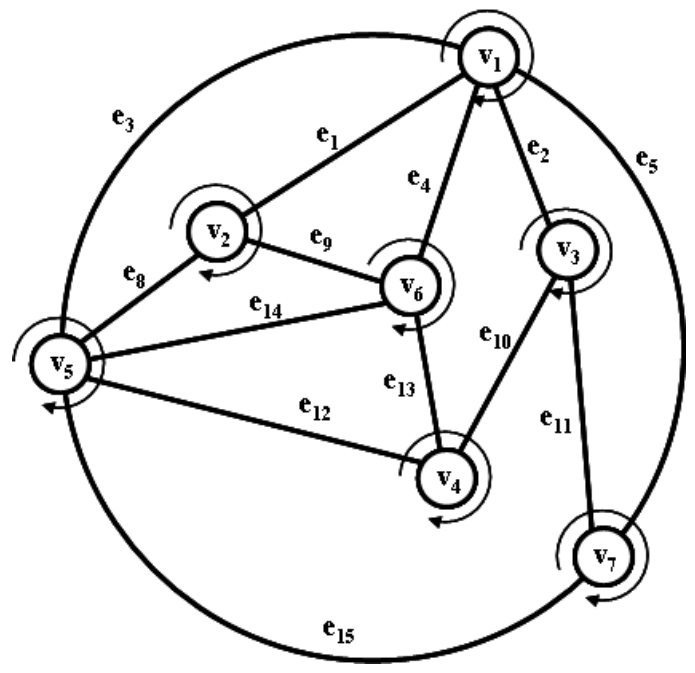

Рис. 8.2. Топологический рисунок плоской части (вариант 1).

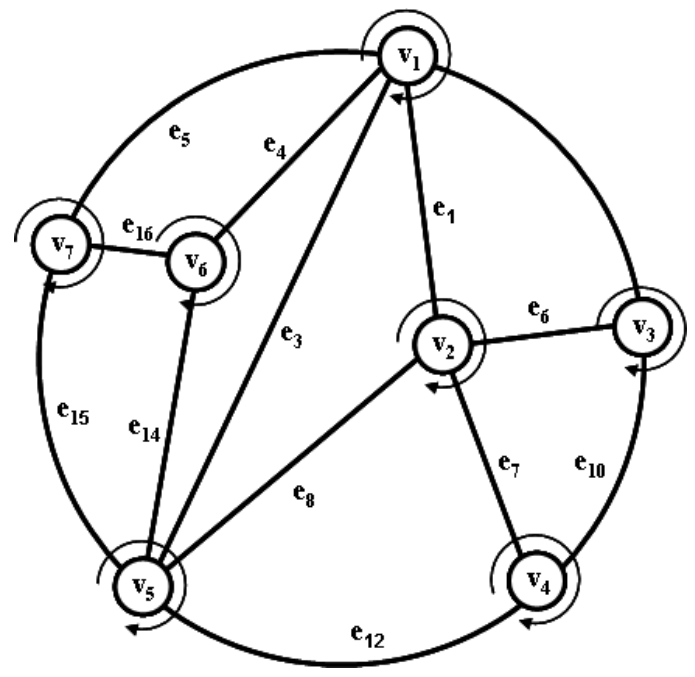

Рис. 8.3. Топологический рисунок плоской части (вариант 2).

**Вариант 2.** Модифицированным алгоритмом Гаусса выделим базис изометрических циклов из следующей фрагментарной перестановки изометрических циклов (циклы базиса подпространства циклов выделены красным цветом).

$f_2 = \langle \mathbf{c_{10},c_{12},c_5,c_7,c_1,c_{19},c_{13},c_2,c_{17},c_9}, c_3,c_{11},c_6,c_{15},c_8,c_{14},c_{18},c_{16},c_4 \rangle$.

$b_{\tau 2} = \{c_{10},c_{12},c_5,c_7,c_1,c_{19},c_{13},c_2,c_{17},c_9\}$.

$c_{10} = \{e_6,e_7,e_{10}\} \rightarrow \{v_2,v_3,v_4\}$;
$c_{12} = \{e_6,e_9,e_{11},e_{16}\} \rightarrow \{v_2,v_3,v_6,v_7\}$;
$c_5 = \{e_2,e_4,e_{10},e_{13}\} \rightarrow \{v_1,v_3,v_4,v_6\}$;
$c_7 = \{e_3,e_4,e_{14}\} \rightarrow \{v_1,v_5,v_6\}$;
$c_1 = \{e_1,e_2,e_6\} \rightarrow \{v_1,v_2,v_3\}$;
$c_{19} = \{e_{14},e_{15},e_{16}\} \rightarrow \{v_5,v_6,v_7\}$.
$c_{13} = \{e_7,e_8,e_{12}\} \rightarrow \{v_2,v_4,v_5\}$;
$c_2 = \{e_1,e_3,e_8\} \rightarrow \{v_1,v_2,v_5\}$;
$c_{17} = \{e_{10},e_{11},e_{13},e_{16}\} \rightarrow \{v_3,v_4,v_6,v_7\}$;
$c_9 = \{e_4,e_5,e_{16}\} \rightarrow \{v_1,v_6,v_7\}$.

Методом градиентного спуска выделяем систему циклов с нулевым значением кубического функционала Маклейна

$$FP(\frac{\partial b_{\tau 2}^3}{\partial c_{12} \partial c_5 \partial c_{17}} = 0.$$

Подмножество независимых изометрических циклов, имеющее нулевое значение кубического функционала Маклейна:

$c_{10} = \{e_6,e_7,e_{10}\} \rightarrow \{v_2,v_3,v_4\}$;
$c_7 = \{e_3,e_4,e_{14}\} \rightarrow \{v_1,v_5,v_6\}$;
$c_1 = \{e_1,e_2,e_6\} \rightarrow \{v_1,v_2,v_3\}$;
$c_{19} = \{e_{14},e_{15},e_{16}\} \rightarrow \{v_5,v_6,v_7\}$;
$c_{13} = \{e_7,e_8,e_{12}\} \rightarrow \{v_2,v_4,v_5\}$;
$c_2 = \{e_1,e_3,e_8\} \rightarrow \{v_1,v_2,v_5\}$;

$c_9 = \{e_4,e_5,e_{16}\} \to \{v_1,v_6,v_7\}$.

Топологический рисунок плоской части графа с удаленными ребрами $e_9,e_{11},e_{13}$ представлен на рис. 8.3.

**Вариант 3**. Модифицированным алгоритмом Гаусса выделим базис изометрических циклов из следующей фрагментарной перестановки изометрических циклов (циклы базиса подпространства циклов выделены красным цветом).

$f_3 = <\mathbf{c_7,c_{16},c_5,c_8,c_{19},c_2,c_{10},c_{17}},c_4,\mathbf{c_{12},c_{13}},c_{15},c_9,c_{18},c_1,c_{11},c_3,c_{14},c_6>$.

$b_{\tau 3} = \{c_7,c_{16},c_5,c_8,c_{19},c_2,c_{10},c_{17},c_{12},c_{13}\}$.

$c_7 = \{e_3,e_4,e_{14}\} \to \{v_1,v_5,v_6\}$;
$c_{16} = \{e_{10},e_{11},e_{12},e_{15}\} \to \{v_3,v_4,v_5,v_7\}$;
$c_5 = \{e_2,e_4,e_{10},e_{13}\} \to \{v_1,v_3,v_4,v_6\}$;
$c_8 = \{e_3,e_5,e_{15}\} \to \{v_1,v_5,v_7\}$;
$c_{19} = \{e_{14},e_{15},e_{16}\} \to \{v_5,v_6,v_7\}$;
$c_2 = \{e_1,e_3,e_8\} \to \{v_1,v_2,v_5\}$;
$c_{10} = \{e_6,e_7,e_{10}\} \to \{v_2,v_3,v_4\}$;
$c_{17} = \{e_{10},e_{11},e_{13},e_{16}\} \to \{v_3,v_4,v_6,v_7\}$;
$c_{12} = \{e_6,e_9,e_{11},e_{16}\} \to \{v_2,v_3,v_6,v_7\}$;
$c_{13} = \{e_7,e_8,e_{12}\} \to \{v_2,v_4,v_5\}$.

Методом градиентного спуска выделяем систему циклов с нулевым значением кубического функционала Маклейна

$$FP(\frac{\partial b_{\tau 3}^3}{\partial c_5 \partial c_8 \partial c_{17}} = 0.$$

Подмножество независимых изометрических циклов, имеющее нулевое значение кубического функционала Маклейна имеет вид:

$c_7 = \{e_3,e_4,e_{14}\} \to \{v_1,v_5,v_6\}$;
$c_{16} = \{e_{10},e_{11},e_{12},e_{15}\} \to \{v_3,v_4,v_5,v_7\}$;
$c_{19} = \{e_{14},e_{15},e_{16}\} \to \{v_5,v_6,v_7\}$;
$c_2 = \{e_1,e_3,e_8\} \to \{v_1,v_2,v_5\}$;
$c_{10} = \{e_6,e_7,e_{10}\} \to \{v_2,v_3,v_4\}$;
$c_{12} = \{e_6,e_9,e_{11},e_{16}\} \to \{v_2,v_3,v_6,v_7\}$;
$c_{13} = \{e_7,e_8,e_{12}\} \to \{v_2,v_4,v_5\}$.

Топологический рисунок плоской части графа с удаленными ребрами $e_2,e_5,e_{13}$ представлен на рис. 8.4.

**Вариант 4**. Модифицированным алгоритмом Гаусса выделим базис изометрических циклов из следующей фрагментарной перестановки изометрических циклов (циклы базиса подпространства циклов выделены красным цветом).

$f_4 = <\mathbf{c_{11},c_{17},c_{10},c_{19},c_{13},c_3},c_{16},\mathbf{c_9,c_2,c_{15}},c_{18},c_8,\mathbf{c_6},c_7,c_4,c_1,c_{12},c_5,c_{14}>$.

$b_{\tau 4} = \{c_{11},c_{17},c_{10},c_{19},c_{13},c_3,c_9,c_2,c_{15},c_6\}$.

$c_{11} = \{e_6,e_8,e_{11},e_{15}\} \to \{v_2,v_3,v_5,v_7\}$;
$c_{17} = \{e_{10},e_{11},e_{13},e_{16}\} \to \{v_3,v_4,v_6,v_7\}$;

$c_{10} = \{e_6, e_7, e_{10}\} \to \{v_2, v_3, v_4\}$;
$c_{19} = \{e_{14}, e_{15}, e_{16}\} \to \{v_5, v_6, v_7\}$;
$c_{13} = \{e_7, e_8, e_{12}\} \to \{v_2, v_4, v_5\}$;
$c_3 = \{e_1, e_4, e_9\} \to \{v_1, v_2, v_6\}$;
$c_9 = \{e_4, e_5, e_{16}\} \to \{v_1, v_6, v_7\}$;
$c_2 = \{e_1, e_3, e_8\} \to \{v_1, v_2, v_5\}$;
$c_{15} = \{e_8, e_9, e_{14}\} \to \{v_2, v_5, v_6\}$;
$c_6 = \{e_2, e_5, e_{11}\} \to \{v_1, v_3, v_7\}$.

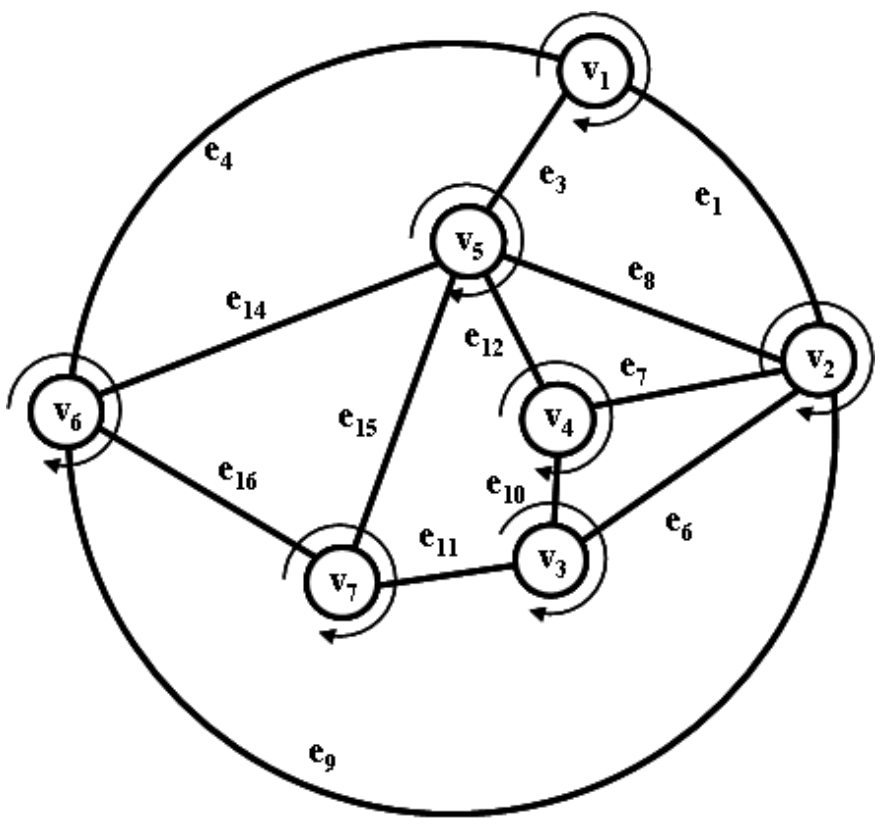

Рис. 8.4. Топологический рисунок плоской части (вариант 3).

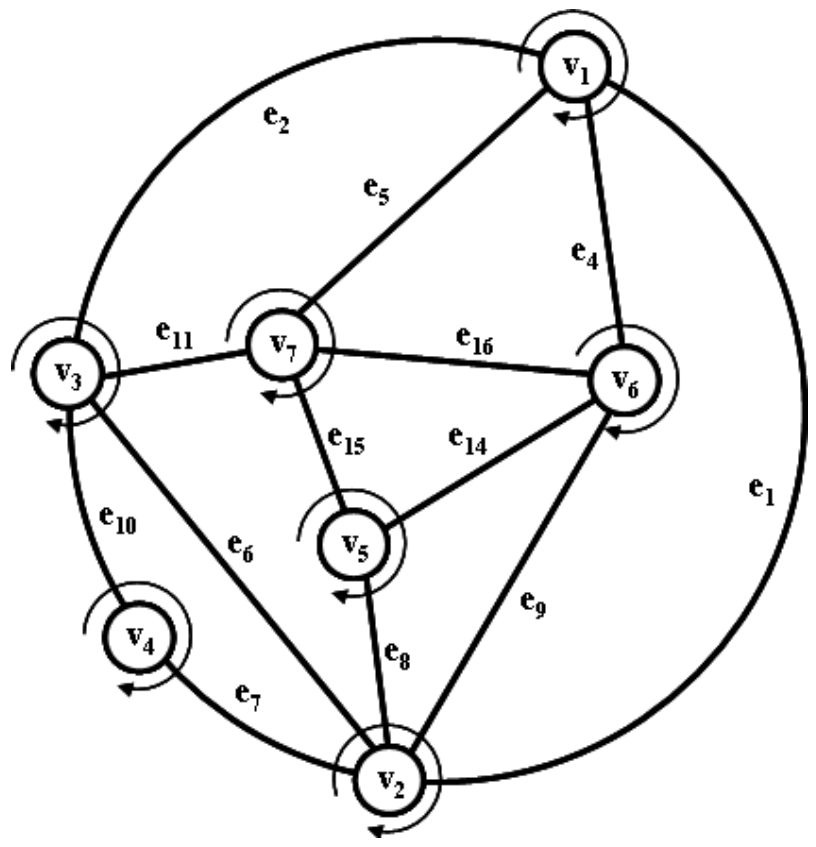

Рис. 8.5. Топологический рисунок плоской части (вариант 4).

Методом градиентного спуска выделяем систему циклов с нулевым значением кубического функционала Маклейна

$$FP(\frac{\partial b_{\tau 4}^3}{\partial c_{17} \partial c_{13} \partial c_2} = 0 .$$

Подмножество независимых изометрических циклов, имеющее нулевое значение кубического функционала Маклейна имеет вид:

$c_{11} = \{e_6, e_8, e_{11}, e_{15}\} \to \{v_2, v_3, v_5, v_7\}$;
$c_{10} = \{e_6, e_7, e_{10}\} \to \{v_2, v_3, v_4\}$;
$c_{19} = \{e_{14}, e_{15}, e_{16}\} \to \{v_5, v_6, v_7\}$;
$c_3 = \{e_1, e_4, e_9\} \to \{v_1, v_2, v_6\}$;
$c_9 = \{e_4, e_5, e_{16}\} \to \{v_1, v_6, v_7\}$;
$c_{15} = \{e_8, e_9, e_{14}\} \to \{v_2, v_5, v_6\}$;
$c_6 = \{e_2, e_5, e_{11}\} \to \{v_1, v_3, v_7\}$.

Топологический рисунок плоской части графа с удаленными ребрами $e_3, e_{12}, e_{13}$ представлен на рис. 8.5.

**Вариант 5.** Модифицированным алгоритмом Гаусса выделим базис изометрических циклов из следующей фрагментарной перестановки изометрических циклов (циклы базиса подпространства циклов выделены красным цветом).

$f_5 = \langle \mathbf{c_{11}, c_{14}, c_6, c_3, c_{16}, c_{19}, c_9, c_{10}, c_1}, c_{17}, \mathbf{c_8}, c_{18}, c_{13}, c_7, c_2, c_{15}, c_5, c_4, c_{12} \rangle$.

$b_{\tau 5} = \{c_{11}, c_{14}, c_6, c_3, c_{16}, c_{19}, c_9, c_{10}, c_1, c_8\}$.

$c_{11} = \{e_6,e_8,e_{11},e_{15}\} \to \{v_2,v_3,v_5,v_7\}$;
$c_{14} = \{e_7,e_9,e_{13}\} \to \{v_2,v_4,v_6\}$;
$c_6 = \{e_2,e_5,e_{11}\} \to \{v_1,v_3,v_7\}$;
$c_3 = \{e_1,e_4,e_9\} \to \{v_1,v_2,v_6\}$;
$c_{16} = \{e_{10},e_{11},e_{12},e_{15}\} \to \{v_3,v_4,v_5,v_7\}$;
$c_{19} = \{e_{14},e_{15},e_{16}\} \to \{v_5,v_6,v_7\}$;
$c_9 = \{e_4,e_5,e_{16}\} \to \{v_1,v_6,v_7\}$;
$c_{10} = \{e_6,e_7,e_{10}\} \to \{v_2,v_3,v_4\}$;
$c_1 = \{e_1,e_2,e_6\} \to \{v_1,v_2,v_3\}$;
$c_8 = \{e_3,e_5,e_{15}\} \to \{v_1,v_5,v_7\}$.

Методом градиентного спуска выделяем систему циклов с нулевым значением кубического функционала Маклейна

$$\text{FP}(\frac{\partial b_{\tau 5}^2}{\partial c_{11}\partial c_8} = 0.$$

Подмножество независимых изометрических циклов, имеющее нулевое значение кубического функционала Маклейна имеет вид:

$c_{14} = \{e_7,e_9,e_{13}\} \to \{v_2,v_4,v_6\}$;
$c_6 = \{e_2,e_5,e_{11}\} \to \{v_1,v_3,v_7\}$;
$c_3 = \{e_1,e_4,e_9\} \to \{v_1,v_2,v_6\}$;
$c_{16} = \{e_{10},e_{11},e_{12},e_{15}\} \to \{v_3,v_4,v_5,v_7\}$;
$c_{19} = \{e_{14},e_{15},e_{16}\} \to \{v_5,v_6,v_7\}$;
$c_9 = \{e_4,e_5,e_{16}\} \to \{v_1,v_6,v_7\}$;
$c_{10} = \{e_6,e_7,e_{10}\} \to \{v_2,v_3,v_4\}$;
$c_1 = \{e_1,e_2,e_6\} \to \{v_1,v_2,v_3\}$.

Данное подмножество циклов описывает топологический рисунок плоской части графа с удалением ребер $e_3,e_8$ представлен на рис. 8.5.

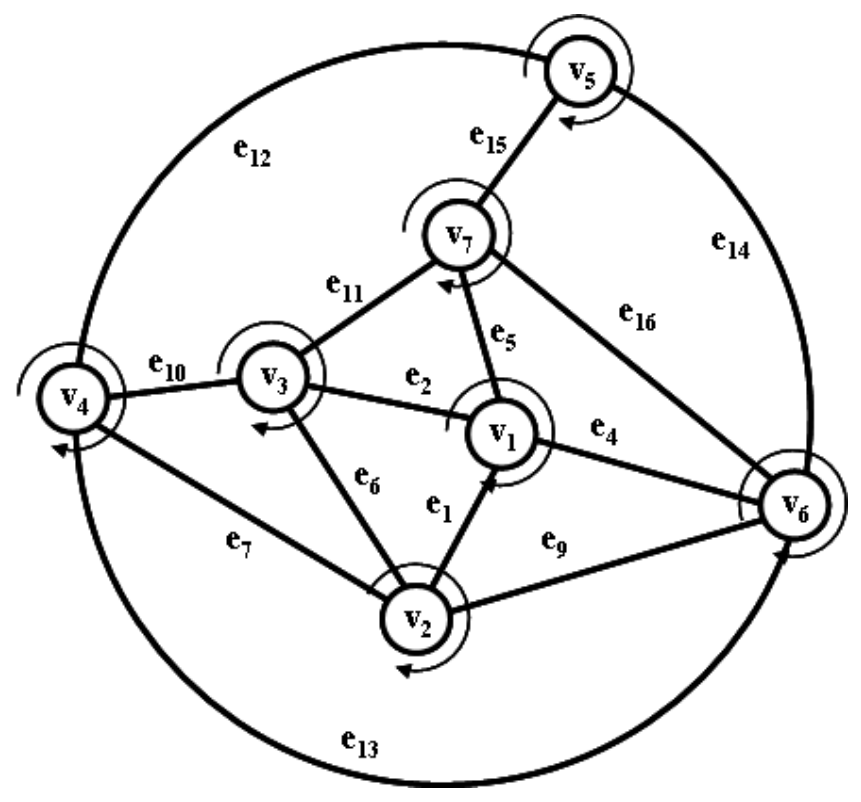

Рис. 8.5. Топологический рисунок плоской части (вариант 5).

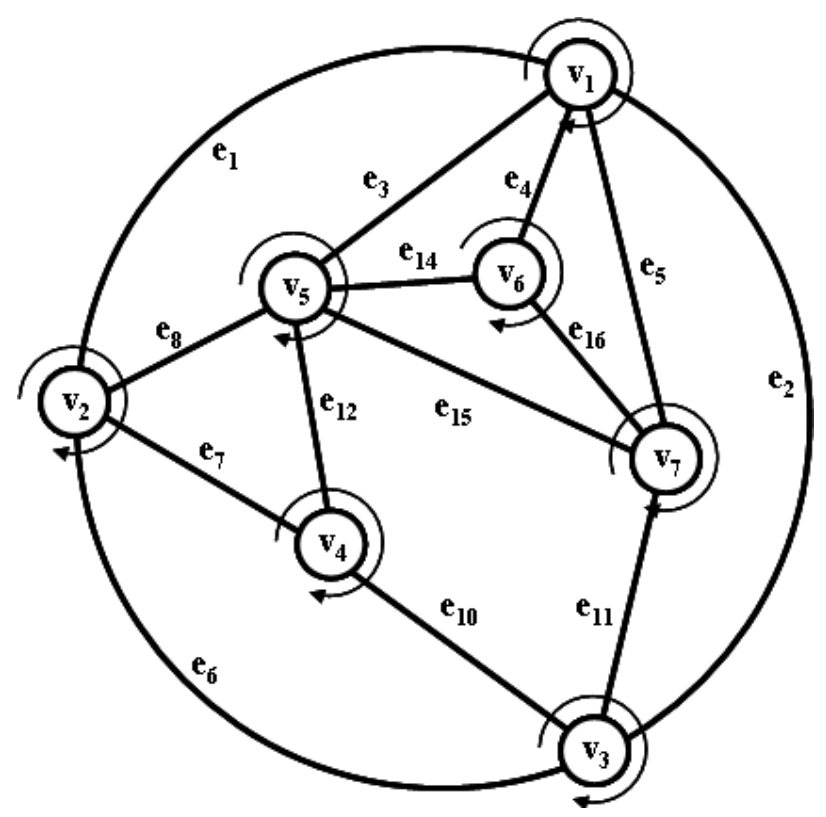

Рис. 8.7. Топологический рисунок плоской части (вариант 6).

**Вариант 5.** Модифицированным алгоритмом Гаусса выделим базис изометрических циклов из следующей фрагментарной перестановки изометрических циклов (циклы базиса подпространства циклов выделены красным цветом).

$f_6 = \langle \mathbf{c_9},\mathbf{c_{13}},\mathbf{c_{19}},\mathbf{c_{15}},\mathbf{c_{18}},\mathbf{c_6},\mathbf{c_{16}},\mathbf{c_5},\mathbf{c_7},c_2,\mathbf{c_{10}},c_4,c_{11},c_8,c_{17},c_1,c_{14},c_{12},c_3\rangle$.

$b_{\tau 6}$ = {$c_9$,$c_{13}$,$c_{19}$,$c_{15}$,$c_{18}$,$c_6$,$c_{16}$,$c_5$,$c_7$,$c_{10}$}.

$c_9 = \{e_4,e_5,e_{16}\} \to \{v_1,v_6,v_7\}$;
$c_{13} = \{e_7,e_8,e_{12}\} \to \{v_2,v_4,v_5\}$;
$c_{19} = \{e_{14},e_{15},e_{16}\} \to \{v_5,v_6,v_7\}$;
$c_{15} = \{e_8,e_9,e_{14}\} \to \{v_2,v_5,v_6\}$;
$c_{18} = \{e_{12},e_{13},e_{14}\} \to \{v_4,v_5,v_6\}$;
$c_6 = \{e_2,e_5,e_{11}\} \to \{v_1,v_3,v_7\}$;
$c_{16} = \{e_{10},e_{11},e_{12},e_{15}\} \to \{v_3,v_4,v_5,v_7\}$;
$c_7 = \{e_3,e_4,e_{14}\} \to \{v_1,v_5,v_6\}$;
$c_2 = \{e_1,e_3,e_8\} \to \{v_1,v_2,v_5\}$;
$c_{10} = \{e_6,e_7,e_{10}\} \to \{v_2,v_3,v_4\}$.

Методом градиентного спуска выделяем систему циклов с нулевым значением кубического функционала Маклейна

$$FP(\frac{\partial b_{\tau 6}^2}{\partial c_{15}\partial c_{18}} = 0.$$

Подмножество независимых изометрических циклов, имеющее нулевое значение кубического функционала Маклейна имеет вид:

$c_9 = \{e_4,e_5,e_{16}\} \to \{v_1,v_6,v_7\}$;
$c_{13} = \{e_7,e_8,e_{12}\} \to \{v_2,v_4,v_5\}$;
$c_{19} = \{e_{14},e_{15},e_{16}\} \to \{v_5,v_6,v_7\}$;
$c_6 = \{e_2,e_5,e_{11}\} \to \{v_1,v_3,v_7\}$;
$c_{16} = \{e_{10},e_{11},e_{12},e_{15}\} \to \{v_3,v_4,v_5,v_7\}$;
$c_7 = \{e_3,e_4,e_{14}\} \to \{v_1,v_5,v_6\}$;
$c_2 = \{e_1,e_3,e_8\} \to \{v_1,v_2,v_5\}$;
$c_{10} = \{e_6,e_7,e_{10}\} \to \{v_2,v_3,v_4\}$.

Данное подмножество циклов описывает топологический рисунок плоской части графа с удалением ребер $e_9$,$e_{13}$ представлен на рис. 8.7.

**Вариант 7.** Модифицированным алгоритмом Гаусса выделим базис изометрических циклов из следующей фрагментарной перестановки изометрических циклов (циклы базиса подпространства циклов выделены красным цветом).

$f_7$ = <**$c_{16}$**,**$c_4$**,**$c_9$**,**$c_7$**,**$c_{19}$**,$c_8$,**$c_{11}$**,**$c_{18}$**,**$c_2$**,**$c_{10}$**,$c_{17}$,$c_1$,**$c_{12}$**,$c_{15}$,$c_5$,$c_{13}$,$c_3$,$c_6$,$c_{14}$>.

$b_{\tau 7}$ = {$c_{16}$,$c_4$,$c_9$,$c_7$,$c_{19}$,$c_{11}$,$c_{18}$,$c_2$,$c_{10}$,$c_{12}$}.

$c_{16} = \{e_{10},e_{11},e_{12},e_{15}\} \to \{v_3,v_4,v_5,v_7\}$;
$c_4 = \{e_2,e_3,e_{10},e_{12}\} \to \{v_1,v_3,v_4,v_5\}$;
$c_9 = \{e_4,e_5,e_{16}\} \to \{v_1,v_6,v_7\}$;
$c_7 = \{e_3,e_4,e_{14}\} \to \{v_1,v_5,v_6\}$;
$c_{19} = \{e_{14},e_{15},e_{16}\} \to \{v_5,v_6,v_7\}$;
$c_{11} = \{e_6,e_8,e_{11},e_{15}\} \to \{v_2,v_3,v_5,v_7\}$;
$c_{18} = \{e_{12},e_{13},e_{14}\} \to \{v_4,v_5,v_6\}$;
$c_2 = \{e_1,e_3,e_8\} \to \{v_1,v_2,v_5\}$;
$c_{10} = \{e_6,e_7,e_{10}\} \to \{v_2,v_3,v_4\}$;
$c_{12} = \{e_6,e_9,e_{11},e_{16}\} \to \{v_2,v_3,v_6,v_7\}$.

Методом градиентного спуска выделяем систему циклов с нулевым значением кубического функционала Маклейна

$$FP(\frac{\partial b_{\tau 7}^{4}}{\partial c_{16}\partial c_{4}\partial c_{18}\partial c_{12}} = 0.$$

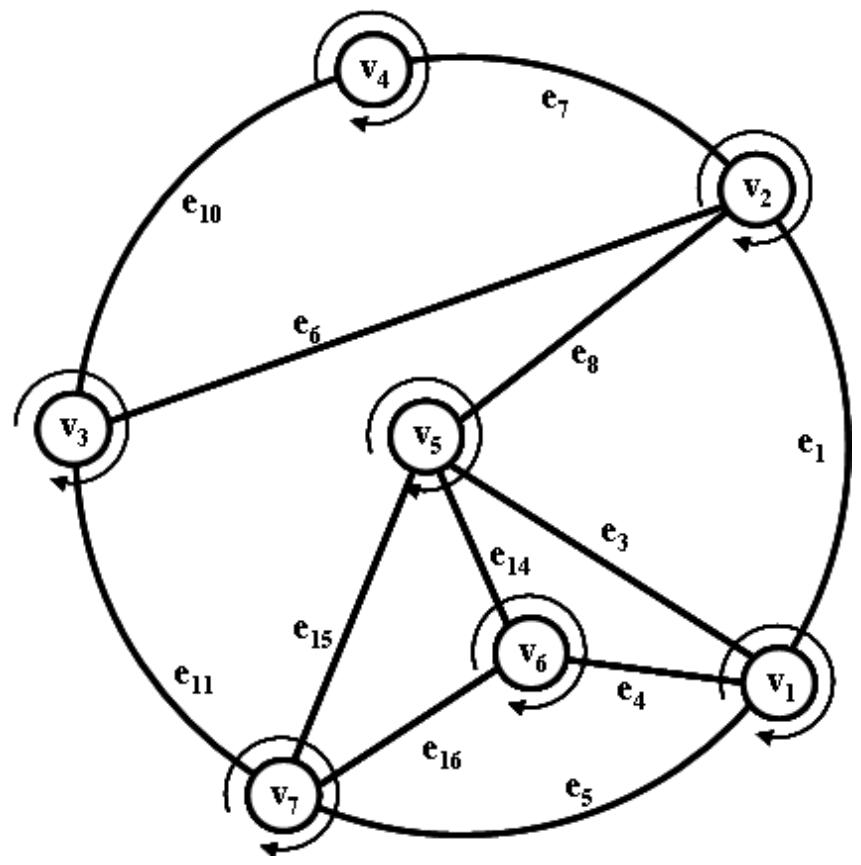

Рис. 8.8. Топологический рисунок плоской части (вариант 7).

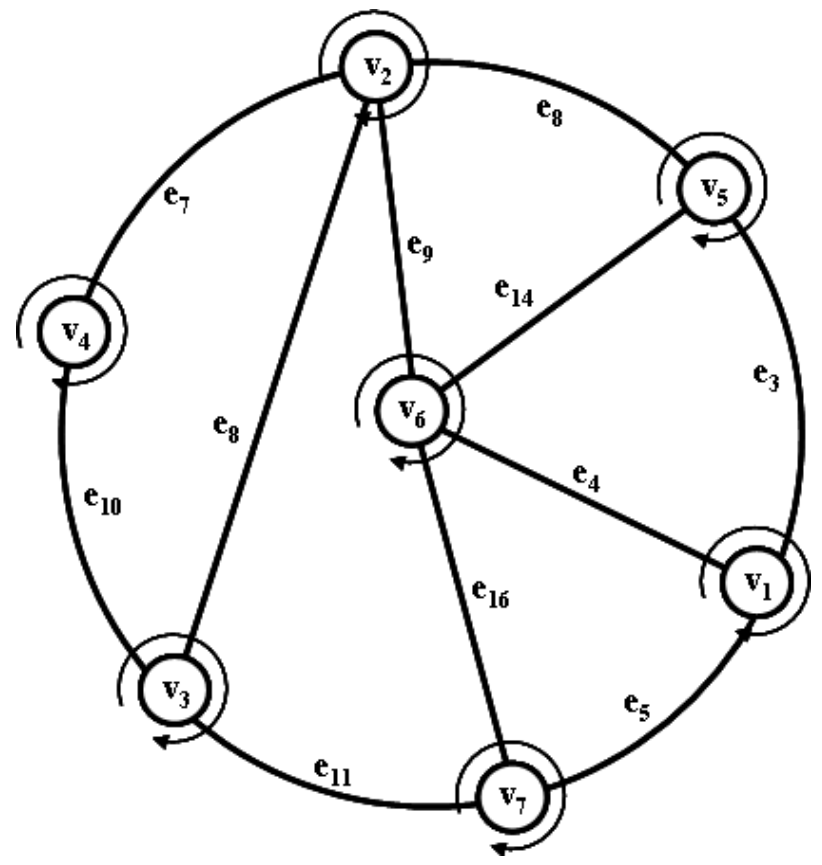

Рис. 8.9. Топологический рисунок плоской части (вариант 8).

Подмножество независимых изометрических циклов, имеющее нулевое значение кубического функционала Маклейна имеет вид:

$c_9 = \{e_4,e_5,e_{16}\} \to \{v_1,v_6,v_7\}$;
$c_7 = \{e_3,e_4,e_{14}\} \to \{v_1,v_5,v_6\}$;
$c_{19} = \{e_{14},e_{15},e_{16}\} \to \{v_5,v_6,v_7\}$;
$c_{11} = \{e_6,e_8,e_{11},e_{15}\} \to \{v_2,v_3,v_5,v_7\}$;
$c_2 = \{e_1,e_3,e_8\} \to \{v_1,v_2,v_5\}$;
$c_{10} = \{e_6,e_7,e_{10}\} \to \{v_2,v_3,v_4\}$.

Данное подмножество циклов описывает топологический рисунок плоской части графа с удалением ребер $e_2,e_9,e_{12},e_{13}$ представлено на рис. 8.8.

**Вариант 8.** Модифицированным алгоритмом Гаусса выделим базис изометрических циклов из следующей фрагментарной перестановки изометрических циклов (циклы базиса подпространства циклов выделены красным цветом).

$f_8 = \langle \mathbf{c_9,c_{15},c_{14},c_{11},c_1,c_{12},c_2,c_{16},c_7,c_{10}},c_3,c_{19},c_6,c_{13},c_{17},c_{18},c_8,c_5,c_4 \rangle$.

$b_{\tau 8} = \{c_9,c_{15},c_{14},c_{11},c_1,c_{12},c_2,c_{16},c_7,c_{10}\}$.

$c_9 = \{e_4,e_5,e_{16}\} \to \{v_1,v_6,v_7\}$;
$c_{15} = \{e_8,e_9,e_{14}\} \to \{v_2,v_5,v_6\}$;
$c_{14} = \{e_7,e_9,e_{13}\} \to \{v_2,v_4,v_6\}$;
$c_{11} = \{e_6,e_8,e_{11},e_{15}\} \to \{v_2,v_3,v_5,v_7\}$;
$c_1 = \{e_1,e_2,e_6\} \to \{v_1,v_2,v_3\}$;
$c_{12} = \{e_6,e_9,e_{11},e_{16}\} \to \{v_2,v_3,v_6,v_7\}$;
$c_2 = \{e_1,e_3,e_8\} \to \{v_1,v_2,v_5\}$;
$c_{16} = \{e_{10},e_{11},e_{12},e_{15}\} \to \{v_3,v_4,v_5,v_7\}$;
$c_7 = \{e_3,e_4,e_{14}\} \to \{v_1,v_5,v_6\}$;
$c_{10} = \{e_6,e_7,e_{10}\} \to \{v_2,v_3,v_4\}$.

Методом градиентного спуска выделяем систему циклов с нулевым значением кубического функционала Маклейна

$$FP(\frac{\partial b_{\tau 8}^{5}}{\partial c_{14}\partial c_{11}\partial c_{1}\partial c_{2}\partial c_{16}} = 0.$$

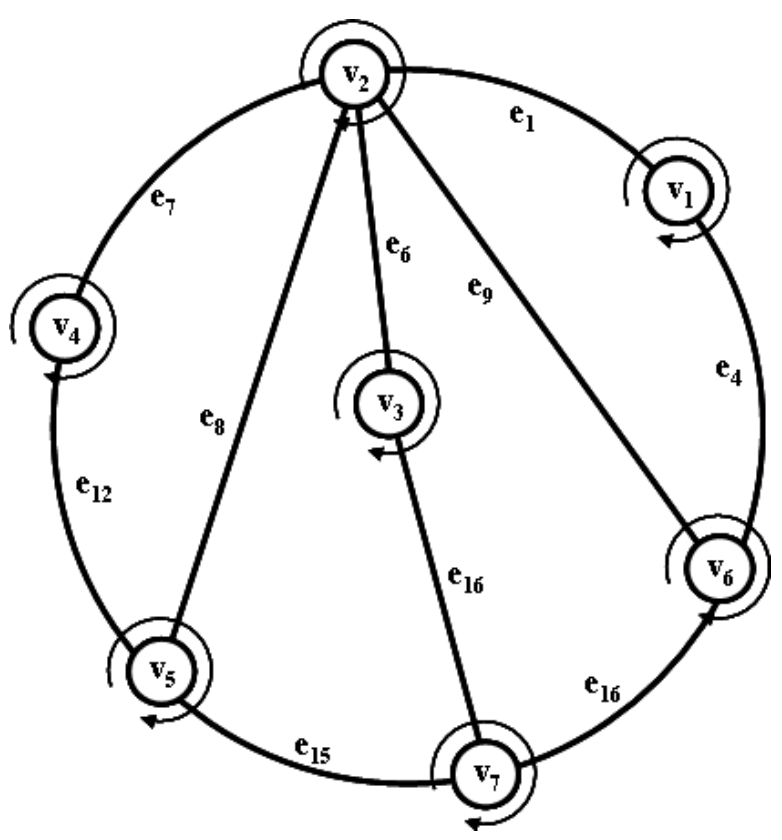

Рис. 8.10. Топологический рисунок плоской части (вариант 9).

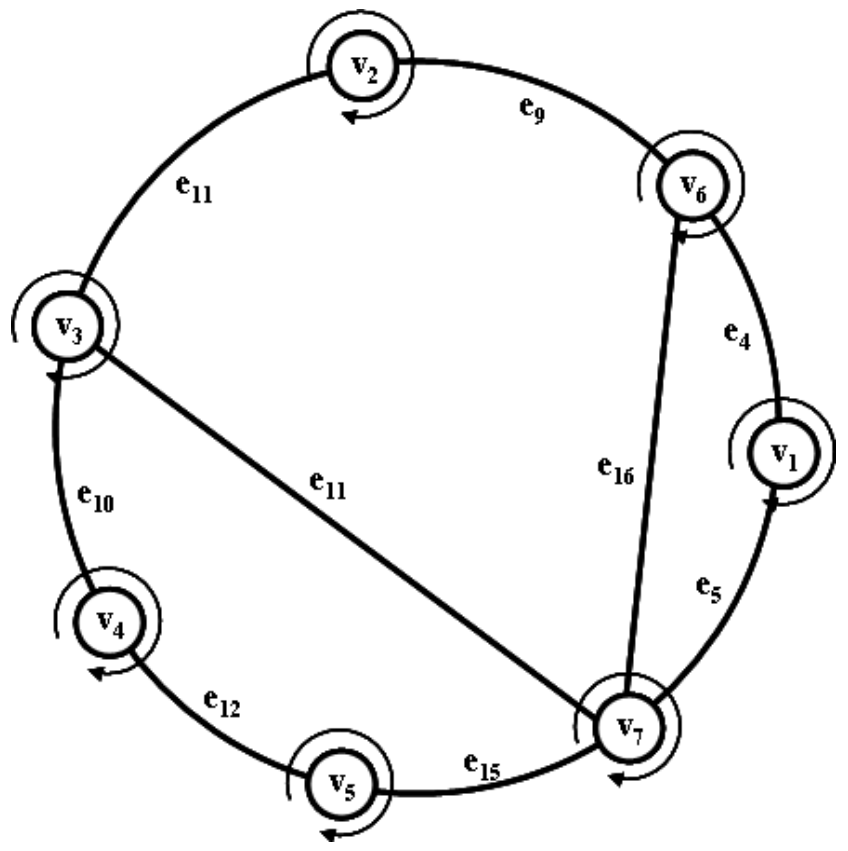

Рис. 8.11. Топологический рисунок плоской части (вариант 10).

Подмножество независимых изометрических циклов, имеющее нулевое значение кубического функционала Маклейна имеет вид:

$c_9 = \{e_4,e_5,e_{16}\} \to \{v_1,v_6,v_7\}$;
$c_{15} = \{e_8,e_9,e_{14}\} \to \{v_2,v_5,v_6\}$;
$c_{12} = \{e_6,e_9,e_{11},e_{16}\} \to \{v_2,v_3,v_6,v_7\}$;
$c_7 = \{e_3,e_4,e_{14}\} \to \{v_1,v_5,v_6\}$;
$c_{10} = \{e_6,e_7,e_{10}\} \to \{v_2,v_3,v_4\}$.

Данное подмножество циклов описывает топологический рисунок плоской части графа с удалением ребер $e_1,e_2,e_{12},e_{13},e_{15}$ представлен на рис. 8.9.

**Вариант 9.** Модифицированным алгоритмом Гаусса выделим базис изометрических циклов из следующей фрагментарной перестановки изометрических циклов (циклы базиса подпространства циклов выделены красным цветом).

$f_9 = <$**$c_3,c_{18},c_4,c_{16},c_{11},c_6,c_{17}$**$,c_8,$**$c_{13},c_{12},c_5$**$,c_2,c_{19},c_{10},c_9,c_{15},c_1,c_{14},c_7>$.

$b_{\tau 9} = \{c_3,c_{18},c_4,c_{16},c_{11},c_6,c_{17},c_{13},c_{12},c_5\}$.

$c_3 = \{e_1,e_4,e_9\} \to \{v_1,v_2,v_6\}$;
$c_{18} = \{e_{12},e_{13},e_{14}\} \to \{v_4,v_5,v_6\}$;
$c_4 = \{e_2,e_3,e_{10},e_{12}\} \to \{v_1,v_3,v_4,v_5\}$;
$c_{16} = \{e_{10},e_{11},e_{12},e_{15}\} \to \{v_3,v_4,v_5,v_7\}$;
$c_{11} = \{e_6,e_8,e_{11},e_{15}\} \to \{v_2,v_3,v_5,v_7\}$;
$c_6 = \{e_2,e_5,e_{11}\} \to \{v_1,v_3,v_7\}$;
$c_{17} = \{e_{10},e_{11},e_{13},e_{16}\} \to \{v_3,v_4,v_6,v_7\}$;
$c_{13} = \{e_7,e_8,e_{12}\} \to \{v_2,v_4,v_5\}$;
$c_{12} = \{e_6,e_9,e_{11},e_{16}\} \to \{v_2,v_3,v_6,v_7\}$;
$c_5 = \{e_2,e_4,e_{10},e_{13}\} \to \{v_1,v_3,v_4,v_6\}$.

Методом градиентного спуска выделяем систему циклов с нулевым значением кубического функционала Маклейна

$$FP(\frac{\partial b_{\tau 9}^{6}}{\partial c_{18}\partial c_{4}\partial c_{16}\partial c_{6}\partial c_{17}\partial c_{5}} = 0.$$

Подмножество независимых изометрических циклов, имеющее нулевое значение кубического функционала Маклейна имеет вид:

$c_3 = \{e_1,e_4,e_9\} \to \{v_1,v_2,v_6\}$;
$c_{11} = \{e_6,e_8,e_{11},e_{15}\} \to \{v_2,v_3,v_5,v_7\}$;
$c_{13} = \{e_7,e_8,e_{12}\} \to \{v_2,v_4,v_5\}$;
$c_{12} = \{e_6,e_9,e_{11},e_{16}\} \to \{v_2,v_3,v_6,v_7\}$.

Данное подмножество циклов описывает топологический рисунок плоской части графа с удалением ребер $e_2,e_3,e_5,e_{10},e_{13},e_{14}$ представлен на рис. 8.10.

**Вариант 10.** Модифицированным алгоритмом Гаусса выделим базис изометрических циклов из следующей фрагментарной перестановки изометрических циклов (циклы базиса подпространства циклов выделены красным цветом).

$f_{10}$ = <**$c_4$,$c_{18}$,$c_{16}$,$c_{15}$,$c_{11}$,$c_7$,$c_3$,$c_9$**,$c_5$,**$c_{13}$**,$c_{14}$,$c_1$,**$c_{12}$**,$c_2$,$c_{10}$,$c_8$,$c_{19}$,$c_6$,$c_{17}$>.

$b_{\tau 10} = \{c_4,c_{18},c_{16},c_{15},c_{11},c_7,c_3,c_9,c_{13},c_{12}\}$.
$c_4 = \{e_2,e_3,e_{10},e_{12}\} \to \{v_1,v_3,v_4,v_5\}$;
$c_{18} = \{e_{12},e_{13},e_{14}\} \to \{v_4,v_5,v_6\}$;
$c_{16} = \{e_{10},e_{11},e_{12},e_{15}\} \to \{v_3,v_4,v_5,v_7\}$;
$c_{15} = \{e_8,e_9,e_{14}\} \to \{v_2,v_5,v_6\}$;
$c_{11} = \{e_6,e_8,e_{11},e_{15}\} \to \{v_2,v_3,v_5,v_7\}$;
$c_7 = \{e_3,e_4,e_{14}\} \to \{v_1,v_5,v_6\}$;
$c_3 = \{e_1,e_4,e_9\} \to \{v_1,v_2,v_6\}$;
$c_9 = \{e_4,e_5,e_{16}\} \to \{v_1,v_6,v_7\}$;
$c_{13} = \{e_7,e_8,e_{12}\} \to \{v_2,v_4,v_5\}$;
$c_{12} = \{e_6,e_9,e_{11},e_{16}\} \to \{v_2,v_3,v_6,v_7\}$.

Методом градиентного спуска выделяем систему циклов с нулевым значением кубического функционала Маклейна

$$FP(\frac{\partial b_{\tau 10}^{7}}{\partial c_{4}\partial c_{18}\partial c_{15}\partial c_{11}\partial c_{7}\partial c_{3}\partial c_{13}} = 0.$$

Подмножество независимых изометрических циклов, имеющее нулевое значение кубического функционала Маклейна имеет вид:

$c_{16} = \{e_{10},e_{11},e_{12},e_{15}\} \to \{v_3,v_4,v_5,v_7\}$;
$c_9 = \{e_4,e_5,e_{16}\} \to \{v_1,v_6,v_7\}$;
$c_{12} = \{e_6,e_9,e_{11},e_{16}\} \to \{v_2,v_3,v_6,v_7\}$.

Данное подмножество циклов описывает топологический рисунок плоской части графа с удалением ребер $e_1,e_2,e_3,e_7,e_8,e_{13},e_{14}$ представлен на рис. 8.11.

Вращение вершин топологического рисунка графа изображенного на рис. 8.5

описывается следующей диаграммой [32]:

$\sigma(v_1)$: $v_2$ $v_6$ $v_7$ $v_3$
$\sigma(v_2)$: $v_1$ $v_4$ $v_3$ $v_5$ $v_6$
$\sigma(v_3)$: $v_1$ $v_7$ $v_2$ $v_4$
$\sigma(v_4)$: $v_2$ $v_3$
$\sigma(v_5)$: $v_7$ $v_6$ $v_2$
$\sigma(v_6)$: $v_1$ $v_2$ $v_5$ $v_7$
$\sigma(v_7)$: $v_1$ $v_6$ $v_5$ $v_3$

Естественно, что в процессе решения возможно получение различного количества удаленных ребер графа. Из множества полученных решений выбирают топологические рисунки плоской части графа с минимальным количеством удаленных ребер [18].

Вычислительная сложность данного метода в основном определяется сложностью выделения независимой системы циклов методом Гаусса $O(p^3)$, где $p$ – количество изометрических циклов. С целью, проверки изложенного метода выделения рисунка плоской части непланарного графа, близкого к максимально плоскому суграфу, разработана экспериментальная программа на языке Паскаль состоящая из следующих частей:

- построение множества изометрических циклов со случайным местоположением элементов;
- выделение независимой системы изометрических циклов модифицированным методом Гаусса;
- выделение подмножества циклов, с нулевым значением кубического функционала Маклейна методом наискорейшего спуска, характеризующее плоскую часть графа;
- создание базы данных характеризующих топологические рисунки плоской части непланарного графа, для проведения анализа и синтеза структур.

## 8.2. Фрагментарно-эволюционный алгоритм с геометрическим оператором кроссовера

Данный алгоритм использует стандартную схему эволюционного алгоритма на перестановках. В качестве базового множества решений выбирается множество $S_n$ перестановок из $n$ элементов [10,11].

Далее с помощью оператора отбора в текущей популяции выбирается множество пар перестановок для кроссовера. К каждой паре из выбранного множества пар применяется оператор кроссовера, а затем к результату кроссовера применяем оператор мутации. Пусть $X=<x_1,x_2,\ldots,x_n>$ и $Z=<z_1,z_2,,,,z_n>$ - две произвольные перестановки. Перестановка потомок строится по алгоритму: на начальном этапе последовательности X и Z пошагово просматривается от начала к концу. На k-ом шаге выбирается наименьшее из первых элементов последовательностей и затем добавляется в новую перестановку – потомок. Затем этот элемент удаляется из двух последовательностей – родителей. Например

$K_r$ = <4,7,3,2,8,1,6,5> <1,4,5,2,6,9,8,7> = <1,4,5,2,6,3,7,8>.

Таким образом, находится множество элементов – потомков. К промежуточной популяции, которая является объединением текущей популяции, и множества потомтов применяется оператор эволюции, выделяющий на этом множестве новую популяцию. Процесс эволюции повторяется до тех пор, пока не будет выполнено условие остановки эволюционного алгоритма.

В качестве примера выберем две перестановки, приводящие к построению топологического рисунка плоской части с максимальным числом удаленных ребер. В нашем случае – это перестановки

$f_9 = \langle c_3,c_{18},c_4,c_{16},c_{11},c_6,c_{17},c_8,c_{13},c_{12},c_5,c_2,c_{19},c_{10},c_9,c_{15},c_1,c_{14},c_7\rangle$

$f_{10} = \langle c_4,c_{18},c_{16},c_{15},c_{11},c_7,c_3,c_9,c_5,c_{13},c_{14},c_1,c_{12},c_2,c_{10},c_8,c_{19},c_6,c_{17}\rangle$

В результате перемешивания получим

$f_{11} = \langle c_3,c_4,c_{18},c_{16},c_{14},c_{11},c_6,c_{15},c_7,c_9,c_5,c_{13},c_1,c_{12},c_2,c_{10},c_8,c_{17},c_{19}\rangle$

Алгоритмом Гаусса выделим базис подпространства циклов C:

$c_3 = \{e_1,e_4,e_9\} \rightarrow \{v_1,v_2,v_6\}$;
$c_4 = \{e_2,e_3,e_{10},e_{12}\} \rightarrow \{v_1,v_3,v_4,v_5\}$;
$c_{18} = \{e_{12},e_{13},e_{14}\} \rightarrow \{v_4,v_5,v_6\}$;
$c_{16} = \{e_{10},e_{11},e_{12},e_{15}\} \rightarrow \{v_3,v_4,v_5,v_7\}$;
$c_{14} = \{e_7,e_9,e_{13}\} \rightarrow \{v_2,v_4,v_6\}$;
$c_{11} = \{e_6,e_8,e_{11},e_{15}\} \rightarrow \{v_2,v_3,v_5,v_7\}$;
$c_6 = \{e_2,e_5,e_{11}\} \rightarrow \{v_1,v_3,v_7\}$;
$c_{15} = \{e_8,e_9,e_{14}\} \rightarrow \{v_2,v_5,v_6\}$;
$c_7 = \{e_3,e_4,e_{14}\} \rightarrow \{v_1,v_5,v_6\}$;
$c_9 = \{e_4,e_5,e_{16}\} \rightarrow \{v_1,v_6,v_7\}$.

Методом градиентного спуска выделяем систему циклов с нулевым значением кубического функционала Маклейна

$$FP(\frac{\partial b^3_{\tau 11}}{\partial c_{14} \partial c_{18} \partial c_{11}} = 0.$$

Подмножество независимых изометрических циклов, имеющее нулевое значение кубического функционала Маклейна имеет вид:

$c_4 = \{e_2,e_3,e_{10},e_{12}\} \rightarrow \{v_1,v_3,v_4,v_5\}$;
$c_{16} = \{e_{10},e_{11},e_{12},e_{15}\} \rightarrow \{v_3,v_4,v_5,v_7\}$;
$c_6 = \{e_2,e_5,e_{11}\} \rightarrow \{v_1,v_3,v_7\}$;
$c_{15} = \{e_8,e_9,e_{14}\} \rightarrow \{v_2,v_5,v_6\}$;
$c_7 = \{e_3,e_4,e_{14}\} \rightarrow \{v_1,v_5,v_6\}$;
$c_9 = \{e_4,e_5,e_{16}\} \rightarrow \{v_1,v_6,v_7\}$.

Строим топологический рисунок графа с диаграммой вращения вершин.

Как мы видим, количество удаленных ребер в топологическом рисунке плоской части потомка уменьшилось на две единицы.

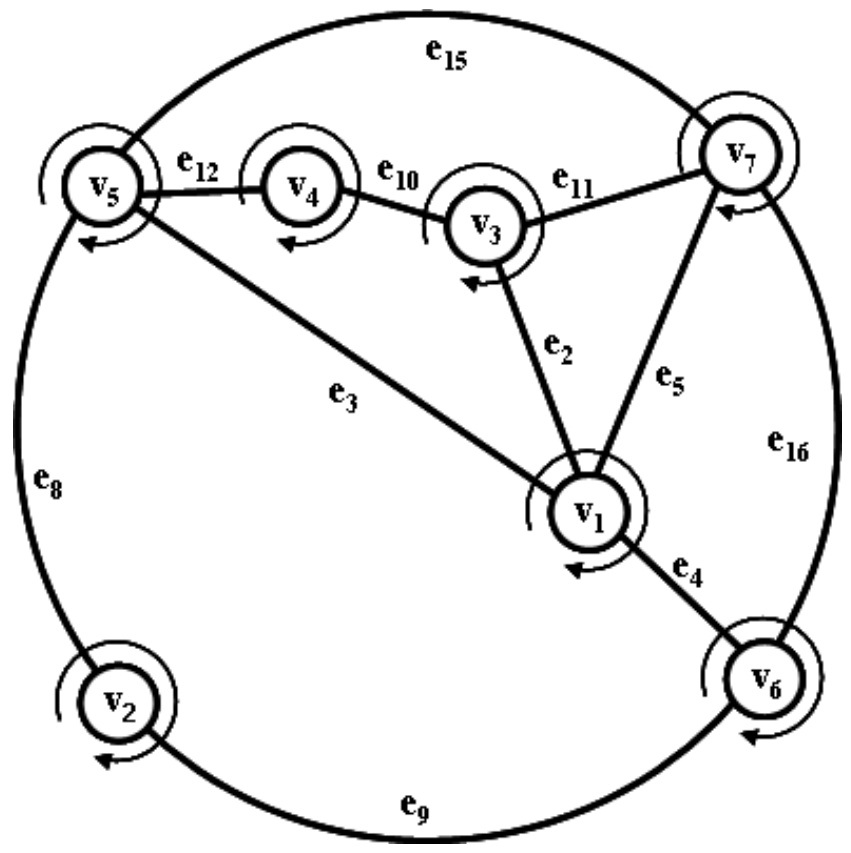

Рис. 8.12. Топологический рисунок графа.

Диаграмма вращения вершин

1: 3 7 6 5
2: 5 6
3: 1 4 7
4: 3 5
5: 1 2 7 4
6: 1 7 2
7: 1 3 5 6

При смешивании элементов в системе перестановок, с целью предотвращения конфликтных ситуаций, связанных с неравномерностью распределения элементов, выбор стоит производить с учетом вероятности его появления в последовательности.

## 8.3. Программа Ra1

Программа Ra1, предназначена для построения множества изометрических циклов графа. Построение изометрических циклов осуществляется на базе построения матрицы совместимостей графа.

```
type
      TMasy = array[1..1000] of integer;
      TMass = array[1..4000] of integer;
var
      F1,F2,F3,F4 : text;
      i,ii,j,jj,K,K1,Np,Nv,Kzikl,M,MakLin : integer;
      Ziklo,Nr,KKK,AB,K9 : integer;
      Masy: TMasy;
      Mass: TMass;
      Masi: TMass;
      MasyT: TMasy;
      MassT : TMass;
      MasMdop : TMasy;
      MasMy1 : TMasy;
      MasMs1 : TMass;
      MasMy2 : TMasy;
      MasMs2 : TMass;
      MasMy3 : TMasy;
      MasMs3 : TMass;
      MasMcg : TMasy;
      MasMcg1 : TMasy;
      MasKol : TMasy;
      Mass1 : TMass;
{*****************************************************************}
 procedure FormVolna(var Nv,Nv1,Nv2 : integer;
                    var My : TMasy;
                    var Ms : TMass;
                    var Mdop : TMasy);
{ Nv -  количество вершин в графе;                    }
{ Nv1 -  номер первой вершины;                        }
```

```
{ Nv2 -  номер второй вершины;                          }
{ My -  массив указателей для маттрицы смежностей графа;  }
{ Ms -  массив элементов матрицы смежностей;              }
{ Mdop -  массив глубины распространения волны.           }
{                                                   }
{    Данная процедура формирует массив распространения    }
{    волны, здесь номер уровень волны.                  }
 var I,Im,J,Kum,Istart,Istop : integer;
 label 1,2,3,4;
 begin
        for I:=1 to Nv do Mdop[I]:=0;
        Mdop[Nv1]:=1;
        Mdop[Nv2]:=2;
        Im:=2;
       1:  Im:=Im+1;
        for I:=1 to Nv do
         if Mdop[I]=0 then goto 2;
        goto 3;
       2:  for J:=1 to Nv do
        begin
           if Mdop[J]<>Im-1 then goto 4;
           Istart:=My[J];
           Istop:=My[J+1]-1;
           for I:=Istart to Istop do
           begin
              Kum:=Ms[I];
              if Mdop[Kum]=0 then Mdop[Kum]:=Im;
           end;
       4:  end;
        goto 1;
       3:;
       end; {FormVolna}
{*********************************************************}
 procedure FormKpris(var Kzikl : integer;
                        var Myy : TMasy;
                        var Mss : TMass;
                        var My : TMasy;
                        var Ms : TMass;
                        var Ms2 : TMass);
{ Kzikl  - количество t-циклов в графе;                  }
{ My  - массив указателей для матрицы смежностей;          }
{ Ms  - массив элементов матрицы смежностей;              }
{ Ms2 -  массив элементов матрицы инциденций;              }
{ Myy - массив указателей для матрицы t-циклов;           }
{ Mss - массив элементов матрицы t-циклов;                }
{                                                   }
{   Процедура переводит запись циклов в виде вершин        }
{   в запись в виде ребер.                            }
{                                                   }
   var I,J,JJ,Ip,Ip1,JJJ,Npn,KK : integer;
   label 1,2;
   begin
       for J:= 1 to Kzikl do
       begin
        for JJ:=Myy[J] to Myy[J+1]-2 do
        begin
          Ip:= Mss[JJ];
          Ip1:=Mss[JJ+1];
          for JJJ:=My[Ip] to My[IP+1]-1 do
           begin
              if Ms[JJJ]<>Ip1 then goto 1;
              Mss[JJ]:=Ms2[JJJ];
1:           end;
```

```
        end;
        Mss[Myy[J+1]-1]:=0;
        end;
        Npn:=Myy[Kzikl+1]-1;
        KK:=0;
        for I:=1 to Npn do
        begin
          if Mss[I]=0 then goto 2;
          KK:=KK+1;
          Mss[KK]:=Mss[I];
2:       end;
        for I:=1 to Kzikl do Myy[I+1]:=Myy[I+1]-I;
 end; {FormKpris}
{************************************************************}
 procedure FormDiz(var Mm1,Mm2,Mm4 : integer;
                   var M1 : TMasy;
                   var M2 : TMasy;
                   var M4 : TMasy);
 { Mm1 - количество элементов в первом цикле;              }
 { Mm2 - количество элементов во втором цикле;              }
 { Mm4 - количество элементов в их пересечении;             }
 { M1 - массив элементов первого цикла;                    }
 { M2 - массив элементов второго цикла;                    }
 { M4 - массив элементов пересечения.                      }
 {                                               }
 {   Процедура определения пересечения двух циклов         }
 {                                               }
 var J,I : integer;
 label 1,2;
       begin
        Mm4:=0;
        for I:=1 to Mm1 do
        begin
          for J:=1 to Mm2 do
          begin
            if M1[I]<>M2[J] then goto 1;
            Mm4:=Mm4+1;
            M4[Mm4]:=M2[J];
            goto 2;
       1:     end;
       2: end;
       end; {FormDiz}
{************************************************************}
 procedure FormSwigug(var Kzikl1,Kzikl2 : integer;
                   var My1 : TMasy;
                   var Ms1 : TMass;
                   var My2 : TMasy;
                   var Ms2 : TMass;
                   var Mcg : TMasy;
                   var Mcg1 : TMasy;
                   var Kol : Tmasy);
{ Kzikl1 - количество первых t-циклов в графе;            }
{ Kzikl2 - количество вторых t-циклов в графе;            }
{ My1 - массив указателей для первых t-циклов в графе;     }
{ Ms1 - массив для первых t-циклов в графе;               }
{ My2 - массив указателей для вторых t-циклов в графе;     }
{ Ms2 - массив для вторых t-циклов в графе;               }
{ Mcg -  вспомогательный массив;                        }
{ Mcg1 -  вспомогательный массив;                       }
{ Kol -  вспомогательный массив.                        }
{                                               }
{   Процедура формирования и сравнения t-циклов          }
{                                               }
```

```
       label 1;
       var I,J,JJ,II,NN1,NN2,NN3 : integer;
       begin
        for I:=1 to Kzikl1 do
        begin
          NN1:=My1[I+1]-My1[I];
          for II:=1 to NN1 do Mcg[II]:=Ms1[My1[I]-1+II];
          for J:=1 to Kzikl2 do
          begin
            NN2:=My2[J+1]-My2[J];
            for JJ:=1 to NN2 do Mcg1[JJ]:=Ms2[My2[J]-1+JJ];
            FormDiz(NN1,NN2,NN3,Mcg,Mcg1,Kol);
            if NN3=NN1 then goto 1;
          end;
          for II:=1 to NN1 do Ms1[My1[I]-1+II]:=0;
       1:  end;
       end; {FormSwigug}
{**************************************************************}
procedure FormDozas(var Kzikl,Kzikl1 : integer;
                    var My1 : TMasy;
                    var Ms1 : TMass;
                    var MyT : TMasy;
                    var MsT : TMass;
                    var Mcg : TMasy;
                    var Mcg1 : TMasy;
                    var Kol : Tmasy);
{ Kzikl  - количество t-циклов в графе;                    }
{ Ms1  - массив указателей для матрицы первых t-циклов;    }
{ Ms1  - массив элементов матрицы первых t-циклов;         }
{ MyT - массив указателей для матрицы t-циклов;            }
{ MsT - массив элементов матрицы t-циклов;                 }
{ Mcg -  вспомогательный массив;                           }
{ Mcg1 -  вспомогательный массив;                          }
{ Kol -  вспомогательный массив.                           }
{                                                          }
{    Процедура формирования t-циклов                       }
{                                                          }
label 1,2,3,4;
var I,J,JJ,II,III,Kp,Ks,NN1,NN2,NN3 : integer;
       begin
        if Kzikl=0  then goto 1;
        Kp:=Kzikl;
        Ks:=MyT[Kzikl+1]-1;
        for I:=1 to Kzikl1 do
        begin
          NN1:=My1[I+1]-My1[I];
          for II:=1 to NN1 do
          begin
            if Ms1[My1[I]-1+II]=0 then goto 2;
            Mcg[II]:=Ms1[My1[I]-1+II];
          end;
          for J:=1 to Kzikl do
          begin
            NN2:=MyT[J+1]-MyT[J];
            for JJ:=1 to NN2 do Mcg1[JJ]:=MsT[MyT[J]-1+JJ];
            FormDiz(NN1,NN2,NN3,Mcg,Mcg1,Kol);
            if NN3=NN1 then goto 2;
          end;
          Kp:=Kp+1;
          MyT[Kp+1]:=MyT[Kp]+NN1;
          for III:=1 to NN1 do
          begin
            Ks:=Ks+1;
```

```
          MsT[Ks]:=Mcg[III];
        end;
      2:  end;
       Kzikl:=Kp;
       goto 3;
      1:  Ks:=0;
       MyT[1]:=1;
       for I:=1 to Kzikl1 do
       begin
         NN1:=My1[I+1]-My1[I];
         for II:=1 to NN1 do
         begin
           if Ms1[My1[I]-1+II]=0 then goto 4;
           Ks:=Ks+1;
           MsT[Ks]:=Ms1[My1[I]-1+II];
         end;
         Kzikl:=Kzikl+1;
         MyT[Kzikl+1]:=MyT[Kzikl]+NN1;
      4:  end;
      3:;
      end; {FormDozas}
{*************************************************************}
 procedure FormSoasda(var Nv1,L,I2 : integer;
                     var Key1 : Boolean;
                     var My1 : TMasy;
                     var Ms1 : TMass;
                     var My : TMasy;
                     var Ms : TMass;
                     var Mcg : TMasy;
                     var Mcg1 : TMasy;
                     var Kol : TMasy);
{ Nv1 - номер вершины в графе;                             }
{ Key1  - признак;                                    }
{ I2 - признак;                                    }
{ L - длина цикла;                                 }
{ My  - массив указателей для матрицы смежностей;          }
{ Ms - массив для элементов мматрицы смежностей;           }
{ My1 - массив указателей для матрицы t-циклов;            }
{ Ms1 - массив элементов матрицы t-циклов;                 }
{ Mcg -  вспомогательный массив;                           }
{ Mcg1 -  вспомогательный массив;                          }
{ Kol -  вспомогательный массив.                           }
{                                               }
{    Процедура формирования цикла заданной длины           }
{                                               }
      label 2,3,14,11,5,6,12;
      var I,J,Ip1,Isu,Ip2 : integer;
      begin
       Key1:=false;
       if I2>0 then goto 2;
       for I:=1 to L do
       begin
         Mcg[I]:=1;
         Kol[I]:=My1[I+1]-My1[I];
       end;
       Mcg[L]:=0;
      2:  Mcg[L]:=Mcg[L]+1;
      14:;
       for I:=1 to L do
       begin
         if Mcg[I]<=Kol[I] then goto 3;
         if Mcg[1]>Kol[1] then goto 11;
         Mcg[I-1]:=Mcg[I-1]+1;
```

```
        for J:=I to L do Mcg[J]:=1;
        goto 14;
      3:  end;
       Mcg1[1]:=Nv1;
       Mcg1[L+1]:=Nv1;
       for I:=2 to L do
       begin
        Ip1:=My1[I]+Mcg[I]-1;
        Isu:=Ms1[Ip1];
        Mcg1[I]:=Isu;
        Ip2:=Mcg1[I-1];
        for J:=My[Ip2] to My[Ip2+1]-1 do
        begin
          if Ms[J]<>Isu then goto 6;
          goto 5;
      6:    end;
        Mcg[I]:=Mcg[I]+1;
        goto 14;
      5:  end;
       Key1:=true;
       I2:=1;
       goto 12;
      11: Key1:=false;
      12:;
      end; {FormSoasda}
{*************************************************************}
 procedure FormWegin(var Nv,Nv1,Nv2,Kzikl : integer;
                     var My : TMasy;
                     var Ms : TMass;
                     var My1 : TMasy;
                     var Ms1 : TMass;
                     var MyT : TMasy;
                     var MsT : TMass;
                     var Mdop : TMasy;
                     var Mcg : TMasy;
                     var Mcg1 : TMasy;
                     var Kol : Tmasy);
 { Nv - количество вершин в графе;                    }
 { Nv1 -  номер первой вершины;                       }
 { Nv2 -  номер второй вершины;                       }
 { Kzikl  - количество t-циклов в графе;              }
 { My  - массив указателей для матрицы смежностей;           }
 { Ms  - массив элементов матрицы смежностей;           }
 { My1 - массив для формирования указателей t-циклов;       }
 { Ms1 -  массив формирования t-циклов;                 }
 { MyT - массив указателей для матрицы t-циклов;          }
 { MsT - массив элементов матрицы t-циклов;             }
 { Mdop -  вспомогательный массив;                    }
 { Mcg -  вспомогательный массив;                     }
 { Mcg1 -  вспомогательный массив;                    }
 { Kol -  вспомогательный массив.                     }
 { Ms2 -  текущая начальная строка элементов;          }
 {                                          }
 {    Процедура построения единичных циклов            }
 {                                          }
 label 1,2,3,22,5,7,8,10,11,13,9;
 var I,J,JJ,II,Kot,Istart,Istop,JJJ,KKK,JI,Imum : integer;
 var Ji1,Kot1,Kot2,Isu,Irr,I2 : integer;
 var Key1 : Boolean;
   begin
      Kot:=0;
      MyT[1]:=1;
      Kzikl:=0;
```

```
for I:=1 to Nv do
begin
  if Mdop[I]<=2 then goto 1;
  for J:=My[I] to My[I+1]-1 do
  begin
    if Ms[J]<>Nv1 then goto 2;
    My1[1]:=1;
    Isu:=Mdop[I];
    My1[2]:=2;
    Ms1[1]:=Nv1;
    My1[3]:=3;
    Ms1[2]:=I;
    Irr:=Isu;
    Imum:=2;
3:      Irr:=Irr-1;
    if Irr=1 then goto 22;
    Istart:=My1[Isu-Irr+1];
    Istop:=My1[Isu-Irr+2]-1;
    for II:=1 to Nv do
    begin
      if Mdop[II]<>Irr then goto 5;
      for JJ:=My[II] to My[II+1]-1 do
      begin
        for JJJ:=Istart to Istop do
        begin
          if Ms[JJ]<>Ms1[JJJ] then goto 7;
          Imum:=Imum+1;
          Ms1[Imum]:=II;
7:          end;
      end;
5:    end;
     Kot1:=Istop+1-Imum;
     if Kot1=0 then goto 9;
     for II:=Istop+1 to Imum-1 do
     begin
       if Ms1[II]=0 then goto 8;
       for JJ:=II+1 to Imum do
        if Ms1[JJ]=Ms1[II] then Ms1[JJ]:=0;
8:    end;
9:    KKK:=Istop;
     for II:=Istop+1 to Imum do
     begin
       if Ms1[II]=0 then goto 10;
       KKK:=KKK+1;
       Ms1[KKK]:=Ms1[II];
10:    end;
     Imum:=KKK;
     My1[Isu-Irr+3]:=Imum+1;
     Kot2:=Isu-Irr+2;
     goto 3;
2:   end;
  goto 1;
22:  I2:=0;
11:  Kzikl:=Kzikl+1;
  FormSoasda(Nv1,Isu,I2,Key1,My1,Ms1,My,Ms,Mcg,Mcg1,Kol);
  if Key1=false then goto 13;
  MyT[Kzikl+1]:=MyT[Kzikl]+Isu+1;
  for Ji:=1 to Isu+1 do
  begin
    Kot:=Kot+1;
    MsT[Kot]:=Mcg1[JI];
  end;
  goto 11;
```

```
      13:   Kzikl:=Kzikl-1;
      1:  end;
    end; {FormWegin}
{**********************************************************}
procedure FormIncide(var Nv : integer;
                     var My : TMasy;
                     var Ms : TMass;
                     var Ms3 : TMass);
{ Nv - количество вершин в графе;                     }
{ My  - массив указателей для матрицы смежностей;          }
{ Ms  - массив элементов матрицы смежностей;               }
{ Ms3 - массив элементов матрицы инциденций.              }
{    Формируется матрица инциденций графа в массиве        }
{    Ms3.                                           }
{                                                 }
      var I,J,K,NNN,P,M,L : integer;
      begin
{     инициализация                               }
      NNN:=My[Nv+1]-1;
      K:=0;
      for J:= 1 to NNN do Ms3[J]:= 0;
{   определение номера элемента                         }
      for I:= 1 to Nv do
       for M:= My[I] to My[I+1]-1 do
        if Ms3[M]=0 then
        begin
          P:=Ms[M];
          K:=K+1;
          Ms3[M]:=K;
          for L:=My[P] to My[P+1]-1 do
           if Ms[L]=I then Ms3[L]:=K;
        end;
      end; {FormIncide}
{********************************************************}
 procedure EinZikle(var Nv,Kzikl,M : integer;
                    var Masy : TMasy;
                    var Mass : TMass;
                    var Mass1 : TMass;
                    var Masi : TMass;
                    var MasyT : TMasy;
                    var MassT : TMass;
                    var MasMy1 : TMasy;
                    var MasMs1 : TMass;
                    var MasMy2 : TMasy;
                    var MasMs2 : TMass;
                    var MasMy3 : TMasy;
                    var MasMs3 : TMass;
                    var MasMdop : Tmasy;
                    var MasMcg : Tmasy;
                    var MasMcg1 : Tmasy;
                    var MasKol : Tmasy);
{Процедура создания множества единичных циклов графа          }
{                                                }
{ Nv - количество вершин в графе;                      }
{ Kzikl  - количество единичных циклов в графе;              }
{ Masy  - массив указателей для матрицы смежностей;           }
{ Mass  - массив элементов матрицы смежностей;                }
{ Mass1 - массив для несмежных элементов строки.              }
{ Masi -  массив элементов матрицы инциденций;               }
{ MasyT - массив указателей для матрицы единичных циклов;    }
{ MassT - массив элементов матрицы единичных циклов;         }
{ MasMy1 -  вспомогательный массив;                     }
{ MasMs1 -  вспомогательный массив;                     }
```

```
{ MasMy2 -  вспомогательный массив;                          }
{ MasMs2 -  вспомогательный массив;                          }
{ MasMy3 -  вспомогательный массив;                          }
{ MasMs3 -  вспомогательный массив;                          }
{ MasMdop -  вспомогательный массив для хранения уровней;    }
{ MasMcg -  вспомогательный массив;                          }
{ MasMcg1 -  вспомогательный массив;                         }
{ MasKol -  вспомогательный массив.                          }
label 4;
var I,J,JJ,Nv1,Nv2,Pr,Kzikl1,Kzikl2,Ip1: integer;
begin
        Pr:= Masy[Nv+1]-1;
        M:= Pr div 2;
        FormIncide(Nv,Masy,Mass,Masi);
        Kzikl:= 0;
        MasyT[1]:= 1;
        for I:= 1 to M do
        begin {1}
        Ip1:= I;
          for J:= 1 to Nv do
          begin {2}
            for JJ:=Masy[J] to Masy[J+1]-1 do
            begin {3}
              if Masi[JJ] = Ip1 then
              begin {4}
                Nv1:= J; {определение первой вершины ребра}
                Nv2:= Mass[JJ]; {определение второй вершины ребра}
                goto 4; {концевые вершины ребра определены}
              end;   {4}
            end;   {3}
          end;   {2}
4:         FormVolna(Nv,Nv1,Nv2,Masy,Mass,MasMdop); {алгоритм поиска в ширину}
           {для ориентированного ребра (Nv1,Nv2)}
          FormWegin(Nv,Nv1,Nv2,Kzikl1,Masy,Mass,MasMy3,MasMs3,
          MasMy1,MasMs1,MasMdop,MasMcg,MasMcg1,MasKol);
{построение кандидатов в единичные циклы проходящих по ребру (Nv1,Nv2)}
          FormKpris(Kzikl1,MasMy1,MasMs1,Masy,Mass,Masi);
          FormVolna(Nv,Nv2,Nv1,Masy,Mass,MasMdop);
{алгоритм поиска в ширину для ориентированного ребра (Nv2,Nv1)}
          FormWegin(Nv,Nv2,Nv1,Kzikl2,Masy,Mass,MasMy3,MasMs3,
          MasMy2,MasMs2,MasMdop,MasMcg,MasMcg1,MasKol);
{построение кандидатов в единичные циклы проходящих по ребру (Nv2,Nv1)}
          FormKpris(Kzikl2,MasMy2,MasMs2,Masy,Mass,Masi);
          FormSwigug(Kzikl1,Kzikl2,MasMy1,MasMs1,MasMy2,MasMs2,
          MasMcg,MasMcg1,MasKol);
          FormDozas(Kzikl,Kzikl1,MasMy1,MasMs1,MasyT,MassT,
          MasMcg,MasMcg1,MasKol);
        end;   {1}
end;{EinZikle}
{**********************************************************}
 procedure FormVerBasis(var Nv,Ziklo : integer;
                     var My : TMasy;
                     var Ms : TMass;
                     var My1 : TMasy;
                     var Ms1 : TMass;
                     var Mass1 : TMass;
                     var Masi : TMass);
{Процедура записи базиса циклов через вершины                  }
{                                                              }
{ Nv - количество вершин в графе;                              }
{ Ziklo - цикломатическое число графа;                         }
{ My - массив указателей для матрицы смежностей;               }
{ Ms - массив элементов матрицы смежностей;                    }
```

```
{ My1 - массив указателей для матрицы базисных циклов             }
{ Ms1 : массив элементов матрицы базисных циклов;                 }
{ Masi - массив элементов матрицы инциденций;                     }
{ Mass1 - элементов матрицы базисных циклов через вершины;        }
{*****************************************************************}
label 1,2;
var i,j,ii,jj,iii,Nr,GG,NvS,NvK :integer;
begin
  GG:= 0;
  for i:= 1 to Ziklo do
    begin{1}
     {writeln(F2,'FormVerBasis: i = ',i);}
     Nr:= Ms1[My1[i]];
     {writeln(F2,'FormVerBasis: Nr = ',Nr);}
     for j:=1 to Nv do
     begin {2}
       for jj:= My[j] to My[j+1]-1 do
       begin {3}
        if Masi[jj]= Nr  then
        begin {4}
          NvS:= Ms[jj];
          NvK:= Ms[jj];
          {writeln(F2,'FormVerBasis: NvS = ',NvS,' NvK = ',NvK);}
          goto 1;
        end; {4}
       end;{3}
     end;{2}
1:      for ii:=My[NvK] to My[NvK+1]-1 do
         begin{3}
         {writeln(F2,'FormVerBasis: ii = ',ii);}
         for iii:=My1[i] to My1[i+1]-1 do
           begin {4}
            {writeln(F2,'FormVerBasis: iii = ',iii);}
            if Ms1[iii]=Masi[ii] then if Ms1[iii]<>Nr then
            begin {2}
             GG:=GG+1;
             Mass1[GG]:=NvK;
             Nr:=Ms1[iii];
             {writeln(F2,'FormVerBasis: Nr = ',Nr);
             writeln(F2,'FormVerBasis: Mass1[GG] = ',Mass1[GG]);
             writeln(F2,'FormVerBasis: NvK = ',NvK);
             writeln(F2,'FormVerBasis: Ms1[iii] = ',Ms1[iii]);
             writeln(F2,'FormVerBasis: Masi[ii] = ',Masi[ii]);}
             NvK:=Ms[ii];
             {writeln(F2,'FormVerBasis: NvK = ',NvK);}
             if NvK=NvS then goto 2 else goto 1;
            end; {2}
           end;{4}
         end;{3}
2:      end; {1}
  {writeln(F2,'Запись базиса через вершины [FormVerBasis]');}
  {for i:= 1 to Ziklo do
     begin
       for j:= My1[i] to My1[i+1]-1 do
       begin
        if j<> My1[i+1]-1 then write(F2, Mass1[j],' ');
        if j= My1[i+1]-1 then writeln(F2, Mass1[j],' ');
       end;
     end;}
end;{FormVerBasis}
{*****************************************************************}
procedure  Shell(var N : integer;
                 var A : TMasy);
```

```
{       процедура Шелла для упорядочивания элементов        }
{                                                 }
{     N - количество элементов в массиве;                    }
{     A - сортируемый массив;                           }
var D,Nd,I,J,L,X : integer;
label 1,2,3,4,5;
begin
  D:=1;
1:D:=2*D;
  if D<=N then goto 1;
2:D:=D-1;
  D:=D div 2;
  if D=0 then goto 5;
  Nd:=N-D;
  for I:=1 to Nd do
  begin
    J:=I;
3:  L:=J+D;
    if A[L]>=A[J] then goto 4;
    X:=A[J];
    A[J]:=A[L];
    A[L]:=X;
    J:=J-D;
    if J>0 then goto 3;
4:end;
  goto 2;
5:end;{Shell}
{**********************************************************}
procedure  ProzYpor(var N : integer;
           var Masy: TMasy;
           var Mass: TMass;
           var A : TMasy);
{Расположение элементов массива в порядке возрастания}
var i,j,K: integer;
begin
    for i:= 1 to N do
    begin
      K:=0;
      for j:= MasY[i] to MasY[i+1]-1 do
      begin
        K:=K+1;
        A[K]:= Mass[j];
      end;
      Shell(K,A);
      for j:= 1 to K do Mass[MasY[i]-1+j]:=A[j];
    end;
end;{ProzYpor}
{****************************************************************}
begin
         assign(F1,'D:\Delphi1\GRF\23.grf');
         reset(F1);
         readln(F1,Nv);
         for I:=1 to Nv+1 do
         begin
           if I<>Nv+1 then read(F1,Masy[I]);
           if I=Nv+1 then readln(F1,Masy[I]);
         end;
         Np:=Masy[Nv+1]-1;
         for I:=1 to Np do
         begin
           if I<>Np then read(F1,Mass[I]);
           if I=Np then read(F1,Mass[I]);
         end;
```

```
close (F1);
{ Создаём новый файл и открываем его в режиме "для чтения и записи"}
Assign(F2,'D:\Delphi1\GR1\23.gr1');
Rewrite(F2);
writeln(F2,Nv);
M:=Np div 2;
for I:=1 to Nv+1 do
begin
  if i<>Nv+1 then write(F2,Masy[i],' ');
  if i=Nv+1 then writeln(F2,Masy[i]);
end;
for I:=1 to Nv do
for j:=Masy[i] to Masy[i+1]-1 do
begin
  if j<>Masy[i+1]-1 then write(F2,Mass[j],' ');
  if j=Masy[i+1]-1 then writeln(F2,Mass[j]);
end;
writeln(F2,M);
FormIncide(Nv,Masy,Mass,Masi);
for I:=1 to Nv do
for j:=Masy[i] to Masy[i+1]-1 do
begin
  if j<>Masy[i+1]-1 then write(F2,Masi[j],' ');
  if j=Masy[i+1]-1 then writeln(F2,Masi[j]);
end;
close (F2);
{ Создаём новый файл и открываем его в режиме "для чтения и записи"}
Assign(F3,'D:\Delphi1\EZI\23.ezi');
Rewrite(F3);
{for I:=1 to Nv do
for j:=Masy[i] to Masy[i+1]-1 do
begin
  if j<>Masy[i+1]-1 then write(F3,Masi[j],' ');
  if j=Masy[i+1]-1 then writeln(F3,Masi[j]);
end;}
EinZikle(Nv,Kzikl,M,Masy,Mass,Mass1,Masi,MasyT,MassT,
        MasMy1,MasMs1,MasMy2,MasMs2,MasMy3,MasMs3,MasMdop,
        MasMcg,MasMcg1,MasKol);
ProzYpor(Kzikl,MasyT,MassT,MasKol);
FormVerBasis(Nv,Kzikl,Masy,Mass,MasyT,MassT,Mass1,Masi);
ProzYpor(Kzikl,MasyT,Mass1,MasKol);
writeln(F3,Kzikl);
for i:=1 to Kzikl+1 do
begin
   if i<> Kzikl+1 then write(F3, MasyT [i],' ');
   if i= Kzikl+1 then writeln(F3, MasyT [i]);
end;
for I:=1 to Kzikl do
begin
  for j:=MasyT[i] to MasyT[i+1]-1 do
  begin
    if j<>MasyT[i+1]-1 then write(F3,MassT[j],' ');
    if j=MasyT[i+1]-1 then writeln(F3,MassT[j]);
  end;
end;
for I:=1 to Kzikl do
begin
  for j:=MasyT[i] to MasyT[i+1]-1 do
  begin
    if j<>MasyT[i+1]-1 then write(F3,Mass1[j],' ');
    if j=MasyT[i+1]-1 then writeln(F3,Mass1[j]);
  end;
end;
```

```
    close (F3);
    writeln('Вес взят!');
end.
```

## 8.4. Входные и выходные файлы программы Ra1

Задается матрица смежностей графа. Разница элементов в массиве указателей определяет локальную степень вершины. Например, разница между третьи и вторым элементом массива указателей (11-6) определяет локальную степень вершины 2 (равно 5), разница между 5 и 4 элементоми массива указателей (19-15) определяет локальную степень вершины 4 (равно 4).

### Входной файл 7.grf

| | |
|---|---|
| 7 | количество вершин в графе |
| 1 6 11 15 19 24 29 33 | массив указателей |
| 2 3 5 6 7 | вершины смежные с вершиной 1 |
| 1 3 4 5 6 | вершины смежные с вершиной 2 |
| 1 2 4 7 | вершины смежные с вершиной 3 |
| 2 3 5 6 | вершины смежные с вершиной 4 |
| 1 2 4 6 7 | вернины смежные с вершиной 5 |
| 1 2 4 5 7 | вернины смежные с вершиной 6 |
| 1 3 5 6 | вернины смежные с вершиной 7 |

### Выходной файл 7.gr1

| | |
|---|---|
| 7 | количество вершин в графе |
| 1 6 11 15 19 24 29 33 | массив указателей |
| 2 3 5 6 7 | вершины смежные с вершиной 1 |
| 1 3 4 5 6 | вершины смежные с вершиной 2 |
| 1 2 4 7 | вершины смежные с вершиной 3 |
| 2 3 5 6 | вершины смежные с вершиной 4 |
| 1 2 4 6 7 | вернины смежные с вершиной 5 |
| 1 2 4 5 7 | вернины смежные с вершиной 6 |
| 1 3 5 6 | вернины смежные с вершиной 7 |
| 16 | количество ребер графа |
| 1 2 3 4 5 | ребра инцидентные вершине 1 |
| 1 6 7 8 9 | ребра инцидентные вершине 2 |
| 2 6 10 11 | ребра инцидентные вершине 3 |
| 7 10 12 13 | ребра инцидентные вершине 4 |
| 3 8 12 14 15 | ребра инцидентные вершине 5 |
| 4 9 13 14 16 | ребра инцидентные вершине 6 |
| 5 11 15 16 | ребра инцидентные вершине 7 |

### Выходной файл 7.ezi

| | |
|---|---|
| 19 | количество изометрических циклов |
| 1 4 7 10 14 18 21 24 27 30 33 | |
| 37 41 44 47 50 54 58 61 64 | массив указателей для изометрических циклов |
| 1 2 6 | изометрический цикл 1 в ребрах |
| 1 3 8 | изометрический цикл 2 в ребрах |
| 1 4 9 | изометрический цикл 3 в ребрах |
| 2 3 10 12 | изометрический цикл 4 в ребрах |
| 2 4 10 13 | изометрический цикл 5 в ребрах |
| 2 5 11 | изометрический цикл 6 в ребрах |
| 3 4 14 | изометрический цикл 7 в ребрах |
| 3 5 15 | изометрический цикл 8 в ребрах |
| 4 5 16 | изометрический цикл 9 в ребрах |
| 6 7 10 | изометрический цикл 10 в ребрах |
| 6 8 11 15 | изометрический цикл 11 в ребрах |
| 6 9 11 16 | изометрический цикл 12 в ребрах |
| 7 8 12 | изометрический цикл 13 в ребрах |
| 7 9 13 | изометрический цикл 14 в ребрах |
| 8 9 14 | изометрический цикл 15 в ребрах |
| 10 11 12 15 | изометрический цикл 16 в ребрах |
| 10 11 13 16 | изометрический цикл 17 в ребрах |
| 12 13 14 | изометрический цикл 18 в ребрах |
| 14 15 16 | изометрический цикл 19 в ребрах |
| 1 2 3 | изометрический цикл 1 в вершинах |
| 1 2 5 | изометрический цикл 2 в вершинах |
| 1 2 6 | изометрический цикл 3 в вершинах |
| 1 3 4 5 | изометрический цикл 4 в вершинах |
| 1 3 4 6 | изометрический цикл 5 в вершинах |
| 1 3 7 | изометрический цикл 6 в вершинах |
| 1 5 6 | изометрический цикл 7 в вершинах |
| 1 5 7 | изометрический цикл 8 в вершинах |
| 1 6 7 | изометрический цикл 9 в вершинах |
| 2 3 4 | изометрический цикл 10 в вершинах |
| 2 3 5 7 | изометрический цикл 11 в вершинах |
| 2 3 6 7 | изометрический цикл 12 в вершинах |
| 2 4 5 | изометрический цикл 13 в вершинах |
| 2 4 6 | изометрический цикл 14 в вершинах |
| 2 5 6 | изометрический цикл 15 в вершинах |
| 3 4 5 7 | изометрический цикл 16 в вершинах |
| 3 4 6 7 | изометрический цикл 17 в вершинах |
| 4 5 6 | изометрический цикл 18 в вершинах |
| 5 6 7 | изометрический цикл 19 в вершинах |

## 8.5. Программа Rebro2

Программа Rebro2 предназначена для выделения плоских конфигураций графа, методом случайного построения последовательности изометрических циклов.

```
program Rebro2;

{$APPTYPE CONSOLE}

type
    TMasy = array[1..100000] of integer;
    TMass = array[1..200000] of integer;
    XMass = array[1..100000] of real;
var
    F1,F2,F3,F4 : text;
    i,ii,j,jj,iii,K,K1,M,Nv,Kzikl,Dl,Tz,Ziklo : integer;
    KKK,Kzikl1,Priz,Makley,Priz1,Nh,KM : integer;
    Nr,Priz2,Priz3,Priz4,L,AB : integer;
    Masy: TMasy;
    Mass: TMass;
    Masi: TMass;
    MasyT: TMasy;
    MassT : TMass;
    Mass1 : TMass;
    MasMdop : TMasy;
    MasKol : TMasy;
    MasMy1 : TMasy;
    MasMs1 : TMass;
    MasMy2 : TMasy;
    MasMs2 : TMass;
    MasMy3 : TMasy;
    MasMs3 : TMass;
    MasMy4 : TMasy;
    MasMs4 : TMass;
    MasMcg : TMasy;
```

```
          MasMy5 : TMasy;
          MasMs5 : TMass;
          MasMcg1 : TMasy;
          MasMcg2 : TMasy;
          MasMcg3 : TMasy;
          MasMcg4 : XMass;
{*****************************************************************}
 procedure Unk22(var N,M,Kzikl : integer;
                        var MasyT : TMasy;
                        var MassT : TMass;
                        var MasMdop : TMasy;
                        var MasKol : TMasy);
 {   Расчет вектора количества циклов по ребру;                      }
 {                                                          }
 { N - количество вершин в графе;                              }
 { M - количество ребер в графе;                              }
 { Kzikl  - количество t-циклов в графе;                         }
 { MasyT - массив указателей для матрицы циклов;                      }
 { MassT - массив элементов матрицы t-циклов;                        }
 { MasKol - массив количества циклов проходящих по ребру               }
 { MasMdop -  массив признаков удаления цикла, 1- неудален, -1 - удален }
   var I,J,III : integer;
begin
{инициализация}
for  i:=1 to M do MasKol[i]:=0;
{определим количество циклов проходящих по каждому ребру для определения функционала Мак-Лейна с учетом удаленных циклов}
    for  i:=1 to Kzikl do
      begin {2}
        {writeln(F4,' i = ',i,' MasMdop[i] = ',MasMdop[i]); }
        if MasMdop[i] > 0 then
        begin
          for  j:= MasyT[i] to MasyT[i+1]-1 do
          begin
           MasKol[MassT[j]]:= MasKol[MassT[j]]+1;
           {writeln(F4,' j = ',j,' MasKol[MassT[j]] = ',MasKol[MassT[j]],
           ' MassT[j] = ',MassT[j]); }
          end;
        end;
      end; {2}
          {writeln(F4, 'Процедура Unk22');
          for i:= 1 to M do
          begin
             if i<> M then write(F4, MasKol[i],' ');
             if i= M then writeln(F4, MasKol[i]);
          end; }
end;{Unk22}
{*****************************************************************}
 procedure  Shell(var N : integer;
                 var A : TMass);
{       процедура Шелла для упорядочивания элементов        }
{                                                   }
{     N - количество элементов в массиве;                   }
{     A - сортируемый массив;                           }
var D,Nd,I,J,L,X : integer;
label 1,2,3,4,5;
begin
  D:=1;
1:D:=2*D;
  if D<=N then goto 1;
2:D:=D-1;
  D:=D div 2;
  if D=0 then goto 5;
```

```
  Nd:=N-D;
  for I:=1 to Nd do
  begin
    J:=I;
3:  L:=J+D;
    if A[L]>=A[J] then goto 4;
    X:=A[J];
    A[J]:=A[L];
    A[L]:=X;
    J:=J-D;
    if J>0 then goto 3;
4:end;
  goto 2;
5:end;{Shell}
{*******************************************************************}
 procedure  Shell1(var N : integer;
                 var A : TMasy);
{       процедура Шелла для упорядочивания элементов        }
{                                                  }
{     N - количество элементов в массиве;                   }
{     A - сортируемый массив;                          }
var D,Nd,I,J,L,X : integer;
label 1,2,3,4,5;
begin
  D:=1;
1:D:=2*D;
  if D<=N then goto 1;
2:D:=D-1;
  D:=D div 2;
  if D=0 then goto 5;
  Nd:=N-D;
  for I:=1 to Nd do
  begin
    J:=I;
3:  L:=J+D;
    if A[L]>=A[J] then goto 4;
    X:=A[J];
    A[J]:=A[L];
    A[L]:=X;
    J:=J-D;
    if J>0 then goto 3;
4:end;
  goto 2;
5:end;{Shell1}
{*******************************************************************}
procedure FormTreeGr(var Nv,Nr : integer;
                            var Masy : TMasy;
                            var Mass : TMass;
                            var MasRd : TMasy);
{    Процедура выделения дерева графа                          }
{                                                        }
{    Nv - количество вершин в графе;                        }
{    Nr - количество ветвей дерева в графе;                    }
{    Masy - массив указателей для матрицы смежностей;              }
{    Mass - массив элементов матрицы смежностей;                }
{    MasRd - массив ребер дерева.                          }
{                                                        }
           var I,P,A,B,K,L,I1 : integer;
           label 1,2,3,4;
           begin
             I:=0;
             P:=1;
             i1:=0;
```

```
      3: A:=Masy[P];
         B:=Masy[P+1]-1;
         for K:=A to B do
         begin
           if I=0 then goto 2;
           for L:=1 to I do
           if Mass[K]=MasRd[L] then goto 1;
           goto 2;
       1: end;
         I1:=I1-2;
         P:=MasRd[I1];
         goto 3;
       2: I:=I+1;
         MasRd[I]:=P;
         I:=I+1;
         MasRd[I]:=Mass[K];
         if I=2*Nv-2 then goto 4;
         I1:=I+1;
         P:=Mass[K];
         goto 3;
       4: Nr:=I div 2;
       end;{FormTreeGr}
{************************************************************}
 procedure FormTreeRGrV1(var Nv,Nr : integer;
                         var Masy : TMasy;
                         var Mass : TMass;
                         var Masi : TMass;
                         var MasRd : TMasy;
                         var MasRd1 : TMasy);
{   Процедура формирования ветвей дерева графа                  }
{   в виде массива состоящего из ребер графа                }
{                                                    }
{   Nv - количество вершин в графе;                          }
{   Nr - количество ветвей дерева в графе;                   }
{   Masy - массив указателей для матрицы смежностей;            }
{   Mass - массив элементов матрицы смежностей;                }
{   Masi - массив элементов матрицы инциденций;                }
{   MasRd - массив ветвей дерева в виде вершин;              }
{   MasRd1 - массив ветвей дерева в виде ребер;              }
{                                                    }
       var I,Ip,J : integer;
       label 1,2,3;
       begin
{          формирование ветвей дерева                   }
         for I:=1 to Nr do
         begin
           Ip:=MasRd[2*I-1];
           for J:=Masy[Ip] to Masy[Ip+1]-1 do
           begin
             if Mass[J]<>MasRd[2*i] then goto 2;
             MasRd1[I]:=Masi[J];
             goto 1;
       2:   end;
       1: end;
 end;{FormTreeRGrV1}
{********************************************************************}
procedure ForMasReber(var Nv,Nr,M : integer;
                      var Masy : TMasy;
                      var Mass : TMass;
                      var Masi : Tmass;
                      var MasRd : TMasy;
                      var MasRd1 : TMasy);
{Процедура формирования массива для определения ветвей дерева и хорд}
```

```
{в массиве MasRd, 2 - ветвь дерева и 1 - хорда                    }
var i,j: integer;
begin
   FormTreeGr(Nv,Nr,Masy,Mass,MasRD);
   FormTreeRGrV1(Nv,Nr,Masy,Mass,Masi,MasRD,MasRd1);
   for i:= 1 to M do MasRd[i]:=1;
   for j:= 1 to Nr do MasRd[MasRd1[j]]:=2;
end;{ForMasReber}
{*************************************************************************}
procedure Zamena(var Kzikl,Numer,L : integer;
                 var My : TMasy;
                 var Ms : TMass;
                 var Ms1 : TMass);
{процедура замены цикла с номером Numer в массиве циклов Ms на новый цикл}
{хранящейся в массиве Ms1 с последующей перенумеровкой массива указателей}
{ Kzikl - количество циклов;                              }
{ Numer - номер цикла для замены;                          }
{ L - длина замещающего нового цикла;                        }
{ My - массив указателей для матрицы циклов                    }
{ Ms - массив элементов матрицы циклов;                       }
{ Ms1 -массив ребер для замещающего цикла;                     }
label 1;
var i,jj,jjj,L1,L2,Star,Stor,PP,PP1,PP2 : integer;
begin
  L1:= My[Numer+1]-My[Numer];
  if L = 0 then
  begin
    PP1:= My[Numer+1];
    PP2:= My[Kzikl+1]-1;
    PP:= PP2-PP1;
    for i:= 1 to PP+1 do
    Ms[My[Numer]+i-1]:= Ms[My[Numer+1]+i-1];
    for jj:= Numer to Kzikl+1 do My[jj+1]:= My[jj+1]-L1;
    goto 1;
  end;
  L2:= L-L1;
  if L2 = 0 then for i:=1 to L do Ms[My[Numer]+ i -1]:= Ms1[i];
  if L2 < 0 then
  begin
    for i:= My[Numer+1] to My[Kzikl +1]-1 do Ms[i+L2]:= Ms[i];
    for jj:= Numer to Kzikl+1 do My[jj+1]:= My[jj+1]+L2;
    for jjj:= 1 to L do Ms[My[Numer]+jjj-1]:= Ms1[jjj];
  end;
  if L2 > 0 then
  begin
   Star:= My[Numer+1];
   Stor:= My[Kzikl +1]-1;
   PP:=Stor-Star+1;
   for i:= 1 to PP do Ms[Stor-i+L2+1]:=Ms[Stor-i+1];
   for jj:= Numer to Kzikl+1 do My[jj+1]:= My[jj+1]+L2;
   for jjj:= 1 to L do Ms[Star-L1+jjj-1]:= Ms1[jjj];
  end;
1:  {writeln(F2,'   Замена отработала');} ;
end; {Zamena}
{***************************************************************}
 procedure MetodGauss(var Ziklo,M,KKK : integer;
                 var MasYkZiD1 : TMasy;
                 var MasZiD1 : TMass;
                 var Mdop : TMasy;
                 var Mdop1 : TMasy;
                 var Mdop2 : TMass;
                 var MasReber : TMasy);
{ Метод Гаусса для определения независимости выбранной системы    }
```

```
{ базисных циклов графа                                   }
{                                                    }
{ Ziklo - количество базисных циклов графа;                      }
{ M - количество ребер графа;                               }
{ MasYkZiD1 - массив указателей для матрицы базисных циклов       }
{ MasZiD1 : массив элементов матрицы базисных циклов;             }
{ Mdop -массив хорд вошедших в главные элементы;                  }
{ Mdop1 - массив для нового замещающего цикла;                    }
{ MasReber : массив признаков для ребер, ветвь дерева или хорда.  }
label 1,2,3,4,5,6;
var KK,KK1,KK2,L,LL3,L11,L22,Pr,i,ii,iii,j,jj,jjj,iiii : integer;
begin
KKK:=0;
KK:=0;
Pr:=0;
for i:=1 to Ziklo-1 do
  begin {1}
     {writeln(F3,'  ');
     writeln(F3,'MetodGauss: Текущий цикл = ',i); }
     L11:= MasYkZiD1[i+1] - MasYkZiD1[i]; {определение длины цикла}
     if L11 = 0 then goto 6 else
    begin {2}
 {нахождение в цикле с номером i ребра принадлежащего хорде}
     for ii:=MasYkZiD1[i] to MasYkZiD1[i+1]-1 do
     begin {3}
        KK1:= MasZiD1[ii];
        if MasReber[KK1] = 1 then
        begin {4}
          KK:= MasZiD1[ii];
          {найден номер ребра для хорды в цикле с номером i}
          KKK:= KKK+1;
          Mdop[KKK]:= KK; {запись очередной хорды}
          {writeln(F3,'MetodGauss: Запись очередной хорды = ',KK);}
          for iii:= 1 to L11 do Mdop1[iii]:= MasZiD1[MasYkZiD1[i]+iii-1];
          {for iii:=1 to L11 do
            begin
             if iii<>L11 then write(F3,Mdop1[iii],' ');
             if iii=L11 then writeln(F3,Mdop1[iii]);
             end; }
          goto 5;
         end;  {4}
     end;   {3}
{формирование массива Mdop1 для сложения двух пересекающихся циклов}
   end; {2};
5:  for j:=i+1 to Ziklo do
   begin {5}
     {writeln(F3,'MetodGauss: Пересекающийся номер цикла = ',j); }
     L22:= MasYkZiD1[j+1]- MasYkZiD1[j];
     if L22 = 0 then goto 4  else
     begin
      Pr:=0;
      for jj:= MasYkZiD1[j] to MasYkZiD1[j+1]-1 do if MasZiD1[jj]=KK
       then  Pr:=1;
     end;
     if Pr<>1 then goto 4 else
       begin
       KK2:=j;
       for iii:= 1 to L22 do Mdop1[L11+iii]:= MasZiD1[MasYkZiD1[j]+iii-1];
       LL3:=L11+L22;
       for iii:= 1 to LL3 do Mdop2[iii]:= Mdop1[iii];
       {writeln(F3,'MetodGauss:  Общая запись двух объединяющих циклов:');
       for iii:=1 to LL3 do
       begin
```

```
            if iii<>LL3 then write(F3,Mdop2[iii],' ');
            if iii=LL3 then writeln(F3,Mdop2[iii]);
        end; }
        for ii:=1 to LL3-1 do
        begin {7}
          for jjj:= ii+1 to LL3 do
          begin {9}
            if Mdop2[jjj] = 0  then  goto 3;
            if Mdop2[jjj] <> Mdop2[ii] then goto 3;
            begin {8}
                Mdop2[jjj]:=0;
                Mdop2[ii]:=0;
            end; {8}
 3:         end; {9}
 1:      end;   {7}
        L:=0;
        for iii:= 1 to LL3 do
        begin {8}
             if Mdop2[iii] = 0 then goto 2;
             L:=L+1;
             Mdop2[L]:=Mdop2[iii];
 2:      end;  {8}
        {writeln(F3,'MetodGauss:  Результат объединяющих циклов:');
        for iii:=1 to L do
        begin
            if iii<>L then write(F3,Mdop2[iii],' ');
            if iii=L then writeln(F3,Mdop2[iii]);
        end; }
        Zamena(Ziklo,KK2,L,MasYkZiD1,MasZiD1,Mdop2);
       end;
 4:   end; {5}
 6:  end; {1}
       Pr:= MasYkZiD1[Ziklo+1]- MasYkZiD1[Ziklo];
        if Pr > 0 then
          begin {1}
          for i:=MasYkZiD1[Ziklo] to MasYkZiD1[Ziklo+1]-1 do
            begin {2}
              KK1:= MasZiD1[i];
              if MasReber[KK1] = 1 then
              begin {3}
                KK:= MasZiD1[i];
                {найден номер ребра для хорды в последнем цикле}
                KKK:= KKK+1;
                Mdop[KKK]:= KK;
              end; {3}
            end; {2}
          end; {1}
end; {MetodGauss}
{*****************************************************************}
 procedure ZamZikl01(var KKK,Kzikl : integer;
                     var MasYkZiD1 : TMasy;
                     var MasZiD1 : TMass;
                     var Mdop : TMasy);
{ Изменение признаков удаления зависимых циклов после проверки
  методом Гаусса                                   }
{                                               }
{ Ziklo - количество базисных циклов графа;                    }
{ Kzikl - количество единичных циклов графа;                    }
{ MasYkZiD1 - массив указателей для матрицы базисных циклов       }
{ MasZiD1 : массив элементов матрицы базисных циклов;            }
{ Mdop -массив признаков удаления ребер;                        }
{ K1 -массив признаков удаления ребер;                      }
label 1;
```

```
var i,j,ii,K2,K3 : integer;
begin
K1:=0;
   for j:= 1 to KKK do if MasYkZiD1[j]>MasYkZiD1[j+1]-1 then
        begin {1}
          K1:=1;
          K2:=j;
          K3:=0;
          for ii:=1 to Kzikl do if Mdop[ii]=1 then
          begin {2}
            K3:=K3+1;
            if K3=K2 then Mdop[ii]:=-1;
          end; {2}
        end; {1}
end; {ZamZikl01}
{*****************************************************************}
 procedure Sbor11b(var Nv,M,Kzikl,Ziklo : integer;
                  var Masy : TMasy;
                  var Mass : TMass;
                  var MasyT : TMasy;
                  var MassT : TMass;
                  var MasMy1 : TMasy;
                  var MasMs1 : TMass;
                  var Masi : TMass;
                  var MasMdop : TMasy;
                  var MasMcg : TMasy;
                  var Mass1 : TMass;
                  var MasMcg1 : TMasy;
                  var MasMcg2 : TMasy;
                  var MasKol : TMasy);
{                                               }
var i,j,K7 : integer;
begin
    ForMasReber(Nv,Nr,M,Masy,Mass,Masi,MasKol,MasMcg);
    writeln(F4, Kzikl);
    for I:=1 to M do MasMcg2[i]:=MasKol[i];
    MetodGauss(Kzikl,M,K7,MasyT,MassT,MasMcg,MasMcg1,Mass1,MasKol);
    writeln(F4, Kzikl);
      for i:=1 to Kzikl+1 do
        begin {07}
          if i<> Kzikl+1 then write(F4, MasyT[i],' ');
          if i= Kzikl+1 then writeln(F4, MasyT[i]);
        end; {07}
    ZamZikl01(K7,Kzikl,MasyT,MassT,MasMdop);
end;{Sbor11}
{*****************************************************************}
procedure Preor111 (var Ziklo,M,L : integer;
                  var MasMy1 : TMasy;
                  var MasReber : TMasy;
                  var MasHord1 : TMasy);
{ Процедура определения хорд для зависимой системы циклов          }
{  результат: L количество зависимых циклов,                   }
{  MasHord1 - множество зависимых хорд.                        }
{                                                   }
{ Ziklo - количество базисных циклов в графе;                  }
{ M - количество ребер в графе;                         }
{ MasReber - массив хорд и ветвей;                         }
{ MasHord1 - массив хорд;                             }
{ MasHord2 - массив независимых хорд.                        }
{******************************************************************}
var i,j,ii,jj,iii,Hord :integer;
begin
  Hord:= 0;
```

```
  for i:= 1 to M do
  begin {1}
    if MasReber[i]=1 then
    begin {2}
      Hord:=Hord+1;
      MasHord1[Hord]:=i;
    end; {2}
  end; {1}
  for j:= 1 to Ziklo do if MasMy1[j]=MasMy1[j+1] then
      MasHord1[j]:=-MasHord1[j];
  L:=0;
  for i:= 1 to Hord do
    if MasHord1[i]<0 then
    begin {3}
      L:=L+1;
      MasHord1[L]:=-MasHord1[i];
    end; {3}
end;{Preor111}
{******************************************************************}
procedure  Shella(var N : integer;
                  var A : XMass;
                  var B : TMasy);
{       процедура Шелла для упорядочивания элементов в
        вещественном массиве А                          }
{                                                   }
{     N - количество элементов в массиве;               }
{     A - сортируемый массив;                           }
{     B - связанный массив;                         }
var D,Nd,I,J,L,Y : integer;
var X : real;
label 1,2,3,4,5;
begin
  D:=1;
1:D:=2*D;
  if D<=N then goto 1;
2:D:=D-1;
  D:=D div 2;
  if D=0 then goto 5;
  Nd:=N-D;
  for I:=1 to Nd do
  begin
    J:=I;
3:  L:=J+D;
    if A[L]>=A[J] then goto 4;
    X:=A[J];
    Y:= B[J];
    A[J]:=A[L];
    B[J]:=B[L];
    A[L]:=X;
    B[l]:=Y;
    J:=J-D;
    if J>0 then goto 3;
4:end;
  goto 2;
5:end;{Shella}
{******************************************************************}
 procedure WedenieN(var Nv,M,Ziklo,Kzikl,Priz,Makley : integer;
                         var MasKol : TMasy);

{    Процедура подготовки информации для базиса циклов          }
{                                                   }
{    Nv - количество вершин в графа;                   }
{    M - количество ребер в графа;                   }
```

```
{    Ziklo - цикломатическое число графа;                                }
{    Kzikl - количество изометрических циклов в графе;                 }
{    Priz - признак существования нулей в векторе циклов по ребрам;   }
{    Makley - значение функционала Маклейна;                           }
{                                                          }
        var i,j,k,kp,i1,i2 : integer;
        begin
            Priz:=1;
            Makley:=0;
            for j:=1 to M do  if MasKol[j]>0 then
                  Makley:=Makley + (MasKol[j]-1)*(MasKol[j]-2);
            {writeln(F3,'Makley = ',Makley); }
            i1:=0;
            for i:=1 to M do if MasKol[i]= 1 then  i1:=i1+1;
            if i1=0 then Priz:=0;
         end;{WedenieN}
{*****************************************************************}
 procedure BasisN(var Nv,M,Kzikl : integer;
                        var MasyT : TMasy;
                        var MassT : TMass;
                        var MasMdop : TMasy;
                        var MasKol : TMasy;
                        var MasMcg : TMasy);
 {   Формирование базиса единичных циклов с min функционала Мак-Лейна   }
 {                                                          }
 { Nv - количество вершин в графе;                                }
 { M - количество ребер в графе;                                 }
 { Kzikl  - количество t-циклов в графе;                           }
 { MasyT - массив указателей для матрицы циклов;                     }
 { MassT - массив элементов матрицы t-циклов;                        }
 { MasKol - массив количества циклов проходящих по ребру              }
 { MasMdop -  массив признаков удаления цикла, 1- неудален, -1 - удален }
 { MasMcg -  массив значений функционалов Мак-Лейна при удалении цикла  }
 label 1;
 var I,J,II,JJ,III,JJJ,IIII,KL,Ziklo1,KML,NUM,Dop,DUB : integer;
begin
{инициализация}
NUM:=0;
for  i:=1 to Kzikl do MasMdop[i]:=1;
{количество удаленных циклов равно нулю}
KL:=0;
{определяем значение цикломатического числа}
Ziklo:= M-Nv+1;
{если количество циклов равна цикломатическому числу, то идти на конец выбора базиса циклов}
while Kzikl-KL > Ziklo do
{определим количество циклов проходящих по каждому ребру для определения функционала Мак-Лейна с учетом удаленных циклов}
  begin {1}
1:  Dop:=1;
    for  j:= 1 to M do MasKol[j]:= 0;
    for  i:=1 to Kzikl do
      begin {2}
        if MasMdop[i] =1 then
        for  j:= MasyT[i] to MasyT[i+1]-1 do
            MasKol[MassT[j]]:= MasKol[MassT[j]]+1;
      end; {2}
    {формируем массив MasMcg для расчета значения функционала Мак-Лейна после удаления I-го цикла}
    for  jjj:=1 to Kzikl do MasMcg[jjj]:=0;
    for  ii:=1 to Kzikl do {удаляем II-ый цикл}
      begin {2}
        if MasMdop[ii] >0 then
        begin {3}
         for  jj:= MasyT[ii] to MasyT[ii+1]-1 do
```

```
          MasKol[MassT[jj]]:= MasKol[MassT[jj]]-1;
    {функционал Мак-Лейна для каждого удаленного цикла}
        for  iii:=1 to M do
        begin
        if MasKol[iii] > 0 then
          MasMcg[ii]:= MasMcg[ii]+((MasKol[iii]-2)*(MasKol[iii]-1));
        end;
        for  jj:= MasyT[ii] to MasyT[ii+1]-1 do
          MasKol[MassT[jj]]:= MasKol[MassT[jj]]+1;
       end; {3}
     end; {2}
    {для всех циклов сформировали массив MasMcg}
    KML:=0;
    {выбираем KML максимальное по значению}
    for  jjj:=1 to Kzikl do
    if MasMdop[jjj] = 1  then
    if MasMcg[jjj] > KML then KML:=MasMcg[jjj];
    {Выбираем минимальный элемент в массиве MasMcg}
    for  jjj:=1 to Kzikl do
    if MasMdop[jjj] = 1 then
    if MasMcg[jjj] <= KML then KML:=MasMcg[jjj];
    {Выбираем среди кандидатов на удаление цикл максимальный по длине}
    IIII:=0;
    for  jjj:=1 to Kzikl do
      if MasMdop[jjj] = 1 then
      if MasMcg[jjj] = KML then
      if (MasyT[jjj+1] - MasyT[jjj]) > IIII then  IIII:= MasyT[jjj+1] - MasyT[jjj];
    for  jjj:=1 to Kzikl do
      if MasMdop[jjj] = 1 then
      if MasMcg[jjj] = KML then
      if (MasyT[jjj+1] - MasyT[jjj]) = IIII then  NUM:=jjj;
  {нужно проверить не принесёт ли это удаление к появлению 0 в массиве MasKol}
    MasMdop[NUM]:=-1;
    for j:=1 to M do MasKol[j]:=0;
    for  i:=1 to Kzikl do
      begin {2}
       if MasMdop[i] >0 then
       for  j:= MasyT[i] to MasyT[i+1]-1 do
           MasKol[MassT[j]]:= MasKol[MassT[j]]+1;
      end; {2}
     for j:=1 to M do if MasKol[j]=0 then
     begin
      MasMdop[NUM]:=2;
      goto 1;
     end;
    {writeln(F2,'NUM = ',NUM); }
    MasMdop[NUM]:= -1;
    KL:=KL+1;
    {writeln(F2,' KL = ',KL); }
  end; {1}
  for i:=1 to M do MasKol[i]:=0;
  for  i:=1 to Kzikl do
      begin
       if MasMdop[i] > 0 then
       for  j:= MasyT[i] to MasyT[i+1]-1 do
           MasKol[MassT[j]]:= MasKol[MassT[j]]+1;
      end;
end;{BasisN}
{*****************************************************************}
procedure FormBasis(var Nv,Kzikl,Ziklo : integer;
                          var MyZ : TMasy;
                          var MsZ : TMass;
                          var MyZ1 : TMasy;
```

```
                    var MsZ1 : TMass;
                    var MasPriz : TMasy);
{Процедура формирования базиса циклов в зависимости от MasPriz;   }
{                                                   }
{ Nv - количество вершин в графе;                          }
{ Kzikl - количество единичных циклов графа;                     }
{ Ziklo - цикломатическое число графа;                          }
{ MyZ - массив указателей для матрицы единичных циклов;        }
{ MsZ - массив единичных циклов графа;                          }
{ MyZ1 - массив указателей для матрицы базисных циклов        }
{ MsZ1 : массив элементов матрицы базисных циклов;             }
{ MasPriz - массив списка удаляемых циклов для получения базиса;  }
var i,FF,SF,j : integer;
begin
  {for i:= 1 to 4000 do MasZiD1[i]:=0; }
  FF:=0;
  SF:= 0;
  MyZ1[1]:= 1;
  for i:=1 to Kzikl do if MasPriz[i] > 0 then
  begin
     SF:= SF+1; {номер базисного цикла}
     FF:= MyZ[i+1] - MyZ[i]; {мощность цикла}
     MyZ1[SF+1]:= MyZ1[SF]+ FF;
     for j:= 1 to FF do MsZ1[MyZ1[SF]+j-1]:= MsZ[MyZ[i]+j-1];
  end;
  Ziklo:=SF;
end;{FormBasis}
{*****************************************************************}
procedure FormVerBasis(var Nv,Ziklo : integer;
                    var My : TMasy;
                    var Ms : TMass;
                    var My1 : TMasy;
                    var Ms1 : TMass;
                    var Mass1 : TMass;
                    var Masi : TMass);
{Процедура записи базиса циклов через вершины                    }
{                                                   }
{ Nv - количество вершин в графе;                          }
{ Ziklo - цикломатическое число графа;                          }
{ My - массив указателей для матрицы смежностей;                 }
{ Ms - массив элементов матрицы смежностей;                     }
{ My1 - массив указателей для матрицы базисных циклов             }
{ Ms1 : массив элементов матрицы базисных циклов;                }
{ Masi - массив элементов матрицы инциденций;                   }
{ Mass1 - элементов матрицы базисных циклов через вершины;        }
{****************************************************************}
label 1,2;
var i,j,ii,jj,iii,Nr,GG,NvS,NvK :integer;
begin
  GG:= 0;
  for i:= 1 to Ziklo do
     begin{1}
      {writeln(F2,'FormVerBasis: i = ',i);}
      Nr:= Ms1[My1[i]];
      {writeln(F4,'FormVerBasis: Nr = ',Nr);}
      for j:=1 to Nv do
      begin {2}
        for jj:= My[j] to My[j+1]-1 do
        begin {3}
         if Masi[jj]= Nr  then
         begin {4}
           NvS:= Ms[jj];
           NvK:= Ms[jj];
```

```
              {writeln(F3,'FormVerBasis: NvS = ',NvS,' NvK = ',NvK);}
              goto 1;
            end; {4}
          end;{3}
        end;{2}
1:      for ii:=My[NvK] to My[NvK+1]-1 do
         begin{3}
         {writeln(F4,'FormVerBasis: ii = ',ii);}
         for iii:=My1[i] to My1[i+1]-1 do
           begin {4}
            {writeln(F4,'FormVerBasis: iii = ',iii);}
            if Ms1[iii]=Masi[ii] then if Ms1[iii]<>Nr then
            begin {2}
              GG:=GG+1;
              Mass1[GG]:=NvK;
              Nr:=Ms1[iii];
              NvK:=Ms[ii];
              {writeln(F4,'FormVerBasis: NvK = ',NvK);}
              if NvK=NvS then goto 2 else goto 1;
            end; {2}
           end;{4}
         end;{3}
2:      end; {1}
end;{FormVerBasis}
{*****************************************************************}
 procedure PodMac (var M,Makl : integer;
                   var MasKol : TMasy);
 { Процедура вычисления функционала Маклейна для заданной системы  }
 { циклов                                                }
 {                                                      }
 { M - количество ребер в графе;                                    }
 { MakLin  - значение функционала Маклейна;                         }
 { MasKol - массив количества циклов проходящих по ребру            }
  var I,J,ii,jj,III,jjj : integer;
begin
   Makl:=0;
   for i:=1 to M do
     if MasKol[i]>0 then
     Makl:=Makl + (MasKol[i]-2)*(MasKol[i]-1);
end;{PodMac}
{*****************************************************************}
 procedure Ydalow (var M,KM,Kzikl : integer;
                   var MasyT : TMasy;
                   var MassT : TMass;
                   var MasMdop : TMasy;
                   var MasMcg : TMasy;
                   var MasMcg1 : TMasy;
                   var MasMcg2 : TMasy;
                   var MasKol : TMasy);
 { Процедура определения количества нулей в векторе циклов по      }
 { ребру после удаления цикла;                                     }
 { циклов                                                }
 {                                                      }
 { M - количество ребер в графе;                                    }
 { MakLin  - значение функционала Маклейна;                         }
 { MasyT - массив указателей для изометрических циклов;             }
 { MassT - массив элементов изометрических циклов;                  }
 { MasKol - массив количества циклов проходящих по ребру;           }
 { MasMcg - массив количества нулей после удаления цикла;           }
 { MasMcg1 - массив количества циклов проходящих по ребру           }
 { MasMcg2 - массив значений функционала Маклейна после  удаления  }
 { цикла;                                                }
 { MasMdop - массив признаков удаления циклов;                      }
```

```
 var I,J,ii,jj,III,jjj,Makl,KKK : integer;
begin
   KM:=0;
   for i:=1 to Kzikl do
   begin  {1}
     if MasMdop[i]>0 then
     begin {2}
       KM:=KM+1;
       for j:=MasyT[i] to MasyT[i+1]-1 do
        MasKol[MassT[j]]:=MasKol[MassT[j]]-1;
        KKK:=0;
        for ii:=1 to M do
         if MasKol[ii]=0 then KKK:=KKK+1;
         begin {3}
          MasMcg[KM]:=i;
          MasMcg1[KM]:=KKK;
          {writeln(F4,' Печать перед PodMac ');
          for iii:=1 to M do
          begin
              if iii<>M then write(F4,MasKol[iii],' ');
              if iii=M then writeln(F4,MasKol[iii]);
          end;}
          PodMac(M,Makl,MasKol);
          {writeln(F4,' i = ',i,' Makl = ',Makl); }
          MasMcg2[KM]:=Makl;
          { Востановление вектора циклов по ребрам }
          for j:=MasyT[i] to MasyT[i+1]-1 do
           MasKol[MassT[j]]:=MasKol[MassT[j]]+1;
         end; {3}
     end; {2}
   end; {1}
         {for I:=1 to KM do
         begin
             if i<>KM then write(F4,MasMcg[i],' ');
             if i=KM then writeln(F4,MasMcg[i]);
         end;
         for I:=1 to KM do
         begin
             if i<>KM then write(F4,MasMcg1[i],' ');
             if i=KM then writeln(F4,MasMcg1[i]);
         end;
         for I:=1 to KM do
         begin
             if i<>KM then write(F4,MasMcg2[i],' ');
             if i=KM then writeln(F4,MasMcg2[i]);
         end;}
end;{Ydalow}
{*****************************************************************}
 procedure Klips (var M,KM,Priz : integer;
                  var MasyT : TMasy;
                  var MassT : TMass;
                  var MasMdop : TMasy;
                  var MasMcg : TMasy;
                  var MasMcg1 : TMasy;
                  var MasMcg2 : TMasy;
                  var MasKol : TMasy);
 { Процедура определения первого цикла на удаление                }
 {                                                      }
 { M - количество ребер в графе;                                  }
 { MakLin  - значение функционала Маклейна;                       }
 { MasKol - массив количества циклов проходящих по ребру          }
 { MasMcg - массив количества нулей после удаления цикла;         }
 { MasMcg1 - массив количества циклов проходящих по ребру         }
```

```
 { MasMcg2 - массив значений функционала Маклейна после  удаления  }
 { цикла;                                              }
 { MasMdop - массив признаков удаления циклов;                    }
  var I,J,ii,jj,III,jjj,Makl,KKK,KK20 : integer;
  var KK30,Min,KT : integer;
  label 1;
begin
   Priz:=0;
   KK20:=0;
   {writeln(F2,' KM = ',KM); }
   for i:=1 to KM do
   begin  {2}
     if MasMcg1[i]=1 then
     begin {1}
       KK20:=KK20+1;
       {writeln(F4,'KK20 = ',KK20);}
       MasMcg[KK20]:=MasMcg[i];
       {writeln(F4,' MasMcg[kk20] = ',MasMcg[KK20],' i = ',i,' MasMcg[i] = ',MasMcg[i]);}
       MasMcg2[KK20]:=MasMcg2[i];
       {writeln(F3,' i = ',i,' MasMcg2[i] = ',MasMcg2[i]);}
     end; {1}
   end; {2}
   if KK20=0 then
   begin {5}
     Priz:=7;
     goto 1;
   end; {5}
   KK30:=0;
   Min:=100000;
   for i:=1 to KK20 do if MasMcg2[i]<Min then Min:=MasMcg2[i];
   for ii:=1 to KK20 do
    if MasMcg2[ii]=Min then
    begin {3}
      KK30:=KK30+1;
      MasMcg[KK30]:=MasMcg[ii];
    end; {3}
   { for iii:=1 to KK30 do
        begin
            if iii<>KK30 then write(F4,MasMcg[iii],' ');
            if iii=KK30 then writeln(F4,MasMcg[iii]);
        end;
    for iii:=1 to KK30 do
        begin
            if iii<>KK30 then write(F4,MasMcg2[iii],' ');
            if iii=KK30 then writeln(F4,MasMcg2[iii]);
        end; }
    KT:= MasMcg[1];
    for j:=MasyT[KT] to MasyT[KT+1]-1 do
    MasKol[MassT[j]]:=MasKol[MassT[j]]-1;
    MasMdop[KT]:=-1;
    {writeln(F4,'Номер удаляемого цикла = ',KT); }
    for j:=1 to M do
    if MasKol[j]=0 then MasKol[j]:=-4;
1:
end;{Klips}
{*****************************************************************}
procedure Sbor16(var M,KM,Priz,Ziklo,Kzikl : integer;
                 var MasyT : TMasy;
                 var MassT : TMass;
                 var MasMdop : TMasy;
                 var MasMcg : TMasy;
                 var MasMcg1 : TMasy;
                 var MasMcg2 : TMasy;
```

```
                 var MasKol : TMasy);
 { Процедура определения первого цикла на удаление                 }
 {                                                   }
 { M - количество ребер в графе;                              }
 { MakLin  - значение функционала Маклейна;                       }
 { MasKol - массив количества циклов проходящих по ребру          }
 { MasMcg - массив количества нулей после удаления цикла;          }
 { MasMcg1 - массив количества циклов проходящих по ребру         }
 { MasMcg2 - массив значений функционала Маклейна после  удаления  }
 { цикла;                                                }
 { MasMdop - массив признаков удаления циклов;                   }
  var I,J,ii,jj,III,jjj,Makl,KKK,KK20 : integer;
  label 1;
begin
   {writeln(F4,'  Входим в процедуру Sbor16 ');
   for I:=1 to M do
        begin
            if i<>M then write(F4,MasKol[i],' ');
            if i=M then writeln(F4,MasKol[i]);
        end; }
   PodMac (M,Makl,MasKol);
  { writeln(F4,'  Процедура Sbor16 ',Makl);   }
   while Makl>0 do
   begin {2}
     {writeln(F4,'  Процедура Sbor16 ',Makl); }
     Ydalow (M,KM,Kzikl,MasyT,MassT,MasMdop,MasMcg,MasMcg1,MasMcg2,MasKol);
     Klips (M,KM,Priz,MasyT,MassT,MasMdop,MasMcg,MasMcg1,MasMcg2,MasKol);
     PodMac (M,Makl,MasKol);
     {writeln(F4,'Sbor16  KM = ',KM);
     for I:=1 to M do
        begin
            if i<>M then write(F4,MasKol[i],' ');
            if i=M then writeln(F4,MasKol[i]);
        end; }
     if Priz>0 then goto 1;
   end; {2}
1:
end;{Sbor16}
{*****************************************************************}
 procedure Pravilo(var Nv,M,Priz,Priz1,Ziklo : integer;
                 var MasKol : TMasy);
 { Процедура вычисления постоянной Эйлера.                       }
 {                                                   }
 { Nv - количество вершин в графе;                             }
 { M - количество ребер в графе;                              }
 { Ziklo  - количество циклов подмножества;                     }
 { MasKol - массив количества циклов проходящих по ребру          }
  var I,J,ii,jj,III,jjj,Kol,Kol1,Kol2,Eir,Eir1 : integer;
  label 1;
begin
   Priz:=0;
   Priz1:=0;
   Kol:=0;
   Kol1:=0;
   Kol2:=0;
   for I:=1 to M do if MasKol[i]=-4 then MasKol[i]:=0;
   for I:=1 to M do
        begin
            if MasKol[i]=0 then Kol:=Kol+1;
            if MasKol[i]=1 then Kol1:=Kol1+1;
            if MasKol[i]=2 then Kol2:=Kol2+1;
        end;
        {writeln(F4, ' Kol = ',kol);
```

```
          writeln(F4, ' Kol1 = ',kol1);
          writeln(F4, ' Kol2 = ',kol2);
          writeln(F4, ' Ziklo = ',Ziklo); }
    Eir:=Ziklo-(M-Kol)+Nv-1;
    Eir1:=4*Kol2+Kol1-3*(2*Kol2+Kol1)+2*(M-Kol);
    if Kol1=0 then  Eir:=(Ziklo-1)-Kol2+Nv-1;
    { writeln(F4, ' Eir = ',Eir);
     writeln(F4, ' Eir1 = ',Eir1);   }
    if Eir<>0 then Priz:=1;
    if Eir1<>0 then Priz1:=1;
end;{Pravilo}
{*********************************************************************}
var kar,kara : integer;
label 1,2,3,4,5,6,7,8;
begin
         assign(F1,'D:\Planar\GR1\14v23r.gr1');
          reset(F1);
          assign(F4,'D:\Planar\BOS\14v23ra.bos');
          Rewrite(F4);
{Вводим количество вершин}
          readln(F1,Nv);
{Вводим массив указателей}
          for I:=1 to Nv+1 do
          begin
            if i<>Nv+1 then read(F1,Masy[i]);
            if i=Nv+1 then readln(F1,Masy[i]);
          end;
{Вводим массив элементов матрицы смежностей}
          for I:=1 to Nv do
          for j:=Masy[i] to Masy[i+1]-1 do
          begin
            if j<>Masy[i+1]-1 then read(F1,Mass[j]);
            if j=Masy[i+1]-1 then readln(F1,Mass[j]);
          end;
{Вводим количество ребер}
          readln(F1,M);
{Вводим элементы матрицы инциденций}
          for I:=1 to Nv do
          for j:=Masy[i] to Masy[i+1]-1 do
          begin
            if j<>Masy[i+1]-1 then read(F1,Masi[j]);
            if j=Masy[i+1]-1 then readln(F1,Masi[j]);
          end;
          close (F1);
          assign(F2,'D:\Planar\EZI\14v23r.ezi');
          reset(F2);
{Вводим количество изометрических циклов}
          readln(F2,Kzikl);
{Вводим массив указателей для изометрических циклов}
          for i:=1 to Kzikl+1 do
          begin
             if i<> Kzikl+1 then read(F2, MasyT [i]);
             if i= Kzikl+1 then readln(F2, MasyT [i]);
          end;
{Вводим множество изометрических циклов в виде ребер}
          for I:=1 to Kzikl do
          begin
            for j:=MasyT[i] to MasyT[i+1]-1 do
            begin
              if j<>MasyT[i+1]-1 then read(F2,MassT[j]);
              if j=MasyT[i+1]-1 then readln(F2,MassT[j]);
            end;
          end;
```

```
{Вводим множество изометрических циклов в виде вершин}
        for I:=1 to Kzikl do
        begin
          for j:=MasyT[i] to MasyT[i+1]-1 do
          begin
             if j<>MasyT[i+1]-1 then read(F2,Mass1[j]);
             if j=MasyT[i+1]-1 then readln(F2,Mass1[j]);
          end;
        end;
        close (F2);
        {Создаём новый файл и открываем его в режиме "для чтения и записи"}
{Создаем случайным образом новую систему изометричеких циклов}
        for iii:=1 to 100 do
        begin {00}
7:       kara:= 44;
        Randomize;
        for i:= 1 to Kzikl do
        begin {01}
            MasMcg1[i]:=i;
            MasMcg4[i]:=Random;
        end; {01}
        Shella(Kzikl,MasMcg4,MasMcg1);
        MasMy2[1]:=1;
        k:=0;
        for i:= 1 to Kzikl do
        begin {02}
          TZ:=MasMcg1[i];
          Dl:= MasyT[Tz+1] - MasyT[Tz];
          for j:= MasyT[Tz] to MasyT[Tz+1]-1 do
          begin {03}
            k:=k+1;
            MasMs2[K]:=MassT[j];
          end; {03}
        MasMy2[i+1]:=MasMy2[i]+Dl;
        end; {02}
        writeln(F4,'       Последовательность циклов');
        for I:=1 to Kzikl do
        begin {05}
          for j:=MasMy2[i] to MasMy2[i+1]-1 do
          begin {06}
             if j<>MasMy2[i+1]-1 then write(F4,MasMs2[j],' ');
             if j=MasMy2[i+1]-1 then writeln(F4,MasMs2[j]);
             end; {06}
        end; {05}
        for i:=1 to Kzikl+1 do MasMy5[i]:=MasMy2[i];
        for i:=1 to MasMy5[Kzikl+1]-1 do MasMs5[i]:=MasMs2[i];
{Новое множество изометрических циклов сформировано и находится в массиве MasMs2}
        Sbor11b(Nv,M,Kzikl,Ziklo,Masy,Mass,MasMy2,MasMs2,
            MasMy3,MasMs3,Masi,MasMdop,MasMcg,Mass1,MasMcg1,
            MasMcg2,MasKol);
        for ii:=1 to Kzikl+1 do
        begin {07}
           if ii<> Kzikl+1 then write(F4, MasMy2[ii],' ');
           if ii= Kzikl+1 then writeln(F4, MasMy2[ii]);
        end; {07}
        for i:=1 to Kzikl do MasMdop[i]:=1;
        for ii:=1 to Kzikl do
          if MasMy2[ii]=MasMy2[ii+1] then MasMdop[ii]:=0;
        writeln(F4,' Массив признаков для базиса  ');
        for i:=1 to Kzikl do
        begin {07}
           if i<> Kzikl then write(F4, MasMdop[i],' ');
           if i= Kzikl then writeln(F4, MasMdop[i]);
```

```
        end; {07}
        FormBasis(Nv,Kzikl,Ziklo,MasMy5,MasMs5,MasMy3,MasMs3,MasMdop);
        FormVerBasis(Nv,Ziklo,Masy,Mass,MasMy3,MasMs3,Mass1,Masi);
        writeln(F4,'  Цикы базиса в ребрах');
        for I:=1 to Ziklo do
        begin {09}
          for j:=MasMy3[i] to MasMy3[i+1]-1 do
          begin {10}
            if j<>MasMy3[i+1]-1 then write(F4,MasMs3[j],' ');
            if j=MasMy3[i+1]-1 then writeln(F4,MasMs3[j]);
          end; {10}
        end; {09}
        writeln(F4,'  Циклы базиса в вершинах');
        for I:=1 to Ziklo do
        begin {09}
          for j:=MasMy3[i] to MasMy3[i+1]-1 do
          begin {10}
            if j<>MasMy3[i+1]-1 then write(F4,Mass1[j],' ');
            if j=MasMy3[i+1]-1 then writeln(F4,Mass1[j]);
          end; {10}
        end; {09}
 {Формируем массив признаков удаления циклов MasMdop}
         for i:=1 to Ziklo do MasMdop[i]:=1;
         Unk22(Nv,M,Ziklo,MasMy3,MasMs3,MasMdop,MasKol);
        for I:=1 to M do
        begin
          if i<>M then write(F4,MasKol[i],' ');
          if i=M then writeln(F4,MasKol[i]);
        end;
{Перезапись множества циклов из массива MasMs3 MasMs1}
        for i:=1 to Ziklo+1 do MasMy1[i]:=MasMy3[i];
        for i:=1 to MasMy1[Ziklo+1]-1 do MasMs1[i]:=MasMs3[i];
{Определение значения функционала Маклейна по MasKol}
        WedenieN(Nv,M,Ziklo,Kzikl,Priz1,Makley,MasKol);
        {writeln(F4,'  после процедуры Wedenie для определения 1 ');
         if Priz1=0 then goto 2;
        {writeln(F4,'   Циклы на удаление по Маклейну ');
        for i:=1 to Ziklo+1 do
        begin
           if i<> Ziklo+1 then write(F4, MasMy1 [i],' ');
           if i= Ziklo+1 then writeln(F4, MasMy1 [i]);
        end;
        for I:=1 to Ziklo do
        begin
          for j:=MasMy1[i] to MasMy1[i+1]-1 do
          begin
            if j<>MasMy1[i+1]-1 then write(F4,MasMs1[j],' ');
            if j=MasMy1[i+1]-1 then writeln(F4,MasMs1[j]);
          end;
        end; }
        Sbor16 (M,KM,Priz2,Ziklo,Kzikl,MasMy1,MasMs1,MasMdop,MasMcg1,
          MasMcg2,MasMcg3,MasKol);
         writeln(F4,'       Вектор циклов по ребрам');
        for i:=1 to M do
        begin
           if i<> M then write(F4, MasKol[i],' ');
           if i= M then writeln(F4, MasKol[i]);
        end;
        writeln(F4,'       Удаленные циклы');
        for i:=1 to Ziklo do
        begin
           if i<> Ziklo then write(F4, MasMdop[i],' ');
           if i= Ziklo  then writeln(F4, MasMdop[i]);
```

```
      end;
      kar:=0;
      for I:=1 to M do if MasKol[i]>kar then kar:=MasKol[i];
      if kar>2 then goto 6;
      k:=0;
      writeln(F4,'       Циклы в ребрах');
      for I:=1 to Ziklo do
      begin {015}
        if MasMdop[i]>0 then
        begin {014}
        k:=k+1;
        write(F4,'цикл ',k,': {');
        for j:=MasMy3[i] to MasMy3[i+1]-1 do
        begin {016}
          if j<>MasMy3[i+1]-1 then write(F4,'e',MasMs3[j],',');
          if j=MasMy3[i+1]-1 then writeln(F4,'e',MasMs3[j],'}');
        end; {016}
       end; {014}
      end; {015}
      k:=0;
      writeln(F4,'       Циклы в вершинах');
      for I:=1 to Ziklo do
      begin {017}
        if MasMdop[i]>0 then
        begin {018}
         k:=k+1;
         write(F4,'цикл ',k,': {');
         for j:=MasMy3[i] to MasMy3[i+1]-1 do
         begin {019}
          if j<>MasMy3[i+1]-1 then write(F4,'v',Mass1[j],',');
          if j=MasMy3[i+1]-1 then writeln(F4,'v',Mass1[j],'}');
         end; {019}
        end; {018}
      end; {017}
      goto 1;
2:  writeln(F4,' Не существует единиц в векторе циклов по ребрам! ');
    goto 1;
4:  writeln(F4,' Не соблюдается правило удаления цикла! ');
    goto 1;
5:  writeln(F4,' Не выполняется правило Эйлера! ');
    goto 1;
6:  writeln(F4,' Функционал Маклейна не равен 0! ');
    goto 7;
1:        writeln(F4,'  Расчет ',iii);
    end;
     close (F4);
     writeln('  Расчет окончен ');
end.
```

## 8.5. Выходной файл программы Rebro2

### Выходной файл 7.bos

| | |
|---|---|
| | изометрические циклы в последовательности |
| 10 11 12 15 | изометрический цикл 1 |
| 2 3 10 12 | изометрический цикл 2 |
| 4 5 16 | изометрический цикл 3 |
| 3 4 14 | изометрический цикл 4 |
| 14 15 16 | изометрический цикл 5 |
| 3 5 15 | изометрический цикл 6 |
| 6 8 11 15 | изометрический цикл 7 |
| 12 13 14 | изометрический цикл 8 |
| 1 3 8 | изометрический цикл 9 |

6 7 10 изометрический цикл 10
10 11 13 16 изометрический цикл 11
1 2 6 изометрический цикл 12
6 9 11 16 изометрический цикл 13
8 9 14 изометрический цикл 14
2 4 10 13 изометрический цикл 15
7 8 12 изометрический цикл 16
1 4 9 изометрический цикл 17
2 5 11 изометрический цикл 18
7 9 13 изометрический цикл 19
19 количество изометрических циклов в последовательности
1 5 9 12 16 19 19 23 26 31 34
34 34 39 39 39 39 39 39 39 массив указателей для независимой системы циклов
Массив признаков для базиса
1 1 1 1 1 0 1 1 1 1 0 0 1 0 0 0 0 0 0 признак 0 – цикл удален, 1 – цикл оставлен
Цикломатическое число графа = 10 (количество независимых циклов)
10 11 12 15 изометрический цикл 1
2 3 10 12 изометрический цикл 2
4 5 16 изометрический цикл 3
3 4 14 изометрический цикл 4
14 15 16 изометрический цикл 5
6 8 11 15 изометрический цикл 7
12 13 14 изометрический цикл 8
1 3 8 изометрический цикл 9
6 7 10 изометрический цикл 10
6 9 11 16 изометрический цикл 13
4 5 7 3 изометрический цикл 1 в вершинах
3 4 5 1 изометрический цикл 2 в вершинах
6 7 1 изометрический цикл 3 в вершинах
5 6 1 изометрический цикл 4 в вершинах
6 7 5 изометрический цикл 5 в вершинах
3 7 5 2 изометрический цикл 7 в вершинах
5 6 4 изометрический цикл 8 в вершинах
2 5 1 изометрический цикл 9 в вершинах
3 4 2 изометрический цикл 10 в вершинах
3 7 6 2 изометрический цикл 13 в вершинах
1 1 3 2 1 3 1 2 1 3 3 3 1 3 3 3 вектор циклов по ребрам
1 -4 2 2 1 2 1 2 -4 1 1 -4 -4 2 2 2 удаленные ребра со знаком минус
-1 -1 1 1 1 1 -1 1 1 -1 -1 удаленные циклы из системы независимых циклов
цикл 1: {e4,e5,e16}
цикл 2: {e3,e4,e14}
цикл 3: {e14,e15,e16}
цикл 4: {e6,e8,e11,e15}
цикл 5: {e1,e3,e8}
цикл 6: {e6,e7,e10}
цикл 1: {v6,v7,v1}
цикл 2: {v5,v6,v1}
цикл 3: {v6,v7,v5}
цикл 4: {v3,v7,v5,v2}
цикл 5: {v2,v5,v1}
цикл 6: {v3,v4,v2}
Расчет 1
…………………………………………………………………………..
изометрические циклы в последовательности
6 7 10 изометрический цикл 1
7 9 13 изометрический цикл 2
3 5 15 изометрический цикл 3
7 8 12 изометрический цикл 4
14 15 16 изометрический цикл 5
1 2 6 изометрический цикл 6
3 4 14 изометрический цикл 7
6 9 11 16 изометрический цикл 8
2 3 10 12 изометрический цикл 9

| | |
|---|---|
| 4 5 16 | изометрический цикл 10 |
| 10 11 13 16 | изометрический цикл 11 |
| 1 4 9 | изометрический цикл 12 |
| 2 4 10 13 | изометрический цикл 13 |
| 2 5 11 | изометрический цикл 14 |
| 12 13 14 | изометрический цикл 15 |
| 8 9 14 | изометрический цикл 16 |
| 10 11 12 15 | изометрический цикл 17 |
| 1 3 8 | изометрический цикл 18 |
| 6 8 11 15 | изометрический цикл 19 |
| 19 | количество изометрических цикло |
| 1 4 8 11 15 18 21 24 28 34 34 | |
| 34 39 39 39 39 39 39 39 39 | массив указателей для независимой системы циклов |
| Массив признаков для базиса | |
| 1 1 1 1 1 1 1 1 1 0 0 1 0 0 0 0 0 0 0 | признак 0 – цикл удален, 1 – цикл |
| Цикломатическое число графа = 10 | (количество циклов в базисе) |
| 6 7 10 | изометрический цикл 1 |
| 7 9 13 | изометрический цикл 2 |
| 3 5 15 | изометрический цикл 3 |
| 7 8 12 | изометрический цикл 4 |
| 14 15 16 | изометрический цикл 5 |
| 1 2 6 | изометрический цикл 6 |
| 3 4 14 | изометрический цикл 7 |
| 6 9 11 16 | изометрический цикл 8 |
| 2 3 10 12 | изометрический цикл 9 |
| 1 4 9 | изометрический цикл 12 |
| 3 4 2 | изометрический цикл 1 в вершинах |
| 4 6 2 | изометрический цикл 2 в вершинах |
| 5 7 1 | изометрический цикл 3 в вершинах |
| 4 5 2 | изометрический цикл 4 в вершинах |
| 6 7 5 | изометрический цикл 5 в вершинах |
| 2 3 1 | изометрический цикл 6 в вершинах |
| 5 6 1 | изометрический цикл 7 в вершинах |
| 3 7 6 2 | изометрический цикл 8 в вершинах |
| 3 4 5 1 | изометрический цикл 9 в вершинах |
| 2 6 1 | изометрический цикл 12 в вершинах |
| 2 2 3 2 1 3 3 1 3 2 1 2 1 2 2 2 | вектор циклов по ребрам |
| 2 2 2 2 -4 2 2 1 1 2 -4 2 -4 2 1 1 | удаленные ребра со знаком минус |
| 1 -1 -1 1 1 1 1 -1 1 1 | -1 удаленные циклы из базиса |

цикл 1: {e6,e7,e10}
цикл 2: {e7,e8,e12}
цикл 3: {e14,e15,e16}
цикл 4: {e1,e2,e6}
цикл 5: {e3,e4,e14}
цикл 6: {e2,e3,e10,e12}
цикл 7: {e1,e4,e9}
цикл 1: {v3,v4,v2}
цикл 2: {v4,v5,v2}
цикл 3: {v6,v7,v5}
цикл 4: {v2,v3,v1}
цикл 5: {v5,v6,v1}
цикл 6: {v3,v4,v5,v1}
цикл 7: {v2,v6,v1}
Расчет 100

## Комментарии

Рассмотрен метод построения плоских конфигураций для случайной последовательности изометрических циклов, с применением модифицированного алгоритма Гаусса. Показана зависимость построения множества плоских конфигураций от порядка расположения элементов в кортеже изометорических циклов.

Приведен текст программы для решения задачи выделения множества плоских конфигураций в графе. Рассмотрены вопросы построения данных для решения задачи построения множества плоских конфигураций.

# Глава 9. ТОПОЛОГИЧЕСКИЙ РИСУНОК НЕПЛАНАРНОГО ГРАФА

## 9.1. Введение соединений в топологический рисунок графа

Рассмотрим процесс укладки ребер удаленных в процессе планаризации графа. Процесс укладки соединений будем рассматривать на примере топологического рисунка графа $G_{10}$ (см. рис. 9.2 или рис. 9.3).

Известно, что топологический рисунок планарного графа можно описать диаграммой вращения его вершин.

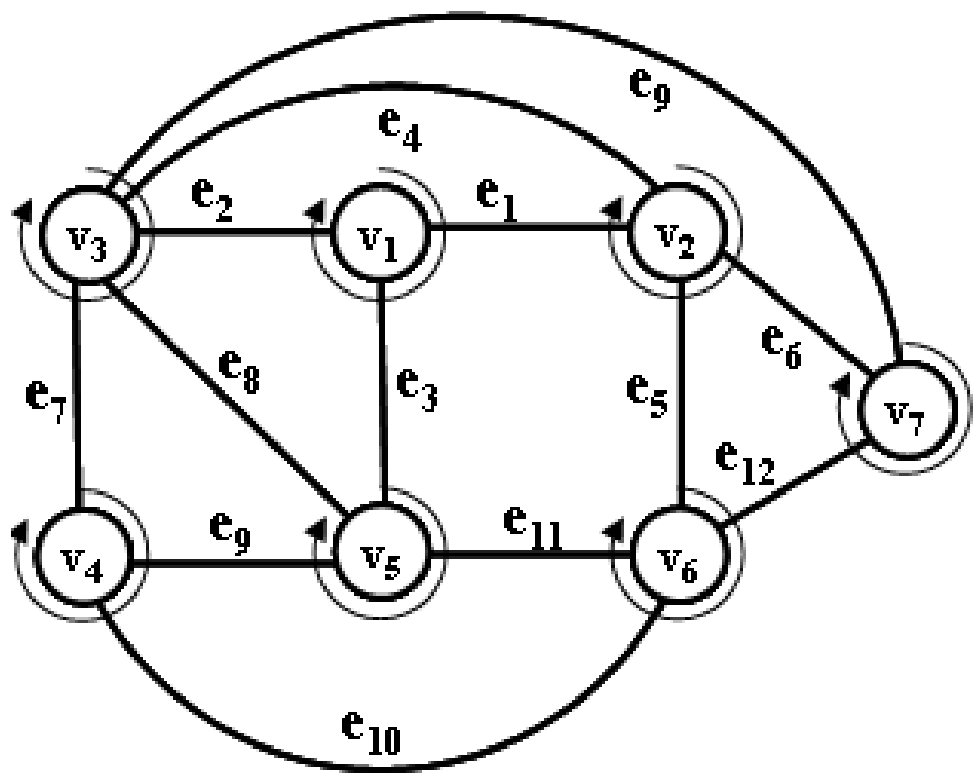

Рис. 9.1. Топологический рисунок плоской части графа $G_{10}$.

Например, топологический рисунок графа $(G_{10},\sigma)$ представленного на рис. 9.1, описывается следующей диаграммой вращения вершин:

$\sigma(v_1)$: $v_3$ $v_2$ $v_5$,
$\sigma(v_2)$: $v_1$ $v_3$ $v_7$ $v_6$,
$\sigma(v_3)$: $v_5$ $v_4$ $v_7$ $v_2$ $v_1$,
$\sigma(v_4)$: $v_3$ $v_5$ $v_6$,
$\sigma(v_5)$: $v_4$ $v_3$ $v_1$ $v_6$,
$\sigma(v_6)$: $v_4$ $v_5$ $v_2$ $v_7$,
$\sigma(v_7)$: $v_2$ $v_3$ $v_5$.

Топологический рисунок может быть записан в виде кортежей вращения вершин графа:

$\sigma(G_{10}) = \{ \sigma(v_1) = \langle v_3,v_2,v_5 \}, \sigma(v_2) = \langle v_1,v_3,v_7,v_6 \rangle, \sigma(v_3) = \langle v_5,v_4,v_7,v_2,v_1 \rangle,$
$\sigma(v_4) = \langle v_3,v_5,v_6 \rangle, \sigma(v_5) = \langle v_4,v_3,v_1,v_6 \rangle, \sigma(v_6) = \langle v_4,v_5,v_2,v_7 \rangle, \sigma(v_7) = \langle v_2,v_3,v_6 \rangle \}.$

Запись вращения вершин трактуется как циклическая запись ненулевых элементов соответствующей строки в матрицы смежностей графа.

В свою очередь, вращение вершин индуцирует ориентированные циклы. Следует заметить, что ориентация ребер в базисных циклов, направлена в противоположную сторону, за исключением обода $c_0$. Запишем циклы в виде циклического кортежа вершин, который можно преобразовать в сумму ориентированных ребер графа (векторов):

$c_1 = \langle v_3,v_4,v_5 \rangle = (v_3,v_4) + (v_4,v_5) + (v_5,v_3);$
$c_2 = \langle v_1,v_3,v_5 \rangle = (v_1,v_3) + (v_3,v_5) + (v_5,v_1);$
$c_3 = \langle v_1,v_5,v_6,v_2 \rangle = (v_1,v_5) + (v_5,v_6) + (v_6,v_2) + (v_2,v_1);$
$c_4 = \langle v_2,v_6,v_7 \rangle = (v_2,v_6) + (v_6,v_7) + (v_7,v_2);$
$c_5 = \langle v_4,v_6,v_5 \rangle = (v_4,v_6) + (v_6,v_5) + (v_5,v_4);$
$c_6 = \langle v_1,v_2,v_3 \rangle = (v_1,v_2) + (v_2,v_3) + (v_3,v_1);$

$c_7$ = <$v_2,v_7,v_3$> = $(v_2,v_7) + (v_7,v_3) + (v_3,v_2)$;
$c_0$ = <$v_3,v_7,v_6,v_4$> = $(v_3,v_7) + (v_7,v_6) + (v_6,v_4) + (v_4,v_3)$.

Данный вид записи нужно рассматривать не как арифметическое сложение ребер, а как последовательность ориентированных ребер в цикле. Например, циклическая запись цикла $c_7$ показывает, что за ребром $(v_2,v_7)$ следует ребро $(v_7,v_3)$, за ребром $(v_7,v_3)$ следует ребро $(v_3,v_2)$ и процесс объявляется замкнутым.

Предположим, что нам нужно в топологический рисунок графа G уложить ребро $\{v_2,v_4\}$. Естественно, что проведение данного соединения обязательно будет пересекать ребра плоской части графа. Будем вводить новые вершины как результат пересечения ребер.

**Определение 9.1.** Будем называть *мнимой вершиной – вершину графа образованную топологическим местоположением пересечения двух ребер.*

Методом поиска в ширину выделим все маршруты, соединяющие выбранные вершины графа. В данном случае, для графа $G_3$ имеем пять маршрутов для проведения соединения $(v_2,v_4)$. Будем отсчитывать маршруты от вершины $v_2$ к вершине $v_4$. Тогда можно записать выбранные маршруты в виде кортежей циклов $h_1$ = <$c_7,c_0$>, $h_2$ = <$c_4,c_0$>, $h_3$ = <$c_3,c_5$>, $h_4$ = <$c_6,c_2,c_1$>, $h_5$ = <$c_3,c_2,c_1$>. В общем случае для проведения соединений выбираются маршруты минимальной длины. Однако в нашем случае, для более детального рассмотрения процесса проведения удаленных соединений, выберем маршрут $h_5$ = <$c_3,c_2,c_1$>. Определим порядок проведения соединения $(v_2,v_4)$. При проведении соединения $\{v_2,v_4\}$ вначале пересекается ребро $e_3 = \{v_1,v_5\}$ принадлежащее циклу $c_3$ и циклу $c_2$. Затем пересекается ребро $e_8 = \{v_3,v_5\}$ принадлежащее циклу $c_2$ и циклу $c_1$. Следовательно, вводятся две мнимые вершины для исходного графа $G_3$. Обозначим эти вершины как $v_8$ и $v_9$.

Введем вершину $v_8$ в ребро $e_3$ циклов $c_3$ и $c_2$

$c_2^*$ = <$v_1,v_3,v_5,v_8$> = $(v_1,v_3) + (v_3,v_5) + (v_5,v_8) + (v_8,v_1)$;
$c_3^*$ = <$v_1,v_8,v_5,v_6,v_2$> = $(v_1,v_8) + (v_8,v_5) + (v_5,v_6) + (v_6,v_2) + (v_2,v_1)$.

Теперь введем вершину $v_9$ в ребро $e_8$ циклов $c_1$ и $c_2$

$c_1^*$ = <$v_3,v_4,v_5,v_9$> = $(v_3,v_4) + (v_4,v_5) + (v_5,v_9) + (v_9,v_3)$;
$c_2^{**}$ = <$v_1,v_3,v_9,v_5,v_8$> = $(v_1,v_3) + (v_3,v_9) + (v_9,v_5) + (v_5,v_8) + (v_8,v_1)$.

Осталось последовательно ввести части соединения $\{v_2,v_4\}$ в соответствующие циклы. Будем последовательно вводить ребро $\{v_2,v_8\}$ в цикл $c_3^*$, ребро $\{v_8,v_9\}$ в цикл $c_2^{**}$ и ребро $\{v_9,v_4\}$ в цикл $c_1^*$.

Вводим ребро $\{v_2,v_8\}$ в цикл $c_3^*$

$c_3^* = (v_1,v_8) + (v_8,v_5) + (v_5,v_6) + (v_6,v_2) + (v_2,v_1) + \ + (v_2,v_8) + (v_8,v_2) =$
$= (v_8,v_5) + (v_5,v_6) + (v_6,v_2) + \ + (v_2,v_8) +$
$+ (v_2,v_1) + (v_1,v_8) + \ + (v_8,v_2).$

В результате образуются два новых цикла:

$c_3^{**} = (v_8,v_5) + (v_5,v_6) + (v_6,v_2) + (v_2,v_8);$
$c_8 = (v_2,v_1) + (v_1,v_8) + (v_8,v_2).$

Вводим ребро $\{v_8,v_9\}$ в цикл $c_2^{**}$

$c_2^{**} = (v_1,v_3) + (v_3,v_9) + (v_9,v_5) + (v_5,v_8) + (v_8,v_1) + \ + (v_9,v_8) + (v_8,v_9) =$
$= (v_8,v_1) + (v_1,v_3) + (v_3,v_9) + \ + (v_9,v_8) +$
$+ (v_9,v_5) + (v_5,v_8) + \ + (v_8,v_9).$

В результате образуются два новых цикла:

$c_2^{***} = (v_8,v_1) + (v_1,v_3) + (v_3,v_9) + (v_9,v_8);$
$c_9 = (v_9,v_5) + (v_5,v_8) + (v_8,v_9).$

Вводим ребро $\{v_9,v_4\}$ в цикл $c_1^{*}$

$c_1^{*} = (v_3,v_4) + (v_4,v_5) + (v_5,v_9) + (v_9,v_3) + \ + (v_9,v_4) + (v_4,v_9) =$
$= (v_4,v_5) + (v_5,v_9) + \ + (v_9,v_4) +$
$+ (v_9,v_3) + (v_3,v_4) + \ + (v_4,v_9).$

В результате образуются два новых цикла:

$c_1^{**} = (v_4,v_5) + (v_5,v_9) + (v_9,v_4);$
$c_{10}= (v_9,v_3) + (v_3,v_4) + (v_4,v_9).$

Новая система циклов имеет вид:

$c_1^{**} = \langle v_4,v_5,v_9\rangle = (v_4,v_5) + (v_5,v_9) + (v_9,v_4);$
$c_2^{***}= \langle v_8,v_1,v_3,v_9\rangle = (v_8,v_1) + (v_1,v_3) + (v_3,v_9) + (v_9,v_8);$
$c_3^{**} = \langle v_8,v_5,v_6,v_2\rangle = (v_8,v_5) + (v_5,v_6) + (v_6,v_2) + (v_2,v_8);$
$c_4 = \langle v_2,v_6,v_7\rangle = (v_2,v_6) + (v_6,v_7) + (v_7,v_2);$
$c_5 = \langle v_4,v_6,v_5\rangle = (v_4,v_6) + (v_6,v_5) + (v_5,v_4);$
$c_6 = \langle v_1,v_2,v_3\rangle = (v_1,v_2) + (v_2,v_3) + (v_3,v_1);$
$c_7 = \langle v_2,v_7,v_3\rangle = (v_2,v_7) + (v_7,v_3) + (v_3,v_2);$
$c_8 = \langle v_2,v_1,v_8\rangle = (v_2,v_1) + (v_1,v_8) + (v_8,v_2);$
$c_9 = \langle v_9,v_5,v_8\rangle = (v_9,v_5) + (v_5,v_8) + (v_8,v_9);$
$c_{10}= \langle v_9,v_3,v_4\rangle = (v_9,v_3) + (v_3,v_4) + (v_4,v_9);$
$c_0 = \langle v_3,v_7,v_6,v_4\rangle = (v_3,v_7) + (v_7,v_6) + (v_6,v_4) + (v_4,v_3).$

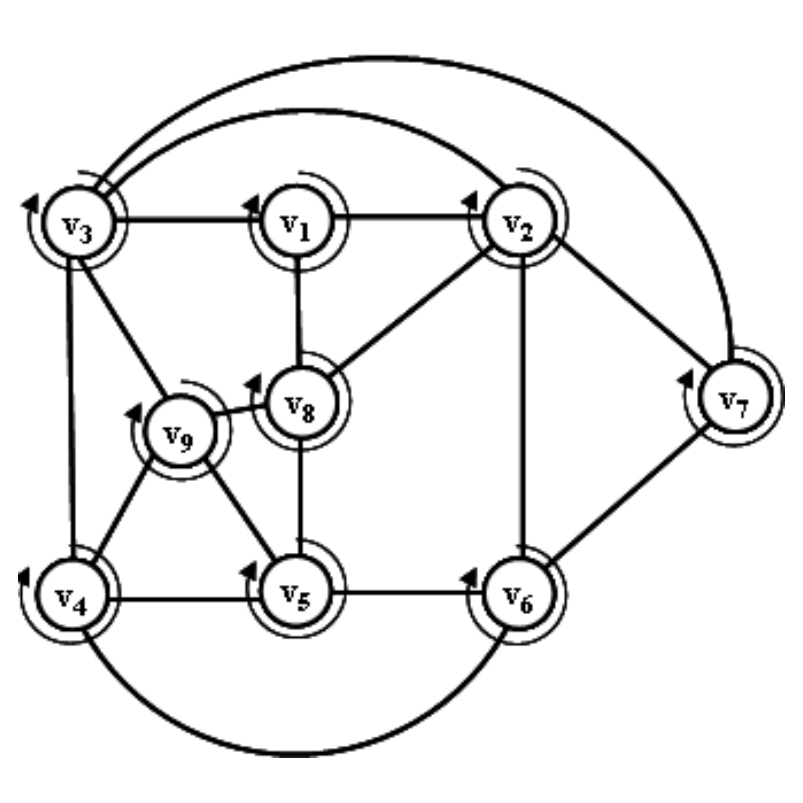

Рис. 9.2. Топологический рисунок графа $G_{10}$ с двумя точками пересечения.

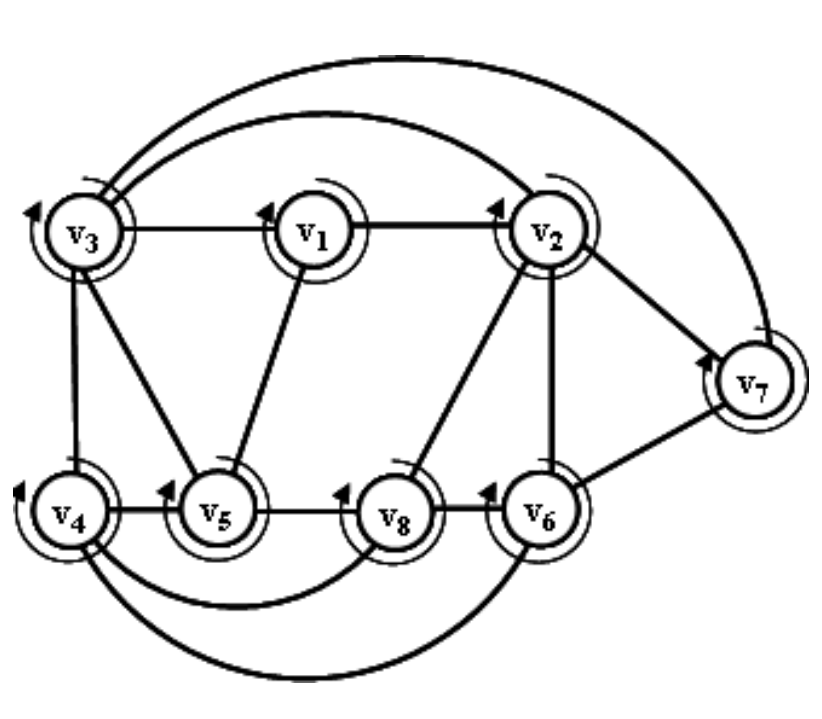

Рис. 9.3. Топологический рисунок графа $G_{10}$ с одной точкой пересечения.

В результате получен топологический рисунок непланарного графа $G_1$ с двумя точками пересечения (см. рис. 9.2). Применяя алгоритм преобразования циклов во вращение вершин, построим диаграмму вращения вершин:

$\sigma(v_1)$: $v_3$ $v_2$ $v_5$,
$\sigma(v_2)$: $v_1$ $v_3$ $v_7$ $v_6$ $v_8$,
$\sigma(v_3)$: $v_5$ $v_4$ $v_7$ $v_2$ $v_9$,
$\sigma(v_4)$: $v_3$ $v_9$ $v_5$ $v_6$,
$\sigma(v_5)$: $v_4$ $v_9$ $v_8$ $v_6$,
$\sigma(v_6)$: $v_4$ $v_5$ $v_2$ $v_7$,
$\sigma(v_7)$: $v_2$ $v_3$ $v_6$,
$\sigma(v_8)$: $v_1$ $v_2$ $v_5$ $v_9$,
$\sigma(v_9)$: $v_4$ $v_3$ $v_8$ $v_5$.

Если воспользоваться маршрутом $h_3$, то полученный топологический рисунок непланарного графа $G_{10}$, будет характеризоваться только одной точкой пересечения (см. рис. 9.3).

## 9.2. Построение рисунка графа с минимальным числом пересечений.

**Определение 9.2.** В теории графов *число пересечений* cr(*G*) графа G — это наименьшее число пересечений рёбер рисунка графа G.

Например, граф является планарным тогда и только тогда, когда его число пересечений равно нулю.

Очень полезное *неравенство для числа пересечений* обнаружили независимо Аитаи, Хватал, Ньюборн (Newborn) и Семереди и Лейтон [27]:

Для неориентированных простых графов G с $n$ вершинами и $m$ рёбрами, таких, что $m > 7n$ имеем:

$$cr(G) \geq \frac{m^3}{29n^2} \tag{9.1}$$

Константа 29 является наилучшей известной. Согласно Акерману, если понизить условие $m > 4n$, то это будет стоить заменой константы 29 на 64.

$$cr(G) \geq \frac{m^3}{64n^2} \tag{9.2}$$

Небольшое изменение приведённой аргументации позволяет заменить 64 на 33,75 для $m > 7{,}5n$.

Так как в процессе укладки ребер, происходит разбиение ребер на части и появление новых циклов, поэтому в целом процесс построения минимального количества пересечений ребер, определяется последовательностью проведения соединений. Данный процесс будем рассматривать с помощью перестановок ребер в последовательности. Каждую перестановку будем описывать последовательностью элементов кортежа, $\varphi_i =< e_1, e_2, ...e_k >$, где $i$ – номер кортежа, $k$ – количество удаленных ребер в процессе планаризации графа.

Выделим в графе $G_{11}$ плоскую часть. В результате планаризации удалены ребра $e_2,e_5,e_{22}$.

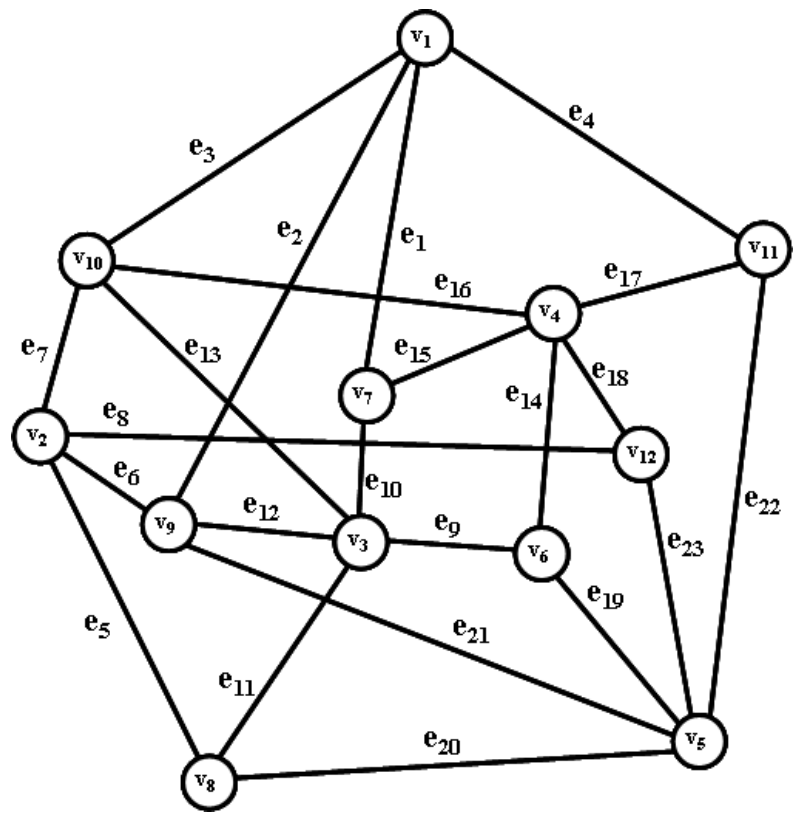

Рис. 9.4. Неориентированный несепарабельный граф $G_{11}$.

Определим плоскую часть графа:

$c_1 = \{e_6,e_7,e_{12},e_{13}\} \to \{v_2,v_3,v_9,v_{10}\}$;
$c_2 = \{e_1,e_3,e_{10},e_{13}\} \to \{v_1,v_3,v_7,v_{10}\}$;
$c_3 = \{e_3,e_4,e_{10},e_{17}\} \to \{v_1,v_4,v_{10},v_{11}\}$;
$c_4 = \{e_1,e_4,e_{15},e_{17}\} \to \{v_1,v_4,v_7,v_{11}\}$;
$c_5 = \{e_9,e_{10},e_{14},e_{15}\} \to \{v_3,v_4,v_6,v_7\}$;
$c_6 = \{e_9,e_{11},e_{19},e_{20}\} \to \{v_3,v_5,v_6,v_8\}$;
$c_7 = \{e_{11},e_{12},e_{20},e_{21}\} \to \{v_3,v_5,v_9,v_8\}$;
$c_8 = \{e_{14},e_{18},e_{19},e_{23}\} \to \{v_4,v_5,v_6,v_{12}\}$;
$c_9 = \{e_6,e_8,e_{21},e_{23}\} \to \{v_2,v_5,v_9,v_{12}\}$;
$c_0 = \{e_7,e_8,e_{16},e_{18}\} \to \{v_2,v_4,v_{10},v_{12}\}$;

Случайным образом выберем перестановку удаленных ребер: $K_r = \langle e_2,e_5,e_{22}\rangle$. Используя алгоритм поиска в ширину, определим проведение первого ребра $e_2$ (см. рис. 9.5).

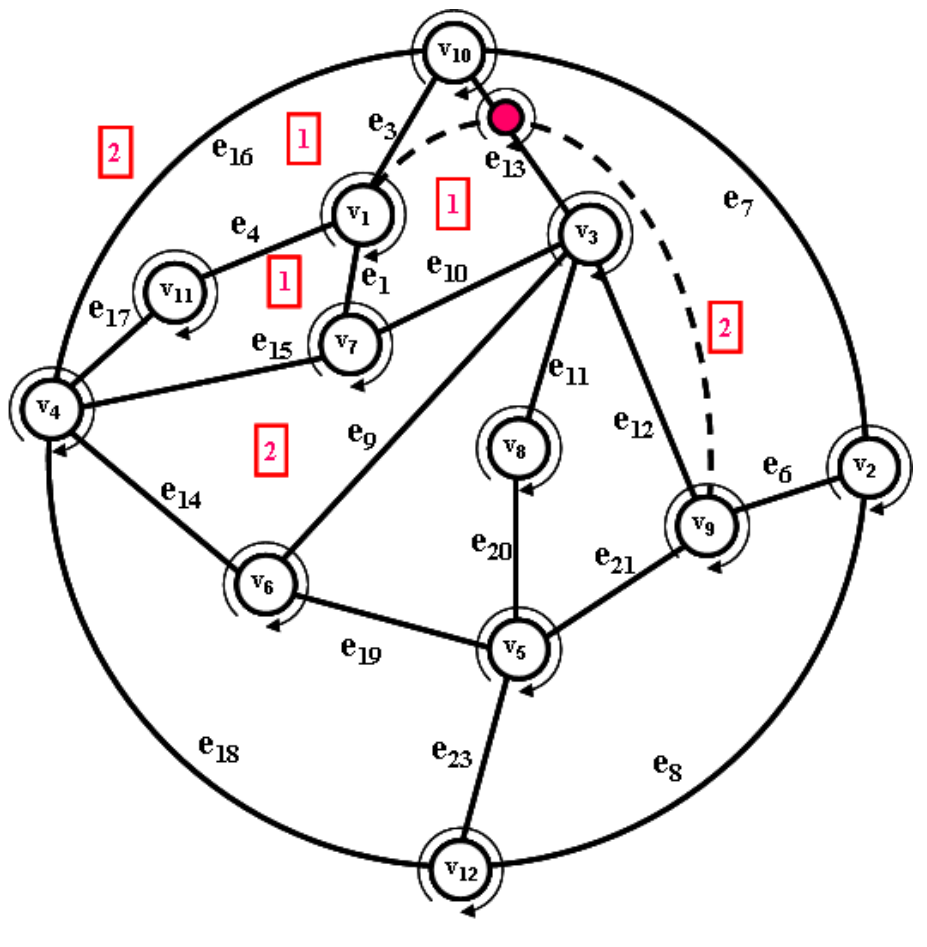

Рис. 9.5. Проведение ребра $e_2$.

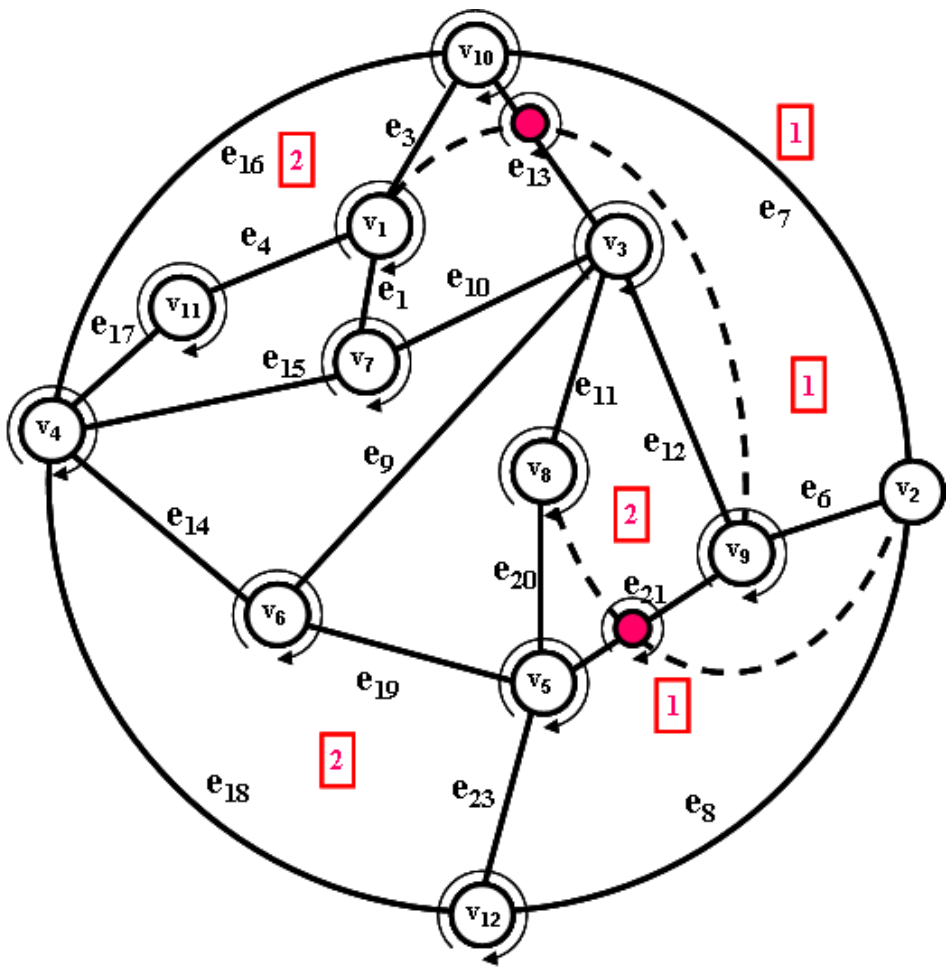

Рис. 9.5. Проведение ребра $e_5$.

Используя алгоритм поиска в ширину [27,36], определим проведение второго ребра $e_5$ (см. рис. 10.6) и третьего ребра $e_{22}$ (см. рис. 9.7).

Рассмотрим другой порядок проведения удаленных ребер. Пусть задан кортеж $\langle e_5,e_{22},e_2\rangle$. Применяя алгоритм поиска в ширину, построим минимальный маршрут для проведения ребра $e_5$ (см. рис. 9.8).

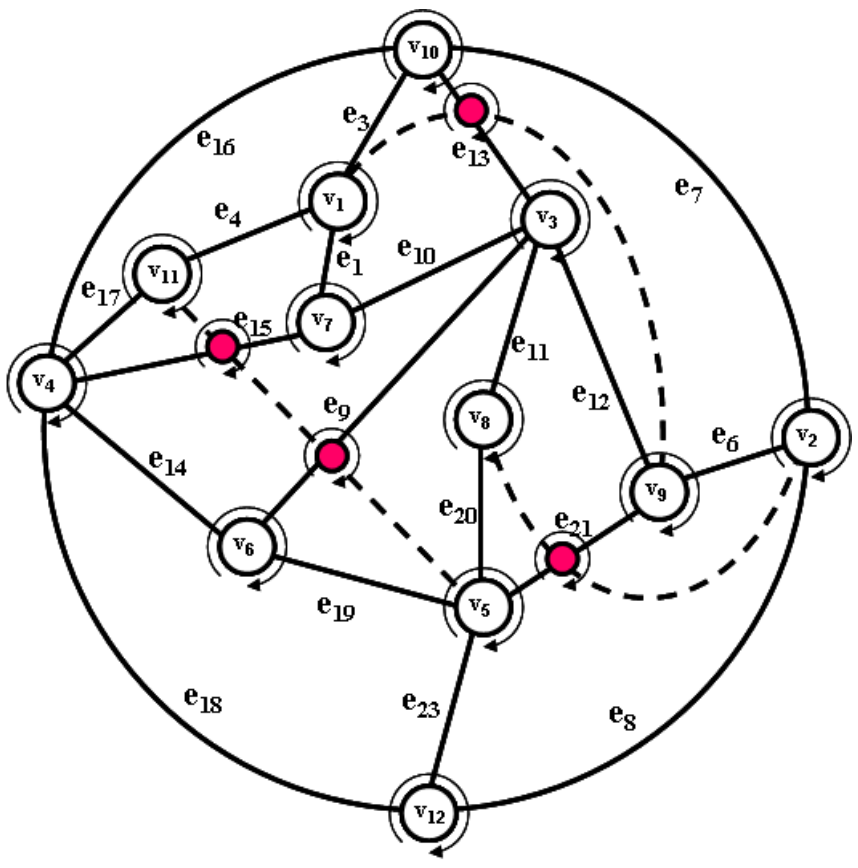

Рис. 9.7. Проведение ребра $e_{22}$.

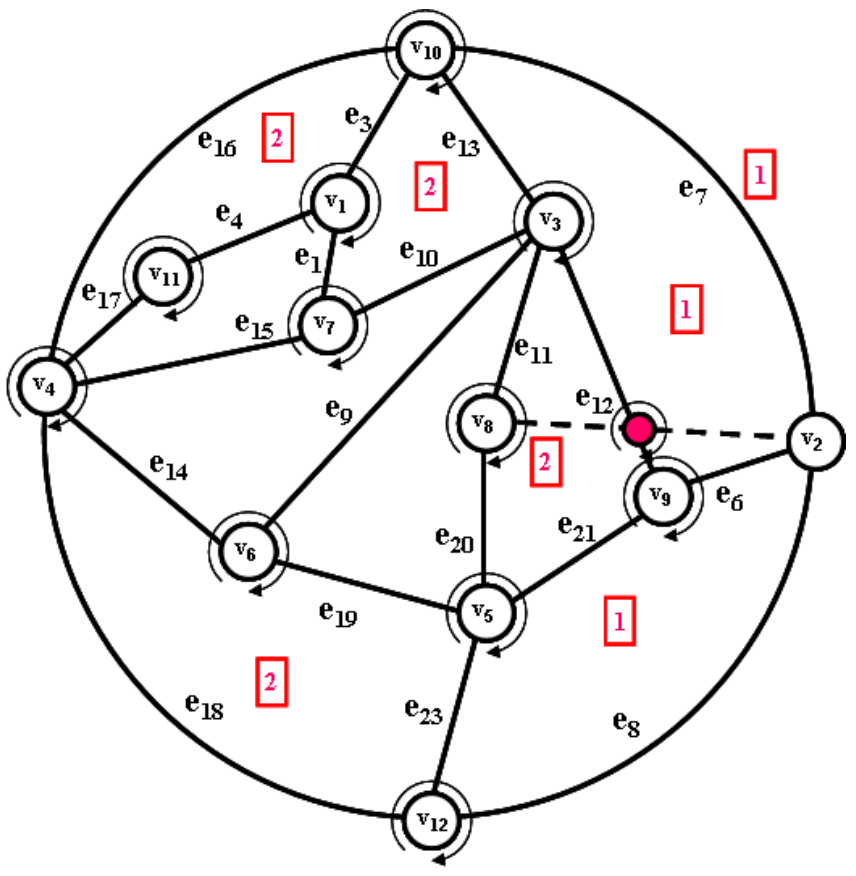

Рис. 9.8. Проведение ребра $e_{25}$.

Тем же способом построим минимальные маршруты для проведениея ребер $e_{22}$ (см. рис. 10.9) и $e_2$ (см. рис. 9.10).

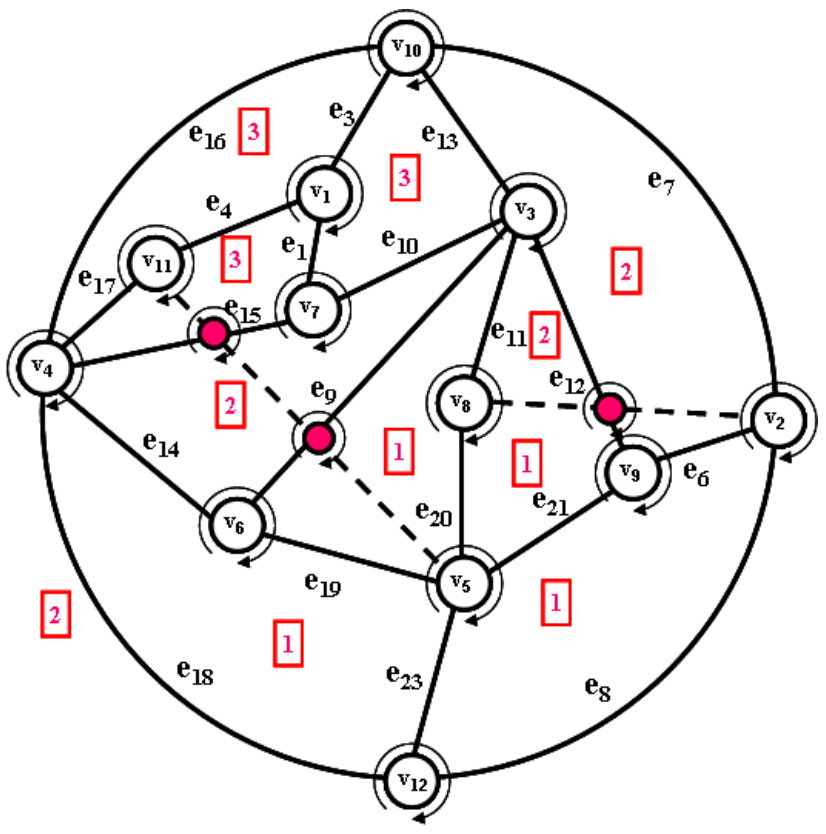

Рис. 9.9. Проведение ребра $e_{22}$.

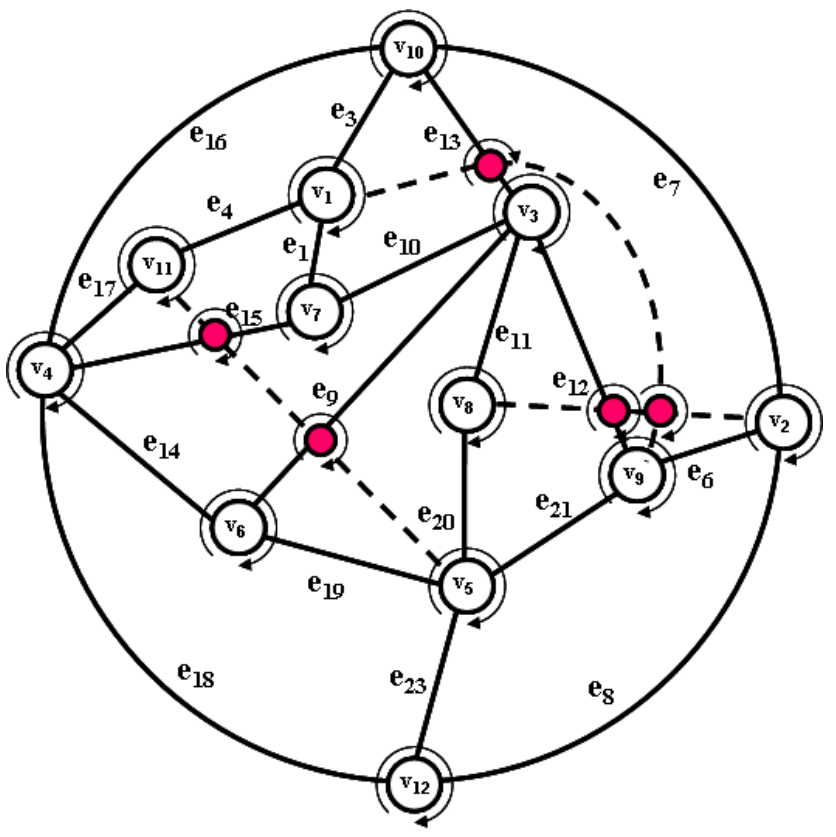

Рис. 9.10. Проведение ребра $e_5$.

В случае проведения соединений относительно первой последовательности, описываемой кортежем перестановок <$e_2$,$e_5$,$e_{22}$>, мы получили четыре пересечения. А в случае проведения соединений согласно второй последовательности, мы получили пять пересечений.

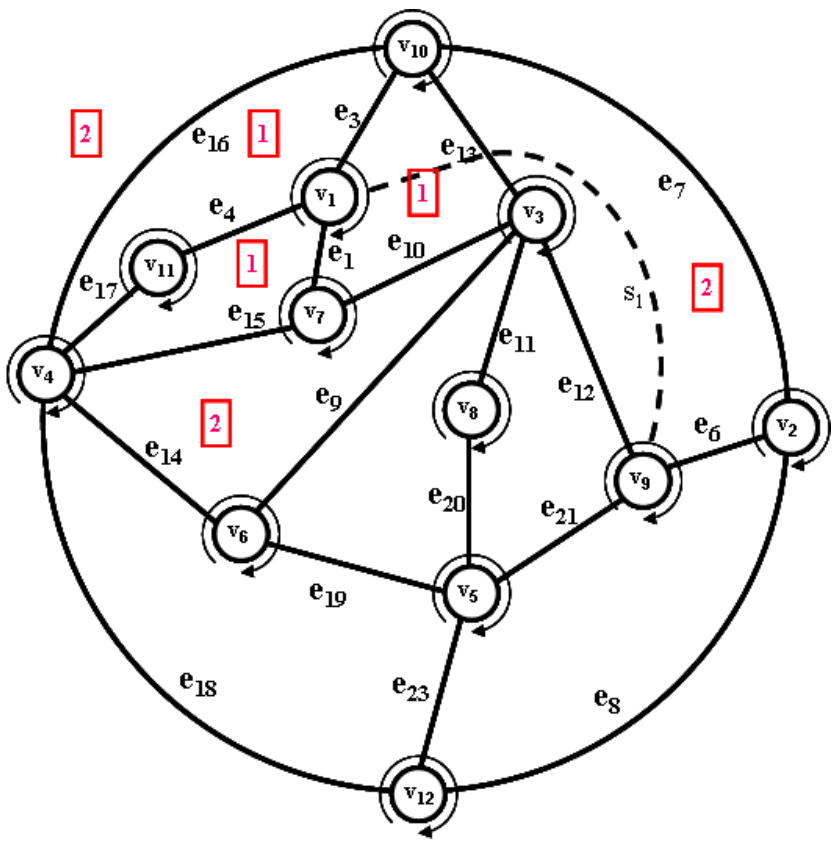

Рис. 9.11. Маршрут для ребра $e_2$.

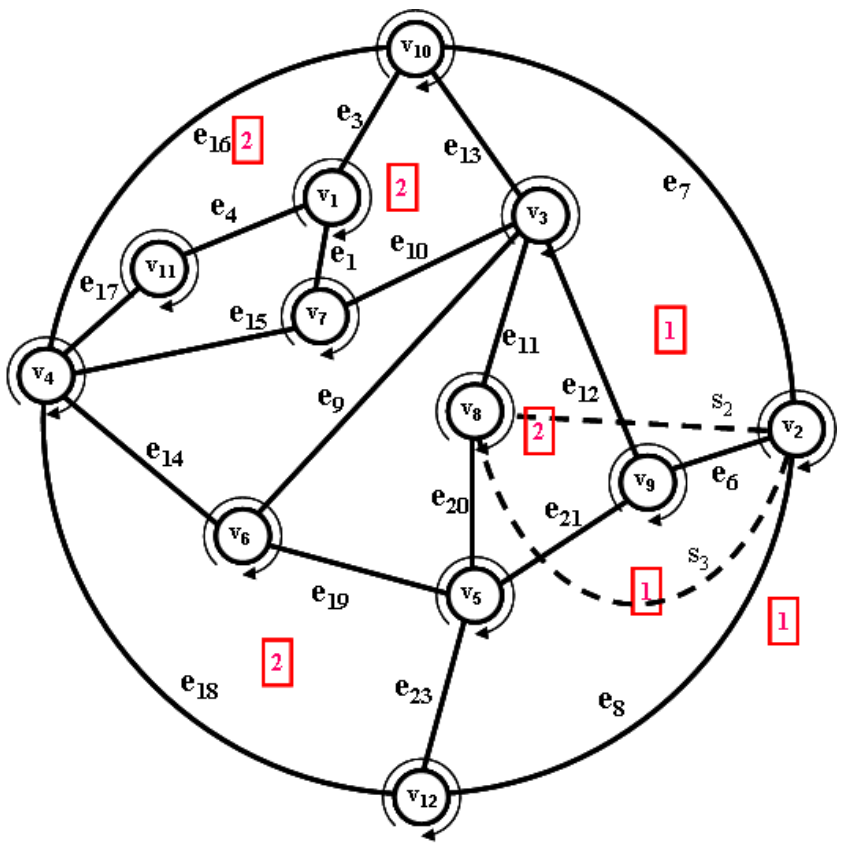

Рис. 9.12. Маршруты для ребра $e_5$.

Оказывается, что количество соединений зависит не только от последовательно проведения соединений, но и от минимального маршрута соединения. Так как одно и то же соединение может иметь несколько минимальных маршрутов.

Соединение $e_2$ имеет единственный минимальный маршрут $s_1$ (см. рис. 9.11). Соединение $e_5$ имеет два минимальных маршрута $s_2$ и $s_3$ (см. рис. 9.12). Соединение $e_{22}$ имеет два минимальных маршрута $s_4$ и $s_5$ (см. рис. 9.13).

Для описания процесса проведения соединений воспользуеся фрагментарной структурой последовательности (перестановки) [10].

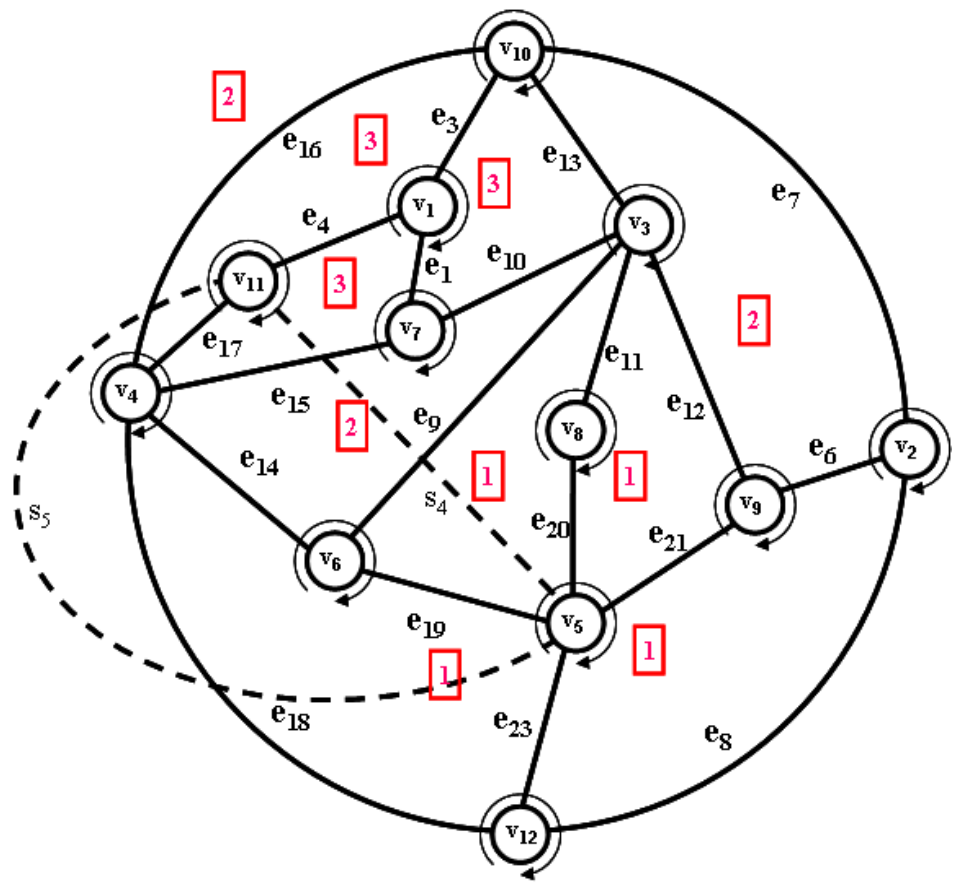

Рис. 9.13. Маршруты для ребра $e_{22}$.

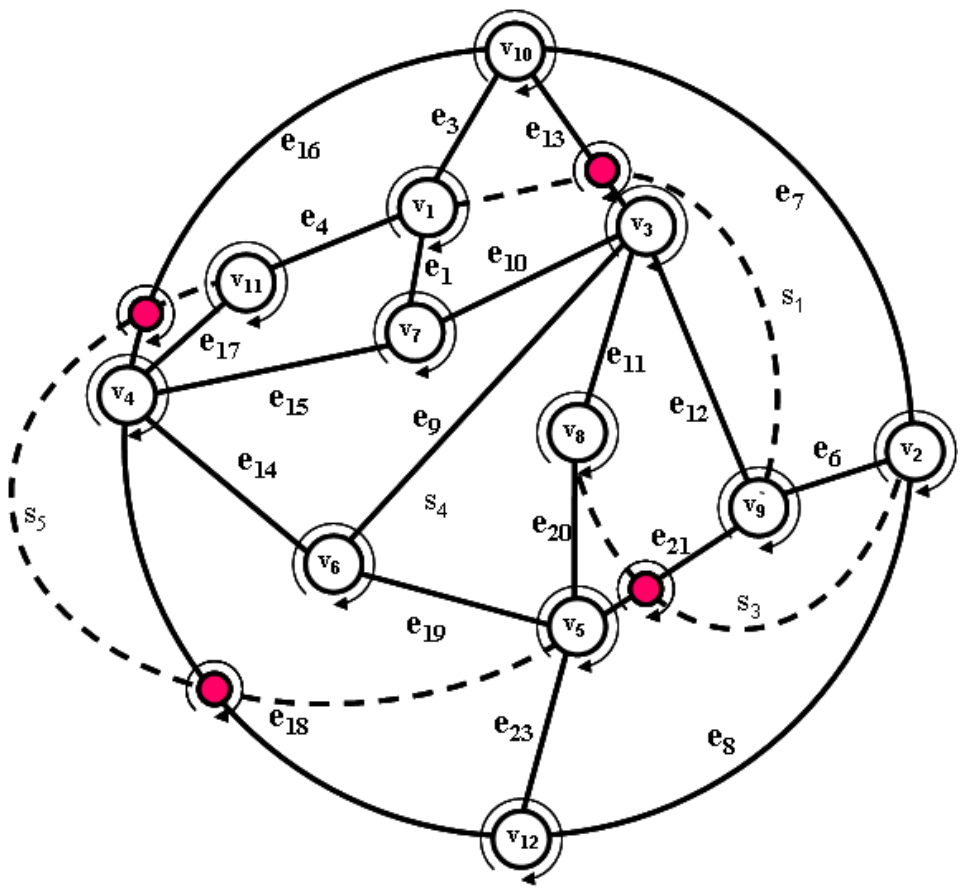

Рис. 9.14. Рисунок графа.

Будем обозначать минимальные маршруты соединением и составом последовательностью циклов принадлежащих этому маршруту. Для данной плоской части имеем пять минимальных маршрутов:

$s_1 = e_2, \langle c_2, c_1 \rangle$;
$s_2 = e_5, \langle c_1, c_7 \rangle$;
$s_3 = e_5, \langle c_9, c_7 \rangle$;
$s_4 = e_{22}, \langle c_6, c_5, c_4 \rangle$;
$s_5 = e_{22}, \langle c_8, c_0, c_3 \rangle$.

Маршрут идентифицируется именем ребра, затем, через запятую, следует последовательность циклов, принадлежащих минимальному маршруту.

С учетом сказанного, рассмотрим следующий кортеж проведения соединений $K_r = \langle s_3, s_1, s_2, s_5, s_4 \rangle$ состоящий из всех маршрутов. Формируем новый кортеж, состоящий из трех минимальных маршрутов равных числу соединений.

$\langle s_3, s_1, s_2, s_5, s_4 \rangle \rightarrow \langle \mathbf{s_3}, s_1, \mathbf{s_2}, s_5, s_4 \rangle \rightarrow \langle \mathbf{s_3}, \mathbf{s_1}, \mathbf{s_2}, s_5, s_4 \rangle \rightarrow \langle \mathbf{s_3}, \mathbf{s_1}, \mathbf{s_2}, \mathbf{s_5}, \mathbf{s_4} \rangle$.

На первом шаге, выбираем маршрут $s_3$ и одновременно удаляем из кортежа, машруты принадлежащие ребру $e_5$. В данном случае это маршрут $s_2$ (обозначим его красным цветом). На втором шаге, выбираем следующий не помесенный маршрут $s_1$, это единственный маршрут. На третьем шаге, выбираем маршрут $s_5$, одновременно удаляя все другие маршруты данного ребра $e_{22}$. Это маршрут $s_4$. Исключив из рассмотрения помеченные маршруты, строим топологический рисунок графа с пересекающимися ребрами.

Топологический рисунок с персечением представлен на рис. 9.14.

Результат расчета метаэвристического алгоритма [10,34,39] основанного на фрагментарных структурах представляет собой прямоугольную матрицу $q_r \times w_k$, где $q_r$ – количество плоских конфигураций с минимальным числом удаленных ребер, а $w_k$ – заданное определенное количество перестановок для проведения соединений.

### 9.3. Рисунок графа минимальной толщины

**Определение 9**.3. *Толщина графа G* — это наименьшее число плоских суграфов, на которые можно разложить рёбра графа G. То есть, если существует набор *k* плоских графов, имеющих одинаковый набор вершин, объединение которых даёт граф G, то толщина графа G не больше *k*.

Таким образом, плоский граф имеет толщину единица [36]. Всесторонний обзор по теме толщины графа (по состоянию на 1998 год) написан Мутцелем, Оденталем и Шарбродтом.

**Определение 9.4**. Будем называть *слоем* топологический рисунок плоского суграфа с максимальным количеством непересекающихся ребер.

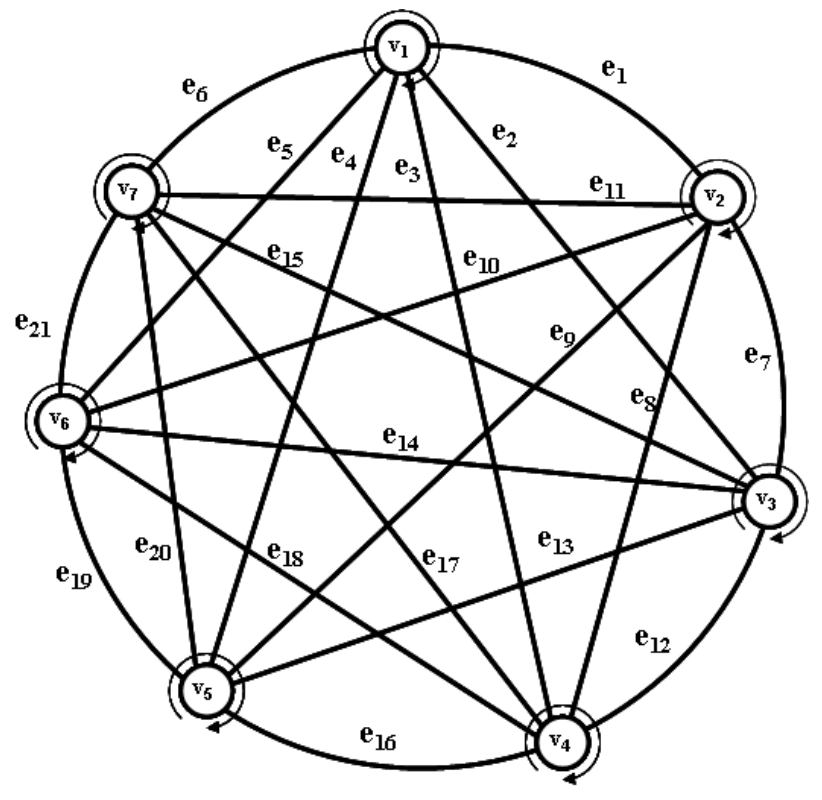

Рис. 9.15. Полный граф $K_7$

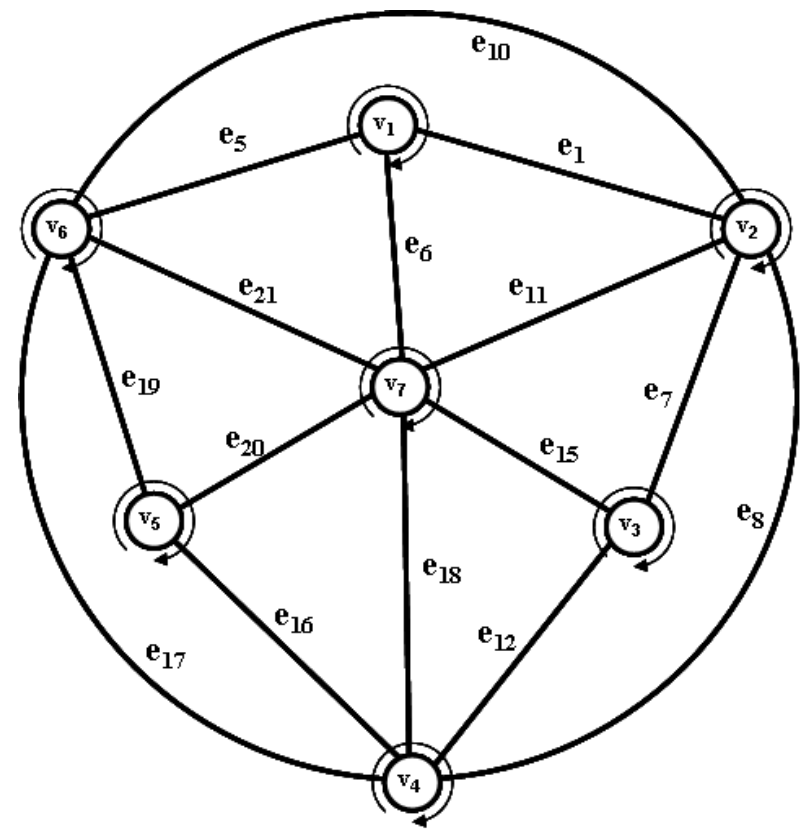

Рис. 9.15. Плоская часть графа $K_7$.

Толщина полного графа с *n* вершинами, $K_n$, равна

$$\left\lfloor \frac{n+7}{6} \right\rfloor \tag{9.3}$$

за исключением случаев $n = 9$, 10, для которых толщина равна трём.

За исключением нескольких случаев, толщина полного двудольного графа $K_{a,b}$ равна

$$\left\lceil \frac{ab}{2(a+b-2)} \right\rceil \tag{9.4}$$

Метод построения топологического рисунка графа с минимальной толщиной, ничем не отличается от принципа построения топологического рисунка с минимальным числом пересечений ребер, кроме выполнения условия невозможности пересечения ребер в одном слое.

Так как проведения соединений (ребер) удаленных в процессе планаризации осуществляется последовательно, то алгоритмом Фишера-Йетса выделим случайный кортеж последовательности вносимых в плоский топологический рисунок соединений.

$K_r$ = <4,2,6,1,5,3> $\rightarrow$ < $e_9,e_3,e_{14},e_2,e_{13},e_4$ >.

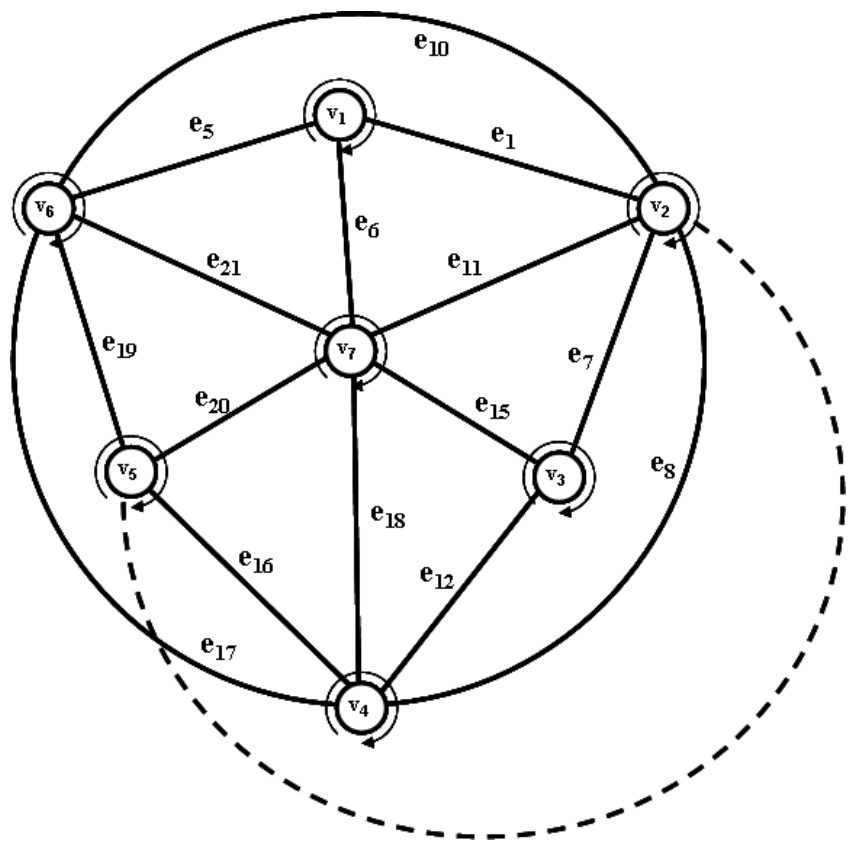

Рис. 9.17. Проведение соединения $e_9$.

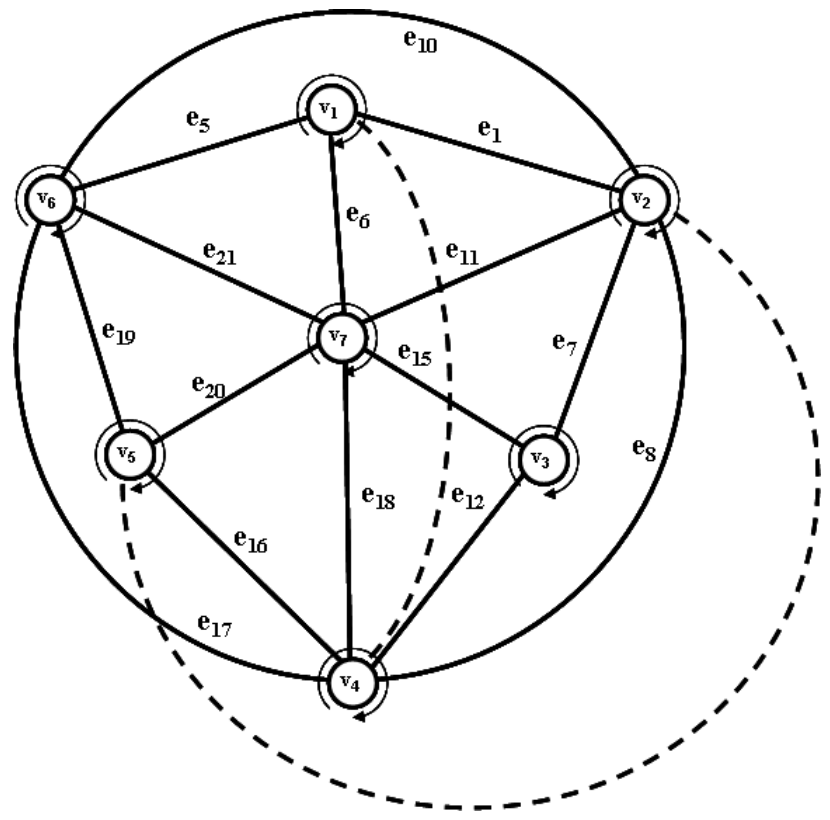

Рис. 9.18. Проведение соединения $e_3$.

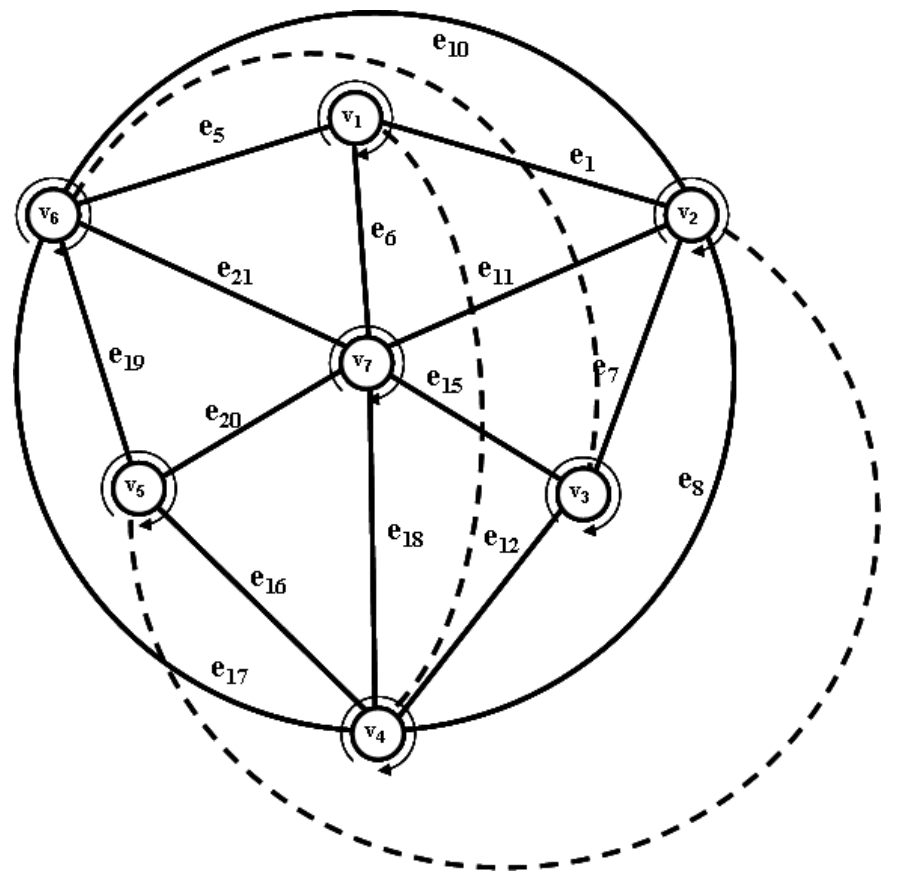

Рис. 9.19. Проведение соединения $e_{14}$.

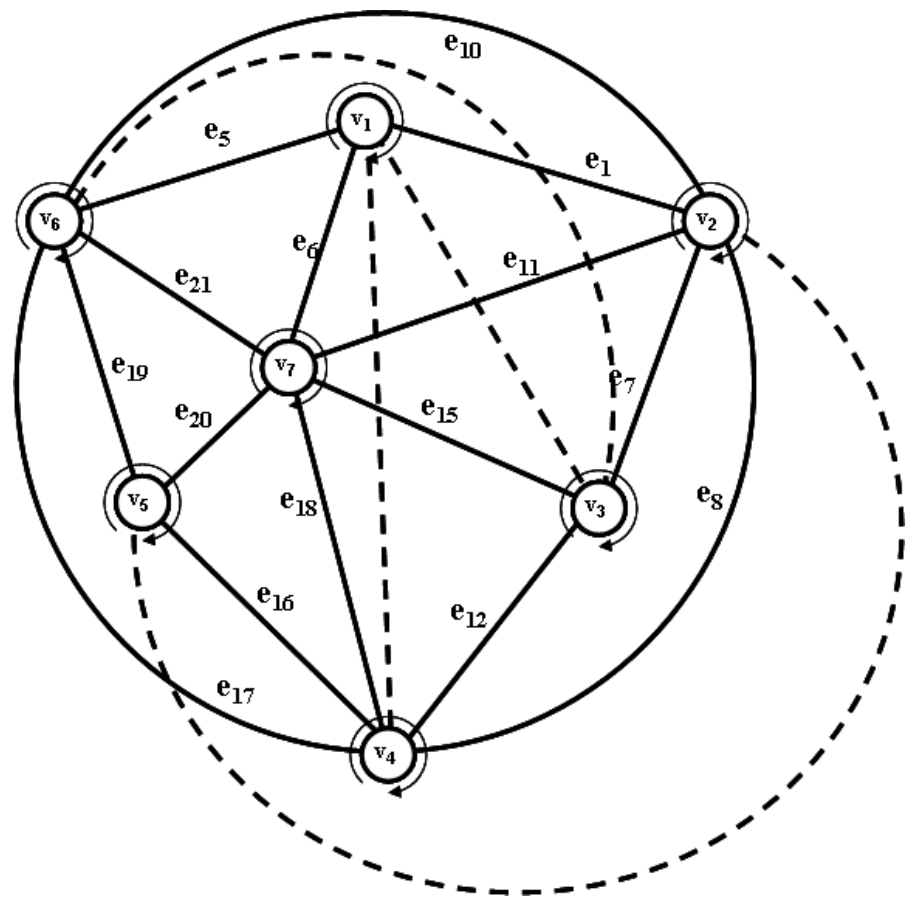

Рис. 9.20. Проведение соединения $e_2$.

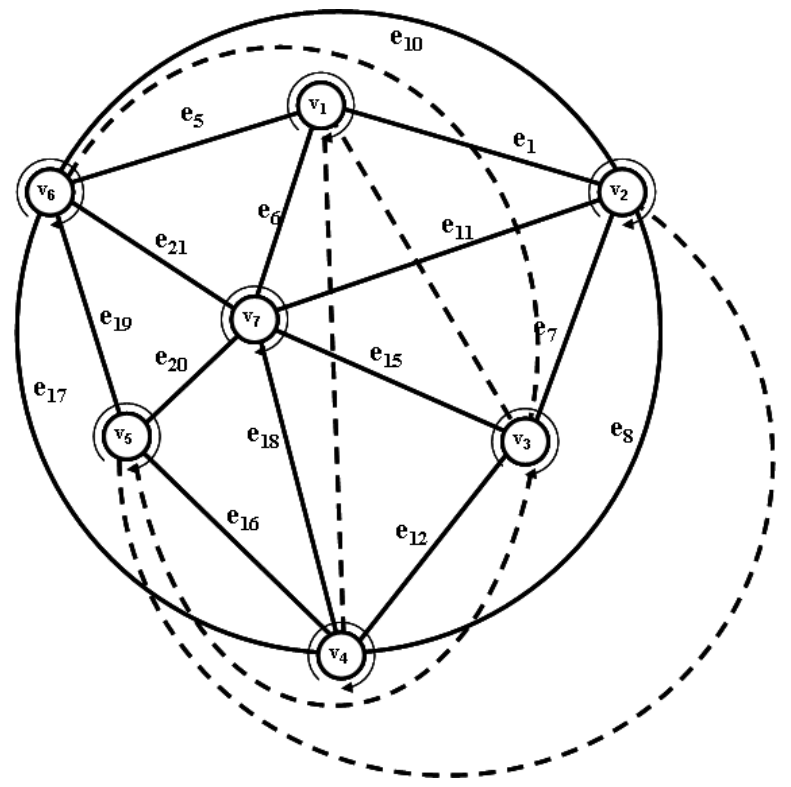

Рис. 9.21. Проведение соединения $e_{13}$.

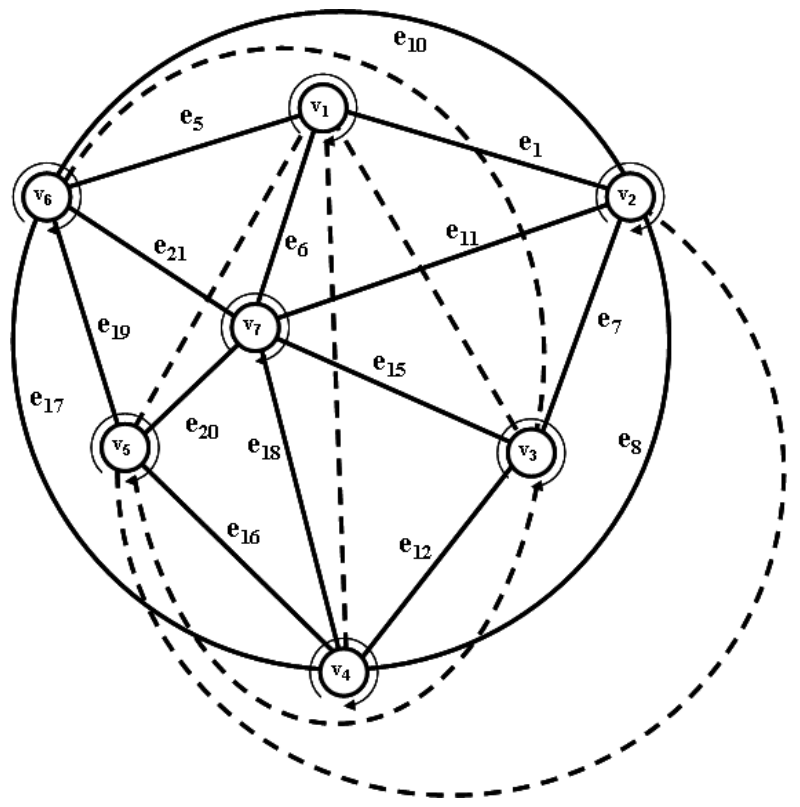

Рис. 9.22. Проведение соединения $e_4$.

Процесс проведения соединений осуществляется по маршруту который определяется алгоритмом поиска в ширину. Процесс проведения соединений представлен на рис. 9.17 – рис. 9.22.

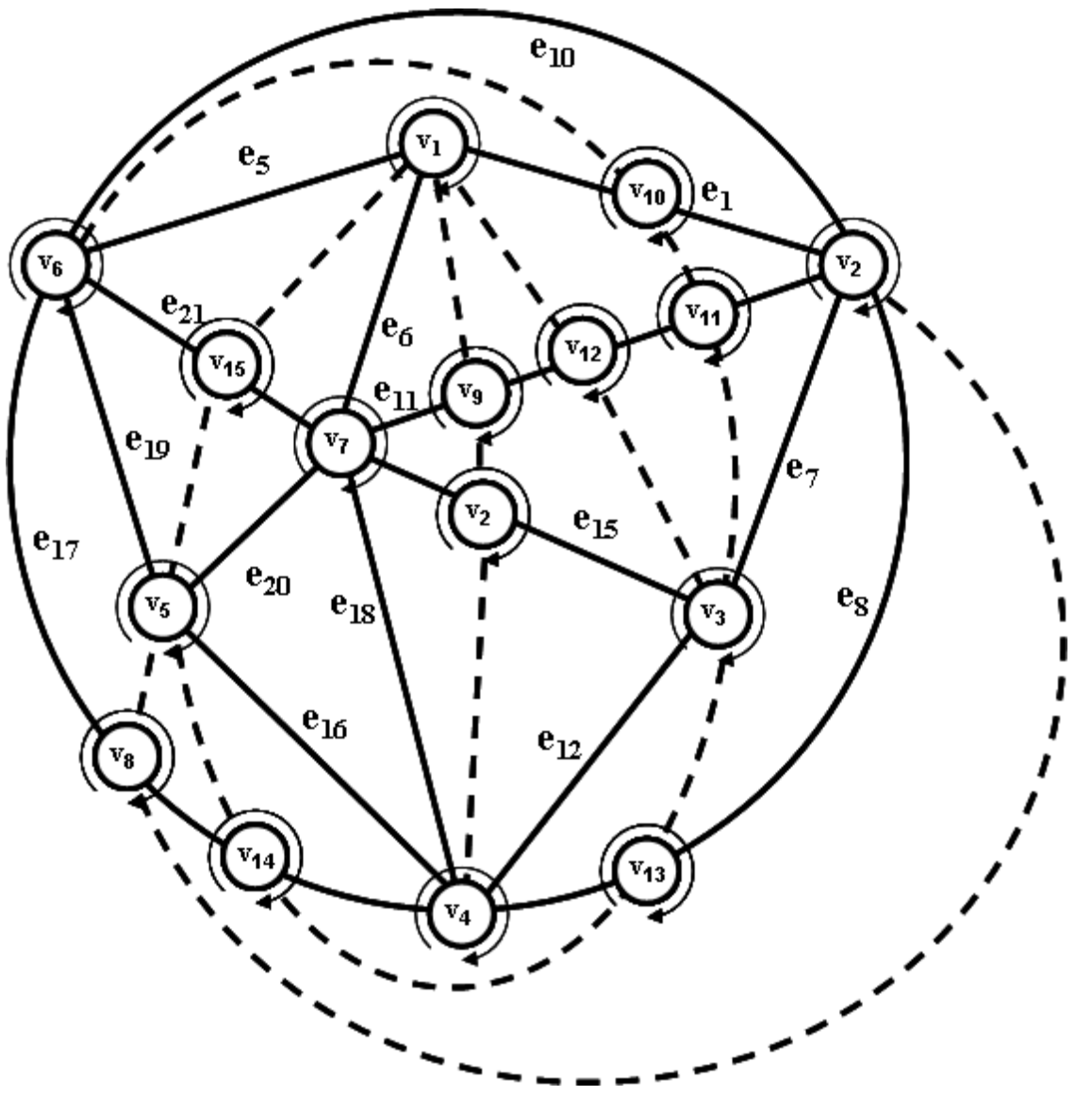

Рис. 9.23. Топологический рисунок графа минимальной толщины.

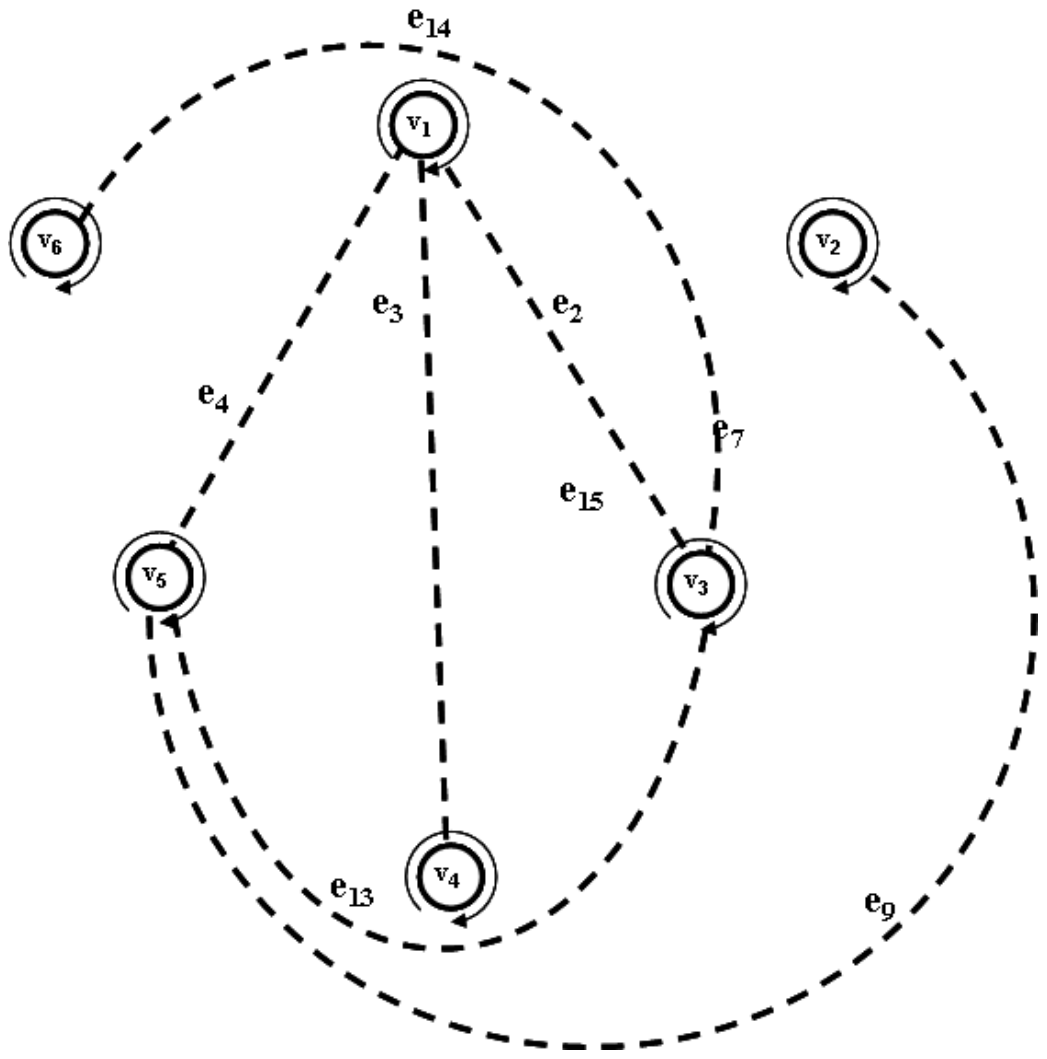

Рис. 9.24. Топологический рисунок второго слоя.

Определим следующую последовательность проведения соединений.

$K_г$ = <5,1,4,6,3,2> $\rightarrow$ <$e_{13}$,$e_2$,$e_9$,$e_{14}$,$e_4$,$e_3$>.

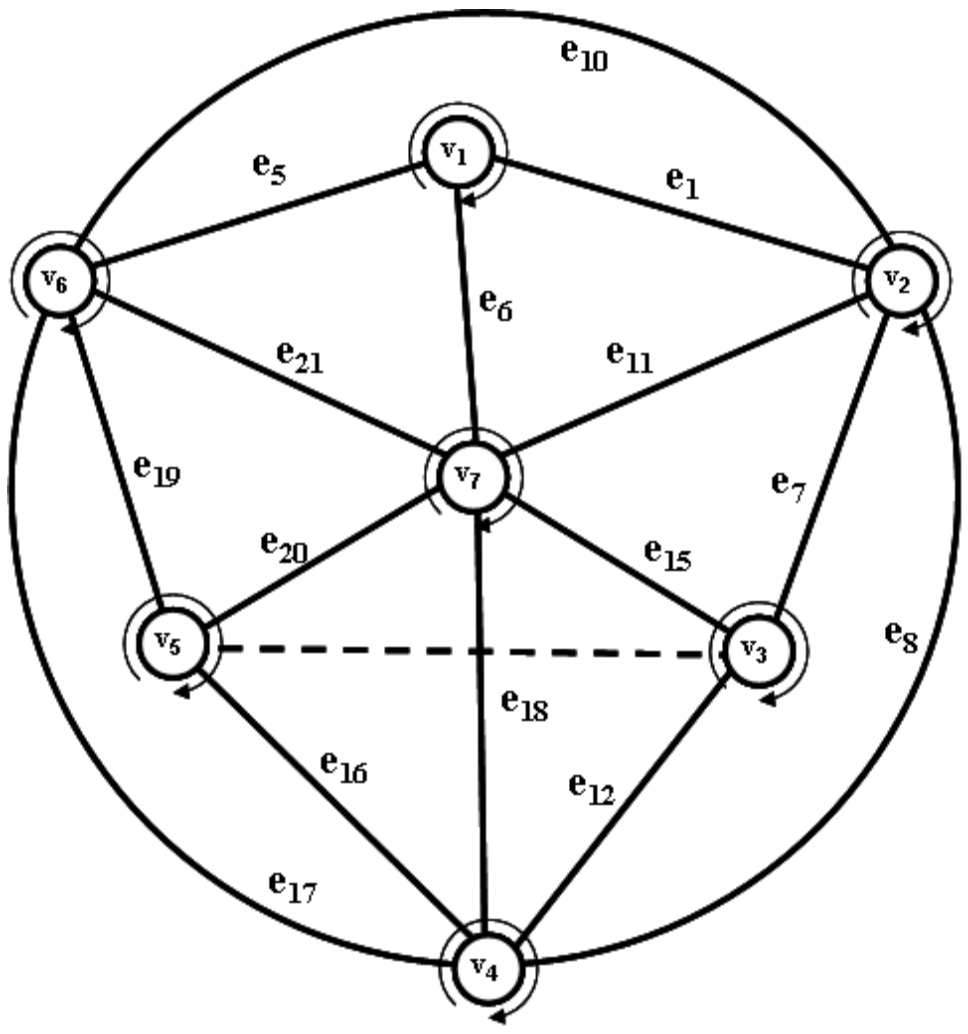

Рис. 9.25. Проведение соединения $e_{13}$.

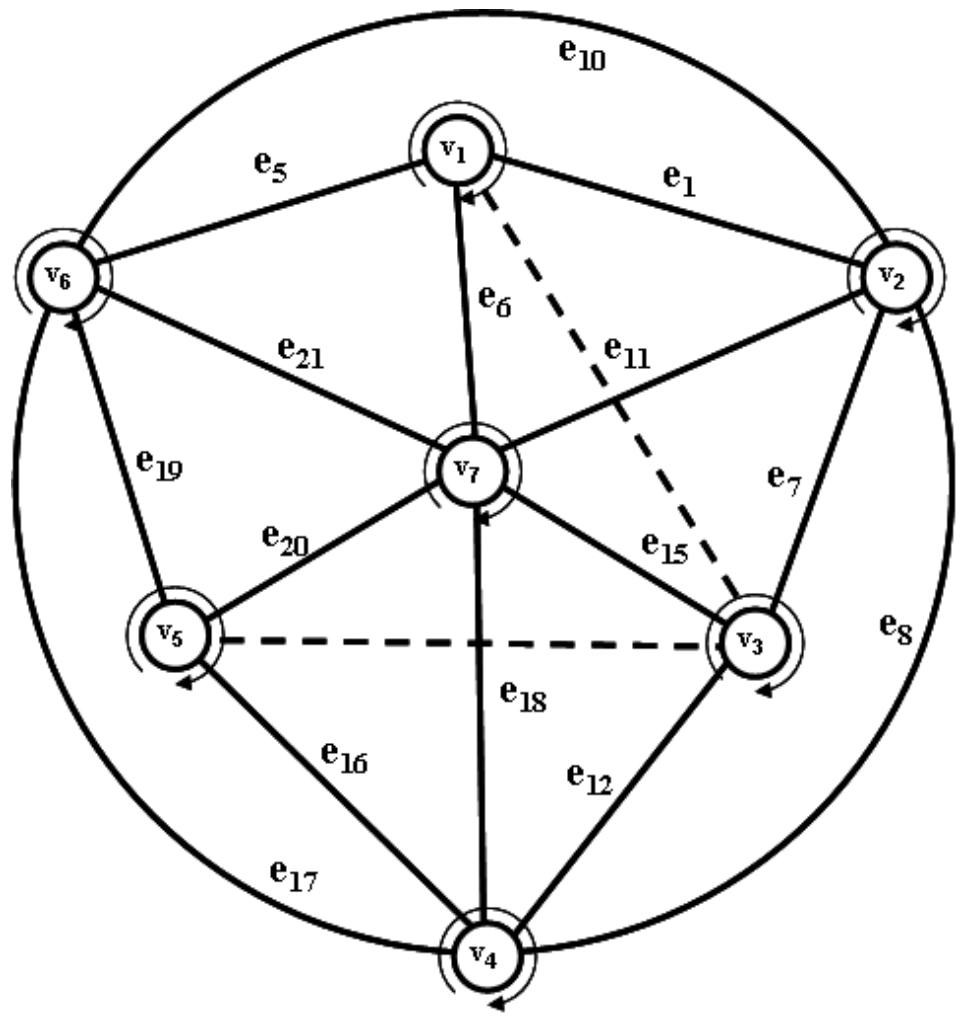

Рис. 9.25. Проведение соединения $e_2$.

Следует заметить, что топологические рисунки слоев, начиная со второго, должны строиться совместно с топологическим рисунком первого слоя, с учетом различения, так как проведение соединений осуществляется в топологическом пространстве определенным вращением вершин первого слоя (см. рис. 9.14).

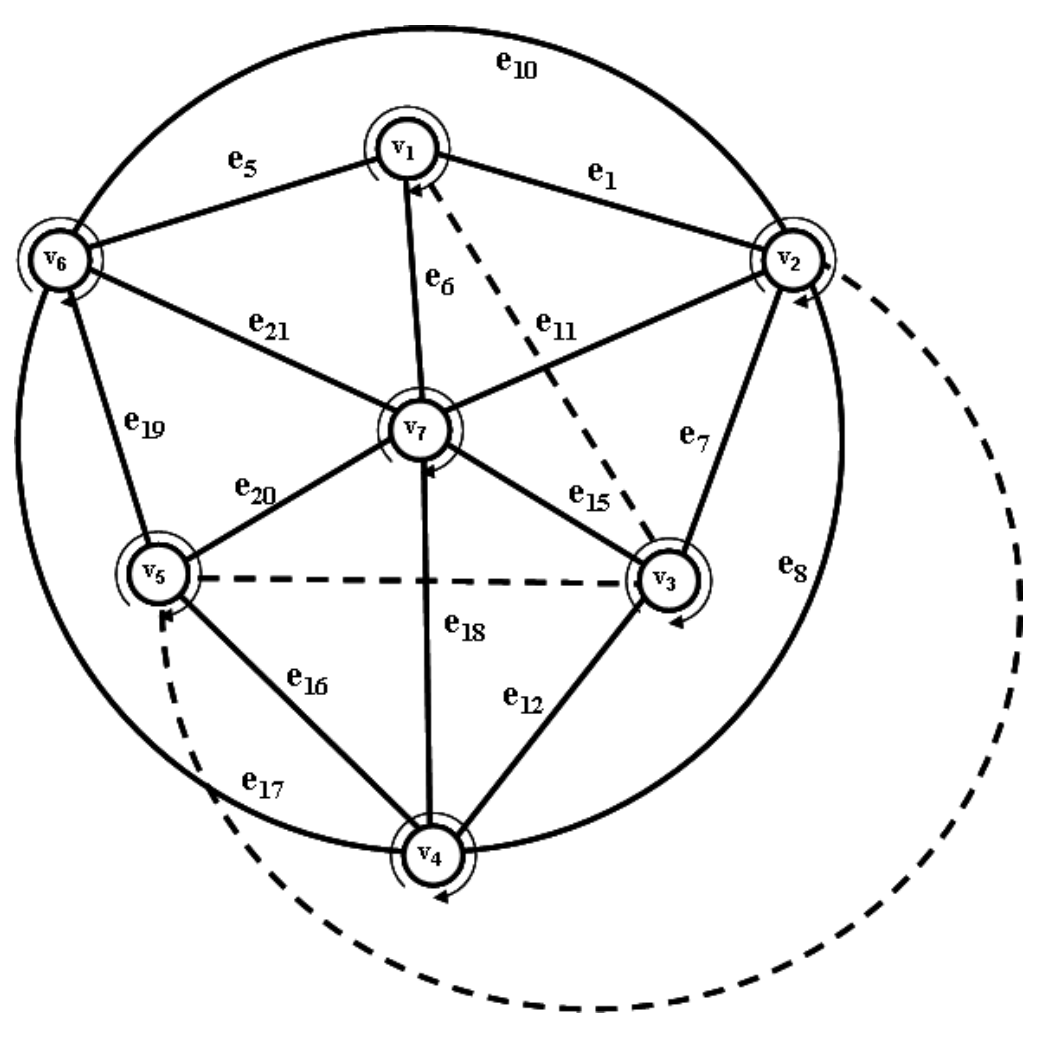

Рис. 9.27. Проведение соединения $e_9$.

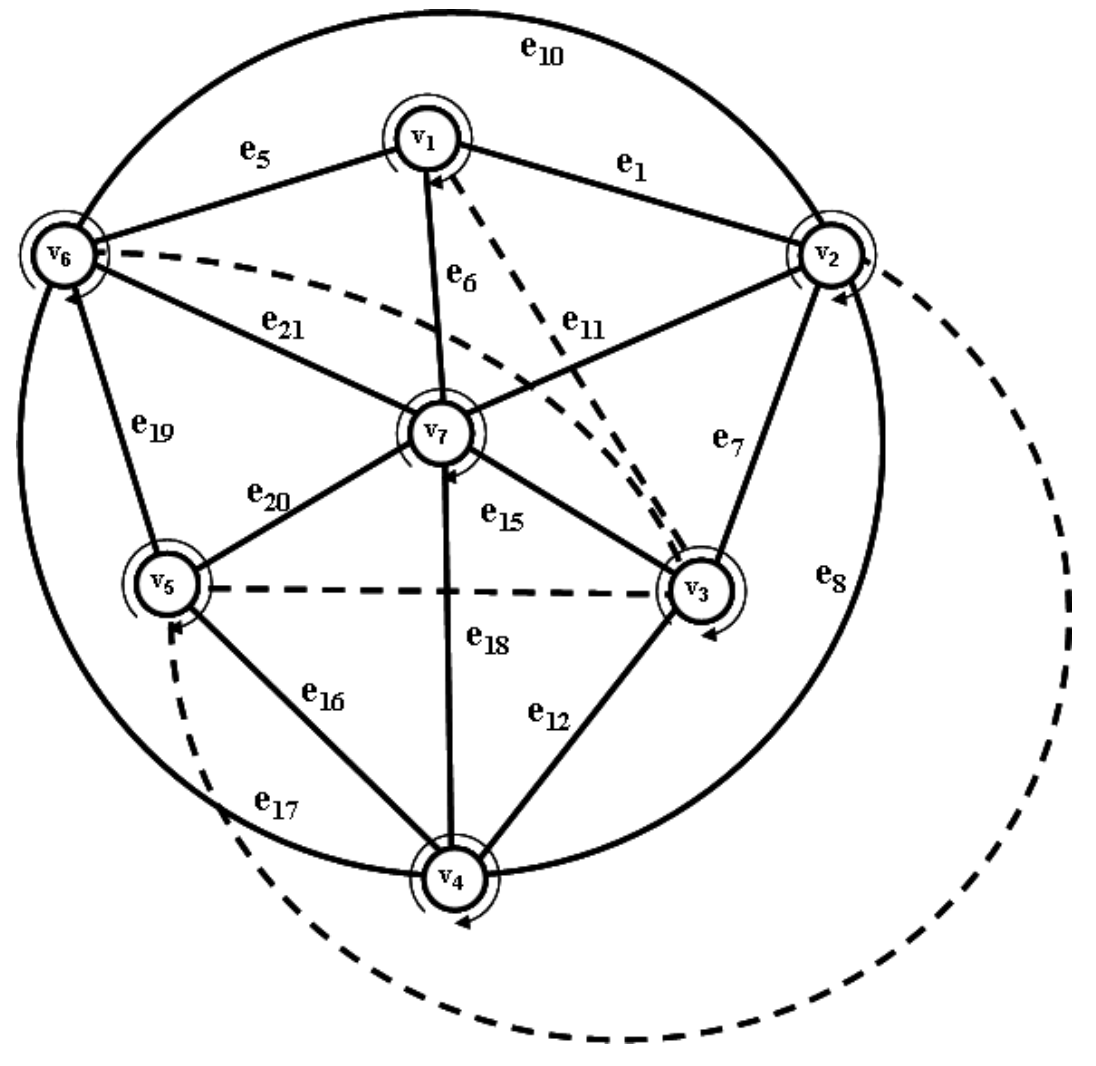

Рис. 9.28. Проведение соединения $e_{14}$.

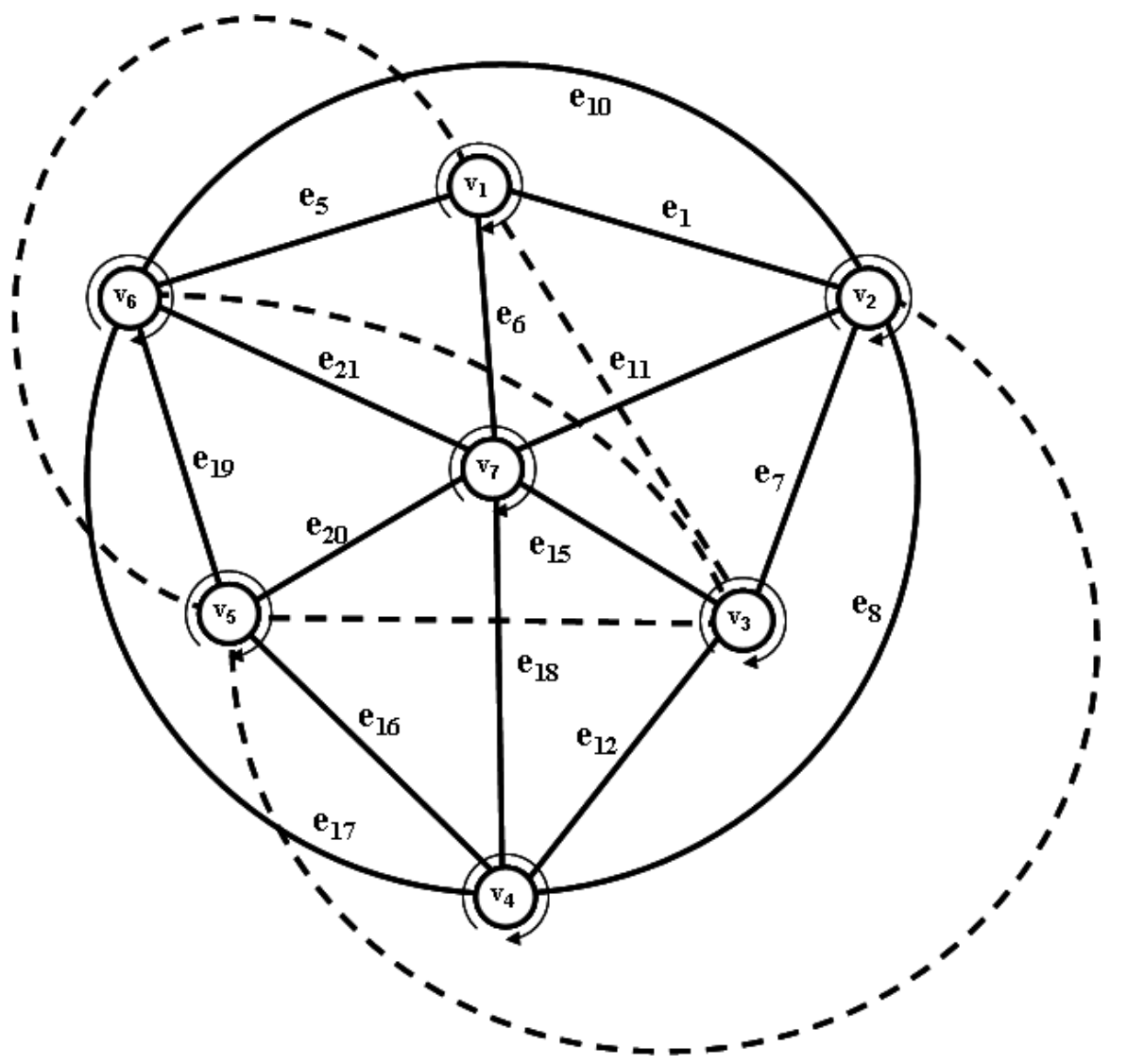

Рис. 9.29. Проведение соединения $e_4$.

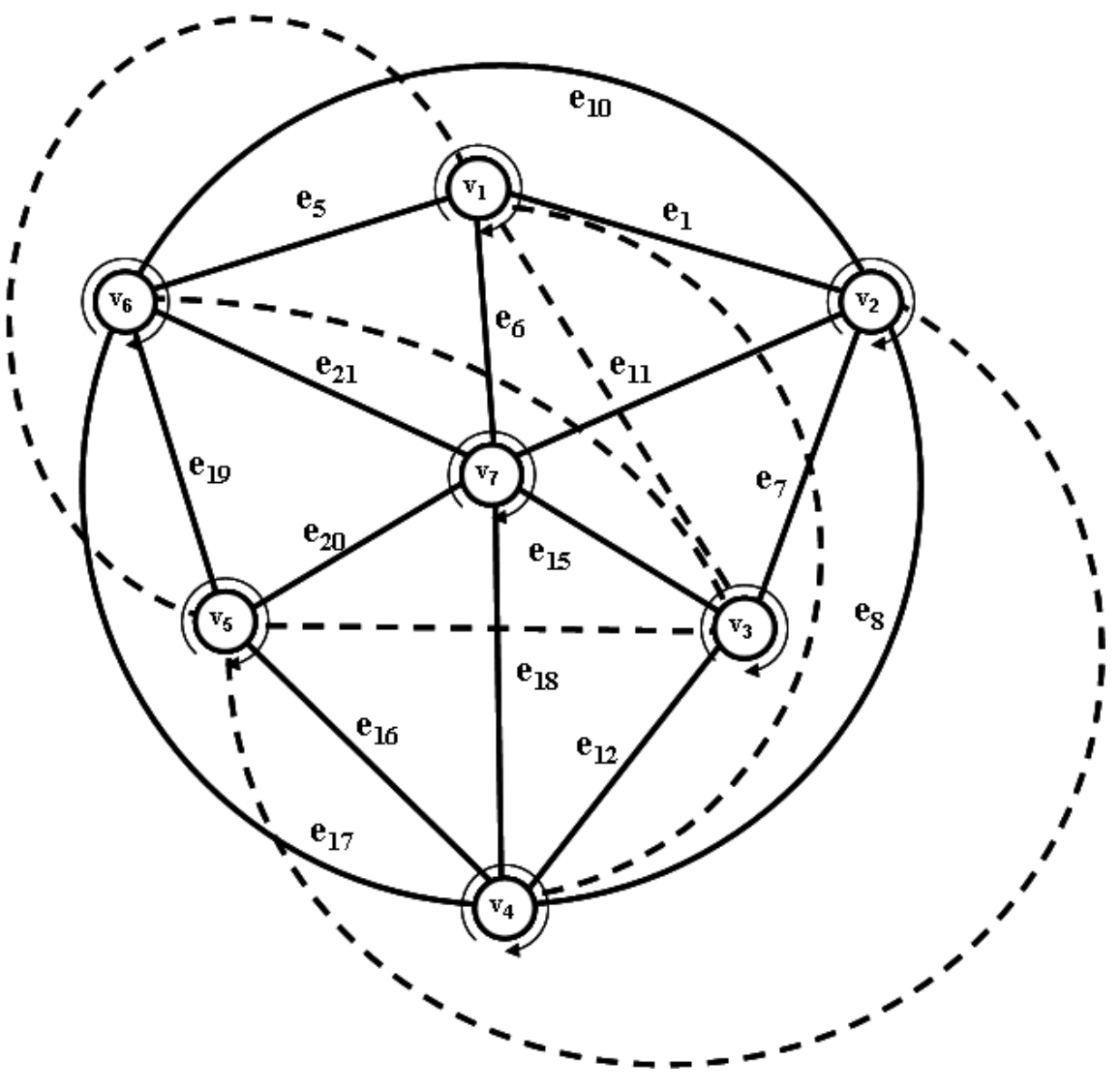

Рис. 9.30. Проведение соединения $e_3$.

Вне зависимости от вида графа, топологический рисунок определяется вращением вершин и индуцированной этим вращением системе циклов.

Например, вращение вершин описывающее топологический рисунок второго слоя графа $K_7$ представленного на рис. 9.33 имеет вид:

$\sigma(v_1) = \langle v_{13},v_8,v_9,v_7,v_{18},v_{14}\rangle$; $\sigma(v_2) = \langle v_6,v_{11},v_4,v_3,v_{16},v_{17}\rangle$;
$\sigma(v_3) = \langle v_{10},v_8,v_{12},v_{16},v_2,v_4\rangle$; $\sigma(v_4) = \langle v_7,v_{10},v_3,v_2,v_{11},v_5\rangle$;
$\sigma(v_5) = \langle v_6,v_{15},v_{18},v_7,v_4,v_{11}\rangle$; $\sigma(v_6) = \langle v_2,v_{17},v_{14},v_{15},v_5,v_{11}\rangle$;
$\sigma(v_7) = \langle v_1,v_9,v_{10},v_4,v_5,v_{18}\rangle$;
$\sigma(v_8) = \langle v_1,v_{12},v_3,v_9\rangle$; $\sigma(v_9) = \langle v_1,v_8,v_{10},v_7\rangle$; $\sigma(v_{10}) = \langle v_7,v_9,v_3,v_{49}\rangle$;
$\sigma(v_{11}) = \langle v_5,v_4,v_2,v_6\rangle$; $\sigma(v_{12}) = \langle v_{13},v_{16},v_3,v_8\rangle$; $\sigma(v_{13}) = \langle v_{14},v_{17},v_{12},v_1\rangle$;
$\sigma(v_{14}) = \langle v_6,v_{13},v_1,v_{15}\rangle$; $\sigma(v_{15}) = \langle v_6,v_{14},v_{18},v_5\rangle$; $\sigma(v_{16}) = \langle v_2,v_3,v_{12},v_{17}\rangle$;
$\sigma(v_{17}) = \langle v_6,v_2,v_{16},v_{13}\rangle$; $\sigma(v_{18}) = \langle v_{15},v_1,v_4,v_5\rangle$.

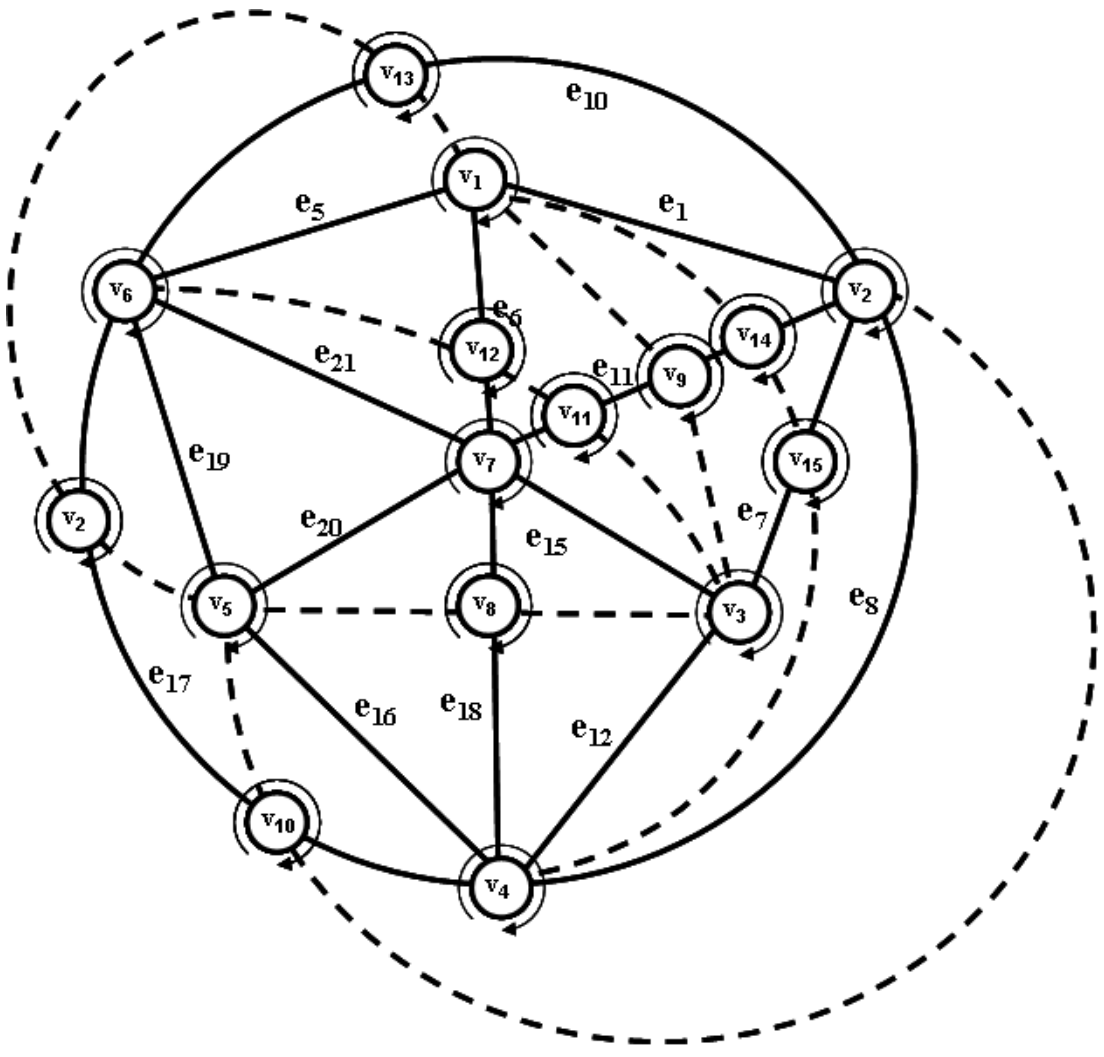

Рис. 9.31. Топологический рисунок графа минимальной толщины.

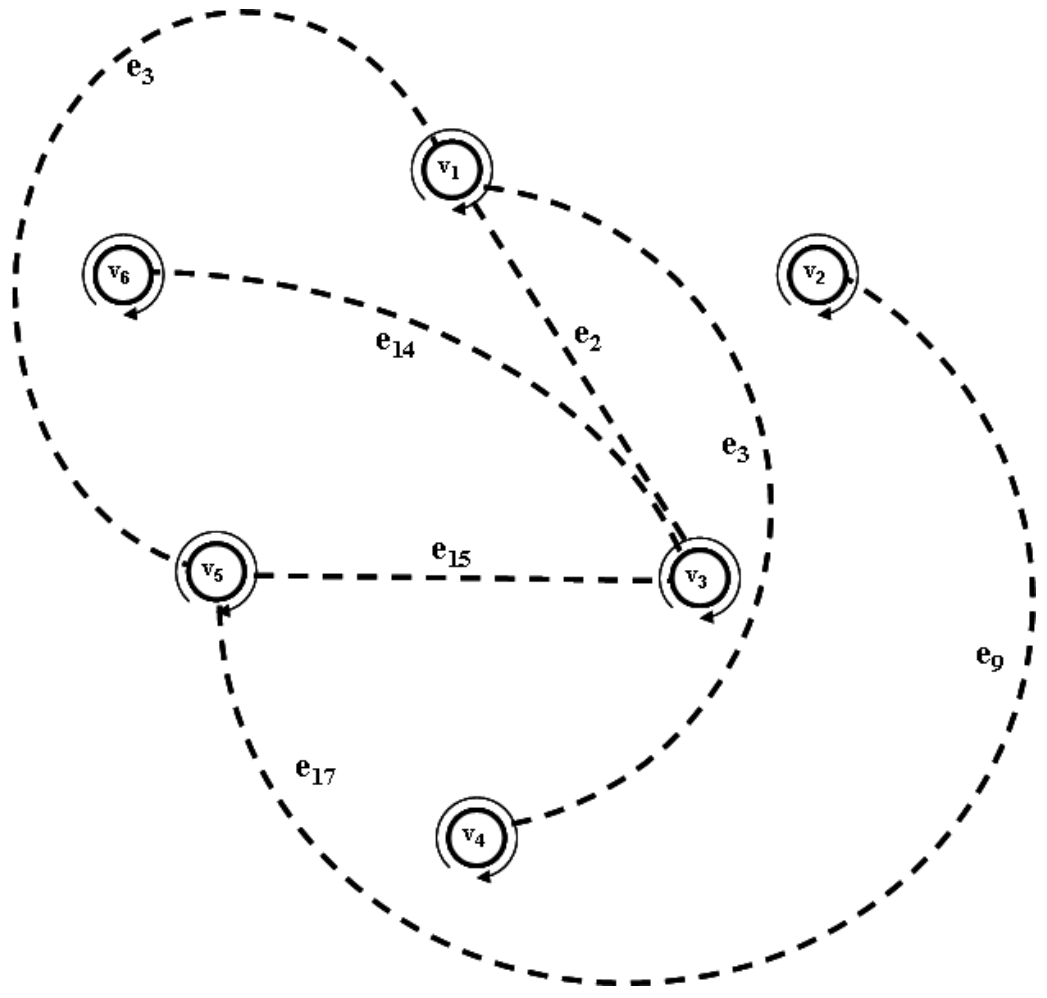

Рис. 9.32. Топологический рисунок второго слоя.

Множество индуцированных циклов данным вращением вершин:

$c_1 = \langle v_2,v_6,v_{17} \rangle = (v_2,v_6) + (v_6,v_{17}) + (v_{17},v_2)$;
$c_2 = \langle v_{17},v_6,v_{14},v_{13} \rangle = (v_{17},v_6) + (v_6,v_{14}) + (v_{14},v_{13}) + (v_{13},v_{17})$;
$c_3 = \langle v_2,v_{17},v_{16} \rangle = (v_2,v_{17}) + (v_{17},v_{16}) + (v_{16},v_2)$;
$c_4 = \langle v_6,v_{15},v_{14} \rangle = (v_6,v_{15}) + (v_{15},v_{14}) + (v_{14},v_6)$;
$c_5 = \langle v_{14},v_1,v_{13} \rangle = (v_{14},v_1) + (v_1,v_{13}) + (v_{13},v_{14})$;
$c_6 = \langle v_{17},v_{13},v_{12},v_{16} \rangle = (v_{17},v_{13}) + (v_{13},v_{12}) + (v_{12},v_{16}) + (v_{16},v_{17})$;
$c_7 = \langle v_{16},v_3,v_2 \rangle = (v_{16},v_3) + (v_3,v_2) + (v_2,v_{16})$;
$c_8 = \langle v_{15},v_{18},v_1,v_{14} \rangle = (v_{15},v_{18}) + (v_{18},v_1) + (v_1,v_{14}) + (v_{14},v_{15})$;
$c_9 = \langle v_{13},v_1,v_8,v_{12} \rangle = (v_{13},v_1) + (v_1,v_8) + (v_8,v_{12}) + (v_{12},v_{13})$;
$c_{10} = \langle v_{12},v_3,v_{16} \rangle = (v_{12},v_3) + (v_3,v_{16}) + (v_{16},v_{12})$;
$c_{11} = \langle v_6,v_5,v_{15} \rangle = (v_6,v_5) + (v_5,v_{15}) + (v_{15},v_6)$;
$c_{12} = \langle v_{15},v_5,v_{18} \rangle = (v_{15},v_5) + (v_5,v_{18}) + (v_{18},v_{15})$;
$c_{13} = \langle v_{18},v_7,v_1 \rangle = (v_{18},v_7) + (v_7,v_1) + (v_1,v_{18})$;
$c_{14} = \langle v_1,v_7,v_9 \rangle = (v_1,v_7) + (v_7,v_9) + (v_9,v_1)$;
$c_{15} = \langle v_1,v_9,v_8 \rangle = (v_1,v_9) + (v_9,v_8) + (v_8,v_1)$;
$c_{16} = \langle v_8,v_3,v_{12} \rangle = (v_8,v_3) + (v_3,v_{12}) + (v_{12},v_8)$;
$c_{17} = \langle v_9,v_{10},v_3,v_8 \rangle = (v_9,v_{10}) + (v_{10},v_3) + (v_3,v_8) + (v_8,v_9)$;
$c_{18} = \langle v_7,v_{10},v_9 \rangle = (v_7,v_{10}) + (v_{10},v_9) + (v_9,v_7)$;
$c_{19} = \langle v_5,v_7,v_{18} \rangle = (v_5,v_7) + (v_7,v_{18}) + (v_{18},v_5)$;
$c_{20} = \langle v_6,v_{11},v_5 \rangle = (v_6,v_{11}) + (v_{11},v_5) + (v_5,v_6)$;
$c_{21} = \langle v_5,v_4,v_7 \rangle = (v_5,v_4) + (v_4,v_7) + (v_7,v_5)$;
$c_{22} = \langle v_7,v_4,v_{10} \rangle = (v_7,v_4) + (v_4,v_{10}) + (v_{10},v_7)$;
$c_{23} = \langle v_{10},v_4,v_3 \rangle = (v_{10},v_4) + (v_4,v_3) + (v_3,v_{10})$;
$c_{24} = \langle v_5,v_{11},v_4 \rangle = (v_5,v_{11}) + (v_{11},v_4) + (v_4,v_5)$;
$c_{25} = \langle v_3,v_4,v_2 \rangle = (v_3,v_4) + (v_4,v_2) + (v_2,v_3)$;
$c_{26} = \langle v_{11},v_2,v_4 \rangle = (v_{11},v_2) + (v_2,v_4) + (v_4,v_{11})$;
$c_{27} = \langle v_{11},v_6,v_2 \rangle = (v_{11},v_6) + (v_6,v_2) + (v_2,v_{11})$.

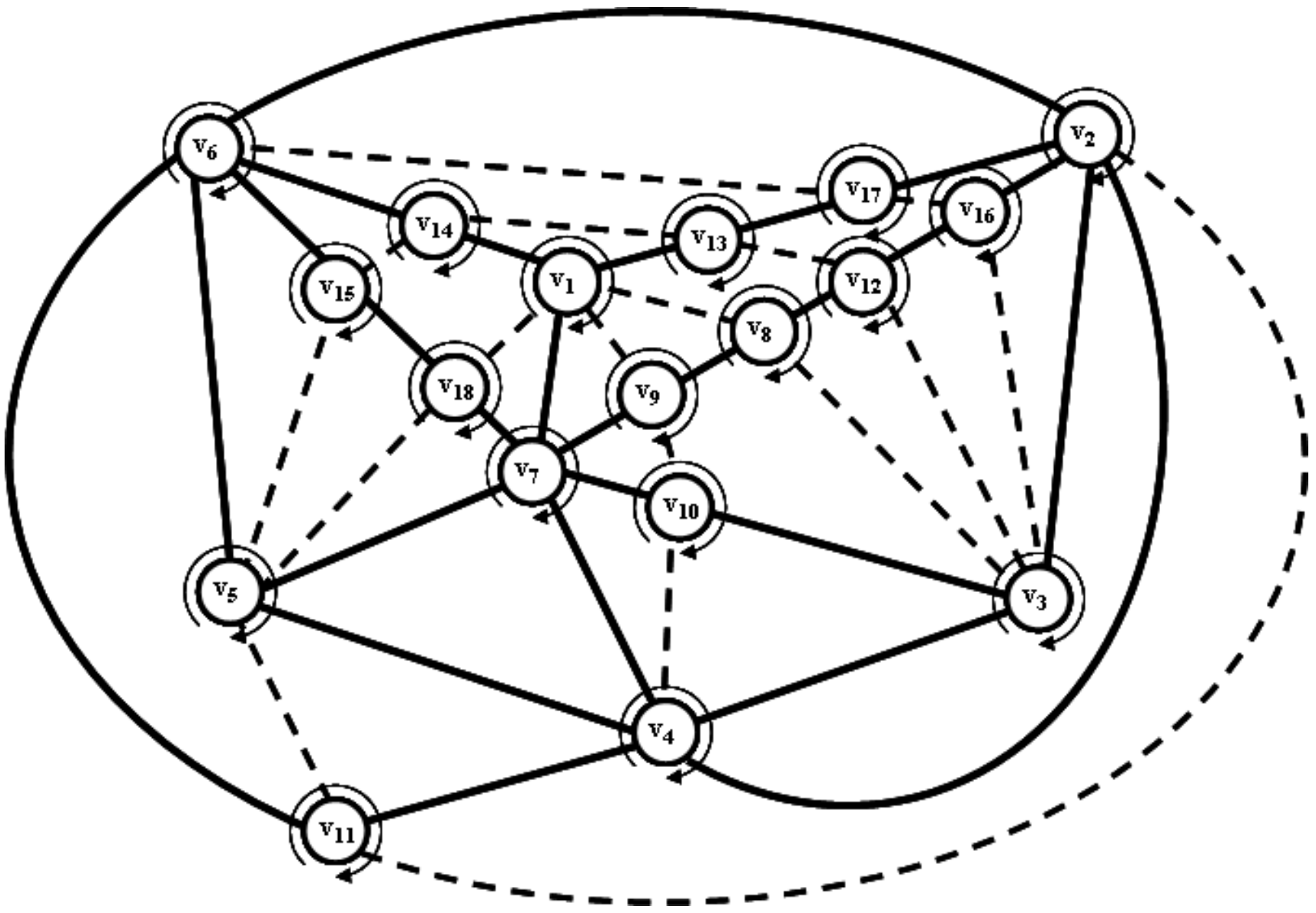

Рис. 9.33. Топологический рисунок графа минимальной толщины.

**Комментарии**

В данной главе рассмотрены методы построения топологического рисунка графа с минимальным числом пересечений и топологического рисунка графа минимальной толщины. Построения осуществляются на базе топологтческого рисунка максимально плоского суграфа. Вычисление маршрутов для проведения соединений производится алгоритмами фрагментарно-эволюционных структур.

# Глава 10. ГЕОМЕТРИЧЕСКИЕ МЕТОДЫ ПОСТРОЕНИЯ ТОПОЛОГИЧЕСКОГО РИСУНКА ГРАФА

## 10.1. Граф уровней и представление топологического рисунка графа

Построение топологического рисунка графа, основанное на построении системы независимых циклов и вращении вершин, является только первой частью визуализации топологического рисунка графа. Важно не только получить математическую модель топологического рисунка графа, но и получить его геометрическое представление. Геометрическое представление позволяет получить визуальное изображение и возможность выбирать оптимальное решение.

Пусть задан топологический рисунок несепарабельного графа $G_{12}$, представленный матрицей смежностей графа и системой независимых циклов с нулевым значением функционала Маклейна. Подмножество изометрических циклов для плоской части графа с учетом ориентации можно представить в векторном виде:

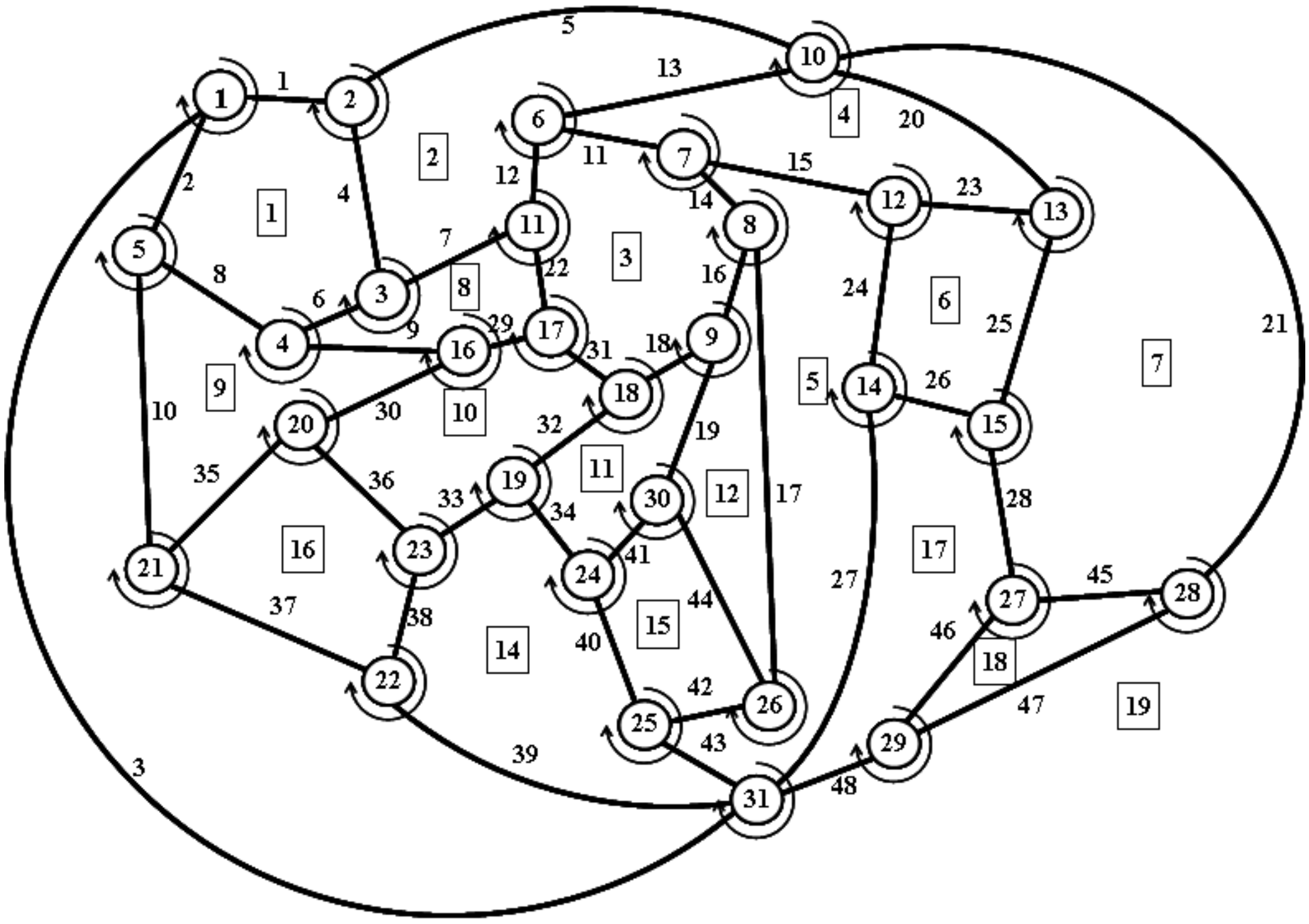

Рис. 10.1. Топологический рисунок графа для заданной системы циклов.

$c_1 = \langle v_2,v_3,v_4,v_5,v_1,v_2\rangle = (v_2,v_3) + (v_3,v_4) + (v_4,v_5) + (v_5,v_1) + (v_1,v_2)$ ;

$c_2 = \langle v_{10},v_6,v_{11},v_3,v_2,v_{10}\rangle = (v_{10},v_6) + (v_6,v_{11}) + (v_{11},v_3) + (v_3,v_2) + (v_2,v_{10})$;

$c_3 = \langle v_6,v_7,v_8,v_9,v_{18},v_{17},v_{11},v_6\rangle= (v_6,v_7) + (v_7,v_8) + (v_8,v_9) + (v_9,v_{18}) + (v_{18},v_{17}) + (v_{17},v_{11}) +$
$+ (v_{11},v_6)$;

$c_4 = \langle v_{10},v_{13},v_{12},v_7,v_6,v_{10}\rangle = (v_{10},v_{13}) + (v_{13},v_{12}) + (v_{12},v_7) + (v_7,v_6) + (v_6,v_{10})$;

$c_5 = \langle v_{31},v_{25},v_{26},v_8,v_7,v_{12},v_{14},v_{31}\rangle= (v_{31},v_{25}) + (v_{25},v_{26}) + (v_{26},v_8) + (v_8,v_7) + (v_7,v_{12}) +$
$+ (v_{12},v_{14}) + (v_{14},v_{31})$;

$c_6 = \langle v_{14},v_{12},v_{13},v_{15},v_{14}\rangle = (v_{14},v_{12}) + (v_{12},v_{13}) + (v_{13},v_{15}) + (v_{15},v_{14})$;

$c_7 = \langle v_{28},v_{27},v_{15},v_{13},v_{10},v_{28}\rangle = (v_{28},v_{27}) + (v_{27},v_{15}) + (v_{15},v_{13}) + (v_{13},v_{10}) + (v_{10},v_{28})$;

$c_8 = \langle v_3,v_{11},v_{17},v_{16},v_4,v_3\rangle = (v_3,v_{11}) + (v_{11},v_{17}) + (v_{17},v_{16}) + (v_{16},v_4) + (v_4,v_3)$;

$c_9 = \langle v_5,v_4,v_{16},v_{20},v_{21},v_5\rangle = (v_5,v_4) + (v_4,v_{16}) + (v_9,v_{20}) + (v_{20},v_{21}) + (v_{21},v_5)$;

$c_{10} = \langle v_{23},v_{20},v_{16},v_{17},v_{18},v_{19},v_{23}\rangle = (v_{23},v_{20}) + (v_{20},v_{16}) + (v_{16},v_{17}) + (v_{17},v_{18}) + (v_{18},v_{19}) + (v_{19},v_{23})$;

$c_{11} = \langle v_{24},v_{19},v_{18},v_9,v_{30},v_{24}\rangle = (v_{24},v_{19}) + (v_{19},v_{18}) + (v_{18},v_9) + (v_9,v_{30}) + (v_{30},v_{24})$;

$c_{12} = \langle v_{26},v_{30},v_9,v_8,v_{26}\rangle = (v_{26},v_{30}) + (v_{30},v_9) + (v_9,v_8) + (v_8,v_{26})$;

$c_{13} = \langle v_{31},v_1,v_5,v_{21},v_{22},v_{31}\rangle = (v_{31},v_1) + (v_1,v_5) + (v_5,v_{21}) + (v_{21},v_{22}) + (v_{22},v_{31})$;

$c_{14} = \langle v_{31},v_{22},v_{23},v_{19},v_{24},v_{25},v_{31}\rangle = (v_{31},v_{22}) + (v_{22},v_{23}) + (v_{23},v_{19}) + (v_{19},v_{24}) + (v_{24},v_{25}) + (v_{25},v_{31})$;

$c_{15} = \langle v_{25},v_{24},v_{30},v_{26},v_{25}\rangle = (v_{25},v_{24}) + (v_{24},v_{30}) + (v_{30},v_{26}) + (v_{26},v_{25})$;

$c_{16} = \langle v_{22},v_{21},v_{20},v_{23},v_{22}\rangle = (v_{22},v_{21}) + (v_{21},v_{20}) + (v_{20},v_{23}) + (v_{23},v_{22})$;

$c_{17} = \langle v_{31},v_{14},v_{15},v_{27},v_{29},v_{31}\rangle = (v_{31},v_{14}) + (v_{14},v_{15}) + (v_{15},v_{27}) + (v_{27},v_{29}) + (v_{29},v_{31})$;

$c_{18} = \langle v_{29},v_{27},v_{28},v_{29}\rangle = (v_{29},v_{27}) + (v_{27},v_{28}) + (v_{28},v_{29})$;

$c_{19} = \langle v_{31},v_{29},v_{28},v_{10},v_2,v_1,v_{31}\rangle = (v_{31},v_{29}) + (v_{29},v_{28}) + (v_{28},v_{10}) + (v_{10},v_2) + (v_2,v_1) + (v_1,v_{31})$.

Для наглядности систему циклов представим в виде рисунка (см. 10.1)

Построим граф уровней для данного рисунка графа относительно выбранного обода. В данном случае в качестве обода выбираем цикл $c_{19}$. Построение графа уровней предпологает применение для каждой вершины обода, алгоритма поиска в ширину. Смежные вершины с вершинами обода помечаются цифрой 2. Смежные вершины с вершинами уровня 2 помечатся цифрой 3. Процесс продолжается до полной разметки всех вершин графа.

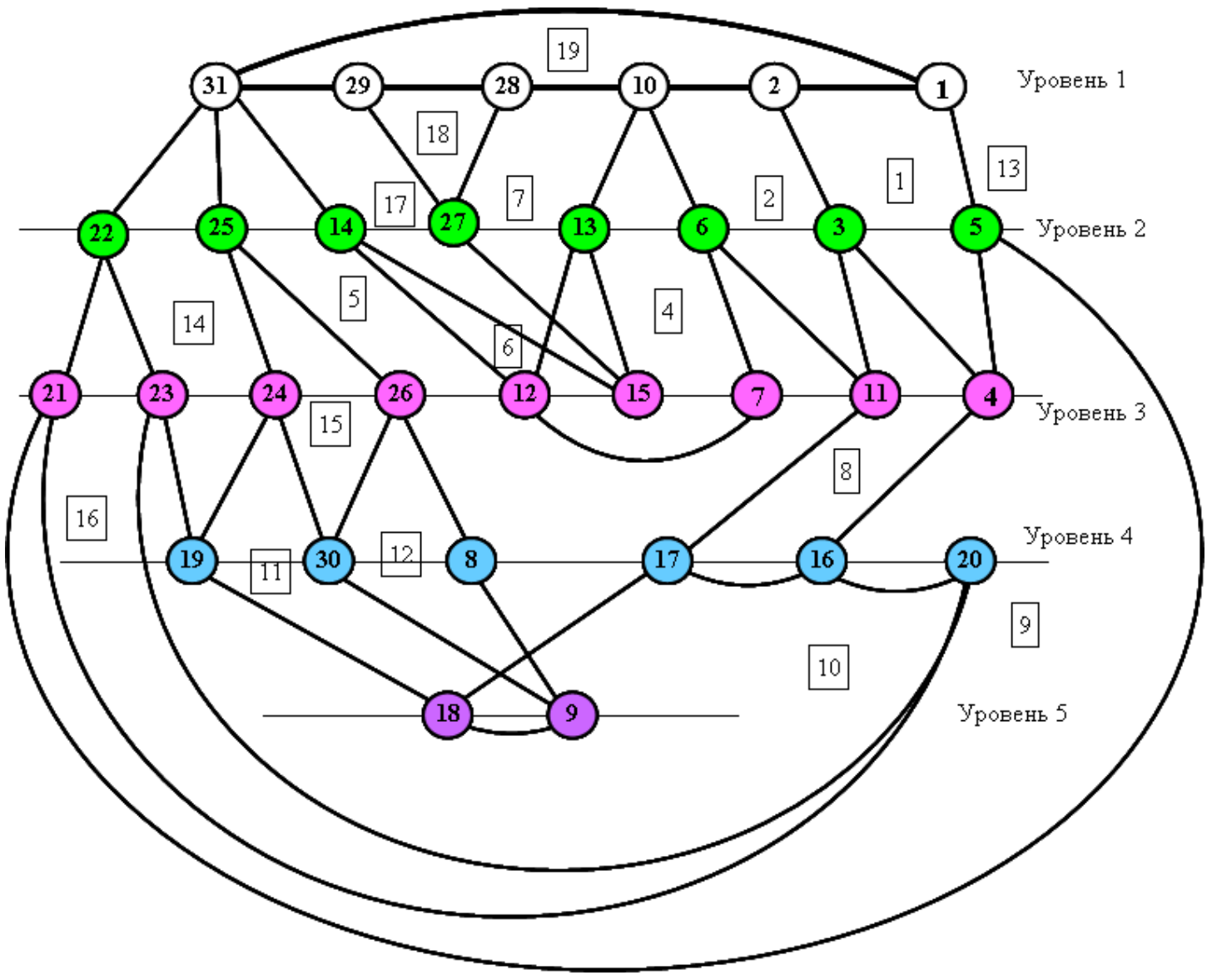

Рис. 10.2. Граф уровней.

Запишем в векторном виде циклы, приписав каждой вершине уровень, записанный в квадратных скобках. Расположим вершины в циклах по возрастанию уровней, тем самым создавая симметричное расположение вершин одинакового уровня.

$c_1 = (v_2[1],v_3[2])+(v_3[2],v_4[3])+(v_4[3],v_5[2])+(v_5[2],v_1[1])+(v_1[1],v_2[1]$;

$c_2 = (v_{10}[1],v_6[2])+(v_6[2],v_{11}[3])+(v_{11}[3],v_3[2])+(v_3[2],v_2[1])+(v_2[1],v_{10}[1])$;

$c_3 = (v_6[2],v_7[3])+(v_7[3],v_8[4])+(v_8[4],v_9[5])+(v_9[5],v_{18}[5])+(v_{18}[5],v_{17}[4])+(v_{17}[4],v_{11}[3])+$
$+(v_{11}[3],v_6[2])$;

$c_4 = (v_{10}[1],v_{13}[2])+(v_{13}[2],v_{12}[3])+(v_{12}[3],v_7[3])+(v_7[3],v_6[2])+(v_6[2],v_{10}[1])$;

$c_5 = (v_{31}[1],v_{25}[2])+(v_{25}[2],v_{26}[3])+(v_{26}[3],v_8[4])+(v_8[4],v_7[3])+(v_7[3],v_{12}[3])+(v_{12}[3],v_{14}[2])+$
$+(v_{14}[2],v_{31}[1])$;

$c_6 = (v_{14}[2],v_{15}[3])+(v_{15}[3],v_{13}[2])+(v_{13}[2],v_{12}[3])+(v_{12}[3],v_{14}[2])$;

$c_7 = (v_{28}[1],v_{27}[2])+(v_{27}[2],v_{15}[3])+(v_{15}[3],v_{13}[2])+(v_{13}[2],v_{10}[1])+(v_{10}[1],v_{28}[1])$;

$c_8 = (v_3[2],v_{11}[3])+(v_{11}[3],v_{17}[4])+(v_{17}[4],v_{16}[4])+(v_{16}[4],v_4[3])+(v_4[3],v_3[2])$;

$c_9 = (v_5[2],v_4[3])+(v_4[3],v_{16}[4])+(v_{16}[4],v_{20}[4])+(v_{20}[4],v_{21}[3])+(v_{21}[3],v_5[2])$;

$c_{10} = (v_{23}[3],v_{20}[4])+(v_{20}[4],v_{16}[4])+(v_{16}[4],v_{17}[4])+(v_{17}[4],v_{18}[5])+(v_{18}[5],v_{19}[4])+$
$+(v_{19}[4],v_{23}[3])$;
$c_{11} = (v_{24}[3],v_{19}[4])+(v_{19}[4],v_{18}[5])+(v_{18}[5],v_{9}[5])+(v_{9}[5],v_{30}[4])+(v_{30}[4],v_{24}[3])$;
$c_{12} = (v_{26}[3],v_{30}[4])+(v_{30}[4],v_{9}[5])+(v_{9}[5],v_{8}[4])+(v_{8}[4],v_{26}[3])$;
$c_{13} = (v_{31}[1],v_{1}[1])+(v_{1}[1],v_{5}[2])+(v_{5}[2],v_{21}[3])+(v_{21}[3],v_{22}[2])+(v_{22}[2],v_{31}[1])$;
$c_{14} = (v_{31}[1],v_{22}[2])+(v_{22}[2],v_{23}[3])+(v_{23}[3],v_{19}[4])+(v_{19}[4],v_{24}[3])+(v_{24}[3],v_{25}[2])+$
$+(v_{25}[2],v_{31}[1])$;
$c_{15} = (v_{25}[2],v_{24}[3])+(v_{24}[3],v_{30}[4])+(v_{30}[4],v_{26}[3])+(v_{26}[3],v_{25}[2])$;
$c_{16} = (v_{22}[2],v_{21}[3])+(v_{21}[3],v_{20}[4])+(v_{20}[4],v_{23}[3])+(v_{23}[3],v_{22}[2])$;
$c_{17} = (v_{31}[1],v_{14}[2])+(v_{14}[2],v_{15}[3])+(v_{15}[3],v_{27}[2])+(v_{27}[2],v_{29}[1])+(v_{29}[1],v_{31}[1])$;
$c_{18} = (v_{29}[1],v_{27}[2])+(v_{27}[2],v_{28}[1])+(v_{28}[1],v_{29}[1])$;
$c_{19} = (v_{31}[1],v_{29}[1])+(v_{29}[1],v_{28}[1])+(v_{28}[1],v_{10}[1])+(v_{10}[1],v_{2}[1])+(v_{2}[1],v_{1}[1])+(v_{1}[1],v_{31}[1])$.

В качестве обода, выделен цикл $c_{19}$. Содержание цикла позволяет построить цепочку последовательного расположения вершин в уровне 1, при заданном репере вершины $v_{31}$. Построение достигается путем простого перечисления расположения вершин в цикле. Последовательность расположения вершин по уровням будем обозначать латинской буквой U и представлять в виде номеров вершин.

$u_1$: 31 29 28 10 02 01
$u_2$:
$u_3$:
$u_4$:
$u_5$:

Будем рассматривать построение последовательности вершин в уровне 2. Построение будем производить в два этапа. На вервом этапе рассматриваются только пары вершин предыдущего уровня. На втором этапе подключаются вершины следущего уровня, порожденные одной вершиной предыдущего уровня.ю

С этой целью рассмотрим цикл $c_{17}$. Вершины $v_{14}$,$v_{27}$ уровня 2 (красный цвет) в цикле $c_{17}$, симметрично расположены между вершинами $v_{31}$ и $v_{29}$ уровня 1 (фиолетовый цвет). Рассмотрение следует начинать с краев цикла.

$c_{17} = (\mathbf{v_{31}[1]},\mathbf{v_{14}[2]})+(v_{14}[2],v_{15}[3])+(v_{15}[3],v_{27}[2])+(\mathbf{v_{27}[2]},\mathbf{v_{29}[1]})+(v_{29}[1],v_{31}[1])$;

Тогда расположение вершин в последовательности U имеет вид:

$u_1$: 31 29 28 10 02 01
$u_2$: 14 27
$u_3$:
$u_4$:
$u_5$:

Рассмотрим следующее расположение вершин $v_{29}$,$v_{28}$ уровня 1 в цикле $c_{18}$. Между этими вершинами находится вершина $v_{27}$. Вершина $v_{27}$ уже включена в последовательность U.

$c_{18} = (\mathbf{v_{29}[1]},\mathbf{v_{27}[2]})+(v_{27}[2],\mathbf{v_{28}[1]})+(v_{28}[1],v_{29}[1])$;

Рассмотрим следующее расположение вершин $v_{28}$,$v_{10}$ уровня 1 в цикле $c_7$. Между вершинами $v_{28}$,$v_{10}$ симметрично находятся вершины $v_{27}$,$v_{13}$ уровня 2.

$c_7 = (\mathbf{v_{28}[1]},\mathbf{v_{27}[2]})+(v_{27}[2],v_{15}[3])+(v_{15}[3],v_{13}[2])+(\mathbf{v_{13}[2]},\mathbf{v_{10}[1]})+(v_{10}[1],v_{28}[1])$;

Тогда расположение вершин в последовательности расположения вершин по уровням U имеет вид:

$u_1$: 31 29 28 10 02 01
$u_2$: 14 27 13
$u_3$:
$u_4$:
$u_5$:

Расположение вершин $v_{10}$,$v_{10}$ уровня 1 в цикле $c_4$, порождает появление пары вершин $v_{13}$,$v_6$ уровня 2.

$c_4 = (v_{10}[1],\mathbf{v_{13}[2]})+(v_{13}[2],v_{12}[3])+(v_{12}[3],v_7[3])+(v_7[3],v_6[2])+(\mathbf{v_6[2]},v_{10}[1]);$

$u_1$: 31 29 28 10 02 01
$u_2$: 14 27 13 06
$u_3$:
$u_4$:
$u_5$:

Расположение вершин $v_{10}$,$v_2$ уровня 1 в цикле $c_2$, порождает появление пары вершин $v_6$,$v_3$ уровня 2.

$c_2 = (\mathbf{v_{10}[1]},\mathbf{v_6[2]})+(v_6[2],v_{11}[3])+(v_{11}[3],v_3[2])+(\mathbf{v_3[2]},\mathbf{v_2[1]})+(v_2[1],v_{10}[1]);$

Тогда расположение вершин в последовательности расположения вершин по уровням U имеет вид:

$u_1$: 31 29 28 10 02 01
$u_2$: 14 27 13 06 03
$u_3$:
$u_4$:
$u_5$:

Расположение вершин $v_2$,$v_1$ уровня 1 в цикле $c_1$, порождает появление пары вершин $v_3$,$v_5$ уровня 2.

$c_1 = (\mathbf{v_2[1]},\mathbf{v_3[2]})+(v_3[2],v_4[3])+(v_4[3],v_5[2])+(\mathbf{v_5[2]},\mathbf{v_1[1]})+(v_1[1],v_2[1];$

Тогда расположение вершин в последовательности расположения вершин по уровням Q имеет вид:

$u_1$: 31 29 28 10 02 01
$u_2$: 14 27 13 06 03 05
$u_3$:
$u_4$:
$u_5$:

Расположение вершин $v_1$,$v_{31}$ уровня 1 в цикле $c_{13}$, порождает появление пары вершин $v_{25}$,$v_{14}$ уровня 2.

$c_{13} = (v_{31}[1],v_1[1])+(\mathbf{v_1[1]},\mathbf{v_5[2]})+(v_5[2],v_{21}[3])+(v_{21}[3],v_{22}[2])+(\mathbf{v_{22}[2]},\mathbf{v_{31}[1]});$

Тогда расположение вершин в последовательности расположения вершин по уровням U имеет вид:

$u_1$: 31 29 28 10 02 01
$u_2$: 14 27 13 06 03 05 22
$u_3$:
$u_4$:
$u_5$:

Закончим процесс рассмотрения пар вершин в предыдущем уровне. Перейдем к рассмотрению второго этапа. Рассмотрим следующее расположение вершин $v_{25}$,$v_{14}$ уровня 2 в цикле $c_5$, между симметричным равположением вершины $v_{31}$ уровня 1.

$c_5$ = (**$v_{31}$[1]**,**$v_{25}$[2]**)+($v_{25}$[2],$v_{26}$[3])+($v_{26}$[3],$v_8$[4])+($v_8$[4],$v_7$[3])+($v_7$[3],$v_{12}$[3])+($v_{12}$[3],$v_{14}$[2])+
+(**$v_{14}$[2]**,**$v_{31}$[1]**);

Тогда расположение вершин в последовательности расположения вершин по уровням Q имеет вид:

$u_1$: 31 29 28 10 02 01
$u_2$: 25 14 27 13 06 03 05 22
$u_3$:
$u_4$:
$u_5$:

Завершает построение последовательности, размещение вершин $v_{22}$,$v_{25}$ во 2-ом уровне цикл $c_{14}$ относительно симметричного расположения вершины $v_{31}$ уровня 1.

$c_{14}$ = (**$v_{31}$[1]**,**$v_{22}$[2]**)+($v_{22}$[2],$v_{23}$[3])+($v_{23}$[3],$v_{19}$[4])+($v_{19}$[4],$v_{24}$[3])+($v_{24}$[3],$v_{25}$[2])+
+ (**$v_{25}$[2]**,**$v_{31}$[1]**);

Будем строить последовательность вершин уровня 3. Рассмотрим первый этап построения, состоящий из выделения пар верншин предыдущего уровня.

Ищем пару $v_{25}$,$v_{14}$ уровня 2. Расположение вершин $v_{25}$,$v_{14}$ уровня 2 в цикле $c_5$, порождает появление пары вершин $v_{26}$,$v_{12}$ уровня 3.

$c_5$ = ($v_{31}$[1],$v_{25}$[2])+(**$v_{25}$[2]**,**$v_{26}$[3]**)+($v_{26}$[3],$v_8$[4])+($v_8$[4],$v_7$[3])+($v_7$[3],$v_{12}$[3])+(**$v_{12}$[3]**,**$v_{14}$[2]**)+
+($v_{14}$[2],$v_{31}$[1]);

Тогда расположение вершин в последовательности расположения вершин по уровням Q имеет вид:

$u_1$: 31 29 28 10 02 01
$u_2$: 25 14 27 13 06 03 05 22
$u_3$: 26 12
$u_4$:
$u_5$:

Расположение вершин $v_{14}$,$v_{27}$ уровня 2 в цикле $c_{17}$, порождает появление вершины $v_{15}$ уровня 3.

$c_{17}$ = ($v_{31}$[1],$v_{14}$[2])+(**$v_{14}$[2]**,**$v_{15}$[3]**)+(**$v_{15}$[3]**,**$v_{27}$[2]**)+($v_{27}$[2],$v_{29}$[1])+($v_{29}$[1],$v_{31}$[1]);

Тогда расположение вершин в последовательности расположения вершин по уровням Q имеет вид:

$u_1$: 31 29 28 10 02 01
$u_2$: 25 14 27 13 06 03 05 22

$u_3$: 26 12 15
$u_4$:
$u_5$:

Расположение вершин $v_{15}$,$v_{12}$ уровня 3 в цикле $c_7$, порождается вершиной $v_{14}$ уровня 2.

$c_6 = (v_{14}[2],\mathbf{v_{15}[3]})+(v_{15}[3],v_{13}[2])+(v_{13}[2],v_{12}[3])+(\mathbf{v_{12}[3]},v_{14}[2]);$

Тогда расположение вершин в последовательности U, остается без изменений:

$u_1$: 31 29 28 10 02 01
$u_2$: 25 14 27 13 06 03 05 22
$u_3$: 26 12 15 12
$u_4$:
$u_5$:

Расположение вершин $v_{27}$,$v_{13}$ уровня 2 в цикле $c_7$, порождает появление тольцо вершины $v_{15}$ уровня 3.

$c_7 = (v_{28}[1],v_{27}[2])+(\mathbf{v_{27}[2]},\mathbf{v_{15}[3]})+(\mathbf{v_{15}[3]},\mathbf{v_{13}[2]})+(v_{13}[2],v_{10}[1])+(v_{10}[1],v_{28}[1]);$

Тогда расположение вершин в последовательности U, остается без изменений:

$u_1$: 31 29 28 10 02 01
$u_2$: 25 14 27 13 06 03 05 22
$u_3$: 26 12 15 12
$u_4$:
$u_5$:

Расположение вершин $v_{13}$,$v_6$ уровня 2 в цикле $c_4$, порождает появление вершин $v_{12}$,$v_7$ уровня 3.

$c_4 = (v_{10}[1],v_{13}[2])+(\mathbf{v_{13}[2]},\mathbf{v_{12}[3]})+(v_{12}[3],v_7[3])+(\mathbf{v_7[3]},\mathbf{v_6[2]})+(v_6[2],v_{10}[1]);$

Тогда расположение вершин в последовательности расположения вершин по уровням Q имеет вид:

$u_1$: 31 29 28 10 02 01
$u_2$: 25 14 27 13 06 03 05 22
$u_3$: 26 12 15 12 07
$u_4$:
$u_5$:

Симметричное расположение вершин $v_6$,$v_3$,$v_5$,$v_{22}$ в циклах $c_3$,$c_8$,$c_9$,$c_{16}$ только подтверждает корректность расположения вершин в последовательности U.

$c_3 = (\mathbf{v_6[2]},\mathbf{v_7[3]})+(v_7[3],v_8[4])+(v_8[4],v_9[5])+(v_9[5],v_{18}[5])+(v_{18}[5],v_{17}[4])+(v_{17}[4],v_{11}[3])+$
$+(\mathbf{v_{11}[3]},\mathbf{v_6[2]});$

Тогда расположение вершин в последовательности расположения вершин по уровням Q имеет вид:

$u_1$: 31 29 28 10 02 01
$u_2$: 25 14 27 13 06 03 05 22
$u_3$: 26 12 15 12 07 11
$u_4$:
$u_5$:

Расположение вершин $v_6,v_3$ уровня 2 в цикле $c_2$, порождает появление вершины $v_{11}$ уровня 3.

$c_2 = (v_{10}[1],v_6[2])+(\mathbf{v_6[2]},\mathbf{v_{11}[3]})+(v_{11}[3],\mathbf{v_3[2]})+(v_3[2],v_2[1])+(v_2[1],v_{10}[1]);$

Тогда расположение вершин в последовательности расположения вершин по уровням Q имеет вид:

| | | | | | | | | |
|---|---|---|---|---|---|---|---|---|
| $u_1$: | 31 | 29 | 28 | 10 | 02 | 01 | | |
| $u_2$: | 25 | 14 | 27 | 13 | 06 | 03 | 05 | 22 |
| $u_3$: | 26 | 12 | 15 | 12 | 07 | 11 | | |
| $u_4$: | | | | | | | | |
| $u_5$: | | | | | | | | |

Далее следует:

$c_8 = (\mathbf{v_3[2]},\mathbf{v_{11}[3]})+(v_{11}[3],v_{17}[4])+(v_{17}[4],v_{16}[4])+(v_{16}[4],v_4[3])+(\mathbf{v_4[3]},\mathbf{v_3[2]});$

$c_9 = (\mathbf{v_5[2]},\mathbf{v_4[3]})+(v_4[3],v_{16}[4])+(v_{16}[4],v_{20}[4])+(v_{20}[4],v_{21}[3])+(\mathbf{v_{21}[3]},\mathbf{v_5[2]});$

Тогда расположение вершин в последовательности расположения вершин по уровням Q имеет вид:

| | | | | | | | | |
|---|---|---|---|---|---|---|---|---|
| $u_1$: | 31 | 29 | 28 | 10 | 02 | 01 | | |
| $u_2$: | 25 | 14 | 27 | 13 | 06 | 03 | 05 | 22 |
| $u_3$: | 26 | 12 | 15 | 12 | 07 | 11 | 04 | 21 |
| $u_4$: | | | | | | | | |
| $u_5$: | | | | | | | | |

Расположение вершин $v_3,v_5$ уровня 2 в цикле $c_1$, порождает появление вершины $v_4$ уровня 3.

$c_1 = (v_2[1],v_3[2])+(\mathbf{v_3[2]},\mathbf{v_4[3]})+(\mathbf{v_4[3]},\mathbf{v_5[2]})+(v_5[2],v_1[1])+(v_1[1],v_2[1];$

Тогда расположение вершин в последовательности расположения вершин по уровням Q имеет вид:

| | | | | | | | | |
|---|---|---|---|---|---|---|---|---|
| $u_1$: | 31 | 29 | 28 | 10 | 02 | 01 | | |
| $u_2$: | 25 | 14 | 27 | 13 | 06 | 03 | 05 | 22 |
| $u_3$: | 26 | 12 | 15 | 12 | 07 | 11 | 04 | 21 |
| $u_4$: | | | | | | | | |
| $u_5$: | | | | | | | | |

Расположение вершин $v_5,v_{22}$ уровня 2 в цикле $c_{13}$, порождает появление вершины $v_{21}$ уровня 3.

$c_{13} = (v_{31}[1],v_1[1])+(v_1[1],v_5[2])+(\mathbf{v_5[2]},\mathbf{v_{21}[3]})+(\mathbf{v_{21}[3]},\mathbf{v_{22}[2]})+(v_{22}[2],v_{31}[1]);$

Тогда расположение вершин в последовательности расположения вершин по уровням Q имеет вид:

| | | | | | | | | |
|---|---|---|---|---|---|---|---|---|
| $u_1$: | 31 | 29 | 28 | 10 | 02 | 01 | | |
| $u_2$: | 25 | 14 | 27 | 13 | 06 | 03 | 05 | 22 |
| $u_3$: | 26 | 12 | 15 | 12 | 07 | 11 | 04 | 21 |
| $u_4$: | | | | | | | | |
| $u_5$: | | | | | | | | |

Расположение вершин $v_{22}$,$v_{25}$ уровня 2 в цикле $c_{14}$, порождает появление вершин $v_{21}$,$v_{23}$,$v_{24}$ уровня 3.

$c_{16}$ = (**$v_{22}$[2]**,**$v_{21}$[3]**)+($v_{21}$[3],$v_{20}$[4])+($v_{20}$[4],$v_{23}$[3])+(**$v_{23}$[3]**,**$v_{22}$[2]**);
$c_{14}$ = ($v_{31}$[1],$v_{22}$[2])+(**$v_{22}$[2]**,**$v_{23}$[3]**)+($v_{23}$[3],$v_{19}$[4])+($v_{19}$[4],$v_{24}$[3])+(**$v_{24}$[3]**,**$v_{25}$[2]**)+
+($v_{25}$[2],$v_{31}$[1]);

Тогда расположение вершин в последовательности вершин по уровням U имеет вид:

$u_1$: 31 29 28 10 02 01
$u_2$: 25 14 27 13 06 03 05 22
$u_3$: 26 12 15 12 07 11 04 21 23 24
$u_4$:
$u_5$:

Построение последовательности вершин уровня 3, завершено. Перейдем к построению второго этапа.

Симметричное расположение вершины $v_{25}$ уровня 2 в цикле $c_{15}$, порождает появление вершин $v_{24}$,$v_{26}$ уровня 3.

$c_{15}$ = (**$v_{25}$[2]**,**$v_{24}$[3]**)+($v_{24}$[3],$v_{30}$[4])+($v_{30}$[4],$v_{26}$[3])+(**$v_{26}$[3]**,**$v_{25}$[2]**);

Тогда расположение вершин в последовательности расположения вершин по уровням Q имеет вид и циклически замыкается:

$u_1$: 31 29 28 10 02 01
$u_2$: 25 14 27 13 06 03 05 22
$u_3$: 24 26 12 15 12 07 11 04 21 23
$u_4$:
$u_5$:

Рассматривая последовательность вершин уровня 3 можно заметить, что существуют две вершины $v_{12}$. Для устранения пересечения ребер, вершины должны быть соединены в последовательности уровня 3, но должны быть соединены с разными вершинами уровня 2.

Будем рассматривать последовательность расположения вершин в уровне 4. Первый этап.

Расположение вершин $v_{24}$,$v_{26}$ уровня 3 в цикле $c_{15}$, порождает появление вершины $v_{30}$ уровня 4.

$c_{14}$ = ($v_{31}$[1],$v_{22}$[2])+(**$v_{22}$[2]**,**$v_{23}$[3]**)+($v_{23}$[3],$v_{19}$[4])+($v_{19}$[4],$v_{24}$[3])+(**$v_{24}$[3]**,**$v_{25}$[2]**)+($v_{25}$[2],$v_{31}$[1]);

Тогда расположение вершин в последовательности вершин по уровням Q имеет вид:

$u_1$: 31 29 28 10 02 01
$u_2$: 25 14 27 13 06 03 05 22
$u_3$: 24 26 12 15 12 07 11 04 21 23
$u_4$: 30
$u_5$:

Расположение вершин $v_{24}$,$v_{26}$ уровня 3 в цикле $c_{15}$, порождает появление вершины $v_{30}$ уровня 4.

$c_{15}$ = ($v_{25}$[2],$v_{24}$[3])+(**$v_{24}$[3]**,**$v_{30}$[4]**)+(**$v_{30}$[4]**,**$v_{26}$[3]**)+($v_{26}$[3],$v_{25}$[2]);

Тогда расположение вершин в последовательности вершин по уровням Q имеет вид:

$u_1$: 31 29 28 10 02 01
$u_2$: 25 14 27 13 06 03 05 22
$u_3$: 24 26 12 15 12 07 11 04 21 23
$u_4$: 30
$u_5$:

Следующей пары $v_{26},v_{12}$ в циклах не существует, а существует пара $v_{26},v_7$ в цикле $c_5$ пропуская $v_{12},v_{15},v_{12}$. Данная пара порождает включение вершины $v_8$ в последовательности уровня 4.

$c_5 = (v_{31}[1],v_{25}[2])+(v_{25}[2],v_{26}[3])+(\mathbf{v_{26}[3],v_8[4]})+(v_8[4],\mathbf{v_7[3]})+(v_7[3],v_{12}[3])+(v_{12}[3],v_{14}[2])+$
$+(v_{14}[2],v_{31}[1]);$

Тогда расположение вершин в последовательности вершин по уровням Q имеет вид:

$u_1$: 31 29 28 10 02 01
$u_2$: 25 14 27 13 06 03 05 22
$u_3$: 24 26 12 15 12 07 11 04 21 23
$u_4$: 30 08
$u_5$:

Расположение вершин $v_7,v_{11}$ уровня 3 в цикле $c_3$, порождает появление вершин $v_8,v_{17}$ уровня 4.

$c_3 = (v_6[2],v_7[3])+(\mathbf{v_7[3],v_8[4]})+(v_8[4],v_9[5])+(v_9[5],v_{18}[5])+(v_{18}[5],v_{17}[4])+(\mathbf{v_{17}[4],v_{11}[3]})+$
$+(v_{11}[3],v_6[2]);$

Тогда расположение вершин в последовательности вершин по уровням Q имеет вид:

$u_1$: 31 29 28 10 02 01
$u_2$: 25 14 27 13 06 03 05 22
$u_3$: 24 26 12 15 12 07 11 04 21 23
$u_4$: 30 08 17
$u_5$:

Расположение вершин $v_{11},v_4$ уровня 3 в цикле $c_8$, порождает появление вершин $v_{17},v_{16}$ уровня 4.

$c_8 = (v_3[2],v_{11}[3])+(\mathbf{v_{11}[3],v_{17}[4]})+(v_{17}[4],v_{16}[4])+(\mathbf{v_{16}[4],v_4[3]})+(v_4[3],v_3[2]);$

Тогда расположение вершин в последовательности вершин по уровням Q имеет вид:

$u_1$: 31 29 28 10 02 01
$u_2$: 25 14 27 13 06 03 05 22
$u_3$: 24 26 12 15 12 07 11 04 21 23
$u_4$: 30 08 17 16
$u_5$:

Расположение вершин $v_4,v_{21}$ уровня 3 в цикле $c_9$, порождает появление вершин $v_{16},v_{20}$ уровня 4.

$c_9 = (v_5[2],v_4[3])+(\mathbf{v_4[3],v_{16}[4]})+(v_{16}[4],v_{20}[4])+(\mathbf{v_{20}[4],v_{21}[3]})+(v_{21}[3],v_5[2]);$

Тогда расположение вершин в последовательности вершин по уровням Q имеет вид:

$u_1$: 31 29 28 10 02 01
$u_2$: 25 14 27 13 06 03 05 22

$u_3$: 24 26 12 15 12 07 11 04 21 23
$u_4$: 30 08 17 16 20
$u_5$:

Расположение вершин $v_{21}, v_{23}$ уровня 3 в цикле $c_{16}$, порождает появление вершины $v_{20}$ уровня 4.

$c_{16} = (v_{22}[2], v_{21}[3]) + (\mathbf{v_{21}[3], v_{20}[4]}) + (\mathbf{v_{20}[4], v_{23}[3]}) + (v_{23}[3], v_{22}[2]);$

Тогда расположение вершин в последовательности вершин по уровням Q имеет вид:

$u_1$: 31 29 28 10 02 01
$u_2$: 25 14 27 13 06 03 05 22
$u_3$: 24 26 12 15 12 07 11 04 21 23
$u_4$: 30 08 17 16 20
$u_5$:

Расположение вершин $v_{23}, v_{26}$ уровня 3 в цикле $c_{14}$, порождает появление вершины $v_{19}$ уровня 4.

$c_{14} = (v_{31}[1], v_{22}[2]) + (v_{22}[2], v_{23}[3]) + (\mathbf{v_{23}[3], v_{19}[4]}) + (v_{19}[4], \mathbf{v_{24}[3]}) + (v_{24}[3], v_{25}[2]) + (v_{25}[2], v_{31}[1]);$

Тогда расположение вершин в последовательности вершин по уровням Q имеет вид:

$u_1$: 31 29 28 10 02 01
$u_2$: 25 14 27 13 06 03 05 22
$u_3$: 24 26 12 15 12 07 11 04 21 23
$u_4$: 30 08 17 16 20 19
$u_5$:

Рассмотрим второй этап построения последовательности уровня 4.

Расположение вершин $v_{23}, v_{26}$ уровня 3 в цикле $c_{14}$, порождает появление вершины $v_{19}$ уровня 4.

$c_{14} = (v_{31}[1], v_{22}[2]) + (v_{22}[2], v_{23}[3]) + (\mathbf{v_{23}[3], v_{19}[4]}) + (\mathbf{v_{19}[4], v_{24}[3]}) + (v_{24}[3], v_{25}[2]) +$
$+ (v_{25}[2], v_{31}[1]);$

Тогда расположение вершин в последовательности вершин по уровням Q имеет вид:

$u_1$: 31 29 28 10 02 01
$u_2$: 25 14 27 13 06 03 05 22
$u_3$: 24 26 12 15 12 07 11 04 21 23
$u_4$: 30 08 17 16 20 19
$u_5$:

Симметрическое расположение вершины $v_{24}$ уровня 3 в цикле $c_{11}$, порождает установку пары вершин $v_{19}, v_{30}$ в последовательности вершин уровня 4.

$c_{11} = (\mathbf{v_{24}[3], v_{19}[4]}) + (v_{19}[4], v_{18}[5]) + (v_{18}[5], v_{9}[5]) + (v_{9}[5], v_{30}[4]) + (\mathbf{v_{30}[4], v_{24}[3]});$

Тогда расположение вершин в последовательности вершин по уровням Q имеет вид:

$u_1$: 31 29 28 10 02 01
$u_2$: 25 14 27 13 06 03 05 22
$u_3$: 24 26 12 15 12 07 11 04 21 23
$u_4$: 19 30 08 17 16 20
$u_5$:

Симметрическое расположение вершины $v_{26}$ уровня 3 в цикле $c_{12}$, порождает установку пары вершин $v_{30}, v_8$ в последовательности вершин уровня 4.

$c_{12} = (\mathbf{v_{26}[3]}, \mathbf{v_{30}[4]}) + (v_{30}[4], v_9[5]) + (v_9[5], v_8[4]) + (\mathbf{v_8[4]}, \mathbf{v_{26}[3]})$;

Расположение вершин в последовательности вершин по уровням Q, не меняется.

Симметрическое расположение вершины $v_{23}$ уровня 3 в цикле $c_{10}$, порождает обратное расположение вершин $v_{20}, v_{16}, v_{17}$ и $v_{19}$, что указывает на пересечение ребер. Расположение вершин в последовательности U остается без изменений.

$c_{10} = (\mathbf{v_{23}[3]}, \mathbf{v_{20}[4]}) + (v_{20}[4], \mathbf{v_{16}[4]}) + (v_{16}[4], \mathbf{v_{17}[4]}) + (v_{17}[4], v_{18}[5]) + (v_{18}[5], \mathbf{v_{19}[4]}) +$
$+ (v_{19}[4], v_{23}[3])$;

| | | | | | | | | | | |
|---|---|---|---|---|---|---|---|---|---|---|
| $u_1$: | 31 | 29 | 28 | 10 | 02 | 01 | | | | |
| $u_2$: | 25 | 14 | 27 | 13 | 06 | 03 | 05 | 22 | | |
| $u_3$: | 24 | 26 | 12 | 15 | 12 | 07 | 11 | 04 | 21 | 23 |
| $u_4$: | 19 | 30 | 08 | 17 | 16 | 20 | | | | |
| $u_5$: | | | | | | | | | | |

Формирование последовательности вершин уровня 4 завершено.

Будем формировать последовательность вершин уровня 5.

Рассмотрим пару вершин $v_{19}, v_{30}$ уровня 4 цикла $c_{11}$ порождает установку вершин $v_{18}, v_9$ уровня 5

$c_{11} = (v_{24}[3], v_{19}[4]) + (\mathbf{v_{19}[4]}, \mathbf{v_{18}[5]}) + (v_{18}[5], v_9[5]) + (\mathbf{v_9[5]}, \mathbf{v_{30}[4]}) + (v_{30}[4], v_{24}[3])$;

Тогда расположение вершин в последовательности вершин по уровням Q имеет вид:

| | | | | | | | | | | |
|---|---|---|---|---|---|---|---|---|---|---|
| $u_1$: | 31 | 29 | 28 | 10 | 02 | 01 | | | | |
| $u_2$: | 25 | 14 | 27 | 13 | 06 | 03 | 05 | 22 | | |
| $u_3$: | 24 | 26 | 12 | 15 | 12 | 07 | 11 | 04 | 21 | 23 |
| $u_4$: | 19 | 30 | 08 | 17 | 16 | 20 | | | | |
| $u_5$: | 18 | 09 | | | | | | | | |

Рассмотрим пару вершин $v_{30}, v_8$ уровня 4 цикла $c_{12}$ порождает установку вершины $v_9$ уровня 5

$c_{12} = (v_{26}[3], v_{30}[4]) + (\mathbf{v_{30}[4]}, \mathbf{v_9[5]}) + (\mathbf{v_9[5]}, \mathbf{v_8[4]}) + (v_8[4], v_{26}[3])$;

Тогда расположение вершин в последовательности вершин по уровням Q имеет вид:

| | | | | | | | | | | |
|---|---|---|---|---|---|---|---|---|---|---|
| $u_1$: | 31 | 29 | 28 | 10 | 02 | 01 | | | | |
| $u_2$: | 25 | 14 | 27 | 13 | 06 | 03 | 05 | 22 | | |
| $u_3$: | 24 | 26 | 12 | 15 | 12 | 07 | 11 | 04 | 21 | 23 |
| $u_4$: | 19 | 30 | 08 | 17 | 16 | 20 | | | | |
| $u_5$: | 18 | 09 | | | | | | | | |

Рассмотрим пару вершин $v_8, v_{17}$ уровня 4 цикла $c_3$ порождает установку вершин $v_9, v_{18}$ уровня 5

$c_3 = (v_6[2], v_7[3]) + (v_7[3], v_8[4]) + (\mathbf{v_8[4]}, \mathbf{v_9[5]}) + (v_9[5], v_{18}[5]) + (\mathbf{v_{18}[5]}, \mathbf{v_{17}[4]}) + (v_{17}[4], v_{11}[3]) +$
$+ (v_{11}[3], v_6[2])$;

Тогда расположение вершин в последовательности вершин по уровням Q имеет вид:

| | | | | | | |
|---|---|---|---|---|---|---|
| $u_1$: | 31 | 29 | 28 | 10 | 02 | 01 |

| | | | | | | | | | | |
|---|---|---|---|---|---|---|---|---|---|---|
| $u_2$: | 25 | 14 | 27 | 13 | 06 | 03 | 05 | 22 | | |
| $u_3$: | 24 | 26 | 12 | 15 | 12 | 07 | 11 | 04 | 21 | 23 |
| $u_4$: | 19 | 30 | 08 | 17 | 16 | 20 | | | | |
| $u_5$: | 18 | 09 | 18 | | | | | | | |

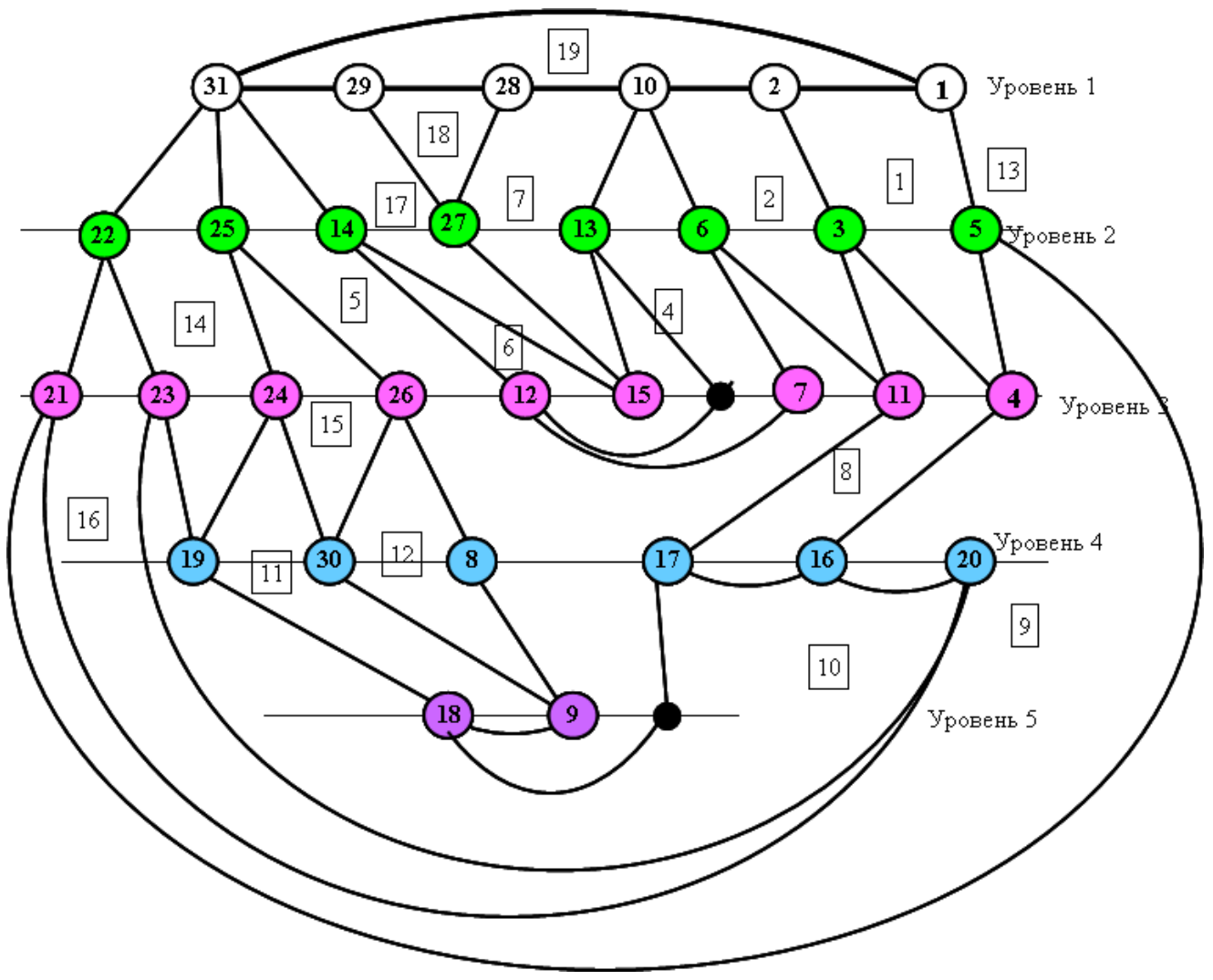

Рис. 10.5. Граф уровней с введением дополнительных вершин.

Построение последовательности верщин в уровнях U завершено. Окончательный вид последовательности U представлен на рис. 10.6:

| | | | | | | | | | | |
|---|---|---|---|---|---|---|---|---|---|---|
| $u_1$: | 31 | 29 | 28 | 10 | 02 | 01 | | | | |
| $u_2$: | 25 | 14 | 27 | 13 | 06 | 03 | 05 | 22 | | |
| $u_3$: | 24 | 26 | 12 | 15 | 12 | 07 | 11 | 04 | 21 | 23 |
| $u_4$: | 19 | 30 | 08 | 17 | 16 | 20 | | | | |
| $u_5$: | 18 | 09 | 18 | | | | | | | |

Наличие двух вершин с одинаковым номером в одном уровне, говорит о наличии их соединения в уровне. Причем, одна из вершин объявляется мнимой (см. рис. 10.6).

## 10.2. Геометрические методы проведения соединений

Построим замкнутую последовательность вершин по уровням U, расположив по окружности самую длинную последовательность. В нашем случае это уровень 3 (см. рис. 10.7).

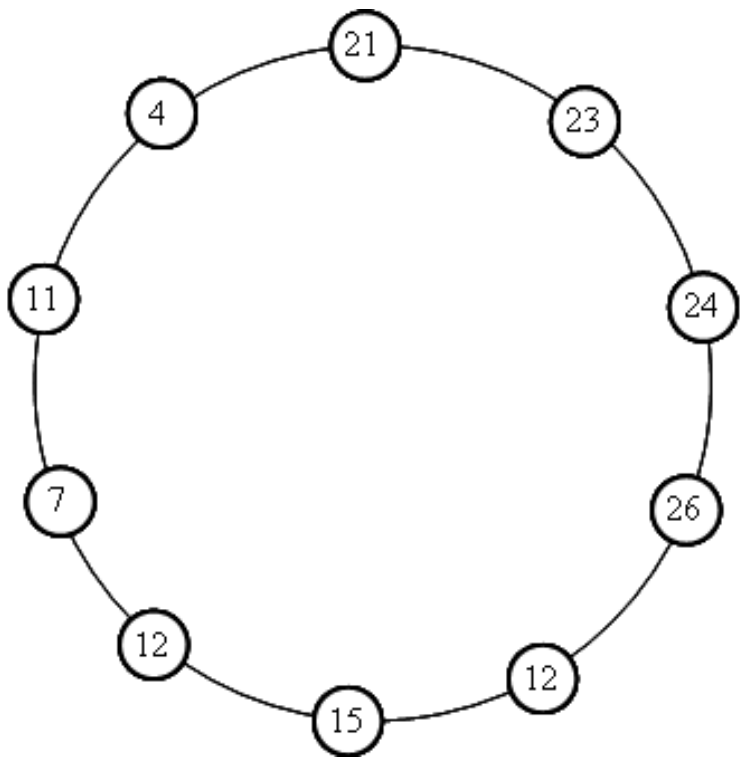

Рис. 10.7. Расположение по часовой стрелке последовательности вершин уровня 3.

**Определение 4.1.** Будем называть топо-метрическим расстоянием L между двумя близлежащими вершинами уровня N, смежные вершины следующего уровня (N+1 или N-1).

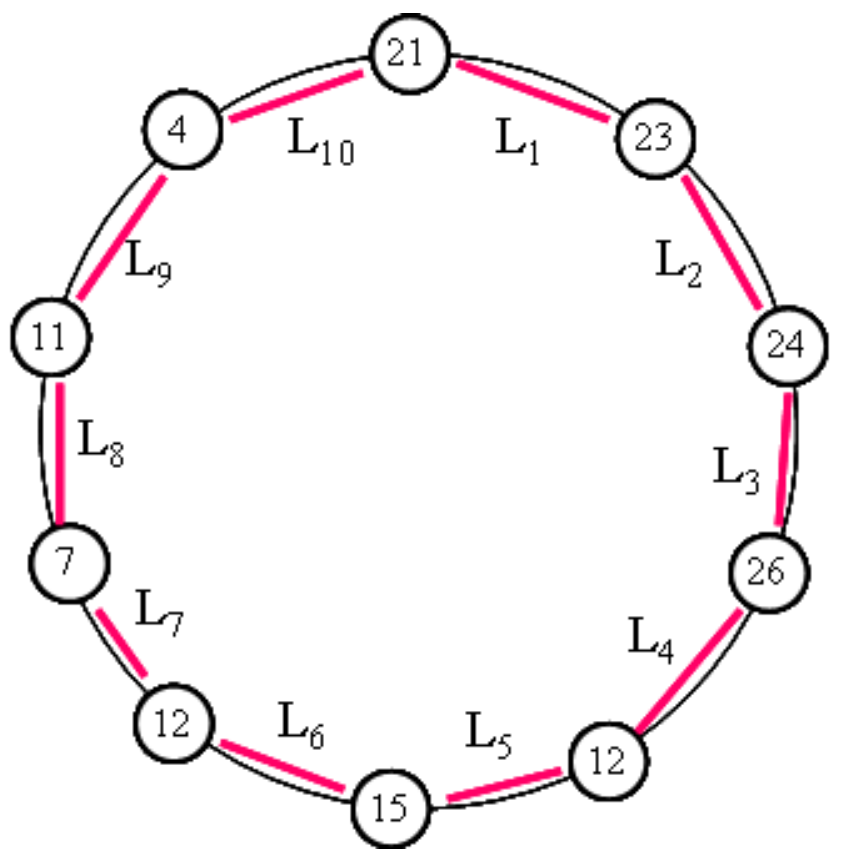

Рис. 10.8. Разметка топо-метрических расстояний L между вершинами уровня 3.

Например, вершина $v_{21}$ уровня 3 смежна с вершиной $v_{22}$ уровня 2, а вершина $v_{23}$ уровня 3 смежна с вершиной $v_{22}$ уровня 2. Этот процесс можно записать в виде:

$$L_1(v_{21},v_{23}) \to (v_{22},v_{22}) \quad (4.1)$$

То есть две близлежащие вершины уровня N порождают концевые вершины уровня N+1 или N-1. Запишем топо-метрические расстояния L для всей последовательности вершин уровня 3 (см. рис. 10.9).

$L_1(v_{21},v_{23}) \to (v_{22},v_{22})$;
$L_2(v_{23},v_{24}) \to (v_{22},v_{25})$;
$L_3(v_{24},v_{26}) \to (v_{25},v_{25})$;
$L_4(v_{26},v_{12}) \to (v_{25},v_{14})$;
$L_5(v_{12},v_{15}) \to (v_{14},v_{27})$;
$L_6(v_{15},v_{12}) \to (v_{27},v_{13})$;
$L_7(v_{12},v_7) \to (v_{13},v_6)$;
$L_8(v_7,v_{11}) \to (v_6,v_3)$;
$L_9(v_{11},v_4) \to (v_3,v_5)$;
$L_{10}(v_4,v_{21}) \to (v_5,v_{22})$.

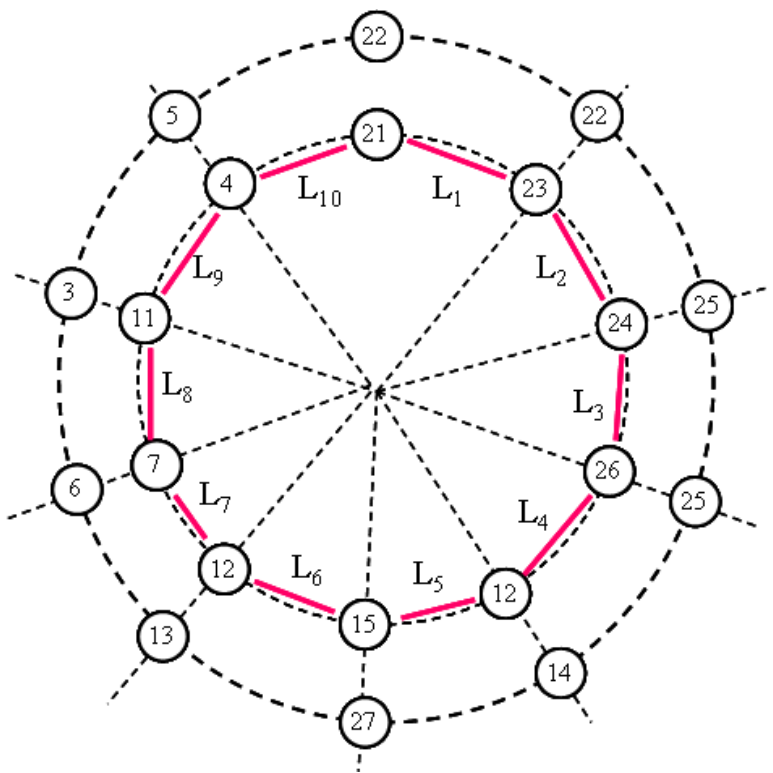

Рис. 10.9. Расположение концевых вершин топо-метрического расстояния L уровня 2.

**Определение 10.2.** Участком топо-метрического расстояния L уровня N+1 или N-1 называется последовательность вершин расположенная между концевыми вершинами топо-метрического расстояния L.

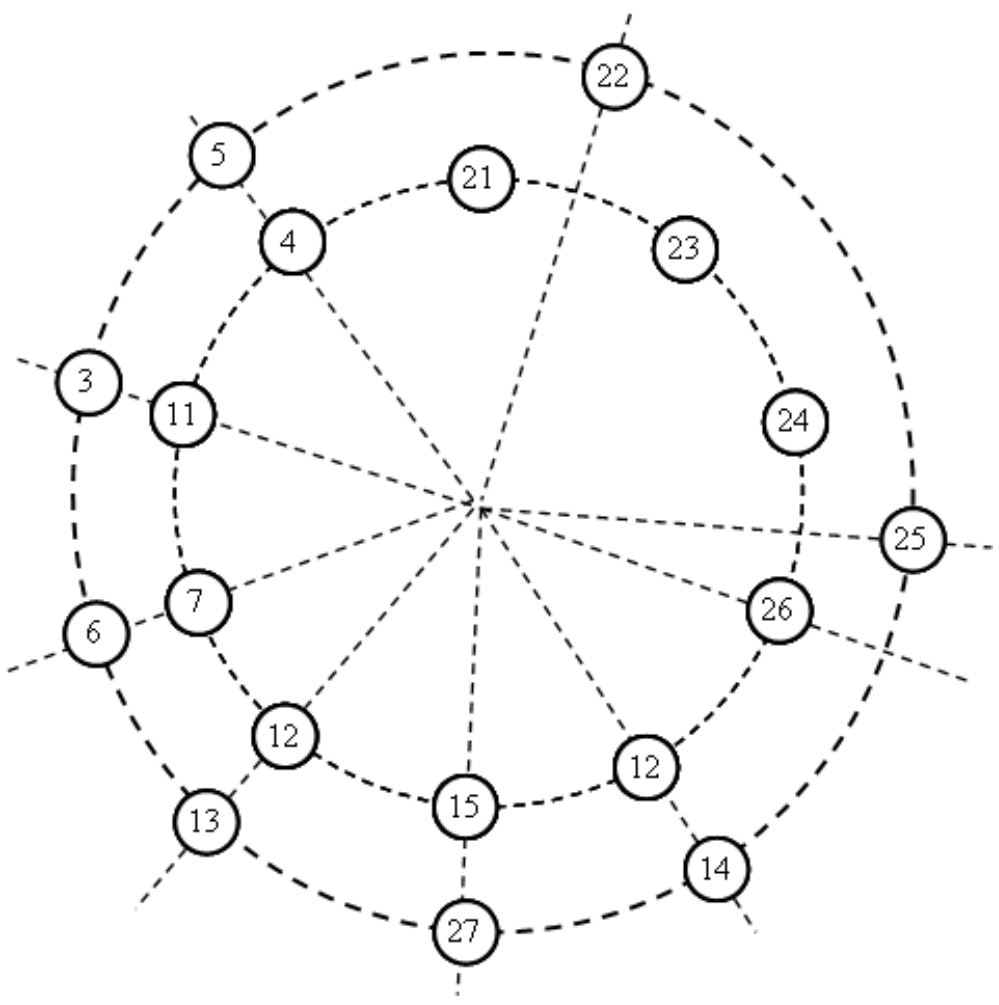

Рис. 10.10. Расположение вершин на участках уровня 2 с учетом их количества.

. Распологаем вершины с учетом их количества на участках уровня 2 (см. рис. 10.10).

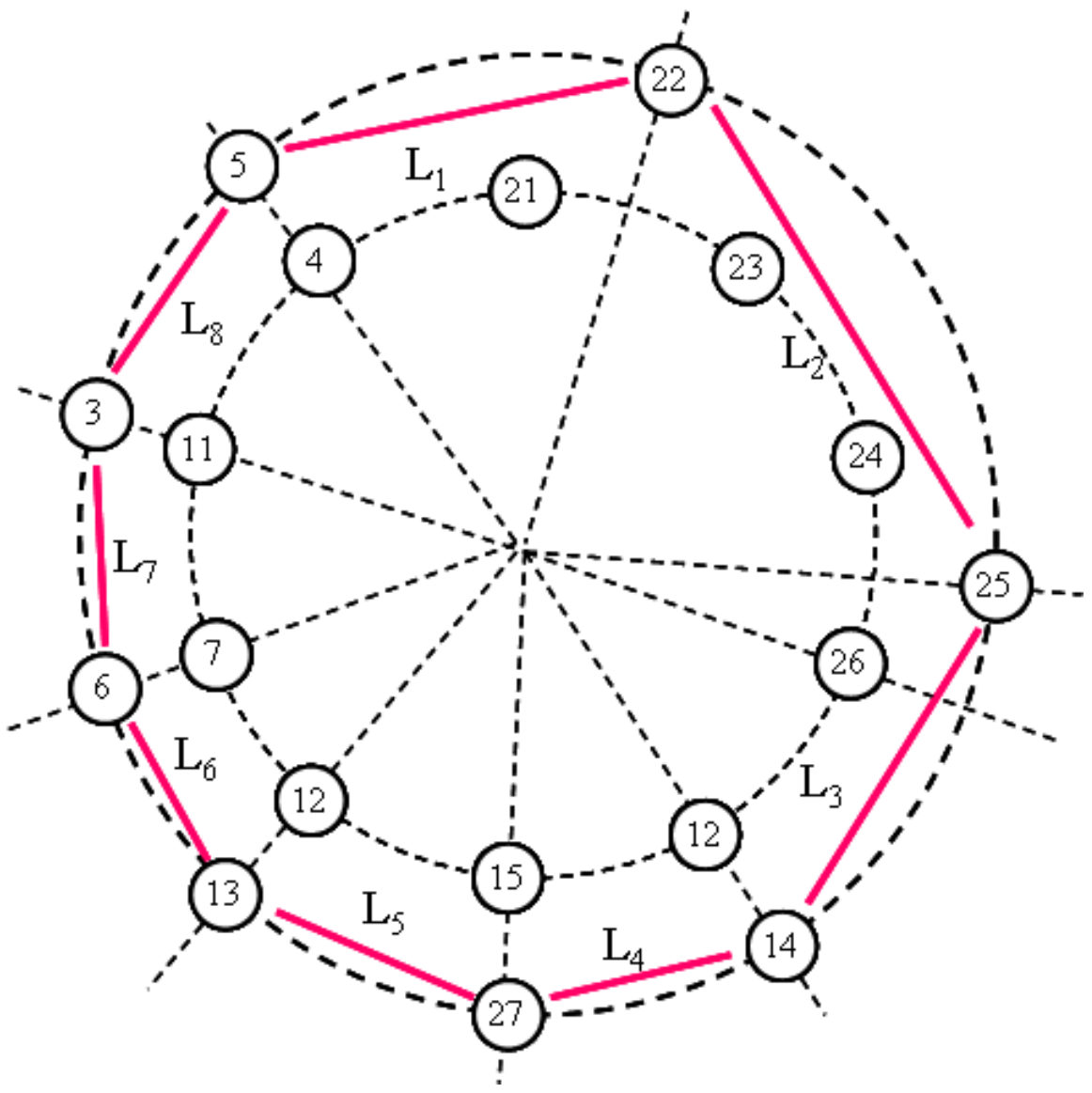

**Рис**. 10.11. Разметка топо-метрических расстояний L между вершинами уровня 2.

На рис. 10.11 представлено расстояние L между вершинами уровня 2.

$L_1(v_{22},v_{25}) \rightarrow (v_{31},v_{31})$;
$L_2(v_{25},v_{14}) \rightarrow (v_{31},v_{31})$;
$L_3(v_{14},v_{27}) \rightarrow (v_{31},v_{29})$;
$L_4(v_{27},v_{13}) \rightarrow (v_{28},v_{10})$;
$L_5(v_{13},v_6) \rightarrow (v_{10},v_2)$;
$L_6(v_6,v_3) \rightarrow (v_2,v_1)$;
$L_7(v_3,v_5) \rightarrow (v_1,v_{31})$;
$L_8(v_5,v_{22}) \rightarrow (v_{31},v_{31})$.

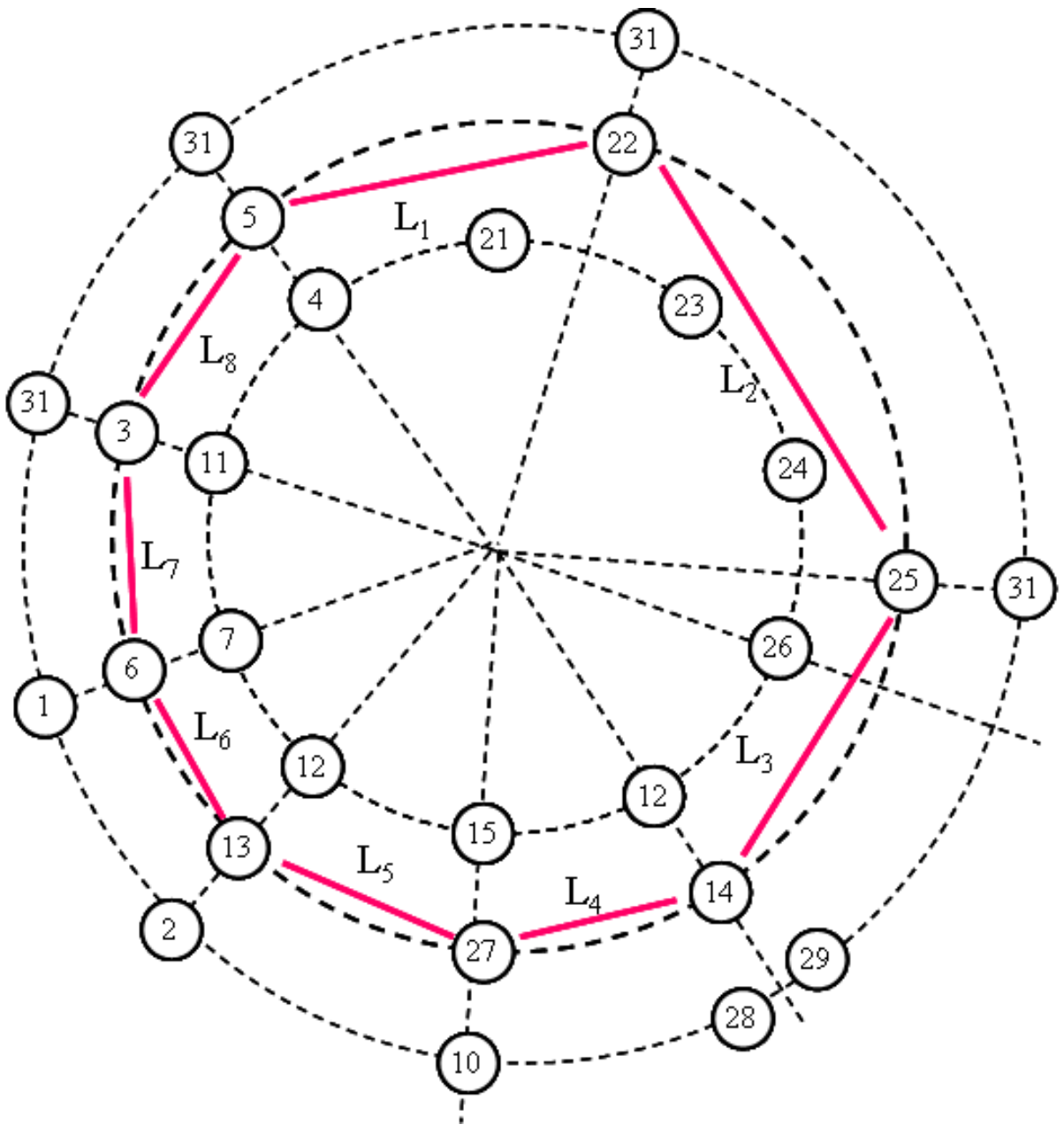

Рис. 10.12 . Расстановка вершин последовательности уровня 1.

На рис. 10.12 представлена расстановка вершин последовательности уровня 1. Расста-новка последовательности вершин для уровня 1 согласно их количеству на участке, представлено на рис. 10.13.

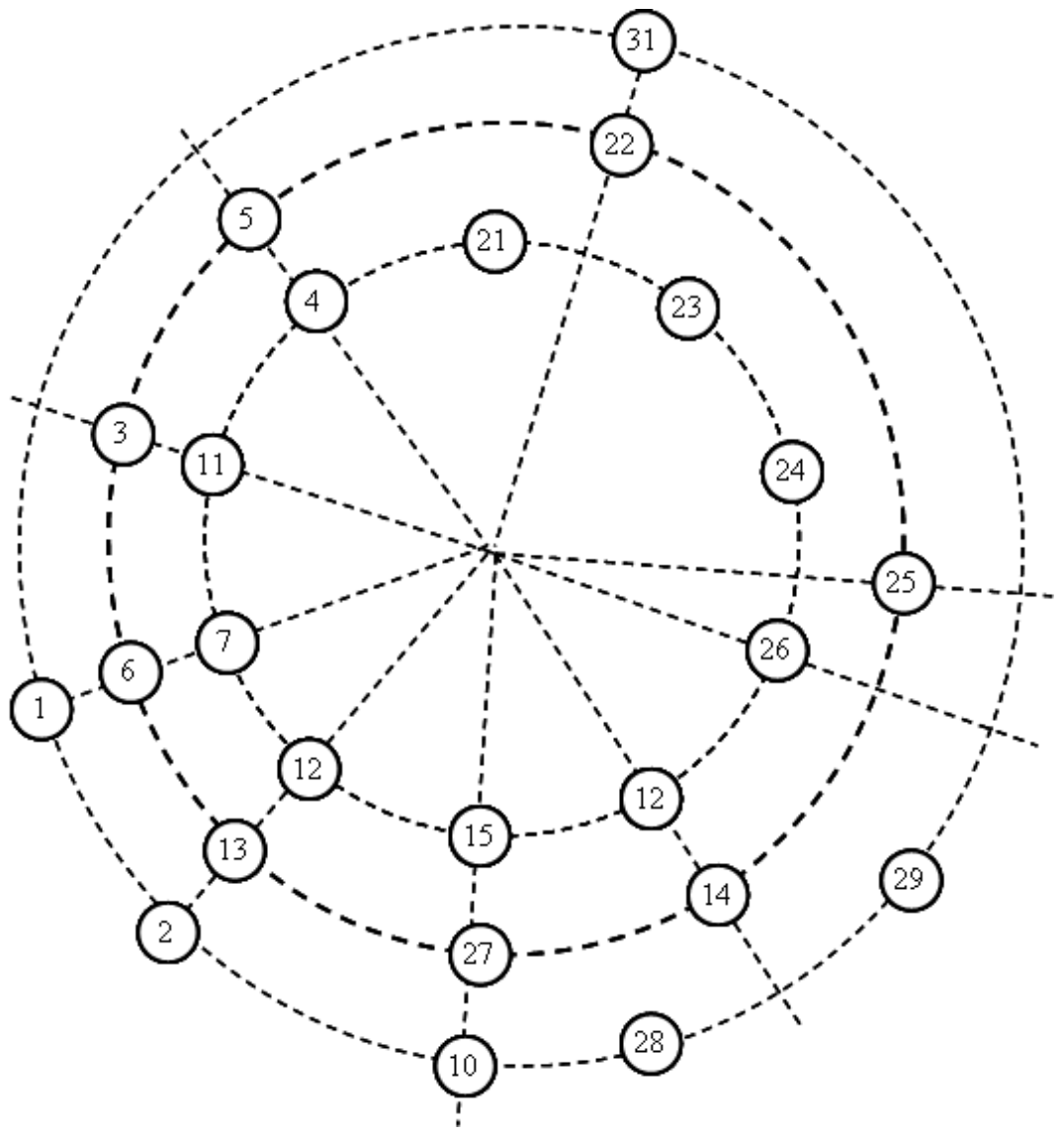

Рис. 10.13. Расстановка вершин последовательности уровня 1 согласно их количеству на участке.

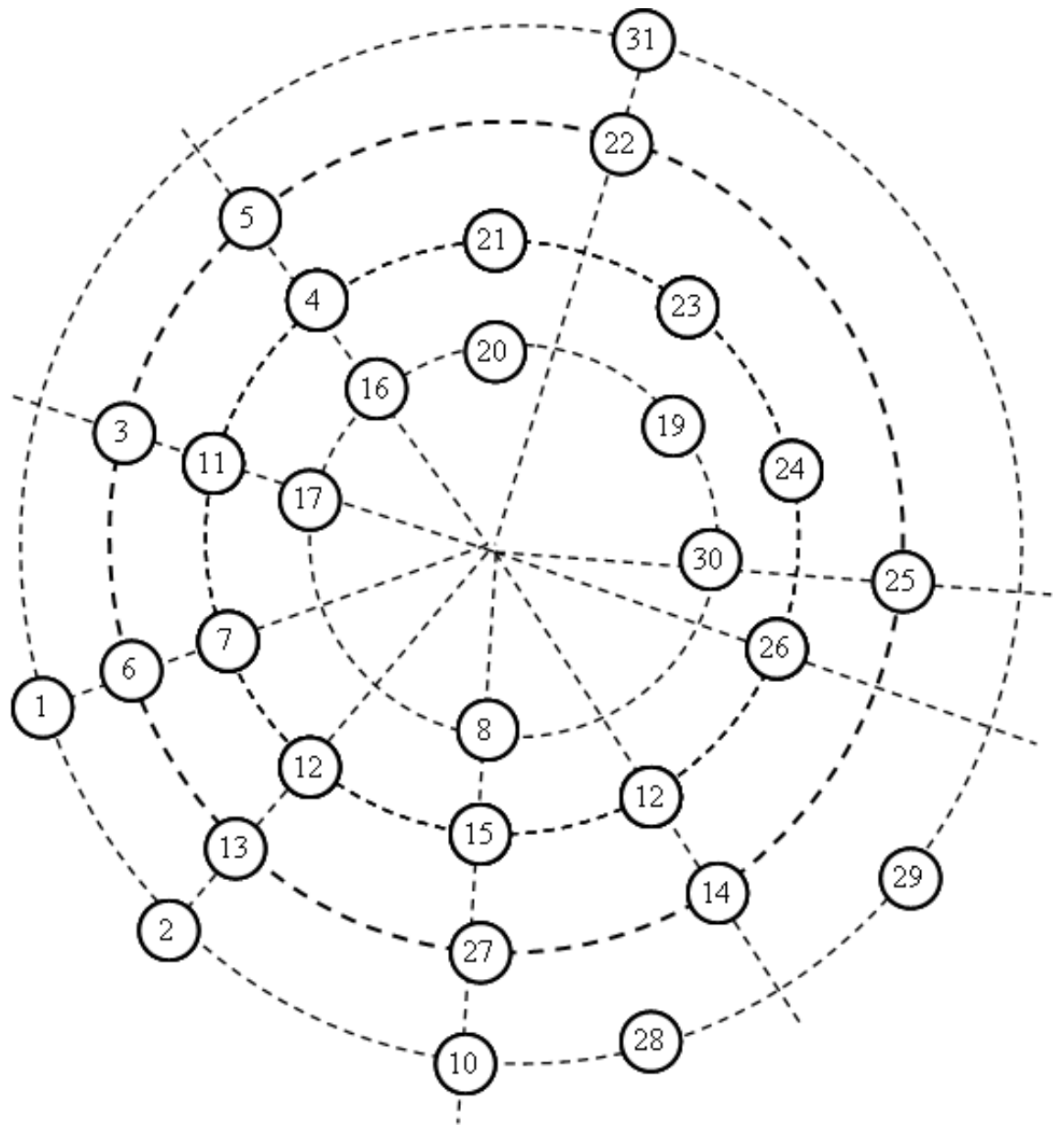

Рис. 10.14. Расстановка вершин в последовательности уровня 4.

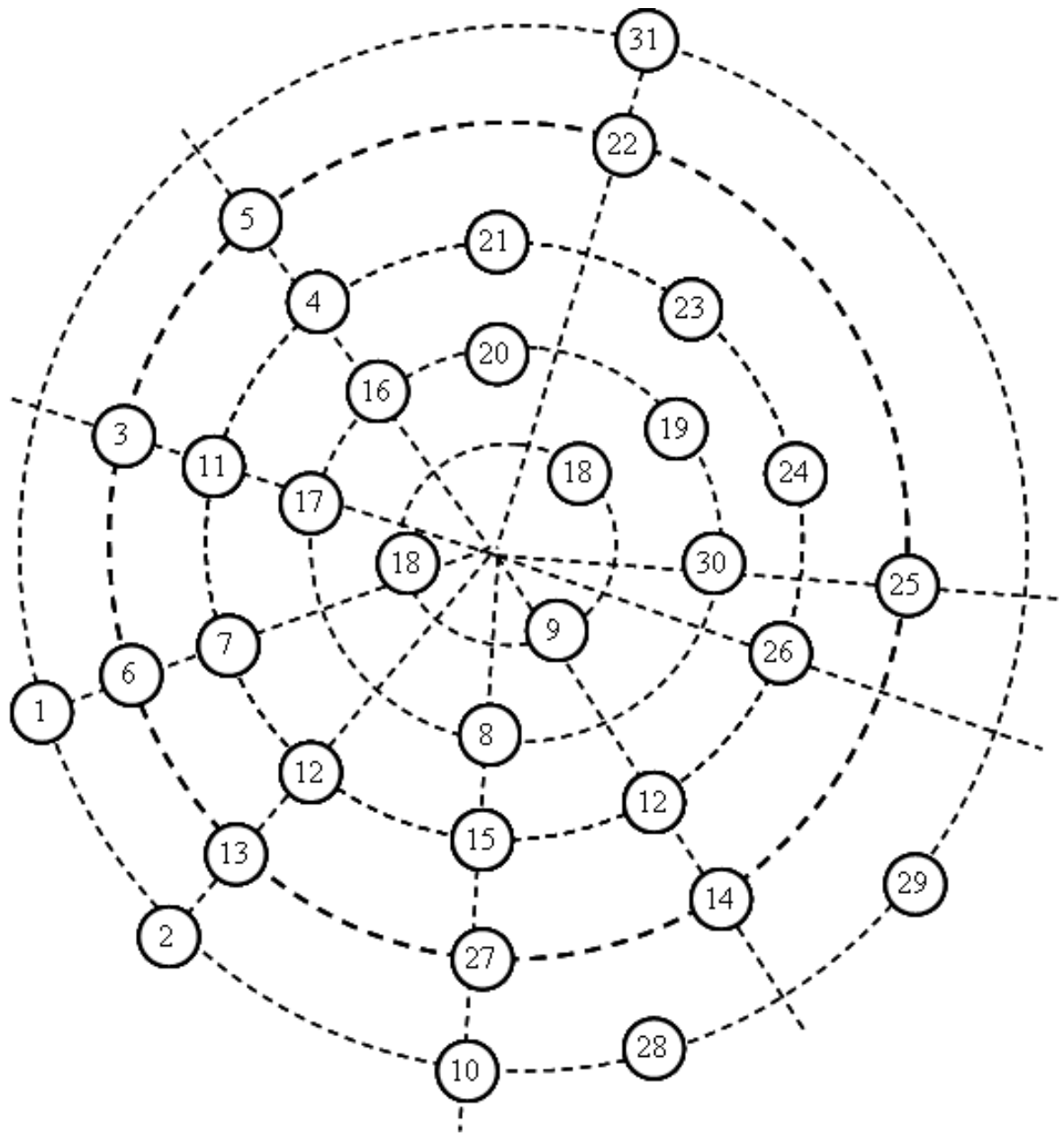

Рис. 10.15. Расстановка вершин в последовательности уровня 5.

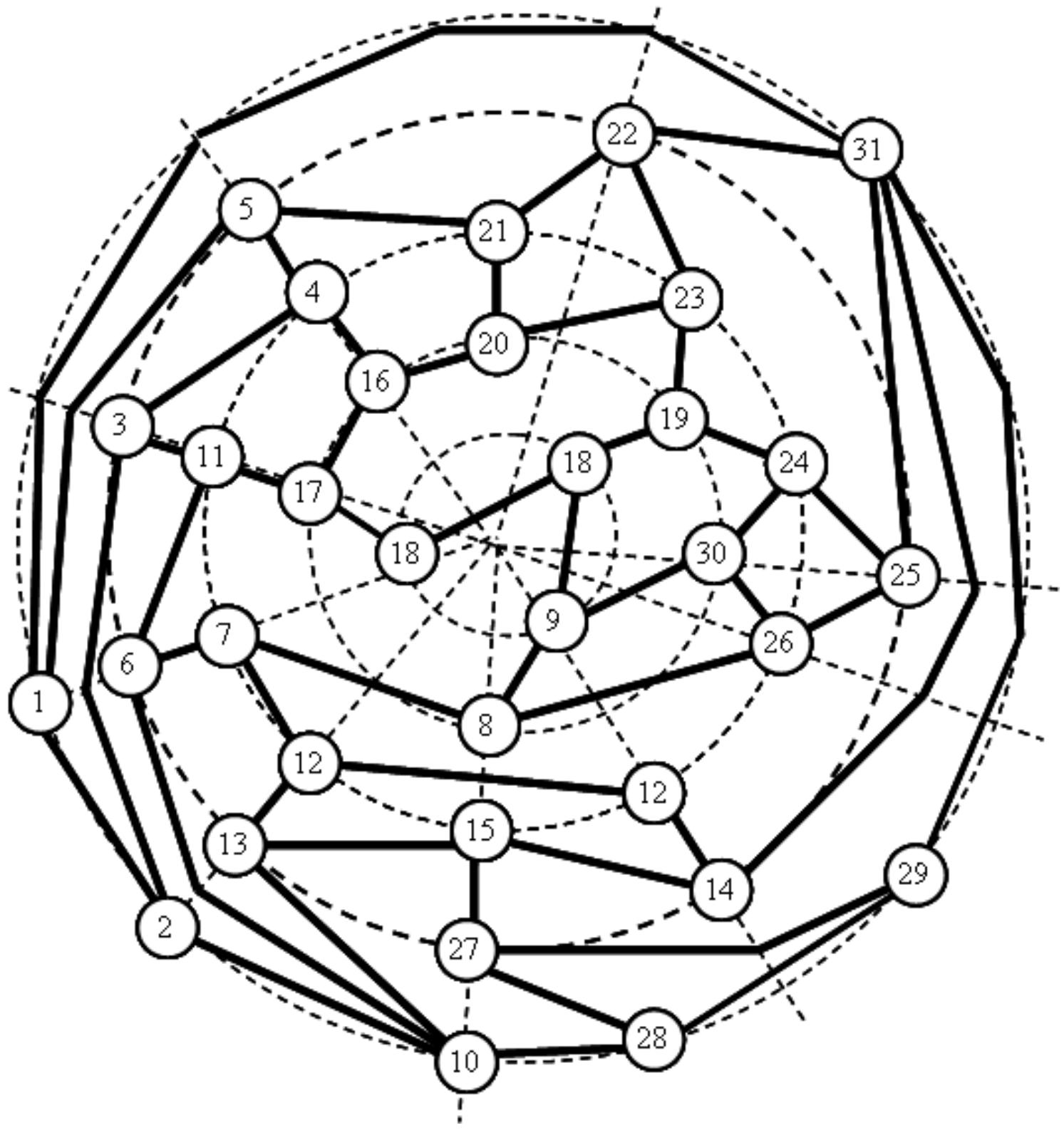

Рис. 10.15. Геометрическое проведение соединений в рисунке графа.

Будем рассматривать расположение вершин в прямоугольных областях.

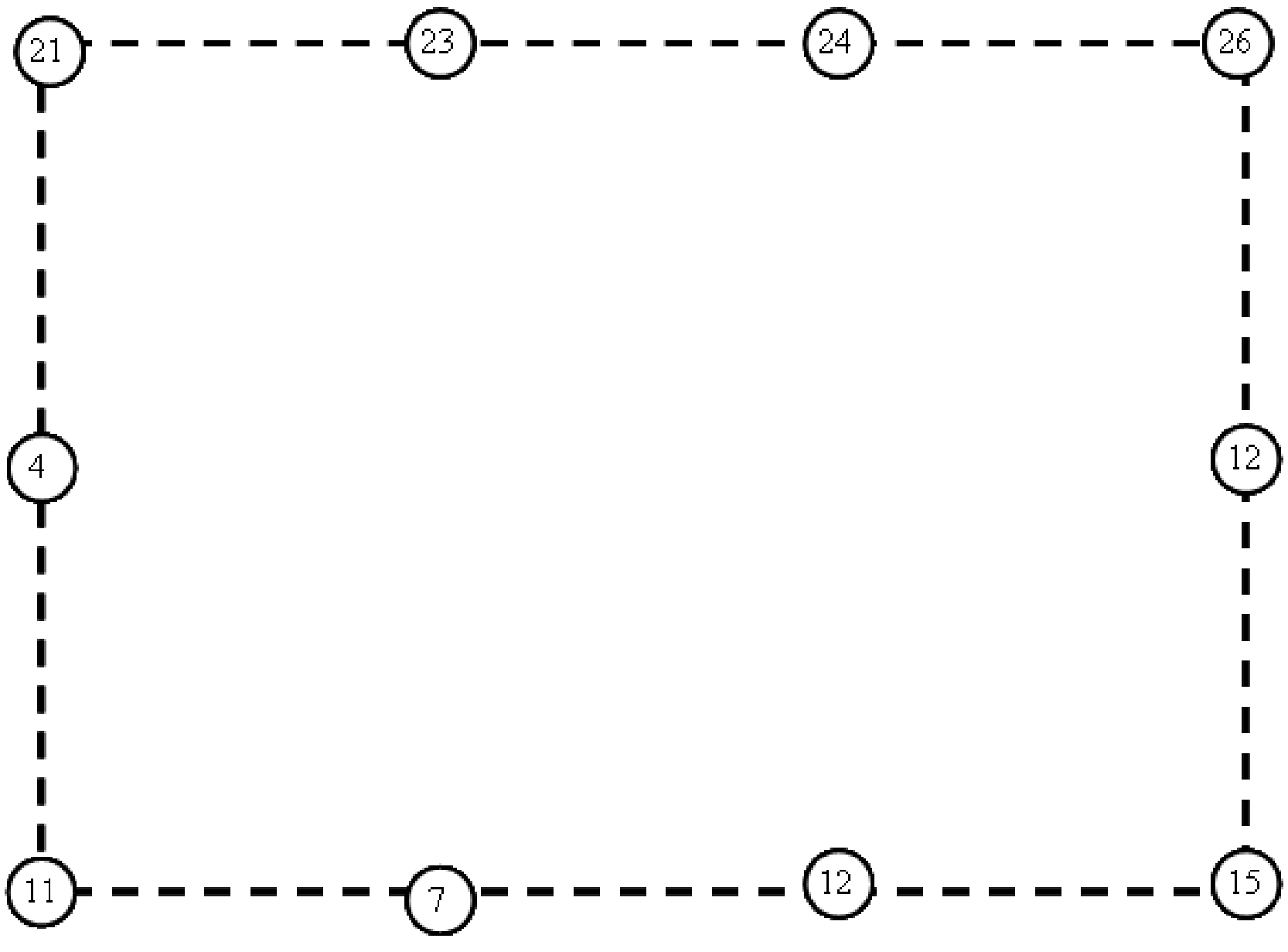

Рис. 10.17. Расстановка вершин в контуре прямоугольника для уровня 3.

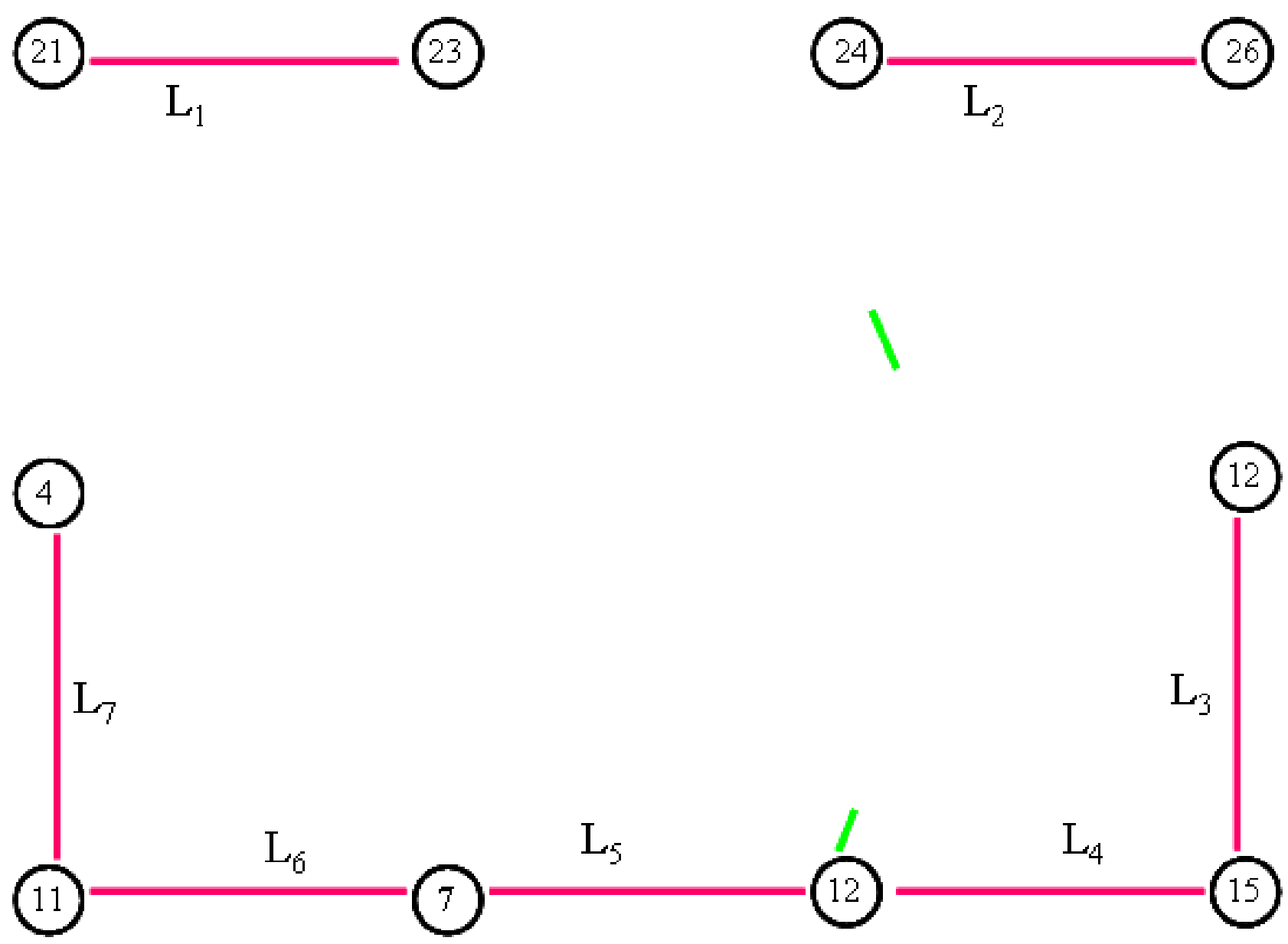

Рис. 10.18. Разметка топо-метрических расстояний L в контуре прямоугольника для уровня 3.

Размещение вершин в прямоугольных областей по сути ни чем не отличается от расположения вершин в круговых областях.

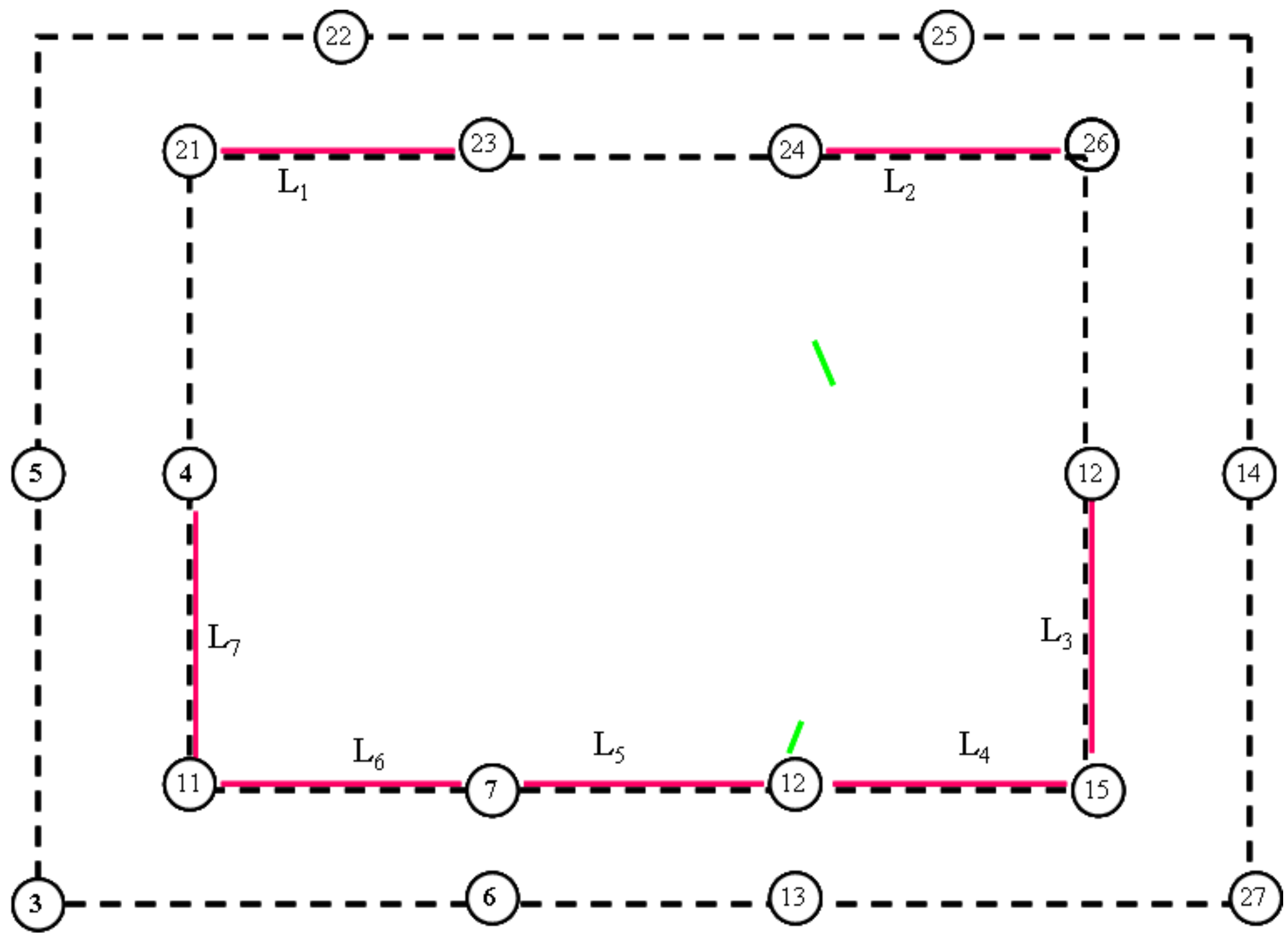

Рис. 10.19. Расстановка вершин в контуре прямоугольника для уровня 2.

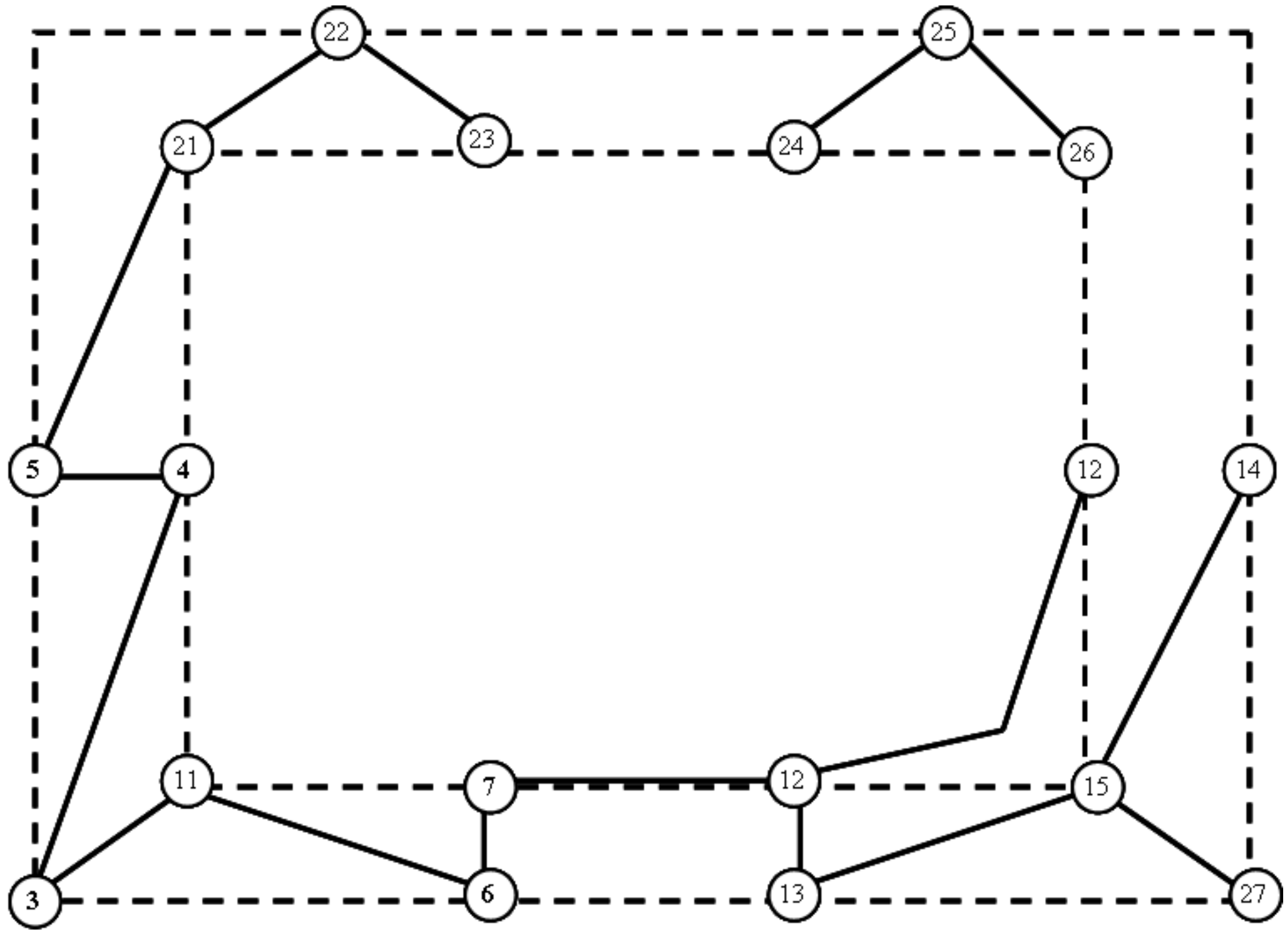

Рис. 10.20. Построение соединений между вершинами уровня 2 и уровня 3 в прямоугольных областях.

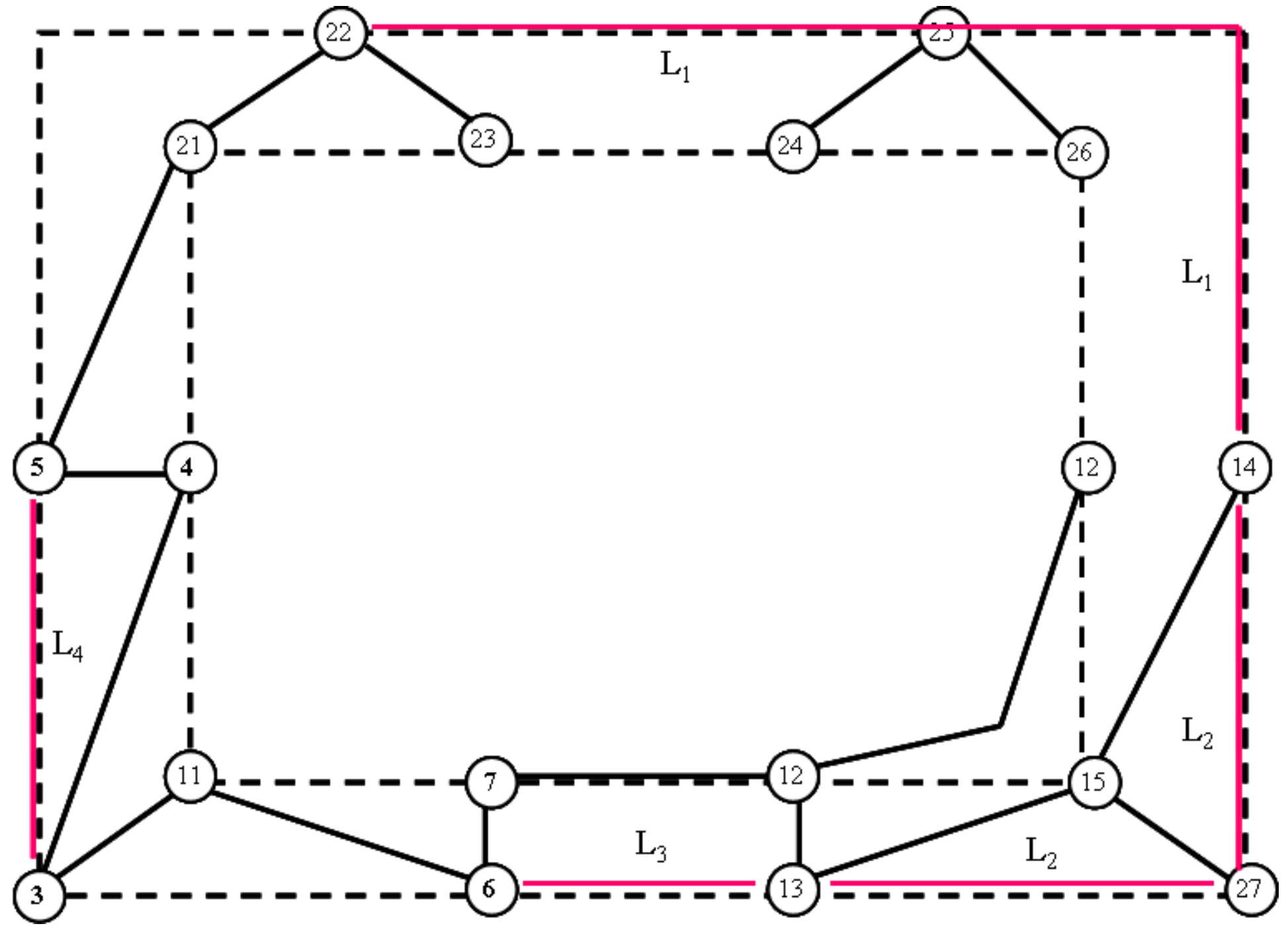

Рис. 10.21. Разметка топо-метрических расстояний L в контуре прямоугольника для уровня 2.

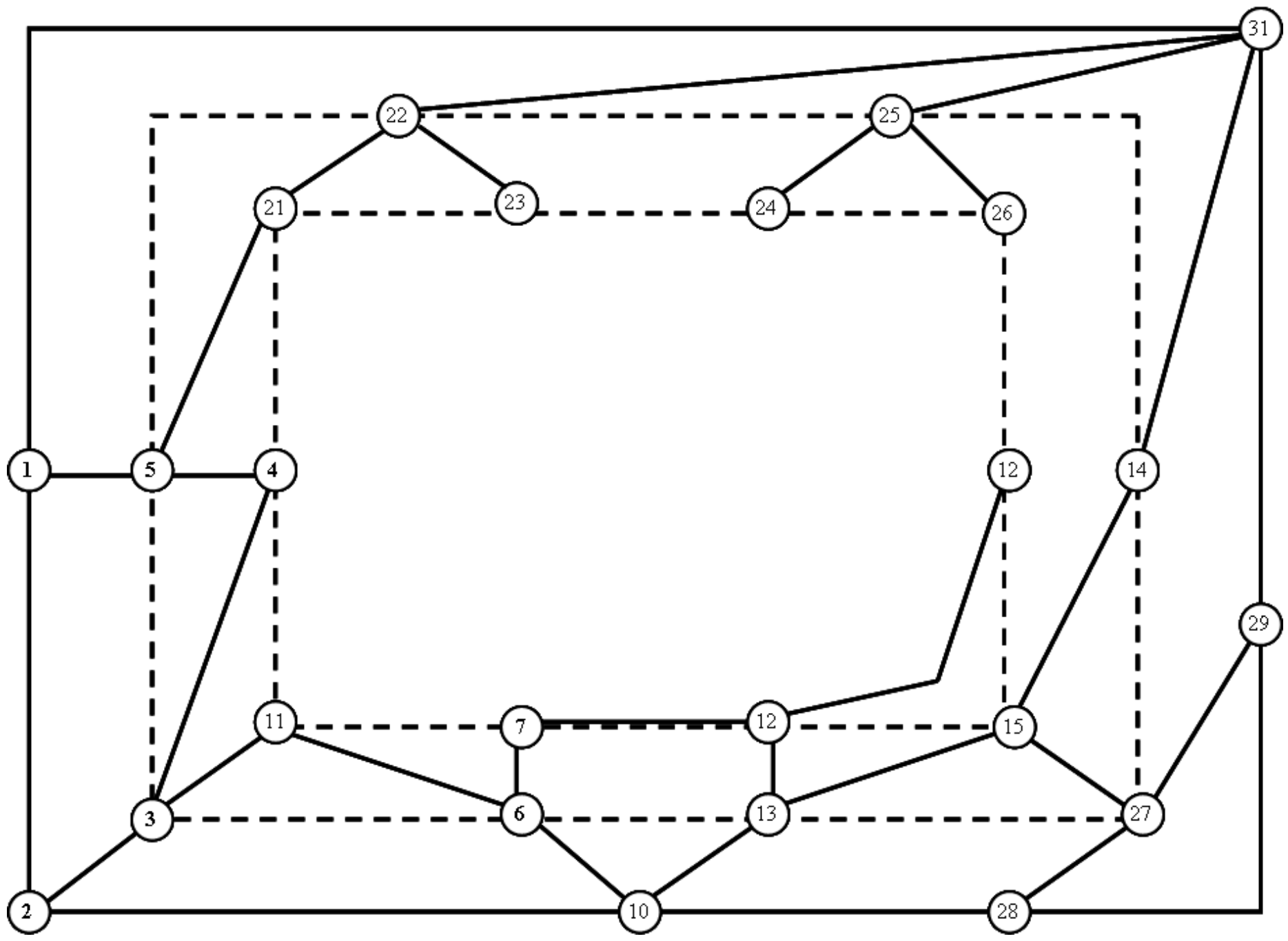

Рис. 10.22. Расстановка вершин в контуре прямоугольника для уровня 1.

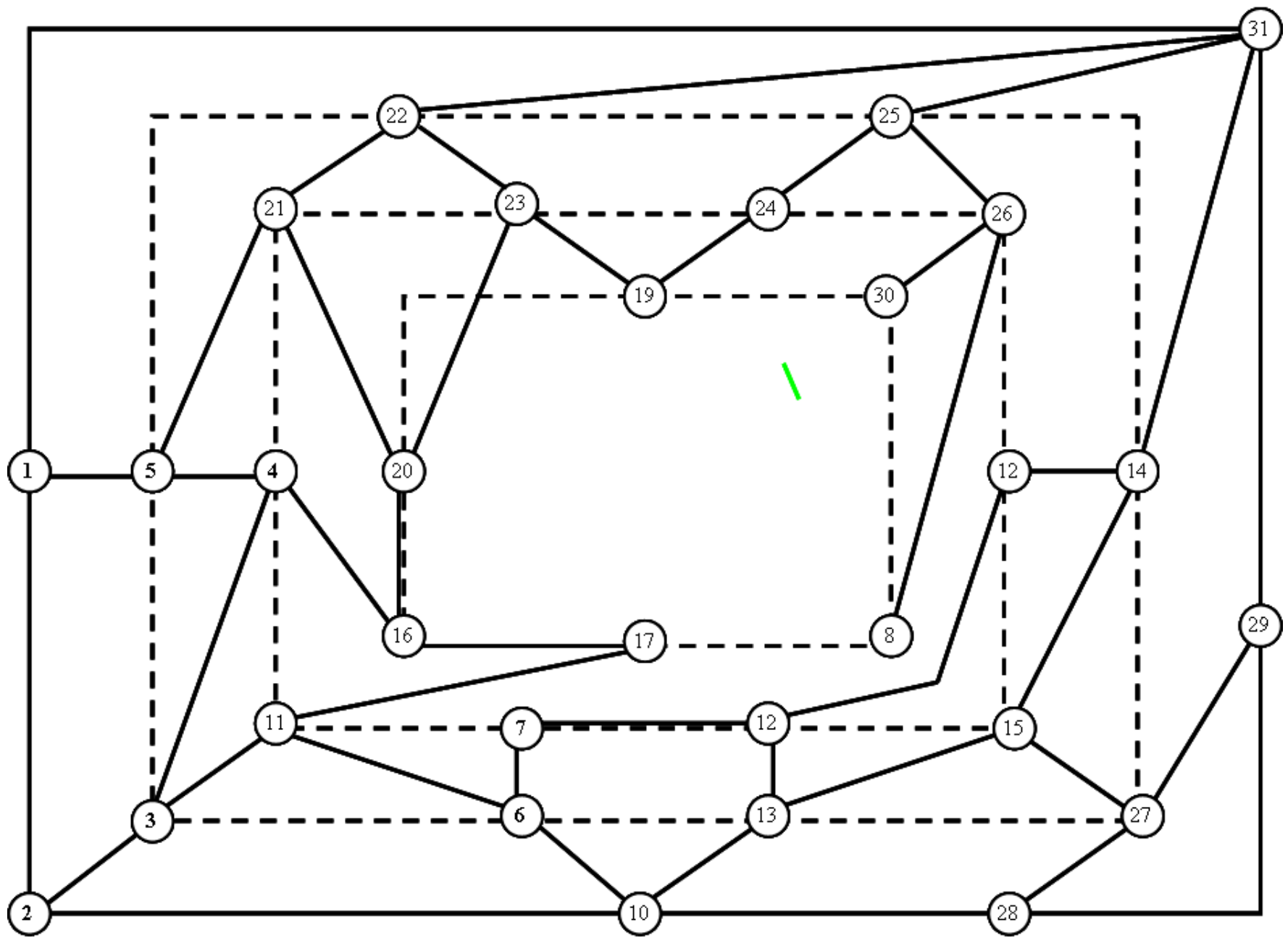

Рис. 10.23. Проведение соединений в уровнях 1,2,3,4 в контуре прямоугольника.

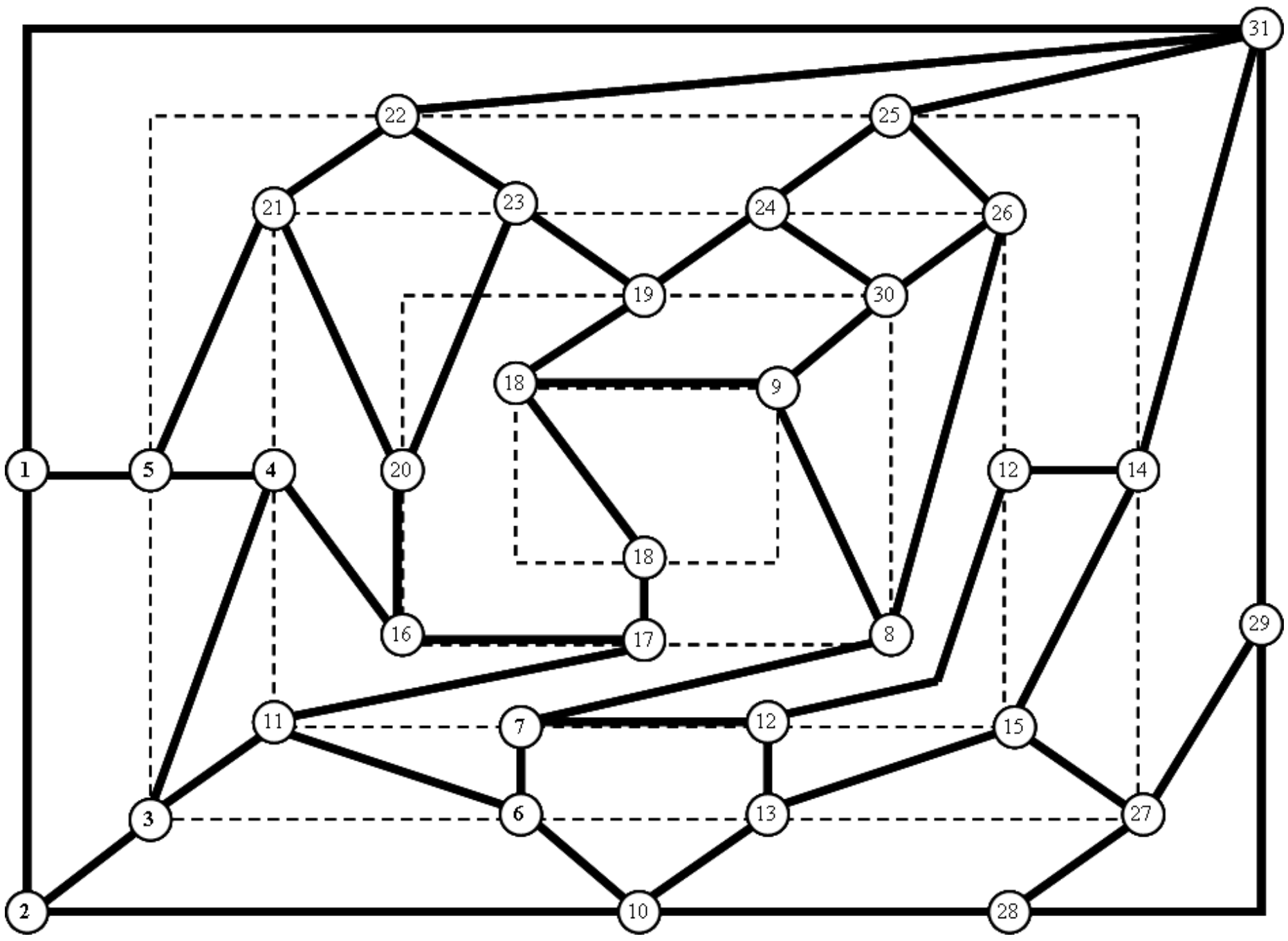

Рис. 10.24. Общий рисунок соединений в прямоугольных областях.

## 10.3. Силовой алгоритм размещения вершин

Математической моделью геометрического рисунка графа может служить силовая модель, представляющая ребра графа как пружины с заданным модулем упругости, причем вершины, принадлежащие выделенному циклу (ободу) жестко закреплены [2]. И тогда каждое ребро графа представляется вектором силы прямо пропорциональным его длине:

$$\overline{\mathrm{F}} = \mathrm{g}\bar{\mathrm{l}} \tag{10.1}$$

Здесь $\overline{\mathrm{F}}$ - вектор силы, g – модуль упругости, $\bar{\mathrm{l}}$ - вектор длины отрезка.

Для силовой модели сумма векторов сил при равновесии в точке, должна быть равна нулю:

$$\sum_{i=1}^{p} \overline{\mathrm{F}_i} = 0 \tag{10.2}$$

здесь *p* – количество сил действующих на точку (локальная степень вершины в графе).

К уравнению (2) можно добавить условие равенства нулю векторной суммы длин отрезков для любого замкнутого контура (цикла в графе):

$$\sum_{i=1}^{k} \overline{l_i} = 0 \tag{10.3}$$

здесь *k* – количество отрезков в контуре (длина цикла в графе).

Уравнения (7.2-7.3) аналогичны первому и второму законам Кирхгофа для электрической цепи, в предположении, что сила пружины F соответствует току ветви I, а длина пружины l соответствует падению напряжения U на ветви. Здесь g проводимость ветви (модуль упругости). Воспользовавшись данной аналогией, можем составить уравнения равновесия относительно координат относительно оси **x**.

$$\begin{cases} x_1\sum\limits_{p=1}^{n_1} g_p - \sum\limits_{p=1}^{n_1-k_1} g_p x_p = \sum\limits_{p=1}^{k_1} g_p x_p; \\ x_2\sum\limits_{p=1}^{n_2} g_p - \sum^{n_2-k_2} g_p x_p = \sum\limits_{p=1}^{k_2} g_p x_p; \\ \dots\dots\dots\dots\dots\dots\dots\dots \\ x_m\sum\limits_{p=1}^{n_m} g_p - \sum\limits_{p=1}^{n_m-k_m} g_p x_p = \sum\limits_{p=1}^{k_m} g_p x_p. \end{cases} \qquad (10.4)$$

И уравнения равновесия относительно координат относительно оси **y**.

$$\begin{cases} y_1\sum\limits_{p=1}^{n_1} g_p - \sum\limits_{p=1}^{n_1-k_1} g_p y_p = \sum\limits_{p=1}^{k_1} g_p y_p; \\ y_2\sum\limits_{p=1}^{n_2} g_p - \sum^{n_2-k_2} g_p y_p = \sum\limits_{p=1}^{k_2} g_p y_p; \\ \dots\dots\dots\dots\dots\dots\dots\dots \\ y_m\sum\limits_{p=1}^{n_m} g_p - \sum\limits_{p=1}^{n_m-k_m} g_p y_p = \sum\limits_{p=1}^{k_m} g_p y_p. \end{cases} \qquad (10.5)$$

Здесь $n_i$ – количество вершин смежных с вершиной *i*, $k_i$ – количество смежных вершин с вершиной *i* принадлежащих ободу.

Пусть задан топологический рисунок графа, представленный вращением вершин σ (см. рис. 10.1). Вращение вершин σ будем представлять в виде множества циклических подмножеств [3].

Например, для нашего рисунка графа (см. рис.10.1), множество подмножеств вращения вершин

$\sigma = \{\sigma(v_1) = \langle v_2,v_{31},v_5\rangle$, $\sigma(v_2) = \langle v_1,v_{10},v_3\rangle$, $\sigma(v_3) = \langle v_2,v_{11},v_4\rangle$, $\sigma(v_4) = \langle v_3,v_{16},v_5\rangle$,

$\sigma(v_5) = \langle v_1,v_4,v_{21}\rangle$, $\sigma(v_6) = \langle v_7,v_{11},v_{10}\rangle$, $\sigma(v_7) = \langle v_6,v_8,v_{12}\rangle$, $\sigma(v_8) = \langle v_7,v_{26},x_9\rangle$,

$\sigma(v_9) = \langle v_8,v_{30},v_{18}\rangle$, $\sigma(v_{10}) = \langle v_2,v_{28},v_{13},v_6\rangle$, $\sigma(v_{11}) = \langle v_3,v_6,v_{17}\rangle$, $\sigma(v_{12}) = \langle v_7,v_{13},v_{14}\rangle$,

$\sigma(v_{13}) = \langle v_{10},v_{15},v_{12}\rangle$, $\sigma(v_{14}) = \langle v_{12},v_{15},v_{31}\rangle$, $\sigma(v_{15}) = \langle v_{13},v_{27},v_{14}\rangle$, $\sigma(v_{16}) = \langle v_4,v_{17},v_{20}\rangle$,

$\sigma(v_{17}) = \langle v_{11},v_{18},v_{16}\rangle$, $\sigma(v_{18}) = \langle v_9,v_{19},v_{17}\rangle$, $\sigma(v_{19}) = \langle v_{18},v_{24},v_{23}\rangle$, $\sigma(v_{20}) = \langle v_{16},v_{23},v_{21}\rangle$,

$\sigma(v_{21}) = \langle v_5,v_{20},v_{22}\rangle$, $\sigma(v_{22}) = \langle v_{21},v_{23},v_{31}\rangle$, $\sigma(v_{23}) = (v_{19},v_{22},v_{20})$, $\sigma(v_{24}) = \langle v_{19},v_{30},v_{25}\rangle$,

$\sigma(v_{25}) = \langle v_{24},v_{26},v_{31}\rangle$, $\sigma(v_{26}) = \langle v_8,v_{25},v_{30}\rangle$, $\sigma(v_{27}) = \langle v_{15},v_{28},v_{29}\rangle$, $\sigma(v_{28}) = \langle v_{10},v_{29},v_{27}\rangle$,

$\sigma(v_{29}) = \langle v_{27}, v_{28}, v_{30} \rangle$, $\sigma(v_{30}) = \langle v_9, v_{26}, v_{24} \rangle$, $\sigma(v_{31}) = \langle v_1, v_{22}, v_{25}, v_{14}, v_{29} \rangle \}$.

Для решения системы линейных уравнений 10.4-10.5, применим метод LU-разложения. В качестве элементов правой части используются геометрические координаты вершин цикла $c_{19}$ равномерно расположенные на контуре квадрата размером 100×100 (см. рис. 10.25).

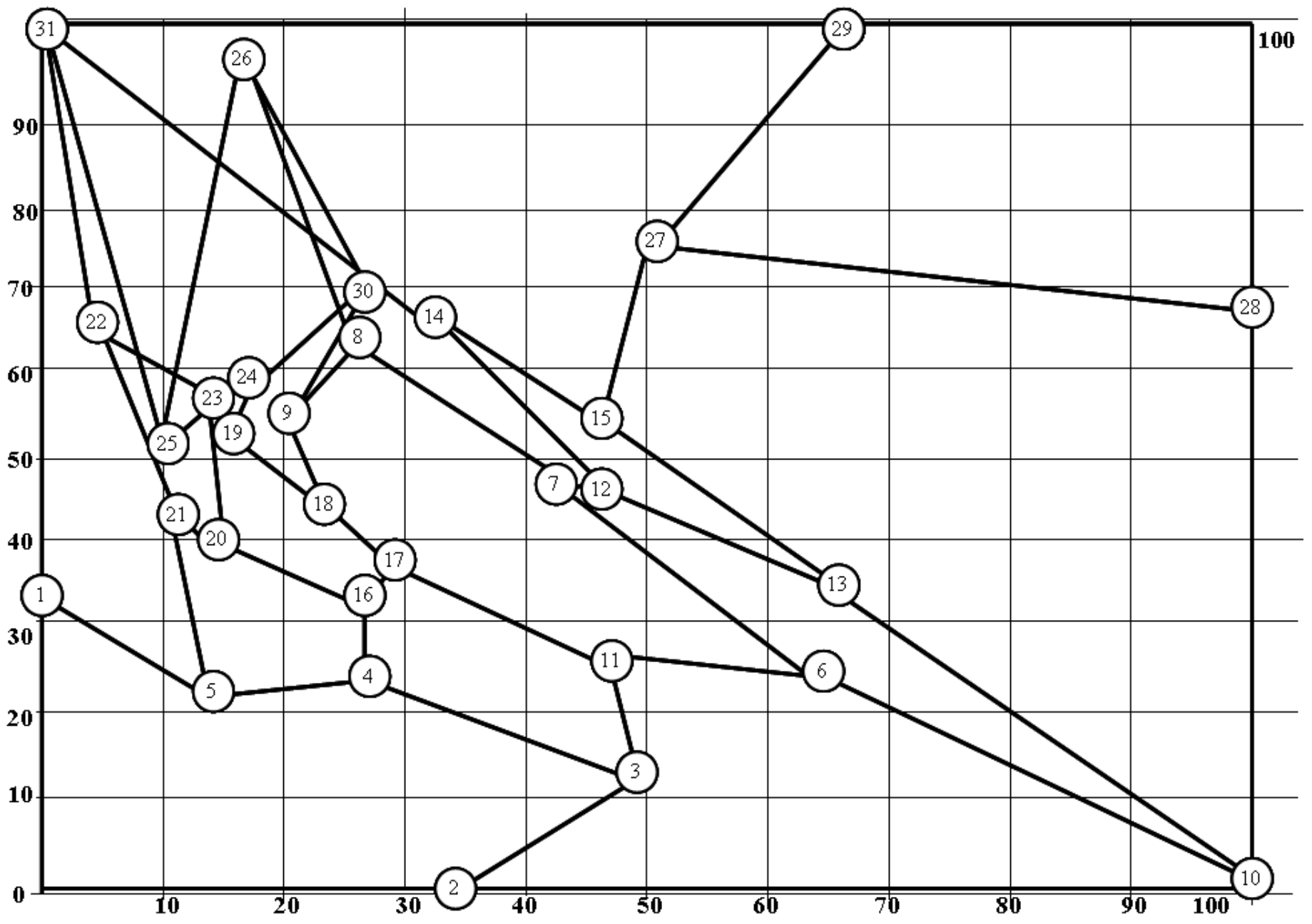

Рис. 10.25. Результат расчета рисунка графа силовым алгоритмом.

Недостатком данного подхода является то, что полученный рисунок как бы «стягивает» вершины, к центру неравномерно распределяя их на плоскости и пересекая соединения. Но, однако, имеется возможность корректировки решения.

С этой целью, заменим каждый цикл кликой графа (см. рис. 10.26).

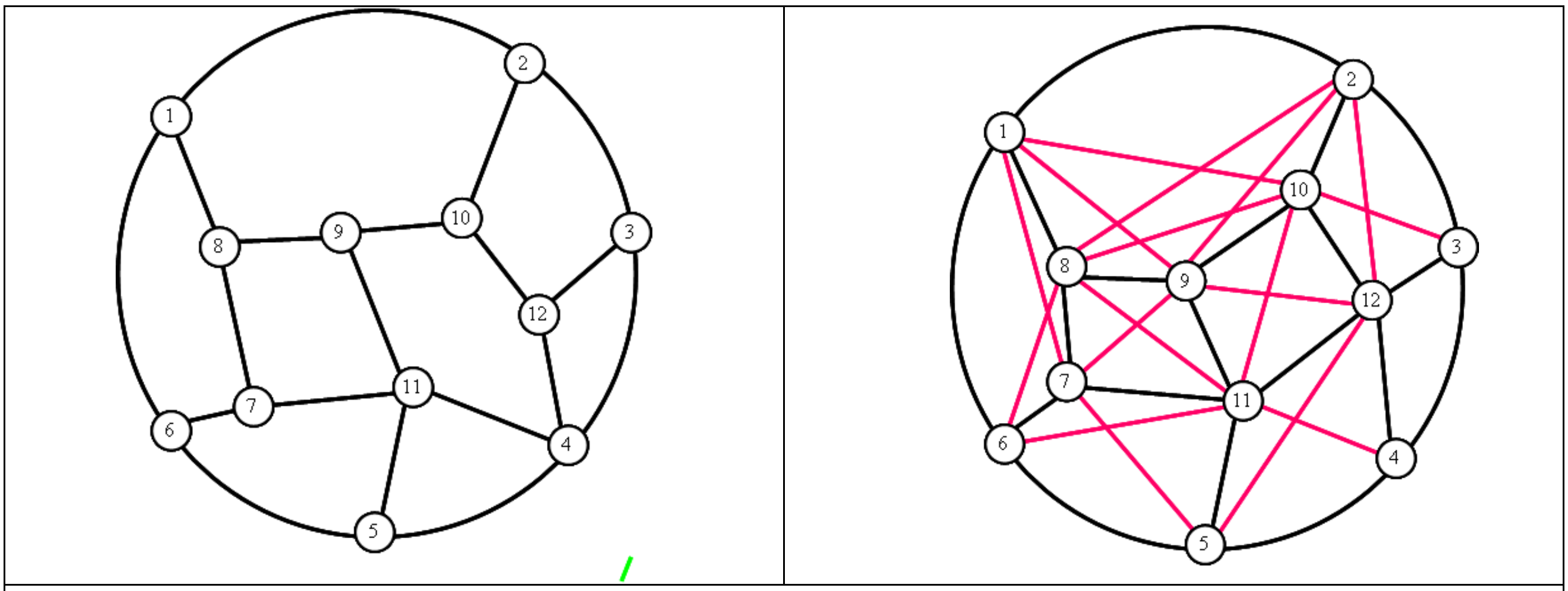

Рис. 10.25. Рисунок графа с представлением циклов кликами.

Рассмотрим итерационный процесс построения геометрического рисунка графа. В качестве итерации будем рассматривать уровневое расположение вершин. В качестве примера, рассмотрим следующий рисунок (см. рис. 10.27).

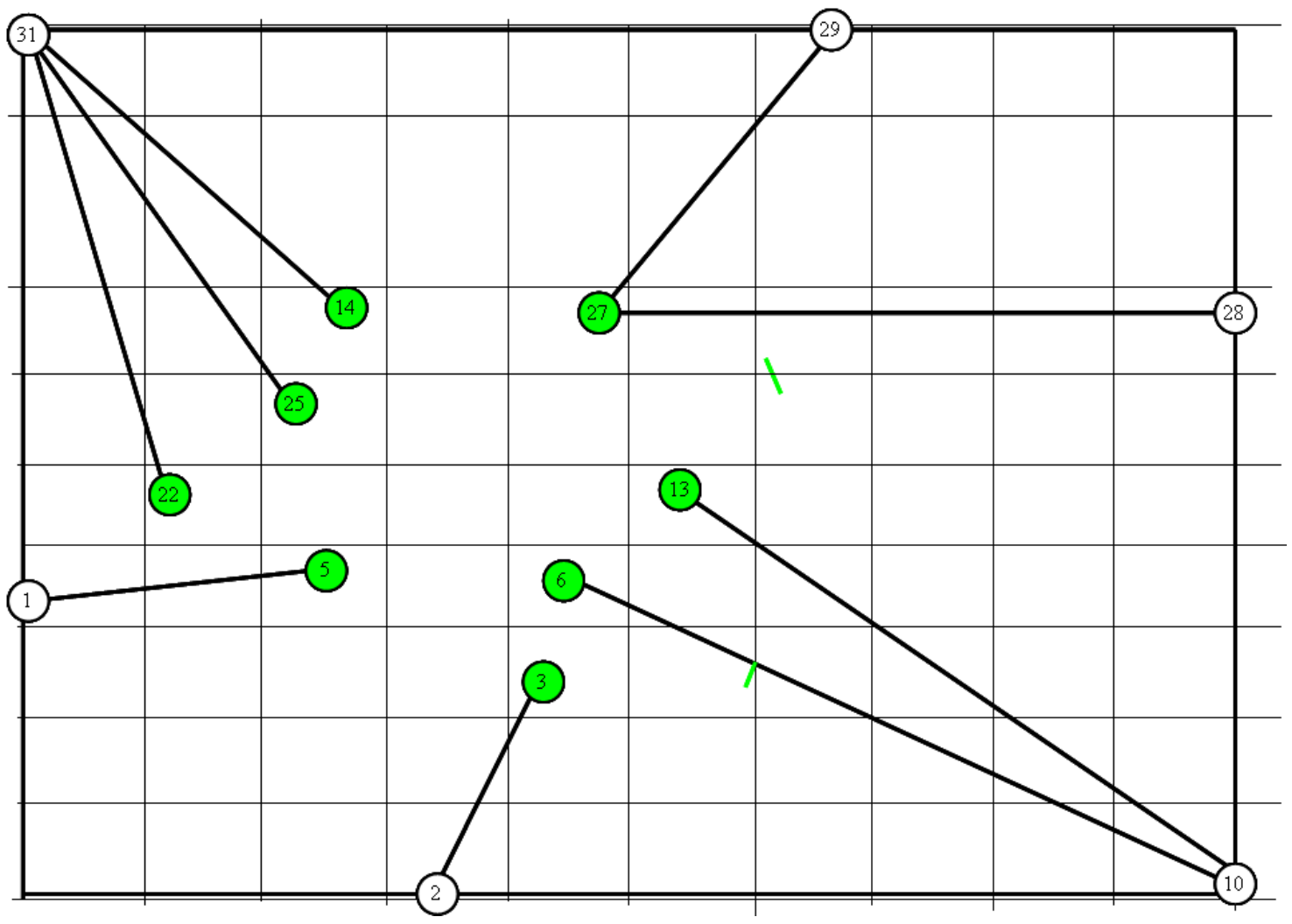

Рис. 10.27. Итерация расположения вершин 1-го и 2-го уровней.

Любая итерация предпологает вычисление координат всех вершин графа. Для первой итерации рассматриваются только вершины 1-го и 2-го уровней, без рассмотрения других вершин. Для второй итерации рассматриваются только расположения вершин первого, второго и третьего уровней. Учитывая явление стягивания вершин к центру, установим положение координат для следующей итерации, уменьшим длину каждого соединения в нестолько раз. На рис. 10.28 представлено сокращение длин соединений в два раза.

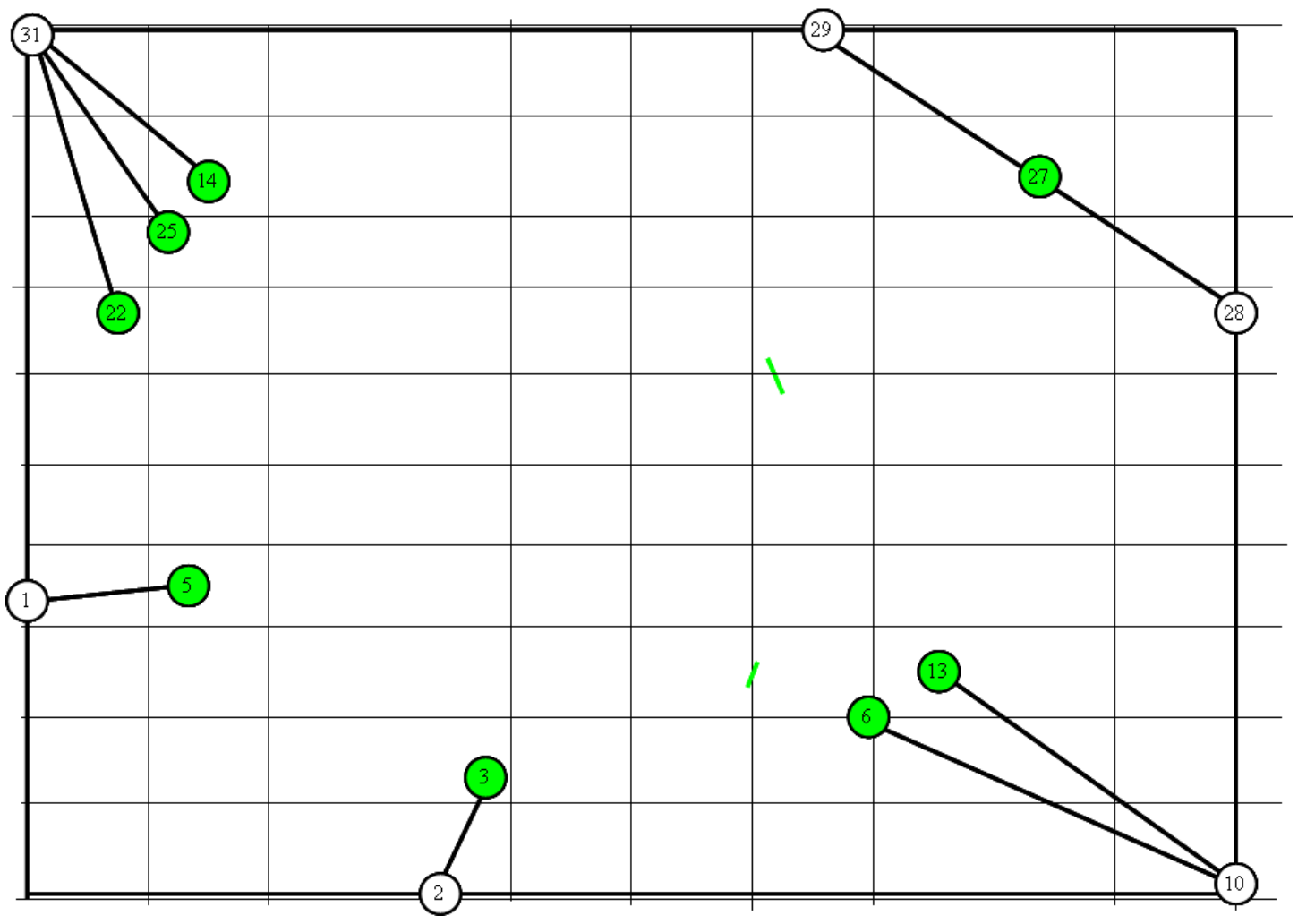

Рис. 10.28. Сокращение длины соединений для вершин 1-го и 2-го уровней.

После сокращения длин соединений вычисляются новые координаты вершин (см. рис. 10.28). Новые координаты вершин подставляются в правую часть уравнений 10.4-10.5. Результат следующей итерации для нового рассположения вершин 2-го уровня представлены на рис. 10.29.

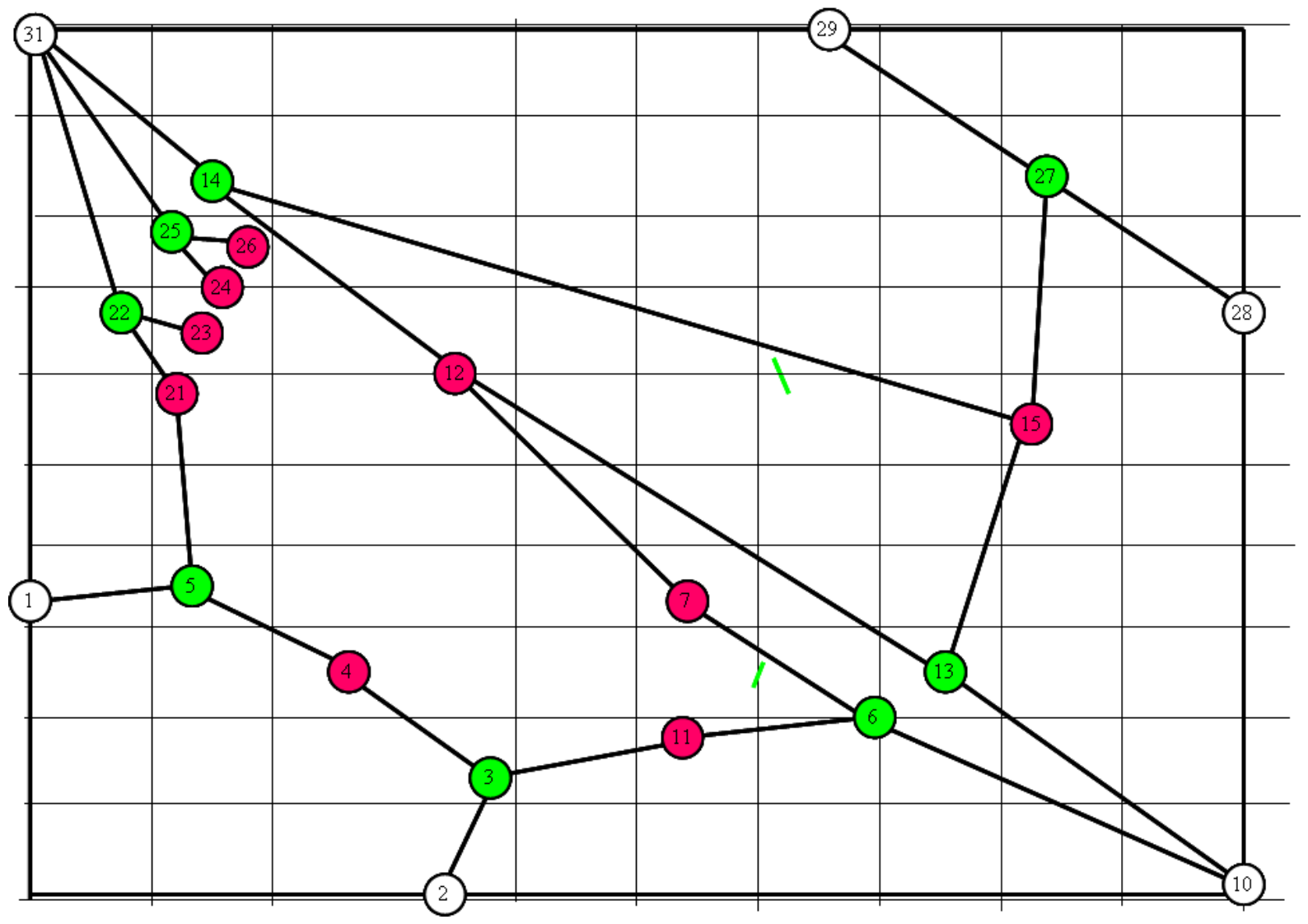

Рис. 10.29. Расположение вершин 1-го, 2-го и 3-го уровней.

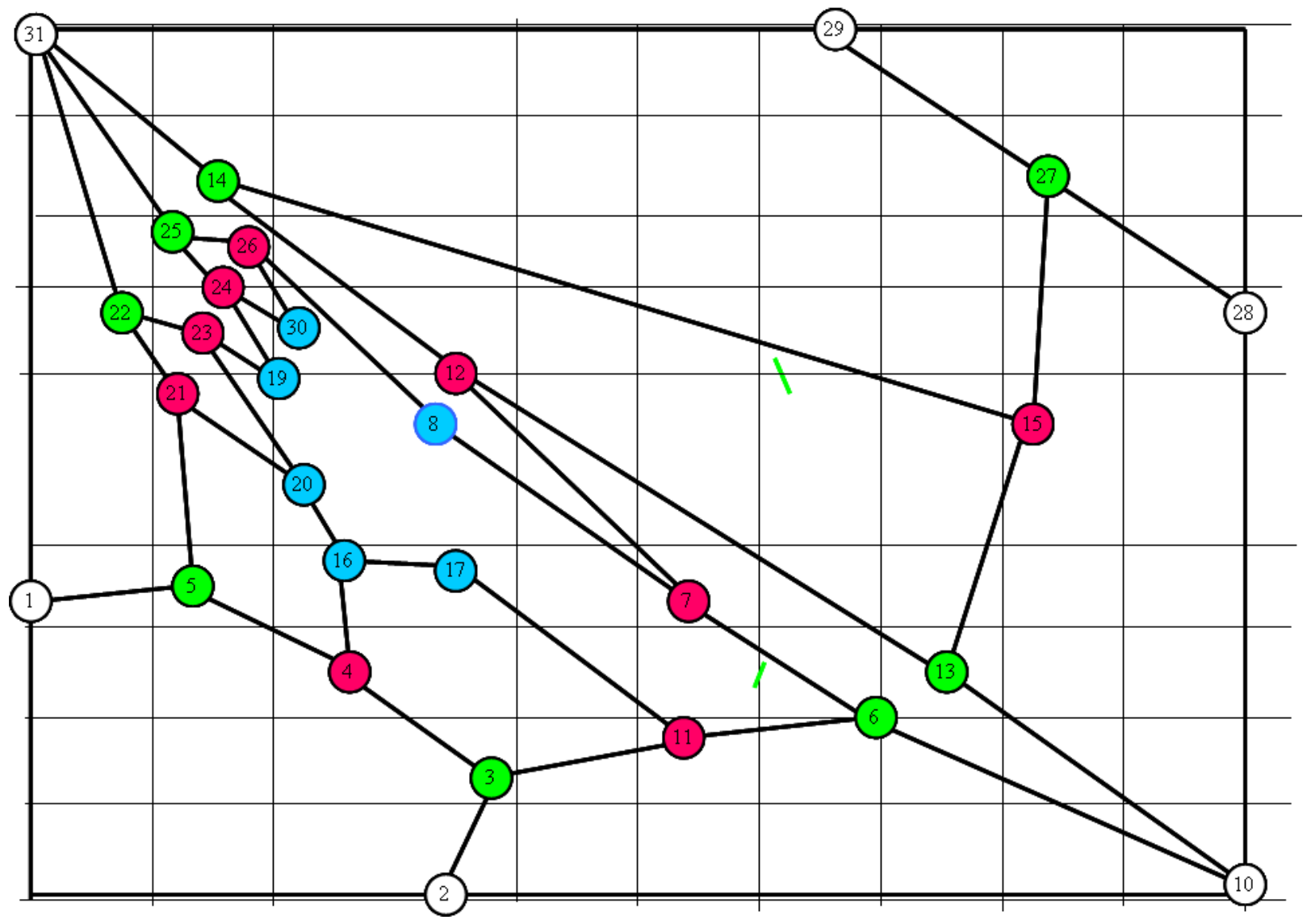

Рис. 10.30. Расположение вершин 1-го, 2-го,3-го и 4-го уровней.

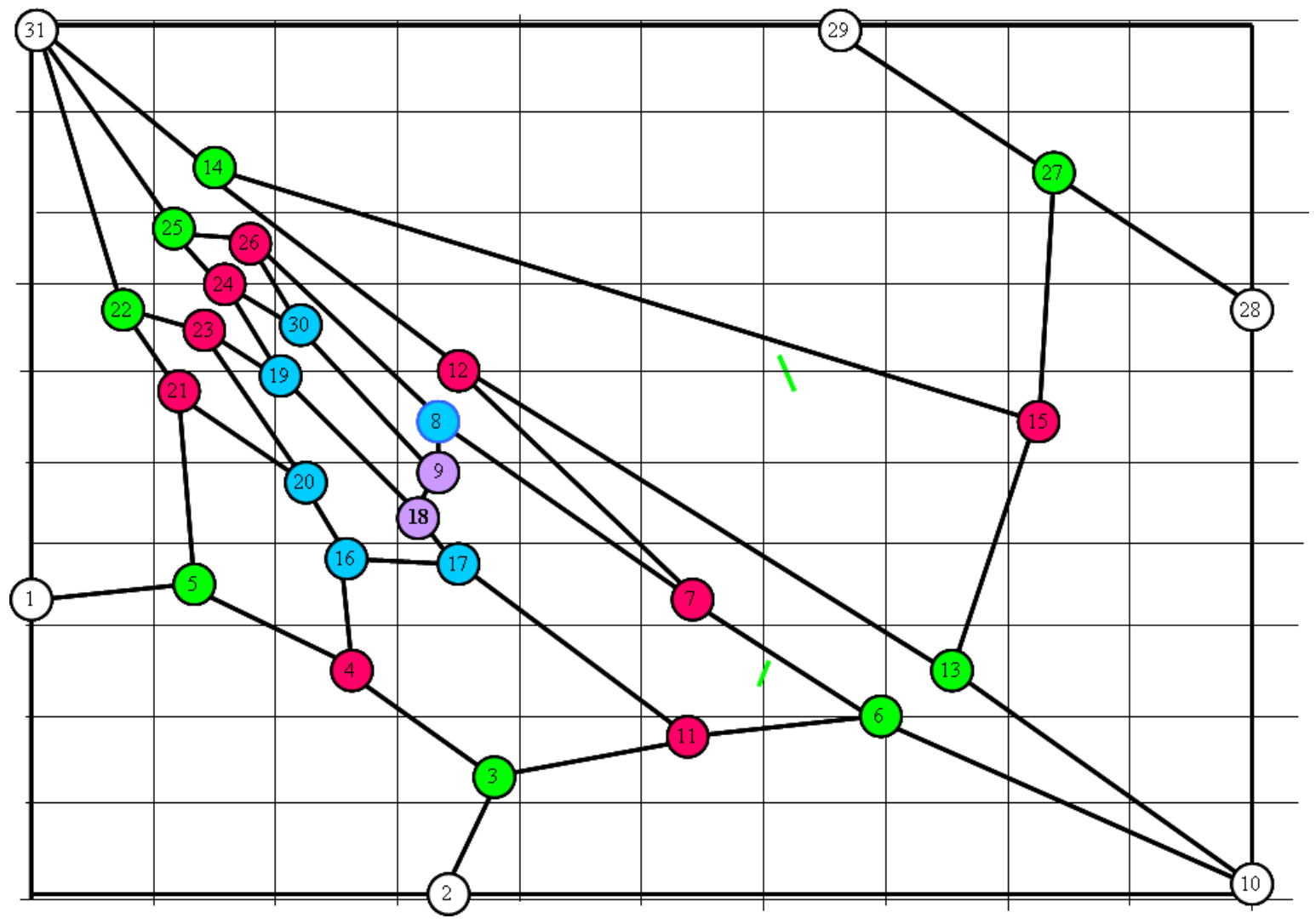

Рис. 10.31. Окончательный вариант расположение вершин плоского графа.

На рисунке 10.29 представлен результат расчета координат для окончания второй итерации. На рисунке 10.30 представлен результат расчета координат для окончания третьей итерации. Окончательный растет представлен на рис. 10.31.

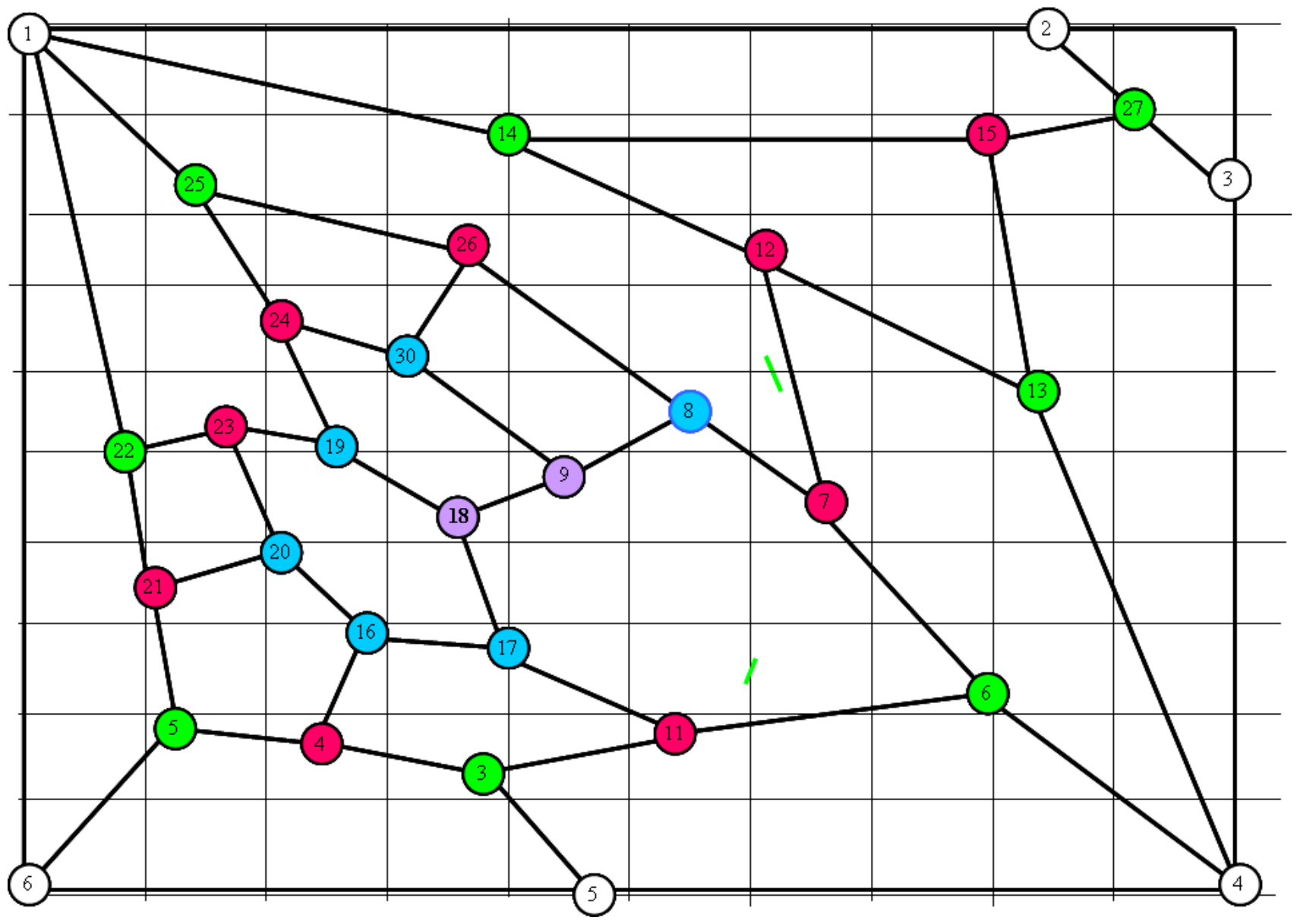

Рис. 10.32. Построение рисунка плоского графа методом равных площадей.

Метод равных площадей предполагает итерационный процесс, размещения вершин в отдельных циклах. Итерационный процесс должет сохранять вращение вершин плоского графа с целью выравнивания площадей циклов.

## 10.4. Программа Raschet9

```
program Raschet9;

type    TMasy = array[1..1000] of integer;
        TMass = array[1..4000] of integer;
        TMasyy = array[1..100] of integer;
        RMasR = array[1..500] of real;
        RMasyy = array[1..100] of real;
 var
   F1,F2 : text;
   Masy: TMasy;
   Mass: TMass;
   MasR: TMass;
   Masi : TMass;
   Kontr: TMasyy;
   Xkontr : RMasyy;
   Ykontr : RMasyy;
   Num : TMasy;
   Most : TMasy;
   Masy1 : TMasy;
   Mass1 : TMass;
   MasR1 : TMass;
   Num1 : TMasy;
   Nu : TMasy;
   Zyd : RMasR;
   Ztx : RMasR;
   Zty : RMasR;
   Zy : RMasR;
```

```
    Zyl : RMasR;
    ZnaprX : RMasR;
    ZnaprY : RMasR;
    KoorX : TMasy;
    KoorY : TMasy;
    KoorXN : TMasy;
    KoorYN : TMasy;
    K1,II,JJ,J,III,Nkontr,Na : Integer;
    MaxX,MaxY,X1,Y1,Np,M,PP,MMM,K11 : Integer;
    XX, YY : real;
    I,M1,K,K2,K22,Nv,Nk,K3,Siz,KKK,KK1,KK2 : integer;
{*****************************************************************}
procedure FormIncide(var Nv : integer;
                     var My : TMasy;
                     var Ms : TMass;
                     var Ms3 : TMass);
{ Nv - количество вершин в графе;                          }
{ My  - массив указателей для матрицы смежностей;          }
{ Ms  - массив элементов матрицы смежностей;               }
{ Ms3 - массив элементов матрицы инциденций.               }
{    Формируется матрица инциденций графа в массиве        }
{    Ms3.                                                  }
{                                                          }
var I,J,K,NNN,P1,M,L,L1 : integer;
label 1;
begin
{      инициализация                                       }
  NNN:=My[Nv+1]-1;
  L1:=NNN div 2;
  K:=0;
  for J:= 1 to L1 do Ms3[J]:= 0;
{   определение номера элемента                            }
  for I:= 1 to Nv do
  begin
   for M:= My[I] to My[I+1]-1 do
    if Ms3[M]=0 then
    begin
      P1:=Ms[M];
      K:=K+1;
      {writeln(F2,' K = ',K:3,' P1 = ',P1:3);}
      Ms3[M]:=K;
       for L:=My[P1] to My[P1+1]-1 do
       begin
       {writeln(F2,' L = ',L:3,' P1 = ',P1:3,' Ms[L] = ',Ms[l]:3);}
        if (Ms[L]=I) and (Ms3[L]=0) then
        begin
          Ms3[L]:=K;
          {writeln(F2,' L = ',L:3,' Ms[L] = ',Ms[L]:3,' Ms3[L]= ',Ms3[L]:3);}
          goto 1;
        end;
      end;
1:   end;
  end;
end; {FormIncide}
{*****************************************************************}
procedure  Shell2(var N : integer;
                  var A : TMasy;
                  var B : TMasy);
{       процедура Шелла для упорядочивания элементов        }
{                                                           }
{     N - количество элементов в массиве;                   }
{     A - сортируемый массив;                               }
{     B - перестаавляемый массив.                           }
```

```
var D,Nd,I,J,L,X,Y : integer;
label 1,2,3,4,5;
begin
  D:=1;
1:D:=2*D;
  if D<=N then goto 1;
2:D:=D-1;
  D:=D div 2;
  if D=0 then goto 5;
  Nd:=N-D;
  for I:=1 to Nd do
  begin
    J:=I;
3:  L:=J+D;
    if A[L]>=A[J] then goto 4;
    X:=A[J];
    Y:=B[J];
    A[J]:=A[L];
    B[J]:=B[L];
    A[L]:=X;
    B[L]:=Y;
    J:=J-D;
    if J>0 then goto 3;
4:end;
  goto 2;
5:end;{Shell2}
{**********************************************************}
procedure FormObr(var Nmyk,Nkontr : integer;
                  var  Myk : TMasy;
                  var  Mas : TMass;
                  var MasR : TMass;
                  var Kontr : TMasyy;
                  var Xkontr : RMasyy;
                  var Ykonyr : RMasyy;
                  var Num1 : TMasy;
                  var Zyd : RMasR;
                  var Ztx : RMasR;
                  var Zty : RMasR);
{      процедура создания массивов Zyd,Ztx,Zty             }
{                                                  }
{  Nmyk -  количество элементов в массиве Myk;              }
{  Myk - новый массив указаателей для матрицы смежностей;    }
{  Mas - новый массив элементов матрицы смежностей;          }
{  MasR - массив коэффициентов упругости;                    }
{  Nkontr - количество внешних вершин;                       }
{  Xkontr - координаты X внешних верших;                     }
{  Ykontr - координаты Y внешних вершин;                     }
{  Kontr - массив номеров внешних вершин;                    }
{  Zy - эквивалент матрицы проводимостей без диагональных    }
{      элементов;                                     }
{  Zyl - эквивалент матрицы проводимостей без диагональных   }
{       элементов;                                    }
{  Zyd - элементы главной диагонали;                        }
{  Ztx - массив правых частей для оси абцисс;               }
{  Zty - массив правых частей для оси ординат;              }
var Nvers,Na,I,I1,K,Ha,Ko,I2,K1,K2 : integer;
label 1,4;
begin
  Nvers:=Nmyk-1;
  Na:=Nvers-Nkontr;
{           инициализация массивов Zyd,Ztx,Zty             }
  for I:=1 to Na do
  begin
```

```
  Zyd[I]:=0.0;
  Ztx[I]:=0.0;
  Zty[I]:=0.0;
 end;
 for I1:=1 to Nvers do
 begin
  K:=Num1[I1];
  if K>Na then goto 1;
  Ha:=Myk[I1];
  Ko:=Myk[I1+1]-1;
  for I2:=Ha to Ko do
   Zyd[K]:=Zyd[K]+MasR[I2];
1:end;
{       формирование массивов Zyd,Ztx,Zty                 }
  for I1:=1 to Nkontr do
  begin
   K:=Kontr[I1];
   Ha:=Myk[K];
   Ko:=Myk[K+1]-1;
   for I2:=Ha to Ko do
   begin
    K1:=Mas[I2];
    K2:=Num1[K1];
    if K2>Na then goto 4;
    Ztx[K2]:=Ztx[K2]+MasR[I2]*Xkontr[I1];
    Zty[K2]:=Zty[K2]+MasR[I2]*Ykontr[I1];
4:   end;
  end;
end;{FormObr}
{**************************************************************}
procedure FormRetsen(var Nmyk1 : integer;
                     var  Myk1 : TMasy;
                     var  Mas1 : TMass;
                     var Zyd : RMasR;
                     var Ztx : RMasR;
                     var Zty : RMasR;
                     var Zy : RMasR;
                     var Zyl : RMasR;
                     var  ZnaprX : RMasR;
                     var ZnaprY : RMasR);
{   решение линейной системы уравнений методом Гаусса       }
{                                                  }
{  Nmyk1 -  количество элементов в массиве Myk1;            }
{  Myk1 - новый массив указаателей для матрицы смежностей;   }
{  Mas1 - новый массив элементов матрицы смежностей;         }
{  Zy - эквивалент матрицы проводимостей без диагональных    }
{      элементов;                                  }
{  Zyl - эквивалент матрицы проводимостей без диагональных   }
{       элементов;                                 }
{  Zyd - элементы главной диагонали;                    }
{  Ztx - массив правых частей для оси абцисс;                }
{  Zty - массив правых частей для оси ординат;               }
{  ZnaprX - массив координат для внутренних точек;           }
{  ZnaprY - массив координат для внутренних точек;           }
var N1,Na,I,J1,J2,Jl,Koh,K1,K2,Jc : integer;
var K,L1,L,J : integer;
var Zcx,Zcy : real;
label 44,47,48,49,50,51,54,53,1,56;
begin
 N1:=Nmyk1-1;
 Na:=Nmyk1;
 for I:=1 to N1 do
 begin
```

```
    J1:=Myk1[I];
    J2:=Myk1[I+1]-1;
    if J1>J2 then goto 44;
    for Jl:=J1 to J2 do
      Zyl[Jl]:=Zyl[Jl]/Zyd[I];
    for J:=J1 to J2 do
    begin
      Koh:=Mas1[J];
      K1:=Myk1[Koh];
      K2:=Myk1[Koh+1]-1;
      for Jc:=J1 to J2 do
      begin
        if Koh>Mas1[Jc] then goto 47;
        if Koh=Mas1[Jc] then goto 48;
        if K1>K2 then goto 47;
        for K:=K1 to K2 do
        begin
          if Mas1[K]<>Mas1[Jc] then goto 49;
          Zy[K]:=Zy[K]-Zy[Jc]*Zyl[J];
          Zyl[K]:=Zyl[K]-Zyl[Jc]*Zy[J];
49:       end;
        goto 47;
48:       Zyd[Koh]:=Zyd[Koh]-Zy[Jc]*Zyl[Jc];
47:     end;
    end;
44:end;
  for I:=1 to N1 do
  begin
    J1:=Myk1[I];
    J2:=Myk1[I+1]-1;
    L1:=I+1;
    if J1>J2 then goto 50;
    for L:=L1 to Na do
    begin
      for J:=J1 to J2 do
      begin
        if L<>Mas1[J] then goto 51;
        begin
          Ztx[L]:=Ztx[L]-Zyl[J]*Ztx[I];
          Zty[L]:=Zty[L]-Zyl[J]*Zty[I];
        end;
51:     end;
    end;
50:end;
   ZnaprX[Na]:=Ztx[Na]/Zyd[Na];
   ZnaprY[Na]:=Zty[Na]/Zyd[Na];
{                обратная прогонка                          }
   for I:=1 to N1 do
   begin
     K:=Na-I;
     J1:=Myk1[K];
     J2:=Myk1[K+1]-1;
     Zcx:=0.0;
     Zcy:=0.0;
     if J1>J2 then goto 53;
     K1:=K+1;
     for L:=K1 to Na do
     begin
       for J:=J1 to J2 do
       begin
         if Mas1[J]<>L then goto 54;
         Zcx:=Zcx+Zy[J]*ZnaprX[L];
         Zcy:=Zcy+Zy[J]*ZnaprY[L];
```

```
54:   end;
    end;
53:  if Zyd[K]=0.0 then goto 1;
    ZnaprX[K]:=(Ztx[K]-Zcx)/Zyd[K];
    ZnaprY[K]:=(Zty[K]-Zcy)/Zyd[K];
    goto 56;
1:   ZnaprX[K]:=0.0;
    ZnaprY[K]:=0.0;
56:end;
end;{FormRetsen}
{***********************************************************}
procedure FormYk(var Nmyk : integer;
                 var Myk1 : TMasy;
                 var Mas1 : TMass;
                 var MasR1 : TMass);
{  Nmyk -  количество элементов в массиве Myk;               }
{  Myk1 - новый массив указаателей для матрицы смежностей;   }
{  Mas1 - новый массив элементов матрицы смежностей;         }
{  MasR1 - новый массив упругих связей;                      }
var Na,I1,J1,J2,I2,Koh,K1,K2,I3,I4,I5,I6 : integer;
var Met,Myh,Ns,Kn,Kr1,Na1 : integer;
label 1,24,26,27;
begin
 Na:=Nmyk;
 for I1:=1 to Na do
 begin
  J1:=Myk1[I1];
  J2:=Myk1[I1+1]-1;
  for I2:=J1 to J2 do
  begin
   Koh:=Mas1[I2];
   K1:=Myk1[Koh];
   K2:=Myk1[Koh+1]+1;
   for I3:=J1 to J2 do
   begin
    if Koh>=Mas1[I3] then goto 26;
    Met:=0;
    if K1>K2 then goto 27;
    for I4:=K1 to K2 do
     if Mas1[I3]=Mas1[I4] then
       Met:=Met+1;
27:   if Met<>0 then goto 26;
     Myh:=Myk1[Na];
     Ns:=Myk1[Koh];
     for I5:=Ns to Myh do
     begin
       Kn:=Myh-I5+Ns;
       Mas1[Kn+1]:=Mas1[Kn];
       MasR1[Kn+1]:=MasR1[Kn];
     end;
     Mas1[Ns]:=Mas1[I3];
     MasR1[Ns]:=0;
     Kr1:=Koh+1;
     Na1:=Na+1;
     for I6:=Kr1 to Na1 do
       Myk1[I6]:=Myk1[I6]+1;
     K2:=K2+1;
26:   end;
   end;
24:end;
end;{FormYk}
{***********************************************************}
procedure FormPostr(var Nmyk,Nkontr : integer;
```

```
          var Myk : TMasy;
          var Mas : TMass;
          var MasR : TMass;
          var Myk1 : TMasy;
          var Mas1 : TMass;
          var MasR1 : TMass;
          var Num : TMasy;
          var Num1 : TMasy;
          var Nu : TMasy;
          var Kontr : TMasyy);
{  построение новой мматрицы смежностей и матрицы упругих    }
{   связей для новой нумерации узлов                          }
{                                                  }
{  Nmyk -  количество элементов в массиве Myk;               }
{  Nkontr - количество элемментов во внешнем контуре;        }
{  Myk - массив указаателей для матрицы смежностей;          }
{  Mas -  массив элементов матрицы смежностей;               }
{  MasR - массив упругих связей;                             }
{  Myk1 - новый массив указаателей для матрицы смежностей;   }
{  Mas1 - новый массив элементов матрицы смежностей;         }
{  MasR1 - новый массив упругих связей;                      }
{  Num -  массив нумерации вершин;                           }
{  Nu -   массив нумерации вершин;                           }
{  Num1 -  массин нумерации новых вершин;                    }
{  Kontr - массив внешних вершин.                            }
var Nvers,Na,I,Bers,H,K,Hhob,K1,Khob : integer;
label 2,3,10;
begin
  Nvers:=Nmyk-1;
  Na:=Nvers-Nkontr;
  for I:=1 to Na do
  begin
    Nu[I]:=Num[I];
    Num1[I]:=I;
  end;
  for I:=1 to Nkontr do
  begin
    Nu[Na+I]:=Kontr[I];
    Num1[Na+i]:=Na+I;
  end;
  Shell2(Nvers,Nu,Num1);
  Na:=Nmyk-1-Nkontr;
  Myk1[1]:=1;
  for I:=1 to Na do
  begin
    Bers:=Num[I];
    H:=Myk[Bers]-1;
    K:=Myk[Bers+1]-1;
    Myk1[I+1]:=Myk1[I];
    if H>=K then goto 10;
    Hhob:=Myk1[I]-1;
2:  H:=H+1;
    K1:=Mas[H];
    Khob:=Num1[K1];
    if (Khob>Na) or (Khob<=I) then goto 3;
    Myk1[I+1]:=Myk1[I+1]+1;
    Hhob:=Hhob+1;
    Mas1[Hhob]:=Khob;
    MasR1[Hhob]:=MasR[H];
3:  if H<K then goto 2;
10:end;
end;{FormPostr}
{*************************************************************}
```

```
procedure FormPernum(var Nmyk,Nkontr : integer;
                   var Myk : TMasy;
                   var Mas : TMass;
                   var Kontr : TMasyy;
                   var Num : TMasy;
                   var Most : TMasy);
{           перенумеровка вершин                    }
var Na,Nvers,K,I,I1,El,Elh,Elk,M,I2 : integer;
label 1,4;
begin
  Na:=Nmyk-1-Nkontr;
  Nvers:=Nmyk-1;
  K:=0;
  for I:=1 to Nvers do
  begin
    for I1:=1 to Nkontr do
     if Kontr[I1]=I then goto 1;
    K:=K+1;
    Num[K]:=I;
1:end;
  for I:=1 to Na do
  begin
    El:=Num[I];
    Elh:=Myk[El];
    Elk:=Myk[El+1]-1;
    M:=0;
    for I1:=Elh  to Elk do
    begin
      for I2:=1 to Nkontr do
        if Kontr[I2]=Mas[I1] then
         goto 4;
      M:=M+1;
4:  end;
    Most[I]:=M;
  end;
  Shell2(Na,Most,Num);
end;{FormPernum}
{*************************************************************}
procedure FormMonitor(var Nv,Nkontr : integer;
                   var Masy: TMasy;
                   var Mass: TMass;
                   var MasR: TMass;
                   var Kontr: TMasyy;
                   var Xkontr : RMasyy;
                   var Ykontr : RMasyy;
                   var Num : TMasy;
                   var  Most : TMasy;
                   var  Masy1 : TMasy;
                   var  Mass1 : TMass;
                   var  MasR1 : TMass;
                   var Num1 : TMasy;
                   var  Nu : TMasy;
                   var Zyd : RMasR;
                   var Ztx : RMasR;
                   var Zty : RMasR;
                   var Zy : RMasR;
                   var Zyl : RMasR;
                   var  ZnaprX : RMasR;
                   var ZnaprY : RMasR);
{            Геометрическое размещение элементов            }
var I,Nmyk,Na,Nvers,N1 : integer;
begin
   for I:=1 to Masy[Nv+1]-1 do
```

```
   MasR[I]:=1;
  Nmyk:=Nv+1;
  FormPernum(Nmyk,Nkontr,Masy,Mass,Kontr,
    Num,Most);
  Na:=Nmyk-1-Nkontr;
  Nvers:=Nmyk-1;
  FormPostr(Nmyk,Nkontr,Masy,Mass,MasR,Masy1,Mass1,
    MasR1,Num,Num1,Nu,Kontr);
  N1:=Masy1[Na]-1;
  FormYk(Na,Masy1,Mass1,MasR1);
  N1:=Masy1[Na]-1;
  FormObr(Nmyk,Nkontr,Masy,Mass,MasR,Kontr,
     Xkontr,Ykontr,Num1,Zyd,Ztx,Zty);
  for I:=1 to N1 do
  begin
   Zy[I]:=-MasR1[I];
   Zyl[I]:=-MasR1[I];
  end;
  FormRetsen(Na,Masy1,Mass1,Zyd,Ztx,Zty,Zy,Zyl,
         ZnaprX,ZnaprY);
end;{FormMonitor}
{************************************************************}
label 1,2,3,4,5,6,12,13,21,22,23,31,32,33,44;
begin
  assign(F1,'D:\Delphi\GM1\30+1b.gm1');
  reset(F1);
  readln(F1,Nv);
  for I:=1 to Nv+1 do
  begin
   if I<>Nv+1 then read(F1,Masy[I]);
   if I=Nv+1 then readln(F1,Masy[I]);
  end;
  Np:=Masy[Nv+1]-1;
  Assign(F2,'D:\Delphi\GM2\30+1b.gm2');
  Rewrite(F2);
  writeln(F2,' Количество элементов = ',Np);
  for I:=1 to Np do
  begin
   if I<>Np then read(F1,Mass[I]);
   writeln(F2,' I = ',I,', Mass[',I,'] = ',Mass[I]);
   if I=Np then readln(F1,Mass[I]);
  end;
  readln(F1,Nkontr);
  writeln(F2,' Nkontr = ',Nkontr);
  for I:=1 to Nkontr do
  begin
   if I<>Nkontr then read(F1,Kontr[I]);
   if I=Nkontr then readln(F1,Kontr[I]);
  end;
  for I:=1 to Nkontr do
  begin
   if I<>Nkontr then read(F1,XKontr[I]);
   if I=Nkontr then readln(F1,XKontr[I]);
  end;
  for I:=1 to Nkontr do
  begin
   if I<>Nkontr then read(F1,YKontr[I]);
   if I=Nkontr then readln(F1,YKontr[I]);
  end;
  close(F1);
  FormMonitor(Nv,Nkontr,Masy,Mass,MasR,Kontr,
  Xkontr,Ykontr,Num,Most,Masy1,Mass1,MasR1,Num1,
  Nu,Zyd,Ztx,Zty,Zy,Zyl,ZnaprX,ZnaprY);
```

```
Na:=Nv-Nkontr;
FormIncide(Nv,Masy,Mass,Masi);
writeln(F2,' ');
writeln(F2,'  Матрица смежностей графа: ');
writeln(F2,' ');
for I:=1 to Nv do
begin
   write(F2,'вершина ',I:3,':');
   Siz:=0;
   for K:=Masy[I] to Masy[I+1]-1 do
   begin
      Siz:=Siz+1;
      Np:=Siz mod 10;
      Nk:=Siz div 10;
      if K<>Masy[I+1]-1 then
      begin
       if (Np=1) and (Nk>0) then
       write(F2,'              ',Mass[K]:3);
       if (Np=1) and (Nk=0) then
       write(F2,' ',Mass[K]:3);
       if (Np>1) and (Np<=9) then
       write(F2,' ',Mass[K]:3);
       if Np=0 then
        writeln(F2,' ',Mass[K]:3);
     end;
     if K=Masy[I+1]-1 then
     begin
       if (Np=1) and (Nk>0) then
        writeln(F2,'              ',Mass[K]:3)
       else
      writeln(F2,' ',Mass[K]:3);
     end;
   end;
end;
 writeln(F2,' ');
 writeln(F2,'  Элементы матрицы инциденций: ');
 writeln(F2,' ');
 M1:=Masy[Nv+1]-1;
 M:=M1 div 2;
 writeln(F2,'M1 = ',M1:3);
 writeln(F2,'M = ',M:3);
  for I:=1 to Nv do
     begin
      for K:=Masy[I] to Masy[I+1]-1 do
      begin
        if K<>Masy[I+1]-1 then write(F2,Masi[K]:3,' ');
        if K=Masy[I+1]-1 then writeln(F2,Masi[K]:3);
      end;
     end;
 for I:=1 to M do
 begin
   for K:=1 to Nv do
   begin
     for M1:=Masy[K] to Masy[K+1]-1 do
     begin
       if Masi[M1]=I then
       begin
         write(F2,' ребро ',I:3,':');
         Np:=Mass[M1];
         write(F2,'  ( ',K:3);
         write(F2,' ',Np:3,' )  или  ');
         write(F2,' ( ',Np:3);
         writeln(F2,'  ',K:3,' )');
```

```
          goto 1;
          end;
        end;
      end;
1:   end;
   writeln(F2,' ');
   writeln(F2,' Количество точек в ободе: ',Nkontr);
   writeln(F2,' ');
   writeln(F2,'  Массив точек входящих в обод графа: ');
   writeln(F2,' ');
   for I:=1 to Nkontr do
   begin
     Siz:=I;
     Np:=Siz mod 10;
     Nk:=Siz div 10;
    if I<>Nkontr then
    begin
     if (Np=1) and (Nk>0) then
      write(F2,' ',Kontr[I]:3);
     if (Np=1) and (Nk=0) then
      write(F2,' ',Kontr[I]:3);
     if (Np>1) and (Np<=9) then
      write(F2,' ',Kontr[I]:3);
     if Np=0 then
       writeln(F2,' ',Kontr[I]:3);
    end;
    if I=Nkontr then
     begin
         if (Np=1) and (Nk>0) then
          writeln(F2,' ',Kontr[I]:3)
         else
        writeln(F2,' ',Kontr[I]:3);
     end;
   end;
    writeln(F2,' ');
    writeln(F2,'  Координаты точек по оси абсцисс в ободе: ');
    writeln(F2,' ');
   for I:=1 to Nkontr do
   begin
     Siz:=I;
     Np:=Siz mod 10;
     Nk:=Siz div 10;
    if I<>Nkontr then
    begin
     if (Np=1) and (Nk>0) then
      write(F2,' ',XKontr[I]:6:3);
     if (Np=1) and (Nk=0) then
      write(F2,' ',XKontr[I]:6:3);
     if (Np>1) and (Np<=9) then
      write(F2,' ',XKontr[I]:6:3);
     if Np=0 then
       writeln(F2,' ',XKontr[I]:6:3);
    end;
    if I=Nkontr then
     begin
         if (Np=1) and (Nk>0) then
          writeln(F2,' ',XKontr[I]:6:3)
         else
        writeln(F2,' ',XKontr[I]:6:3);
     end;
   end;
    writeln(F2,' ');
    writeln(F2,'  Координаты точек по оси ординат в ободе: ');
```

```
 writeln(F2,' ');
for I:=1 to Nkontr do
begin
   Siz:=I;
   Np:=Siz mod 10;
   Nk:=Siz div 10;
  if I<>Nkontr then
  begin
   if (Np=1) and (Nk>0) then
    write(F2,' ',YKontr[I]:6:3);
   if (Np=1) and (Nk=0) then
    write(F2,' ',YKontr[I]:6:3);
   if (Np>1) and (Np<=9) then
    write(F2,' ',YKontr[I]:6:3);
   if Np=0 then
     writeln(F2,' ',YKontr[I]:6:3);
  end;
  if I=Nkontr then
   begin
        if (Np=1) and (Nk>0) then
         writeln(F2,' ',YKontr[I]:6:3)
        else
       writeln(F2,' ',YKontr[I]:6:3);
   end;
end;
writeln(F2,' ');
writeln(F2,'  Массив точек не входящих в обод графа: ');
writeln(F2,' ');
for I:=1 to Na do
begin
   Siz:=I;
   Np:=Siz mod 10;
   Nk:=Siz div 10;
  if I<>Na then
  begin
   if (Np=1) and (Nk>0) then
    write(F2,' ',Num[I]:3);
   if (Np=1) and (Nk=0) then
    write(F2,' ',Num[I]:3);
   if (Np>1) and (Np<=9) then
    write(F2,' ',Num[I]:3);
   if Np=0 then
     writeln(F2,' ',Num[I]:3);
  end;
  if I=Na then
   begin
        if (Np=1) and (Nk>0) then
         writeln(F2,' ',Num[I]:3)
        else
       writeln(F2,' ',Num[I]:3);
   end;
end;
writeln(F2,' ');
writeln(F2,
'  Массив точек по оси абсцисс не входящих в обод графа: ');
writeln(F2,' ');
for I:=1 to Na do
begin
   Siz:=I;
   Np:=Siz mod 10;
   Nk:=Siz div 10;
  if I<>Na then
  begin
```

```
        if (Np=1) and (Nk>0) then
         write(F2,' ',ZnaprX[I]:6:3);
        if (Np=1) and (Nk=0) then
         write(F2,' ',ZnaprX[I]:6:3);
        if (Np>1) and (Np<=9) then
         write(F2,' ',ZnaprX[I]:6:3);
        if Np=0 then
          writeln(F2,' ',ZnaprX[I]:6:3);
       end;
       if I=Na then
        begin
            if (Np=1) and (Nk>0) then
             writeln(F2,' ',ZnaprX[I]:6:3)
            else
            writeln(F2,' ',ZnaprX[I]:6:3);
        end;
    end;
    writeln(F2,' ');
    writeln(F2,
      '  Массив точек по оси ординат не входящих в обод графа: ');
    writeln(F2,' ');
    for I:=1 to Na do
    begin
        Siz:=I;
        Np:=Siz mod 10;
        Nk:=Siz div 10;
       if I<>Na then
       begin
        if (Np=1) and (Nk>0) then
         write(F2,' ',ZnaprY[I]:6:3);
        if (Np=1) and (Nk=0) then
         write(F2,' ',ZnaprY[I]:6:3);
        if (Np>1) and (Np<=9) then
         write(F2,' ',ZnaprY[I]:6:3);
        if Np=0 then
          writeln(F2,' ',ZnaprY[I]:6:3);
       end;
       if I=Na then
        begin
            if (Np=1) and (Nk>0) then
             writeln(F2,' ',ZnaprY[I]:6:3)
            else
            writeln(F2,' ',ZnaprY[I]:6:3);
        end;
    end;
close(F2);
writeln('Вес взят');
end.
```

## 10.5. Входные и выходные данные

Смежность графа:

вершина $v_1$: $v_2$ $v_3$ $v_4$ $v_5$ $v_{21}$ $v_{22}$ $v_{31}$
вершина $v_2$: $v_1$ $v_3$ $v_4$ $v_5$ $v_6$ $v_{10}$ $v_{11}$
вершина $v_3$: $v_1$ $v_2$ $v_4$ $v_5$ $v_6$ $v_{10}$ $v_{11}$ $v_{16}$ $v_{17}$
вершина $v_4$: $v_1$ $v_2$ $v_3$ $v_5$ $v_{11}$ $v_{16}$ $v_{17}$ $v_{20}$ $v_{21}$
вершина $v_5$: $v_1$ $v_2$ $v_3$ $v_4$ $v_{16}$ $v_{20}$ $v_{21}$ $v_{22}$ $v_{31}$
вершина $v_6$: $v_2$ $v_3$ $v_7$ $v_8$ $v_9$ $v_{10}$ $v_{11}$ $v_{12}$ $v_{13}$ $v_{17}$ $v_{18}$
вершина $v_7$: $v_6$ $v_8$ $v_9$ $v_{10}$ $v_{11}$ $v_{12}$ $v_{13}$ $v_{14}$ $v_{17}$ $v_{18}$ $v_{25}$ $v_{26}$ $v_{31}$
вершина $v_8$: $v_6$ $v_7$ $v_9$ $v_{11}$ $v_{12}$ $v_{14}$ $v_{17}$ $v_{18}$ $v_{25}$ $v_{26}$ $v_{30}$ $v_{31}$
вершина $v_9$: $v_6$ $v_7$ $v_8$ $v_{11}$ $v_{17}$ $v_{18}$ $v_{19}$ $v_{24}$ $v_{26}$ $v_{30}$
вершина $v_{10}$: $v_2$ $v_3$ $v_6$ $v_7$ $v_{11}$ $v_{12}$ $v_{13}$ $v_{15}$ $v_{27}$ $v_{28}$
вершина $v_{11}$: $v_2$ $v_3$ $v_4$ $v_6$ $v_7$ $v_8$ $v_9$ $v_{10}$ $v_{16}$ $v_{17}$ $v_{18}$

вершина $v_{12}$: $v_6$ $v_7$ $v_8$ $v_{10}$ $v_{13}$ $v_{14}$ $v_{15}$ $v_{25}$ $v_{26}$ $v_{31}$
вершина $v_{13}$: $v_6$ $v_7$ $v_{10}$ $v_{12}$ $v_{14}$ $v_{15}$ $v_{27}$ $v_{28}$
вершина $v_{14}$: $v_7$ $v_8$ $v_{12}$ $v_{13}$ $v_{15}$ $v_{25}$ $v_{26}$ $v_{27}$ $v_{29}$ $v_{31}$
вершина $v_{15}$: $v_{10}$ $v_{12}$ $v_{13}$ $v_{14}$ $v_{27}$ $v_{28}$ $v_{29}$ $v_{31}$
вершина $v_{16}$: $v_3$ $v_4$ $v_5$ $v_{11}$ $v_{17}$ $v_{18}$ $v_{19}$ $v_{20}$ $v_{21}$ $v_{23}$
вершина $v_{17}$: $v_3$ $v_4$ $v_6$ $v_7$ $v_8$ $v_9$ $v_{11}$ $v_{16}$ $v_{18}$ $v_{19}$ $v_{20}$ $v_{23}$
вершина $v_{18}$: $v_6$ $v_7$ $v_8$ $v_9$ $v_{11}$ $v_{16}$ $v_{17}$ $v_{19}$ $v_{20}$ $v_{23}$ $v_{24}$ $v_{30}$
вершина $v_{19}$: $v_9$ $v_{16}$ $v_{17}$ $v_{18}$ $v_{20}$ $v_{22}$ $v_{23}$ $v_{24}$ $v_{25}$ $v_{30}$ $v_{31}$
вершина $v_{20}$: $v_4$ $v_5$ $v_{16}$ $v_{17}$ $v_{18}$ $v_{19}$ $v_{21}$ $v_{22}$ $v_{23}$
вершина $v_{21}$: $v_1$ $v_4$ $v_5$ $v_{16}$ $v_{20}$ $v_{22}$ $v_{23}$ $v_{31}$
вершина $v_{22}$: $v_1$ $v_5$ $v_{19}$ $v_{20}$ $v_{21}$ $v_{23}$ $v_{24}$ $v_{25}$ $v_{31}$
вершина $v_{23}$: $v_{16}$ $v_{17}$ $v_{18}$ $v_{19}$ $v_{20}$ $v_{21}$ $v_{22}$ $v_{24}$ $v_{25}$ $v_{31}$
вершина $v_{24}$: $v_9$ $v_{18}$ $v_{19}$ $v_{22}$ $v_{23}$ $v_{25}$ $v_{26}$ $v_{30}$ $v_{31}$
вершина $v_{25}$: $v_7$ $v_8$ $v_{12}$ $v_{14}$ $v_{19}$ $v_{22}$ $v_{23}$ $v_{24}$ $v_{26}$ $v_{30}$ $v_{31}$
вершина $v_{26}$: $v_7$ $v_8$ $v_9$ $v_{12}$ $v_{14}$ $v_{24}$ $v_{25}$ $v_{30}$ $v_{31}$
вершина $v_{27}$: $v_{10}$ $v_{13}$ $v_{14}$ $v_{15}$ $v_{28}$ $v_{29}$ $v_{31}$
вершина $v_{28}$: $v_{10}$ $v_{13}$ $v_{15}$ $v_{27}$ $v_{29}$
вершина $v_{29}$: $v_{14}$ $v_{15}$ $v_{27}$ $v_{28}$ $v_{31}$
вершина $v_{30}$: $v_8$ $v_9$ $v_{18}$ $v_{19}$ $v_{24}$ $v_{25}$ $v_{26}$
вершина $v_{31}$: $v_1$ $v_5$ $v_7$ $v_8$ $v_{12}$ $v_{14}$ $v_{15}$ $v_{19}$ $v_{21}$ $v_{22}$ $v_{23}$ $v_{24}$ $v_{25}$ $v_{26}$ $v_{27}$ $v_{29}$

Инцидентность графа:

ребро $e_1$: $(v_1,v_2)$ или $(v_2,v_1)$;
ребро $e_2$: $(v_1,v_3)$ или $(v_3,v_1)$;
ребро $e_3$: $(v_1,v_4)$ или $(v_4,v_1)$;
ребро $e_4$: $(v_1,v_5)$ или $(v_5,v_1)$;
ребро $e_5$: $(v_1,v_{21})$ или $(v_{21},v_1)$;
ребро $e_6$: $(v_1,v_{22})$ или $(v_{22},v_1)$;
ребро $e_7$: $(v_1,v_{31})$ или $(v_{31},v_1)$;
ребро $e_8$: $(v_2,v_3)$ или $(v_3,v_2)$;
ребро $e_9$: $(v_2,v_4)$ или $(v_4,v_2)$;
ребро $e_{10}$: $(v_2,v_5)$ или $(v_5,v_2)$;
ребро $e_{11}$: $(v_2,v_6)$ или $(v_6,v_2)$;
ребро $e_{12}$: $(v_2,v_{10})$ или $(v_{10},v_2)$;
ребро $e_{13}$: $(v_2,v_{11})$ или $(v_{11},v_2)$;
ребро $e_{14}$: $(v_3,v_4)$ или $(v_4,v_3)$;
ребро $e_{15}$: $(v_3,v_5)$ или $(v_5,v_3)$;
ребро $e_{16}$: $(v_3,v_6)$ или $(v_6,v_3)$;
ребро $e_{17}$: $(v_3,v_{10})$ или $(v_{10},v_3)$;
ребро $e_{18}$: $(v_3,v_{11})$ или $(v_{11},v_3)$;
ребро $e_{19}$: $(v_3,v_{16})$ или $(v_{16},v_3)$
ребро $e_{20}$: $(v_3,v_{17})$ или $(v_{17},v_3)$;
ребро $e_{21}$: $(v_4,v_5)$ или $(v_5,v_4)$;
ребро $e_{22}$: $(v_4,v_{11})$ или $(v_{11},v_4)$;
ребро $e_{23}$: $(v_4,v_{16})$ или $(v_{16},v_4)$;
ребро $e_{24}$: $(v_4,v_{17})$ или $(v_{17},v_4)$;
ребро $e_{25}$: $(v_4,v_{20})$ или $(v_{20},v_4)$;
ребро $e_{26}$: $(v_4,v_{21})$ или $(v_{21},v_4)$;
ребро $e_{27}$: $(v_5,v_{16})$ или $(v_{16},v_5)$;
ребро $e_{28}$: $(v_5,v_{20})$ или $(v_{20},v_5)$;
ребро $e_{29}$: $(v_5,v_{21})$ или $(v_{21},v_5)$;
ребро $e_{30}$: $(v_5,v_{22})$ или $(v_{22},v_5)$;
ребро $e_{31}$: $(v_5,v_{31})$ или $(v_{31},v_5)$;
ребро $e_{32}$: $(v_6,v_7)$ или $(v_7,v_6)$;
ребро $e_{33}$: $(v_6,v_8)$ или $(v_8,v_6)$;
ребро $e_{34}$: $(v_6,v_9)$ или $(v_9,v_6)$;
ребро $e_{35}$: $(v_6,v_{10})$ или $(v_{10},v_6)$;
ребро $e_{36}$: $(v_6,v_{11})$ или $(v_{11},v_6)$;
ребро $e_{37}$: $(v_6,v_{12})$ или $(v_{12},v_6)$;
ребро $e_{38}$: $(v_6,v_{13})$ или $(v_{13},v_6)$;
ребро $e_{39}$: $(v_6,v_{17})$ или $(v_{17},v_6)$;
ребро $e_{40}$: $(v_6,v_{18})$ или $(v_{18},v_6)$;
ребро $e_{41}$: $(v_7,v_8)$ или $(v_8,v_7)$;
ребро $e_{42}$: $(v_7,v_9)$ или $(v_9,v_7)$;
ребро $e_{43}$: $(v_7,v_{10})$ или $(v_{10},v_7)$;
ребро $e_{44}$: $(v_7,v_{11})$ или $(v_{11},v_7)$;
ребро $e_{45}$: $(v_7,v_{12})$ или $(v_{12},v_7)$;
ребро $e_{46}$: $(v_7,v_{13})$ или $(v_{13},v_7)$;
ребро $e_{47}$: $(v_7,v_{14})$ или $(v_{14},v_7)$;
ребро $e_{48}$: $(v_7,v_{17})$ или $(v_{17},v_7)$;
ребро $e_{49}$: $(v_7,v_{18})$ или $(v_{18},v_7)$;
ребро $e_{50}$: $(v_7,v_{25})$ или $(v_{25},v_7)$;
ребро $e_{51}$: $(v_7,v_{26})$ или $(v_{26},v_7)$;
ребро $e_{52}$: $(v_7,v_{31})$ или $(v_{31},v_7)$;
ребро $e_{53}$: $(v_8,v_9)$ или $(v_9,v_8)$;
ребро $e_{54}$: $(v_8,v_{11})$ или $(v_{11},v_8)$;
ребро $e_{55}$: $(v_8,v_{12})$ или $(v_{12},v_8)$;
ребро $e_{56}$: $(v_8,v_{14})$ или $(v_{14},v_8)$;
ребро $e_{57}$: $(v_8,v_{17})$ или $(v_{17},v_8)$;
ребро $e_{58}$: $(v_8,v_{18})$ или $(v_{18},v_8)$;
ребро $e_{59}$: $(v_8,v_{25})$ или $(v_{25},v_8)$;
ребро $e_{60}$: $(v_8,v_{26})$ или $(v_{26},v_8)$;
ребро $e_{61}$: $(v_8,v_{30})$ или $(v_{30},v_8)$;
ребро $e_{62}$: $(v_8,v_{31})$ или $(v_{31},v_8)$;
ребро $e_{63}$: $(v_9,v_{11})$ или $(v_{11},v_9)$;
ребро $e_{64}$: $(v_9,v_{17})$ или $(v_{17},v_9)$;
ребро $e_{65}$: $(v_9,v_{18})$ или $(v_{18},v_9)$;
ребро $e_{66}$: $(v_9,v_{19})$ или $(v_{19},v_9)$;
ребро $e_{67}$: $(v_9,v_{24})$ или $(v_{24},v_9)$;
ребро $e_{68}$: $(v_9,v_{26})$ или $(v_{26},v_9)$;
ребро $e_{69}$: $(v_9,v_{30})$ или $(v_{30},v_9)$;
ребро $e_{70}$: $(v_{10},v_{11})$ или $(v_{11},v_{10})$;
ребро $e_{71}$: $(v_{10},v_{12})$ или $(v_{12},v_{10})$;
ребро $e_{72}$: $(v_{10},v_{13})$ или $(v_{13},v_{10})$;
ребро $e_{73}$: $(v_{10},v_{15})$ или $(v_{15},v_{10})$;
ребро $e_{74}$: $(v_{10},v_{27})$ или $(v_{27}, v_{10})$;
ребро $e_{75}$: $(v_{10},v_{28})$ или $(v_{28},v_{10})$;
ребро $e_{76}$: $(v_{11},v_{16})$ или $(v_{16},v_{11})$;
ребро $e_{77}$: $(v_{11},v_{17})$ или $(v_{17},v_{11})$;
ребро $e_{78}$: $(v_{11},v_{18})$ или $(v_{18},v_{11})$;
ребро $e_{79}$: $(v_{12},v_{13})$ или $(v_{13},v_{12})$;
ребро $e_{80}$: $(v_{12},v_{14})$ или $(v_{14},v_{12})$;

ребро $e_{81}$: $(v_{12},v_{15})$ или $(v_{15},v_{12})$;
ребро $e_{82}$: $(v_{12},v_{25})$ или $(v_{25},v_{12})$;
ребро $e_{83}$: $(v_{12},v_{26})$ или $(v_{26},v_{12})$;
ребро $e_{84}$: $(v_{12},v_{31})$ или $(v_{31},v_{12})$;
ребро $e_{85}$: $(v_{13},v_{14})$ или $(v_{14},v_{13})$;
ребро $e_{86}$: $(v_{13},v_{15})$ или $(v_{15},v_{13})$;
ребро $e_{87}$: $(v_{13},v_{27})$ или $(v_{27},v_{13})$;
ребро $e_{88}$: $(v_{13},v_{28})$ или $(v_{28},v_{13})$;
ребро $e_{89}$: $(v_{14},v_{15})$ или $(v_{15},v_{14})$;
ребро $e_{90}$: $(v_{14},v_{25})$ или $(v_{25},v_{14})$;
ребро $e_{91}$: $(v_{14},v_{26})$ или $(v_{26},v_{14})$;
ребро $e_{92}$: $(v_{14},v_{27})$ или $(v_{27},v_{14})$;
ребро $e_{93}$: $(v_{14},v_{29})$ или $(v_{29},v_{14})$;
ребро $e_{94}$: $(v_{14},v_{31})$ или $(v_{31},v_{14})$;
ребро $e_{95}$: $(v_{15},v_{27})$ или $(v_{27},v_{15})$;
ребро $e_{96}$: $(v_{15},v_{28})$ или $(v_{28},v_{15})$;
ребро $e_{97}$: $(v_{15},v_{29})$ или $(v_{29},v_{15})$;
ребро $e_{98}$: $(v_{15},v_{31})$ или $(v_{31},v_{15})$;
ребро $e_{99}$: $(v_{16},v_{17})$ или $(v_{17},v_{16})$;
ребро $e_{100}$: $(v_{16},v_{18})$ или $(v_{18},v_{16})$;
ребро $e_{101}$: $(v_{16},v_{19})$ или $(v_{19},v_{16})$;
ребро $e_{102}$: $(v_{16},v_{20})$ или $(v_{20},v_{16})$;
ребро $e_{103}$: $(v_{16},v_{21})$ или $(v_{21},v_{16})$;
ребро $e_{104}$: $(v_{16},v_{23})$ или $(v_{23},v_{16})$;
ребро $e_{105}$: $(v_{17},v_{18})$ или $(v_{18},v_{17})$;
ребро $e_{106}$: $(v_{17},v_{19})$ или $(v_{19},v_{17})$;
ребро $e_{107}$: $(v_{17},v_{20})$ или $(v_{20},v_{17})$;
ребро $e_{108}$: $(v_{17},v_{23})$ или $(v_{23},v_{17})$;
ребро $e_{109}$: $(v_{18},v_{19})$ или $(v_{19},v_{18})$;
ребро $e_{110}$: $(v_{18},v_{20})$ или $(v_{20},v_{18})$;
ребро $e_{111}$: $(v_{18},v_{23})$ или $(v_{23},v_{18})$;
ребро $e_{112}$: $(v_{18},v_{24})$ или $(v_{24},v_{18})$;
ребро $e_{113}$: $(v_{18},v_{30})$ или $(v_{30},v_{18})$;
ребро $e_{114}$: $(v_{19},v_{20})$ или $(v_{20},v_{19})$;
ребро $e_{115}$: $(v_{19},v_{22})$ или $(v_{22},v_{19})$;
ребро $e_{116}$: $(v_{19},v_{23})$ или $(v_{23},v_{19})$;
ребро $e_{117}$: $(v_{19},v_{24})$ или $(v_{24},v_{19})$;
ребро $e_{118}$: $(v_{19},v_{25})$ или $(v_{25},v_{19})$;
ребро $e_{119}$: $(v_{19},v_{30})$ или $(v_{30},v_{19})$;
ребро $e_{120}$: $(v_{19},v_{31})$ или $(v_{31},v_{19})$;
ребро $e_{121}$: $(v_{20},v_{21})$ или $(v_{21},v_{20})$;
ребро $e_{122}$: $(v_{20},v_{22})$ или $(v_{22},v_{20})$;
ребро $e_{123}$: $(v_{20},v_{23})$ или $(v_{23},v_{20})$;
ребро $e_{124}$: $(v_{21},v_{22})$ или $(v_{22},v_{21})$;
ребро $e_{125}$: $(v_{21},v_{23})$ или $(v_{23},v_{21})$;
ребро $e_{126}$: $(v_{21},v_{31})$ или $(v_{31},v_{21})$;
ребро $e_{127}$: $(v_{22},v_{23})$ или $(v_{23},v_{22})$;
ребро $e_{128}$: $(v_{22},v_{24})$ или $(v_{24},v_{22})$;
ребро $e_{129}$: $(v_{22},v_{25})$ или $(v_{25},v_{22})$;
ребро $e_{130}$: $(v_{22},v_{31})$ или $(v_{31},v_{22})$;
ребро $e_{131}$: $(v_{23},v_{24})$ или $(v_{24},v_{23})$;
ребро $e_{132}$: $(v_{23},v_{25})$ или $(v_{25},v_{23})$;
ребро $e_{133}$: $(v_{23},v_{31})$ или $(v_{31},v_{23})$;
ребро $e_{134}$: $(v_{24},v_{25})$ или $(v_{25},v_{24})$;
ребро $e_{135}$: $(v_{24},v_{26})$ или $(v_{26},v_{24})$;
ребро $e_{136}$: $(v_{24},v_{30})$ или $(v_{30},v_{24})$;
ребро $e_{137}$: $(v_{24},v_{31})$ или $(v_{31},v_{24})$;
ребро $e_{138}$: $(v_{25},v_{26})$ или $(v_{26},v_{25})$;
ребро $e_{139}$: $(v_{25},v_{30})$ или $(v_{30},v_{25})$;
ребро $e_{140}$: $(v_{25},v_{31})$ или $(v_{31},v_{25})$;
ребро $e_{141}$: $(v_{26},v_{30})$ или $(v_{30},v_{26})$;
ребро $e_{142}$: $(v_{26},v_{31})$ или $(v_{31},v_{26})$;
ребро $e_{143}$: $(v_{27},v_{28})$ или $(v_{28},v_{27})$;
ребро $e_{144}$: $(v_{27},v_{29})$ или $(v_{29},v_{27})$;
ребро $e_{145}$: $(v_{27},v_{31})$ или $(v_{31},v_{27})$;
ребро $e_{146}$: $(v_{28},v_{29})$ или $(v_{29},v_{28})$;
ребро $e_{147}$: $(v_{29},v_{31})$ или $(v_{31},v_{29})$.

**Итерация 1:**

**файл 31.gm1**

| | |
|---|---|
| 6 | количество вершин в ободе |
| 31 29 28 10 2 1 | номера вершин |
| 0.0 0.67 100.0 100.0 67.0 0.0 | координаты вершин по оси абцисс |
| 100.0 100.0 67.0 0.0 0.0 0.00 | координаты вершин по оси ординат |

**файл 31.gm2**

Количество точек в ободе: 6
Массив вершин входящих в обод графа:
31 29 28 10 2 1

Координаты вершин по оси абсцисс в ободе:
0.000 0.670 100.000 100.000 67.000 0.000

Координаты вершин по оси ординат в ободе:
100.000 100.000 67.000 0.000 0.000 0.000

Массив вершин не входящих в обод графа:
27 15 3 5 21 13 4 22 30 14
26 24 12 11 23 6 20 19 9 25
16 8 7 17 18

Массив вершин по оси абсцисс не входящих в обод графа:
47.712 45.792 40.839 24.928 18.030 55.860 30.500 17.346 25.266 30.665

25.447 21.885 39.426 42.498 21.393 45.679 24.252 22.513 30.045 23.110
25.390 30.021 35.607 30.558 29.842

Массив вершин по оси ординат не входящих в обод графа:
62.948 61.688 22.107 34.887 41.894 48.438 25.921 47.673 51.183 63.507
57.662 55.664 51.611 31.019 50.709 34.154 40.400 50.577 44.271 54.924
34.986 51.202 45.598 35.994 43.980

**Итерация 2:**

**файл 31.gm1**

| | |
|---|---|
| 14 | количество жестко установленных вершин |
| 31 29 28 10 2 1 22 25 14 27<br>13 6 3 5 | номера жестко установленных вершин |
| 0.0 0.67 100.0 100.0 67.0 0.0 08.0 11.0 15.0 83.0<br>75.0 69.0 38.0 13.0 | координаты вершин по оси абцисс |
| 100.0 100.0 67.0 0.0 0.0 0.0 68.0 75.0 82.0 84.0<br>25.0 19.0 11.0 34.0 | координаты вершин по оси ординат |

**файл 31.gm2**

Количество точек в ободе: 14
Массив точек входящих в обод графа:
31 29 28 10 2 1 22 25 14 27
13 6 3 5

Координаты точек по оси абсцисс в ободе:
0.000 0.670 100.000 100.000 67.000 0.000 8.000 11.000 15.000 83.000
75.000 69.000 38.000 13.000

Координаты точек по оси ординат в ободе:
100.000 100.000 67.000 0.000 0.000 0.000 68.000 75.000 82.000 84.000
25.000 19.000 11.000 34.000

Массив точек не входящих в обод графа:
15 12 21 4 30 24 26 11 7 23
20 19 8 16 9 17 18

Массив точек по оси абсцисс не входящих в обод графа:
52.036 41.617 14.525 28.827 24.179 18.499 23.460 45.263 39.276 18.301
21.533 20.628 30.394 25.539 33.541 33.586 31.728

Массив точек по оси ординат не входящих в обод графа:
63.941 53.529 45.913 27.929 59.332 63.767 64.540 31.001 49.554 58.164
45.820 59.123 55.259 40.390 48.922 41.240 47.714

## Комментарии

В данной главе рассматриваются вопросы построения геометрического изображения плоской части графа. Показана тесная связь построения геометрического изображения с построением графа уровней для топологического рисунка графа. Построение графа уровней позволяет создать последовательность расположения вершин в геометрическом пространстве. Линия уровня графа в топологическом рисунке - это замкнутая линия с определенной последовательностью расположения вершин. Пространственное расположение образов линий уровня в виде окружности, эллипса, контура прямоугольника позволяет определить геометрические координаты вершин. Представлен модифицированный метод и программа для построения плоского рисунка графа силовым алгоритмом.

## Выводы

Следует различать следующие три вида топологических рисунков графов:

- топологический рисунок максимально плоского суграфа;
- топологический рисунок графа с минимальным числом пересечений;
- топологический рисунок графа минимальной толщины.

Последние два вида топологических рисунков строятся с участием мнимых вершин.

Вне зависимости от вида графа, топологический рисунок определяется вращением вершин и индуцированной этим вращением системе циклов.

Основой построения любого вида графа является множество изометрических циклов графа. В работе рассматриваются методы решения задачи построения топологического рисунка плоской части непланарного графа методами дискретной оптимизации[]. Указаны три метода выделения базиса изометрических циклов:

- выделение независимой системы циклов алгоритмом Гаусса;
- метод градиентного спуска;
- методы алгебры структурных чисел.

Рассмотрены недостатки этих методов.

Существование нескольких методов выделения базиса подпространства циклов, позволяет создавать комбинаторные решения для практического применения математических моделей [5,26]. Однако каждая сложная задача дискретной оптимизации специфична. Это находит отражение в методах разбиения задачи на подзадачи, определения последовательности этапов решения. Требует увязки нескольких математических моделей их представления и описания.

Решение задачи построения топологического рисунка плоской части непланарного графа с минимальным числом удаленных ребер относится к сложно решаемым задачам.

На первом этапе производится выделение базиса состоящего из изометрических циклов графа различными методами.

Второй этап предназначен для выделения из базиса подмножество циклов с нулевым значением кубического функционала Маклейна и соблюдения условий Эйлера.

На третьем этапе, методами векторной алгебры пересечений производится введение дополнительных ребер и формирование дополнительного подмножества циклов.

Последний четвертый этап осуществляет преобразование подмножества циклов во вращение вершин графа, тем самым создавая описание топологического рисунка плоского суграфа.

Основой решения задачи построения топологического рисунка с минимальным числом пересечений и задачи построения топологического рисунка графа минимальной толщины является топологический рисунок максимально плоского суграфа. Так как в графе может существовать несколько топологических рисунков максимально плоского суграфа, то все подмножество таких рисунков должно быть внесено в базу данных рассматриваемого графа.

Окончательным этапом решения, служит выбор оптимального варианта из матрицы решений.

В работе показано, что существует общее связывающее звено всех методов – это последовательность выбора элемента. Последовательнось выбора элемента определяется кортежем фрагментарной структуры.

Показано применение метаэвристических алгоритмов основанных на фрагментарно-эволюционных структурах для решения задач указанного типа [34,39].

Рассмоьрен алгоритм решения задачи построения последовательности расположения вершин в уровнях для построения геометрического изображения рисунка графа на плоскости.

Представлены экспериментальные программы Ra1 и Rebro2 для генерации плоских конфигураций. Для построения геометрического образа плоского суграфа представлена программа Raschet9 решения линейной системы алгебраических уравнений.

## СПИСОК ИСПОЛЬЗУЕМОЙ ЛИТЕРАТУРЫ